%% file: K-holonomic-smf.tex
\documentclass[10pt,english]{smfbook}
\usepackage{amsmath,amscd,amssymb}
\usepackage[english,francais]{babel}

\author{Takuro Mochizuki}
\address{Research Institute for Mathematical Sciences,
Kyoto University, Kyoto 606-8502, Japan}
\email{takuro@kurims.kyoto-u.ac.jp}
\title{Holonomic $\nbigd$-modules with
 Betti structure}
\alttitle{}

\input{notation}
\input{new_theorem}

\makeindex

\begin{document}

\frontmatter


\begin{abstract}
We define the notion of Betti structure
for holonomic $\nbigd$-modules which are not necessarily 
regular singular.
We establish the fundamental functorial properties.
We also give auxiliary analysis 
of holomorphic functions of various types
on the real blow up.
\end{abstract}

\begin{altabstract}
Nous d\'{e}finissons la notion de la structure Betti
pour les $\nbigd$-modules holonome
qui ne sont pas n\'{e}cessairement 
singuli\`{e}re r\'{e}guli\`{e}re.
Nous \'{e}tablissons les propri\'et\'e fondamentaux.
Nous aussi donnons l'analyse suppl\'{e}mentaire
pour les fonctions holomorphes diverses
sur l'\'{e}clatement r\'{e}el.
\end{altabstract}

\subjclass{14F10, 32C38}
\keywords{holonomic D-modules, Betti structure,
Stokes structure}
\altkeywords{}

\maketitle
\tableofcontents

\mainmatter

\chapter{Introduction}

\input{1}

\chapter{Preliminary}
\section{Notation and words}

\input{2.1}

\section{Beilinson's construction}
\label{subsection;09.10.19.32}
\input{2.2}

\chapter{Good holonomic $\nbigd$-modules
and their de Rham complexes}

\section{Good holonomic $\nbigd$-modules}
\input{3.1}

\section{De Rham complexes}
\input{3.2}

\chapter{Some sheaves on the real blow up}

\section{Holomorphic functions}
\input{4.1}

\section{$C^{\infty}$-functions}
\input{4.2}

\section{Nilsson type functions}
\input{4.3}

\section{Push-forward}
\input{4.4}

\section{Characterization by growth order}
\label{subsection;14.1.19.30}
\input{4.5}

\section[Flatness]{Flatness of the sheaf of holomorphic functions
with moderate growth}
\input{4.6}

\section{Push-forward of good $\nbigd$-modules
and real blow up}
\label{subsection;14.1.18.1}
\input{4.7}

\chapter[Complexes]{Complexes on the real blow up
associated to good meromorphic flat bundles}

\section{De Rham complexes}
\input{5.1}

\section{Duality}
\input{5.2}

\section{Functoriality}
\label{subsection;09.10.28.13}
\input{5.3}

\section{A rigidity property (Appendix)}
\input{5.4}

\chapter{Good $K$-structure}

\section{Good meromorphic flat bundles}
\input{6.1}

\section{Good holonomic $\nbigd$-modules with
good $K$-structure (Local case)}
\input{6.2}

\section{Good pre-$K$-holonomic $\nbigd$-modules}
\label{subsection;13.4.27.21}
\input{6.3}

\section{Meromorphic flat connections with good $K$-structure}
\label{subsection;13.4.25.200}
\input{6.4}

\section{Preliminary for push-forward}
\input{6.5}

\chapter{$K$-holonomic $\nbigd$-modules}

\section{Preliminary}
\input{7.1}

\section{$K$-Betti structure}
\label{subsection;09.11.11.2}
\input{7.2}

\section{$K(\ast D)$-Betti structure}
\input{7.3}

\chapter{Functoriality properties}
\label{section;09.11.11.3}

\section{Statements}
\label{subsection;13.4.27.300}
\input{8.1}

\section{Step 1}
\label{subsection;09.12.5.30}
\input{8.2}

\section{Step 2}
\label{subsection;09.12.5.31}
\input{8.3}

\section{Some resolutions}
\input{8.4}

\section{Step 3}
\label{subsection;09.10.17.150}
\input{8.5}

\chapter{Derived category of
algebraic $K$-holonomic $\nbigd$-modules}
\input{9}

\section{Standard exact functors}
\input{9.1}

\section{Push-forward and pull-back}
\input{9.2}

\section{Tensor product and inner homomorphism}
\input{9.3}

\section{$K$-structure of the space of morphisms}
\input{9.4}

\backmatter

\input{K-holonomic_ref}
\printindex

\end{document}

%% file: notation.tex
\newcommand{\nbiga}{\mathcal{A}}
\newcommand{\nbigb}{\mathcal{B}}
\newcommand{\nbigc}{\mathcal{C}}
\newcommand{\nbigd}{\mathcal{D}}
\newcommand{\nbige}{\mathcal{E}}
\newcommand{\nbigf}{\mathcal{F}}
\newcommand{\nbigg}{\mathcal{G}}
\newcommand{\nbigh}{\mathcal{H}}
\newcommand{\nbigi}{\mathcal{I}}
\newcommand{\nbigj}{\mathcal{J}}
\newcommand{\nbigk}{\mathcal{K}}
\newcommand{\nbigl}{\mathcal{L}}
\newcommand{\nbigm}{\mathcal{M}}
\newcommand{\nbign}{\mathcal{N}}
\newcommand{\nbigo}{\mathcal{O}}
\newcommand{\nbigp}{\mathcal{P}}
\newcommand{\nbigq}{\mathcal{Q}}
\newcommand{\nbigr}{\mathcal{R}}
\newcommand{\nbigs}{\mathcal{S}}
\newcommand{\nbigt}{\mathcal{T}}
\newcommand{\nbigu}{\mathcal{U}}
\newcommand{\nbigv}{\mathcal{V}}

\newcommand{\proj}{\mathbb{P}}
\newcommand{\seisuu}{{\mathbb Z}}
\newcommand{\rnum}{{\mathbb Q}}

\newcommand{\cnum}{{\mathbb C}}
\newcommand{\real}{{\mathbb R}}

\newcommand{\gbigh}{\mathfrak H}
\newcommand{\gbigi}{\mathfrak I}

\newcommand{\gbigp}{\mathfrak P}

\newcommand{\gminia}{\mathfrak a}
\newcommand{\gminib}{\mathfrak b}

\newcommand{\gminii}{\mathfrak i}
\newcommand{\gminij}{\mathfrak j}

\newcommand{\vece}{{\boldsymbol e}}

\newcommand{\vecw}{{\boldsymbol w}}

\newcommand{\veczero}{{\boldsymbol 0}}
\newcommand{\vecalpha}{{\boldsymbol \alpha}}
\newcommand{\veca}{{\boldsymbol a}}
\newcommand{\vecb}{{\boldsymbol b}}
\newcommand{\vecbeta}{{\boldsymbol \beta}}
\newcommand{\vecdelta}{{\boldsymbol \delta}}

\newcommand{\veck}{{\boldsymbol k}}
\newcommand{\vecm}{{\boldsymbol m}}

\newcommand{\vecL}{{\boldsymbol L}}

\newcommand{\vecn}{{\boldsymbol n}}

\newcommand{\veczeta}{{\boldsymbol \zeta}}
\newcommand{\vecz}{{\boldsymbol z}}

\newcommand{\larr}{\leftarrow}
\newcommand{\llarr}{\longleftarrow}

\newcommand{\rarr}{\rightarrow}
\newcommand{\lrarr}{\longrightarrow}
\newcommand{\darr}{\downarrow}

\newcommand{\pf}{{\bf Proof}\hspace{.1in}}

\def\Hom{\mathop{\rm Hom}\nolimits}

\def\Ext{\mathop{\rm Ext}\nolimits}

\def\Cok{\mathop{\rm Cok}\nolimits}

\def\Image{\mathop{\rm Im}\nolimits}

\def\Re{\mathop{\rm Re}\nolimits}

\def\Gr{\mathop{\rm Gr}\nolimits}

\def\Tot{\mathop{\rm Tot}\nolimits}
\def\Cone{\mathop{\rm Cone}\nolimits}
\def\rank{\mathop{\rm rank}\nolimits}
\def\Spec{\mathop{\rm Spec}\nolimits}

\def\Ker{\mathop{\rm Ker}\nolimits}

\def\Gr{\mathop{\rm Gr}\nolimits}

\def\Res{\mathop{\rm Res}\nolimits}

\def\ord{\mathop{\rm ord}\nolimits}

\def\can{\mathop{\rm can}\nolimits}
\def\var{\mathop{\rm var}\nolimits}
\def\id{\mathop{\rm id}\nolimits}

\def\gcd{\mathop{\rm g.c.d.}\nolimits}
\def\codim{\mathop{\rm codim}\nolimits}

\def\Supp{\mathop{\rm Supp}\nolimits}

\def\Irr{\mathop{\rm Irr}\nolimits}

\newcommand{\del}{\partial}
\newcommand{\delbar}{\overline{\del}}

\newcommand{\nhom}{{\mathcal Hom}}

\newcommand{\next}{{\mathcal Ext}}

\newcommand{\jbar}{\underline{j}}
\newcommand{\mbar}{\underline{m}}

\newcommand{\ibar}{\underline{i}}

\newcommand{\barz}{\overline{z}}
\newcommand{\zbar}{\barz}

\newcommand{\fbar}{\overline{f}}

\newcommand{\lefttop}[1]{{}^{#1}\!}

\def\reg{\mathop{\rm reg}\nolimits}

\def\nil{\mathop{\rm nil}\nolimits}
\def\Nil{\mathop{\rm Nil}\nolimits}

\newcommand{\tildepsi}{\widetilde{\psi}}
\newcommand{\psitilde}{\tildepsi}

\newcommand{\closedopen}[2]{[#1,#2[}
\newcommand{\openclosed}[2]{]#1,#2]}
\newcommand{\openopen}[2]{]#1,#2[}

\newcommand{\Gtilde}{\widetilde{G}}
\newcommand{\rhotilde}{\widetilde{\rho}}

\newcommand{\Vhat}{\widehat{V}}
\newcommand{\nablahat}{\widehat{\nabla}}
\newcommand{\Vtilde}{\widetilde{V}}

\newcommand{\Dhat}{\widehat{D}}

\newcommand{\Sbar}{\overline{S}}
\newcommand{\nbigehat}{\widehat{\nbige}}

\newcommand{\ellsitabar}{\underline{\ell}}
\newcommand{\Utilde}{\widetilde{U}}
\newcommand{\Dtilde}{\widetilde{D}}
\newcommand{\Xtilde}{\widetilde{X}}

\newcommand{\Zhat}{\widehat{Z}}

\newcommand{\nbigmtilde}{\widetilde{\nbigm}}

\def\ord{\mathop{\rm ord}\nolimits}

\def\moderate{\mathop{\rm mod}\nolimits}

\def\Hol{\mathop{\rm Hol}\nolimits}

\def\Glue{\mathop{\rm Glue}\nolimits}

\def\DR{\mathop{\rm DR}\nolimits}
\def\ob{\mathop{\rm ob}\nolimits}
\def\Mor{\mathop{\rm Mor}\nolimits}
\def\Ch{\mathop{\rm Ch}\nolimits}
\def\good{\mathop{\rm good}\nolimits}

\def\hol{\mathop{\rm hol}\nolimits}
\def\Loc{\mathop{\rm Loc}\nolimits}
\def\Per{\mathop{\rm Per}\nolimits}
\def\rapid{\mathop{\rm rapid}\nolimits}
\def\Cat{\mathop{\rm Cat}\nolimits}
\def\Tor{\mathop{\rm Tor}\nolimits}
\def\Forget{\mathop{\rm Forget}\nolimits}
\def\pre{\mathop{\rm pre}\nolimits}
\def\Gd{\mathop{\rm Gd}\nolimits}

\newcommand{\Ibar}{\overline{I}}

\newcommand{\vecj}{{\boldsymbol j}}
\newcommand{\vecD}{{\boldsymbol D}}

\newcommand{\nbigutilde}{\widetilde{\nbigu}}

\newcommand{\Phat}{\widehat{P}}

\newcommand{\gtilde}{\widetilde{g}}
\newcommand{\Ztilde}{\widetilde{Z}}

\newcommand{\varphitilde}{\widetilde{\varphi}}

\newcommand{\Deltatilde}{\widetilde{\Delta}}

\newcommand{\Ytilde}{\widetilde{Y}}
\newcommand{\Ybar}{\overline{Y}}

\newcommand{\itilde}{\widetilde{i}}

\newcommand{\nbigitilde}{\widetilde{\nbigi}}

\newcommand{\vecnbigi}{{\boldsymbol \nbigi}}

\newcommand{\nrhom}{R{\mathcal Hom}}
\newcommand{\DDD}{\boldsymbol D}
\newcommand{\jtilde}{\widetilde{j}}
\newcommand{\vecH}{{\boldsymbol H}}

\newcommand{\Kbar}{\overline{K}}

\newcommand{\Hhat}{\widehat{H}}

\newcommand{\gminiatilde}{\widetilde{\gminia}}
\newcommand{\gminibtilde}{\widetilde{\gminib}}
\newcommand{\cnumtilde}{\widetilde{\cnum}}

\newcommand{\phitilde}{\widetilde{\phi}}

\newcommand{\Xbar}{\overline{X}}

\newcommand{\Gammatilde}{\widetilde{\Gamma}}

\newcommand{\piinverseDhat}{\widehat{\pi^{-1}(D)}}

\newcommand{\veci}{\boldsymbol i}

\newcommand{\mnuleq}{\prec}

\newcommand{\nbigctilde}{\widetilde{\nbigc}}
\newcommand{\iotatilde}{\widetilde{\iota}}
\newcommand{\nbiggtilde}{\widetilde{\nbigg}}
\newcommand{\Jbar}{\overline{J}}
\newcommand{\nbigmbar}{\overline{\nbigm}}

\newcommand{\nbigcbar}{\overline{\nbigc}}
\newcommand{\nbigptilde}{\widetilde{\nbigp}}

\newcommand{\nbigqtilde}{\widetilde{\nbigq}}
\newcommand{\bikkuri}{!}
\newcommand{\nbigotilde}{\widetilde{\nbigo}}
\newcommand{\lambdatilde}{\widetilde{\lambda}}
\newcommand{\Lhat}{\widehat{L}}

%% file: new_theorem.tex
\newtheorem{thm}{Theorem}[section]
\newtheorem{cor}[thm]{Corollary}

\newtheorem{rem}[thm]{Remark}
\newtheorem{lem}[thm]{Lemma}
\newtheorem{prop}[thm]{Proposition}
\newtheorem{df}[thm]{Definition}

\newtheorem{notation}[thm]{Notation}

%% file: 1.tex
In this paper, we introduce the notion of Betti structure for
holonomic $\nbigd$-modules,
motivated by a question in 
\cite{Deligne-Malgrange-Ramis}.
For regular holonomic $\nbigd$-modules,
it is clearly defined  
by the Riemann-Hilbert correspondence,
which is a basis of the theory of mixed Hodge modules
(\cite{saito1}--\cite{saito4}).
Namely, a Betti structure of 
a regular holonomic $\nbigd_X$-module $\nbigm$
is defined to be a $\rnum$-perverse sheaf $\nbigf$
with an isomorphism
$\alpha:\nbigf\otimes\cnum\simeq \DR_X\nbigm$.
It has a nice functorial property
for some of standard functors such as
pull back, push-forward, dual etc.,
in the algebraic situation.

As for the non-regular case,
there has been a significant progress toward
a generalized Riemann-Hilbert correspondence
between holonomic $\nbigd$-modules
and some topological objects,
a kind of perverse sheaves
equipped with ``Stokes structure''
in some sense.
The asymptotic analysis for good meromorphic flat bundles
(\cite{majima}, \cite{sabbah4} and \cite{mochi7})
and the existence of resolution of turning points
(\cite{kedlaya}, \cite{kedlaya2}, \cite{mochi7})
lead us a rather satisfactory understanding of 
the structure of meromorphic flat bundles.
Moreover, the recent work
of A. D'Agnolo and M. Kashiwara
\cite{DAgnolo-Kashiwara},
\cite{DAgnolo-Kashiwara2}
based on the theory of Ind-sheaves
\cite{Kashiwara-Schapira-Ind-sheaves}
gives us a description of holonomic $\nbigd$-modules
in terms of some topological objects.
It should also lead us to a thorough
theory of Betti structure of holonomic $\nbigd$-modules.

However,
except in the one dimensional case,
it turned out that
a rather complicated machinery is necessary 
for the complete description of
generalized Riemann-Hilbert correspondence.
(See \cite{DAgnolo-Kashiwara2} and 
\cite{Kashiwara-Schapira-Ind-sheaves}.
See also \cite{sabbah_lecture_Stokes}.)
In this study,
we shall directly define the notion of
``Betti structure'' for holonomic $\nbigd$-modules
with functorial property 
by using only the classical machinery
of holonomic $\nbigd$-modules
and perverse sheaves.
It still requires non-trivial tasks,
and provides us with non-trivial consequences
on the compatibility of the Stokes structure
and the $\rnum$-structure.
We hope that it would be useful
for direct understanding of Betti structures
and for a further study toward
the generalized Riemann-Hilbert correspondence,
at least temporarily.

\section{Pre-Betti structure}

To define the notion of Betti structure of
a holonomic $\nbigd_X$-module $\nbigm$,
it is a most naive idea to consider 
a pair of $\rnum$-perverse sheaf $\nbigf$
and an isomorphism 
$\alpha:\nbigf\otimes\cnum\simeq\DR_X(\nbigm)$
as above,
which is called a pre-Betti structure of $\nbigm$
in this paper.
\index{pre-Betti structure}
A holonomic $\nbigd_X$-module with 
a pre-Betti structure is called
a pre-$\rnum$-holonomic $\nbigd_X$-module.
We should say that pre-Betti structure
is too naive for the following reasons:
\begin{itemize}
\item
It is not so intimately related
with Stokes structure.
\item
Although pre-Betti structures have
nice functoriality with respect to dual and 
proper push-forward,
they are not functorial with respect to
the push-forward for open immersion,
the pull back, 
the nearby cycle and vanishing cycle functors.
Recall that the de Rham functor is
not compatible with the latter class of functors,
when irregular singularities are present.
\end{itemize}
It is the main goal in this paper
to introduce a condition
for a pre-Betti structure 
to be a ``Betti structure''.
We use an inductive way on the dimension
of the support,
which was a strategy of M. Saito
to define his mixed and pure Hodge modules
\cite{saito1} and \cite{saito2}.

In the following,
a $\rnum$-structure of
a $\cnum$-perverse sheaf $\nbigf_{\cnum}$
is a $\rnum$-perverse sheaf $\nbigf_{\rnum}$
with an isomorphism
$\nbigf_{\rnum}\otimes_{\rnum}\cnum
\simeq
 \nbigf_{\cnum}$.

\section{Betti structure in the one dimensional case}

We explain our condition for Betti structure
in the one dimensional case.

\subsection{The generalized Riemann-Hilbert correspondence
in the one dimensional case}

We know the well established theory
on the general structure of holonomic $\nbigd$-modules
on curves
(the generalized Riemann-Hilbert correspondence).
Namely, in the one dimensional case,
we have a natural bijective correspondence
between meromorphic flat bundles
and local systems with Stokes structure,
and any holonomic $\nbigd$-modules 
are described as the gluing of
meromorphic flat bundles
and skyscraper $\nbigd$-modules.
We shall review it very briefly.
For simplicity,
we consider holonomic $\nbigd$-modules
on $X=\Delta=\{|z|<1\}$ which may have 
a singularity at the origin $D=\{O\}$.

\subsubsection{The Stokes structure of meromorphic flat bundles}

Let $V$ be a meromorphic flat bundle
on $(X,D)$.
Let $\pi:\Xtilde(D)\lrarr X$ be the real blow up
along $D$.
Let $\nbigl$ be the local system on $\Xtilde(D)$
associated to the flat bundle $V_{|X-D}$.
Let $P$ be any point of $\pi^{-1}(D)$.
According to the classical asymptotic analysis,
we have the Stokes filtration $\nbigf^P$
of the stalk $\nbigl_{P}$
given by the growth order of flat sections
with respect to any meromorphic frame of $V$.
The meromorphic flat bundle
$V$ can be 
reconstructed from
the flat bundle $V_{|X-D}$ and
the system of filtrations
$\bigl\{\nbigf^P\,\big|\,
 P\in \pi^{-1}(D)\bigr\}$,
which is the Riemann-Hilbert-Birkhoff correspondence
for meromorphic flat bundles on curves.

Let $V^{\lor}$ be the dual of $V$
as a meromorphic flat bundle,
and let $V_!:=\DDD_X V^{\lor}$ be the dual of
$V^{\lor}$ as a $\nbigd_X$-module.
Let us recall that the de Rham complexes
$\DR_X(V)$ and $\DR_X (V_!)$ can be described
in terms of Stokes filtrations.
Let $\nbigl^{\leq D}$ and $\nbigl^{<D}$ be 
the constructible subsheaves of $\nbigl$
such that 
$\nbigl^{\leq D}_P
=\nbigf_{\leq 0}^P(\nbigl_P)$
and 
$\nbigl^{<D}_{P}=
 \nbigf_{<0}^P(\nbigl_P)$.
Then, we have natural isomorphisms:
\begin{equation}
\label{eq;09.11.9.1}
\DR(V)\simeq
 R\pi_{\ast}\nbigl^{\leq D}[1],
\quad\quad
 \DR(V_!)\simeq
 R\pi_{\ast}\nbigl^{<D}[1].
\end{equation}

\subsubsection{Gluing of holonomic $\nbigd$-modules}

Let us very briefly recall 
a key construction due to A. Beilinson
\cite{beilinson2} on the gluing of
holonomic $\nbigd$-modules,
which we will review in 
\S\ref{subsection;09.10.19.32} in more details.
(See also \cite{MacPherson-Vilonen}
and \cite{Verdier} for the other formalisms
for gluing.)
Let $\nbigm$ be any holonomic $\nbigd_X$-module
such that 
$V:=\nbigm(\ast D)$ is a meromorphic flat bundle
on $(X,D)$.
We have the natural morphisms
$V_!\stackrel{a_0}{\lrarr}
 \nbigm\stackrel{b_0}{\lrarr} V$.
According to \cite{beilinson2},
we have the $\nbigd$-modules
$\Xi_z(V)$ and $\psi_z(V)$
associated to $V$,
with morphisms
\begin{equation}
 \label{eq;09.11.10.1}
 \psi_z(V)\stackrel{a_1}{\lrarr}
 \Xi_z(V)\stackrel{b_1}{\lrarr} \psi_z(V),
\quad\quad
 V_!\stackrel{a_2}{\lrarr}
 \Xi_z(V)\stackrel{b_2}{\lrarr} V.
\end{equation}
It can be shown that
$b_0\circ a_0=b_2\circ a_2$.
We also have $b_2\circ a_1=0$
and $b_1\circ a_2=0$.
We obtain the $\nbigd$-module
$\phi_z(\nbigm)$ as the cohomology of
the naturally associated complex:
\begin{equation}
\label{eq;09.11.9.3}
 V_!\lrarr \Xi_z(V)\oplus\nbigm
 \lrarr V
\end{equation}
We have the naturally induced morphisms
$\psi_z(V)\stackrel{\can}{\lrarr}
 \phi_z(\nbigm)\stackrel{\var}{\lrarr}
 \psi_z(V)$.
Then, $\nbigm$ is reconstructed 
as the cohomology of the complex:
\begin{equation}
 \label{eq;09.11.9.2}
 \psi_z(V)\lrarr 
 \Xi_z(V)\oplus\phi_z(\nbigm)\lrarr
 \psi_z(V)
\end{equation}
Recall that
$\Xi_z(V)$,
$\psi_z(V)$,
and $\phi_z(\nbigm)$ are called
the maximal extension,
the nearby cycle sheaf,
and the vanishing cycle sheaf 
of $\nbigm$.

\subsection{Betti structure
of holonomic $\nbigd$-modules on curves}

We explain when
a pre-Betti structure of holonomic $\nbigd$-modules
seems eligible to be called a Betti structure
in the one dimensional case.
Essentially, the condition describes
a compatibility with the Stokes structure.

\subsubsection{Good $\rnum$-structure
of meromorphic flat bundles}

Let $V$ be a meromorphic flat bundle
on $(X,D)$,
and let $\nbigl$ denote the associated
local system on $\Xtilde(D)$
with the Stokes structure.
A $\rnum$-structure of $V$
is a $\rnum$-structure
of the associated local system on $X\setminus D$,
which is equivalent to a $\rnum$-structure of $\nbigl$.
It is called a good $\rnum$-structure of $V$
if the Stokes filtrations $\nbigf^P$
$(P\in\pi^{-1}(D))$
are defined over $\rnum$,
with respect to the induced $\rnum$-structure
of $\nbigl$.
By the isomorphisms (\ref{eq;09.11.9.1}),
we obtain the pre-Betti structures of
$V$ and $V_!$.
Moreover, it is easy to observe that
$\psi_z(V)$ and $\Xi_z(V)$ are also naturally
equipped with pre-Betti structures
such that the morphisms
$a_i$ and $b_i$ $(i=1,2)$ 
are compatible with pre-Betti structures.

\subsubsection{Betti structure of holonomic
$\nbigd$-modules on curves}

Let $\nbigm$ be a holonomic $\nbigd$-module
on $(X,D)$
such that $V:=\nbigm(\ast D)$
is a meromorphic flat bundle.
Let $(\nbigf,\alpha)$ be a pre-Betti structure of
$\nbigm$.
We call it a Betti structure
if the following holds:
\begin{itemize}
\item
The induced $\rnum$-structure
on $\DR(V_{|X-D})$ induces
a good $\rnum$-structure of $V$.
As remarked above,
we have the induced pre-Betti structures
on $V$ and $V_!$.
\item
The natural morphisms
$a_0$ and $b_0$ are compatible with
the pre-Betti structures.
\end{itemize}
Note that we obtain a pre-Betti structure on $\phi_z(\nbigm)$
from the expression as the cohomology
of the complex (\ref{eq;09.11.9.3}),
and the morphisms
$\var$ and $\can$ are compatible
with the pre-Betti structures.
The pre-Betti structure of $\nbigm$
can be reconstructed from
the pre-Betti structure of $\phi_z(\nbigm)$
and the good $\rnum$-structure of $V$.

\section{Betti structure in the higher dimensional case}

We would like to generalize 
the notion of Betti structure
in the higher dimensional case.

\subsection{Good meromorphic flat bundle
and good $\rnum$-structure}

Let $X$ be any complex manifold
with a simple normal crossing hypersurface $D$.
It is fundamental to understand the structure of
good meromorphic flat bundles on $(X,D)$,
which is now well established
after the work of H. Majima, C. Sabbah
and the author.
(See \cite{majima},
\cite{mochi7}, \cite{mochi8},
\cite{sabbah4} and
\cite{sabbah_lecture_Stokes}.
See \cite{mochi10} for a survey.)
Very briefly,
the asymptotic analysis for meromorphic flat bundles
on curves can be naturally generalized
for good meromorphic flat bundles
in the higher dimensional case,
and we obtain 
the Riemann-Hilbert-Birkhoff correspondence,
which is a natural correspondence
between good meromorphic flat bundles
and local systems with Stokes structure.

Let us recall it very briefly.
Let $(V,\nabla)$ be a good meromorphic flat bundle.
Let $\pi:\Xtilde(D)\lrarr X$ be the real blow up
along $D$,
which means in this paper
the fiber product of the real blow up
along the irreducible components of $D$
taken over $X$.
Let $\nbigl$ be the local system on
$\Xtilde(D)$ associated to $V_{|X-D}$.
For any point $P\in \pi^{-1}(D)$,
we have the Stokes filtration
$\nbigf^P$ of the stalk $\nbigl_{P}$.
It satisfies a compatibility condition
with the Stokes filtrations
$\nbigf^Q$ for $Q$ which are close to $P$.
We can reconstruct $V$
from $V_{|X-D}$
and the system of filtrations
$\bigl\{
 \nbigf^P\,\big|\,P\in \pi^{-1}(D)
 \bigr\}$.
Moreover, if we are given 
a local system with the family of Stokes filtrations
$\{\nbigf^P\,|\,P\in\pi^{-1}(D)\}$
satisfying the compatibility condition,
we have the corresponding 
good meromorphic flat bundle $V$.
This is the Riemann-Hilbert-Birkhoff correspondence
for good meromorphic flat bundles.

As in the one dimensional case,
the de Rham complexes of $V$ and $V_!$
are described in terms of
the local system $\nbigl$ with the Stokes structure.
We obtain the constructible subsheaf
$\nbigl^{\leq D}$ of $\nbigl$ 
which consists of flat sections
with the moderate growth.
It is described as 
$\nbigl^{\leq D}_{P}=
 \nbigf^P_{\leq 0}(\nbigl_P)$
$(P\in\pi^{-1}(D))$
in terms of the Stokes filtrations.
Let $\nbigl^{<D}$ be the 
constructible subsheaf of $\nbigl$,
which consists of flat sections 
with rapid decay along $D$.
It is also described in terms of the Stokes filtration
(see \S\ref{subsection;09.10.28.1}).
Then,
we have
$\DR_X(V)\simeq  R\pi_{\ast}\nbigl^{\leq D}[\dim X]$
and 
$\DR_X(V_{!})\simeq  R\pi_{\ast}\nbigl^{< D}[\dim X]$
as in (\ref{eq;09.11.9.1}).

For any holomorphic function $g$ on $X$
such that $g^{-1}(0)=D$,
we obtain $\nbigd_X$-modules
$\psi_g(V)$ and $\Xi_g(V)$
with morphisms as in (\ref{eq;09.11.10.1})
by using the formalism of Beilinson.
Their de Rham complexes are also 
described in terms of the local system
$\nbigl$ with the Stokes filtrations.

As in the one dimensional case,
a $\rnum$-structure of $V$
is a $\rnum$-structure of the associated local system
on $X\setminus D$,
which is equivalent to a $\rnum$-structure
of $\nbigl$.
It is called a good $\rnum$-structure of $V$
if the Stokes filtrations are defined over $\rnum$.
If $V$ is equipped with a good $\rnum$-structure,
the $\nbigd_X$-modules
$V$, $V_!$, $\Xi_g(V)$ and $\psi_g(V)$
are naturally equipped with pre-Betti structures,
and the natural morphisms 
as in (\ref{eq;09.11.10.1})
are compatible with the pre-Betti structures.

\subsection{Good $\rnum$-structure
of meromorphic flat connections}
\label{subsection;13.4.27.200}

In the higher dimensional case,
not all meromorphic flat bundles are good,
which is one of the main difficulties.
Let us recall local resolutions of turning points
due to K. Kedlaya 
\cite{kedlaya}, \cite{kedlaya2}.
(See \cite{sabbah4} for the original conjecture.
See also \cite{mochi6} and \cite{mochi7}
for the algebraic case.)

Let $X$ be a complex manifold
with a hypersurface $D$.
Let $V$ be a reflexive $\nbigo_X(\ast D)$-module
with a flat connection,
which is called a meromorphic flat connection
\cite{Malgrange-Lattice}.
For any $P\in X$,
there exist a neighbourhood $X_P$ of $P$ in $X$
and a projective birational morphism
$\lambda_P:\check{X}_P\lrarr X_P$
such that 
(i) $\check{X}_P$ is smooth and 
 $\check{D}_P:=\lambda_P^{-1}(D)$ is normal crossing,
(ii) $\check{X}_P\setminus \check{D}_P
\simeq
 X_P\setminus D$,
(iii) $\check{V}_P:=
 \lambda_P^{\ast}V$
is a good meromorphic flat bundle
on $(\check{X}_P,\check{D}_P)$.
(See Theorem 8.2.2 of \cite{kedlaya2}.)
Such $(X_P,\lambda_P)$
is called a local resolution of $V$
in this paper.
If $X$ and $V$ are algebraic,
we have a global resolution.
(See Theorem 8.1.3 of \cite{kedlaya2}
or Theorem 16.2.1 of \cite{mochi7}.)

Then, the notion of good $\rnum$-structure
is generalized for meromorphic flat connections
which are not necessarily good.
Namely,
a $\rnum$-structure of $V$ is called good
if the induced $\rnum$-structure
of good meromorphic flat bundles $\check{V}_P$
are good for any local resolutions
$(X_P,\lambda_P)$.
Even in this case,
the de Rham complexes
$\DR_X(V)$ and $\DR_X(V_!)$
have naturally induced $\rnum$-structures.
Moreover, 
if we are given a holomorphic function $g$ 
on $X$ such that $g^{-1}(0)=D$,
the holonomic $\nbigd_X$-modules
$\psi_g(V)$ and $\Xi_g(V)$
are naturally equipped with
pre-Betti structures,
with which the morphisms in (\ref{eq;09.11.10.1})
are compatible.

\subsection{Cells and gluing}

Let us recall that
any holonomic $\nbigd$-module $\nbigm$
can be described as the gluing of
a ``cell'' and a holonomic $\nbigd$-module $\nbigm'$
whose support $\Supp\nbigm'$
is strictly smaller than $\Supp\nbigm$.
Namely,
for any $P\in\Supp\nbigm$,
there exists a tuple 
$\nbigc=(Z,U,\varphi,V)$ as follows:
\begin{description}
\item[(Cell 1)]
 $\varphi:Z\lrarr X$ is a morphism of
 complex manifolds
 such that $P\in\varphi(Z)$ and 
that $\dim Z$ is equal to the dimension of
 $\Supp\nbigm$ at $P$.
 We impose that 
 there exists a neighbourhood $X_P$
 of $P$ in $X$ such that
 $\varphi:Z\lrarr X_P$ is projective.
\item[(Cell 2)]
 $U\subset Z$ is the complement 
 of a hypersurface $D_Z$.
 We impose that 
 the restriction $\varphi_{|U}$ is an immersion,
 and that there exists a hypersurface $H$ of $X_P$
 such that $\varphi^{-1}(H)=D_Z$.
\item[(Cell 3)]
 $V$ is a good meromorphic flat bundle
 on $(Z,D_Z)$.
 We impose 
 $\nbigm(\ast H)=\varphi_{\dagger}V$
 for a hypersurface $H$ as in (Cell 2).
 Note that we obtain the natural morphisms
 $\varphi_{\dagger}V_!\lrarr
 \nbigm\lrarr\varphi_{\dagger}V$.
\end{description}
Such $\nbigc$ is called a cell of $\nbigm$
at $P$.
A holomorphic function $g$ on $X$
is called a cell function for $\nbigc$
if $\varphi(U)=\Supp\nbigm\setminus g^{-1}(0)$.
We set $g_Z:=g\circ\varphi$.
We have natural isomorphisms
$\varphi_{\dagger}\Xi_{g_Z}(V)
\simeq
 \Xi_g\varphi_{\dagger}(V)$
and
$\varphi_{\dagger}\psi_{g_Z}(V)
\simeq
 \psi_g\varphi_{\dagger}(V)$.
By the formalism of Beilinson,
the $\nbigd_X$-module $\phi_g(\nbigm)$
is obtained as the cohomology
of the complex:
\begin{equation}
\label{eq;09.11.10.2}
 \varphi_{\dagger}V_!
\lrarr
 \Xi_g\varphi_{\dagger}(V)
\oplus
 \nbigm
\lrarr
 \varphi_{\dagger}V
\end{equation}
We have the description of $\nbigm$ around $P$
as the cohomology of the complex:
\[
\psi_g(\varphi_{\dagger}V)
\lrarr
 \Xi_g(\varphi_{\dagger}V)
\oplus 
 \phi_g(\nbigm)
\lrarr
 \psi_g(\varphi_{\dagger}V).
\]
In other words,
$\nbigm$ is described as 
the gluing of the cell $\nbigc$
and $\phi_g(\nbigm)$.

\subsection{Betti structure}

\subsubsection{Compatibility of cell and pre-Betti structure}

We introduce the compatibility condition
of a cell $\nbigc$
and a pre-Betti structure $\nbigf$
of $\nbigm$.
We say that $\nbigf$ and $\nbigc$
are compatible
if the following holds:
\begin{itemize}
\item
Note that the flat bundle $V_{|U}$ has 
an induced $\rnum$-structure.
We suppose that it is a good $\rnum$-structure
in the sense of \S\ref{subsection;13.4.27.200}.
\item
By the first condition,
$\varphi_{\dagger}V$,
$\varphi_{\dagger}V_!$,
$\Xi_g\varphi_{\dagger}V$
and $\psi_g\varphi_{\dagger}V$
are equipped with the induced
pre-Betti structures.
Then, we impose that
the morphisms
$\varphi_{\dagger}V_!\lrarr\nbigm
\lrarr\varphi_{\dagger}V$
are compatible with pre-Betti structures.
\end{itemize}
Such a cell $\nbigc$ is called
a $\rnum$-cell of $\nbigm$ at $P$.
Since $\phi_g(\nbigm)$ is the cohomology
of the complex (\ref{eq;09.11.10.2}),
it is equipped with the induced pre-Betti structure.

\subsubsection{Inductive definition of Betti structure}

Let us define the notion of Betti structure
of $\nbigm$ at $P$,
inductively on the dimension of $\Supp\nbigm$.
If $\dim_P\Supp\nbigm=0$,
a Betti structure is defined to be
a pre-Betti structure.
Let us consider the case $\dim_P\Supp\nbigm\leq n$.
We say that a pre-Betti structure of $\nbigm$
is a Betti structure at $P$
if there exists an $n$-dimensional $\rnum$-cell
$\nbigc=(Z,\varphi,U,V)$ at $P$
with the following property:
\begin{itemize}
\item
 $\dim_P\Bigl(
 \bigl(\Supp\nbigm\cap X_P\bigr)
 \setminus
 \varphi(Z)\Bigr)
 <n$
 for some neighbourhood $X_P$ of $P$ in $X$.
\item
For a cell function $g$ for $\nbigc$, 
the induced pre-Betti structure of $\phi_g(\nbigm)$
is a Betti structure at $P$.
Note that
$\dim\Supp\phi_g(\nbigm)<n$
by the first condition.
\end{itemize}
A holonomic $\nbigd$-module with
Betti structure is called
a $\rnum$-holonomic $\nbigd$-module.
Morphisms of $\rnum$-holonomic $\nbigd_X$-modules
are defined to be morphisms of 
pre-$\rnum$-holonomic $\nbigd_X$-modules.

\begin{rem}
The above is not exactly the same
as the definition in {\S\ref{subsection;09.11.11.2}},
but they give equivalent objects.
\hfill\qed
\end{rem}

\section{Main goal}

\subsection{The category of $\rnum$-holonomic
$\nbigd$-modules}

Besides giving the details on the above arguments,
it is our main purpose to show that
our notion of Betti structure is nice.
The category of $\rnum$-holonomic $\nbigd$-modules
should contain the holonomic $\nbigd$-modules
naturally induced from any meromorphic flat connections
with a good $\rnum$-structure,
for which we have the following theorem.

\begin{thm}
\label{thm;13.4.27.301}
Let $X$ be any complex manifold
with a hypersurface $D$.
Let $V$ be any meromorphic flat connection on $(X,D)$
with a good $\rnum$-structure.
Then, the natural pre-Betti structures
of $V$ and $V_!$ are Betti structures.
\end{thm}

See Theorem \ref{thm;13.4.23.20}
for a refined result.
Some of the functors for holonomic $\nbigd$-modules
should be enriched with Betti structures,
as in the following theorems.

\begin{thm}[Theorem \ref{thm;09.10.16.5}]
Let $F:X\lrarr Y$ be any projective morphism
of complex manifolds.
For any $\rnum$-holonomic $\nbigd_X$-module $\nbigm$,
the push-forward $F_{\dagger}^i\nbigm$
are also naturally $\rnum$-holonomic for any $i$.
\end{thm}

\begin{thm}[Theorem \ref{thm;09.10.18.400}]
Let $X$ be any complex manifold with a hypersurface $D$.
Let $\nbigm$ be any $\rnum$-holonomic $\nbigd_X$-module.
Then, $\nbigm\otimes\nbigo_X(\ast D)$ has 
a unique Betti structure, 
for which
$\nbigm\lrarr\nbigm\otimes\nbigo_X(\ast D)$
is compatible with the Betti structures.
\end{thm}

\begin{thm}[Proposition \ref{prop;09.12.6.1}]
Let $X$ be any complex manifold
with a hypersurface $D$.
Let $\nbigm$ be any $\rnum$-holonomic
$\nbigd_X$-module.
Let $V$ be any meromorphic connection on $(X,D)$
with a good $\rnum$-structure.
Then, $\nbigm\otimes V$ is naturally 
a $\rnum$-holonomic $\nbigd_X$-module.
\end{thm}

The following is an easier result.
\begin{thm}
\mbox{{}}
\begin{itemize}
\item
The category of $\rnum$-holonomic
$\nbigd_X$-modules is abelian.
\item
The dual of $\rnum$-holonomic $\nbigd_X$-modules
are naturally $\rnum$-holonomic.
\item
Let $\nbigm$ be a $\rnum$-holonomic $\nbigd_X$-module.
Let $\nbigm'\subset\nbigm$ be 
a subobject in the category of
pre-$\rnum$-holonomic $\nbigd_X$-modules.
Then, $\nbigm'$ is also $\rnum$-holonomic.
We have a similar claim for quotients.
\end{itemize}
\end{thm}

By using the theorems,
we obtain that 
the category of $\rnum$-holonomic 
$\nbigd$-modules contains
expected objects.
For example,
it contains the holonomic $\nbigd$-modules
obtained from the structure sheaf of any algebraic variety
by successive use of
the pull back and the push-forward by algebraic morphisms,
and the exponential twist by algebraic functions.
(This type of holonomic $\nbigd$-modules are 
closely related with extended exponential-motivic
$\nbigd$-modules in \cite{Kontsevich-JJM}.)
It implies the compatibility of 
the $\rnum$-structure and the Stokes structure
for some naturally obtained
meromorphic flat bundles.
Such phenomena are expected
in the non-commutative Hodge theory
\cite{Katzarkov-Kontsevich-Pantev}.

\vspace{.1in}

In the algebraic case,
the derived category of $\rnum$-holonomic
$\nbigd$-modules is equipped with standard functoriality,
so called $6$-operations.
\begin{thm}
The category of $\rnum$-holonomic 
algebraic $\nbigd$-modules
is equipped with the standard functors
such as dual, push-forward, pull-back,
tensor product,
and inner homomorphism,
compatible with those for the category of
holonomic algebraic $\nbigd$-modules
with respect to the forgetful functor.
\end{thm}

\subsection{Analysis on real blow up}

We also give some analysis on the real blow up,
which is a complement to \cite{sabbah_lecture_Stokes}.
Very briefly,
we can capture the Stokes structure
by considering the de Rham complex
on the real blow up,
at least in the case of good meromorphic flat bundles.
We have several useful classes of functions on the real blow up,
the moderate growth,
the rapid decay,
and the Nilsson type.
We study or review the fundamental property
of the sheaves of such functions
and the corresponding de Rham complexes.
We will not restrict ourselves
to our main purpose, i.e., the study on Betti structure.
For example,
we shall prove that the sheaf of holomorphic functions
of moderate growth is flat over 
the sheaf of holomorphic functions on the underlying space
(Theorem \ref{thm;14.1.22.1}).
Although we will not use it in this paper,
it is quite basic,
and the author expects that it would be useful
for a further study.

\begin{rem}
G. Morando informed the author 
that the theory of ind-sheaves 
{\rm\cite{Kashiwara-Schapira-Ind-sheaves}}
provides us with a powerful method
to study analysis on the real blow up.
(See also the recent work by A. D'Agnolo 
and M. Kashiwara \cite{DAgnolo-Kashiwara}.)
While the author hopes that it would make
the subject more transparent,
he also hopes that
his direct way would also be significant 
for our understanding at this moment.
\hfill\qed
\end{rem}

\section{Acknowledgement}

I am grateful to C. Sabbah for discussions
and for his kindness.
He sent an earlier version of
his monograph \cite{sabbah_lecture_Stokes},
which invited me to this study.
He also pointed an error in an earlier version
of this paper.
I thank H. Esnault,
who attracted my attention
to Betti structure of holonomic $\nbigd$-modules.
This study much owes to 
the foundational works 
on holonomic $\nbigd$-modules
and perverse sheaves
due to many people,
among all
A. Beilinson, 
J. Bernstein,
P. Deligne,
M. Kashiwara,
Z. Mebkhout,
and M. Saito.
I am grateful to G. Morando for useful discussions.
Special thanks goes to K. Vilonen.
I thank P. Schapira for a valuable discussion.
I am grateful to A. Ishii and Y. Tsuchimoto
for their constant encouragement.
I am grateful to 
M. Hien,
M-H. Saito, C. Simpson and Sz. Szabo
for stimulating discussions.
The author thanks the referee for his thorough reading
and a lot of of useful and valuable suggestions 
to improve this monograph.

This work was supported by
the Grant-in-Aid for Scientific Research (C)
(No. 22540078),
Japan Society for the Promotion of Science.

%% file: 2.1.tex
\subsection{Dual, push-forward and de Rham functor}

We prepare some notation.
See very useful text books \cite{hotta-tanisaki}
and \cite{kashiwara_text} for more details and precisions
on $\nbigd$-modules.
Let $X$ be a complex manifold
with $\dim X=d_X$.
Let $\nbigd_X$ denote the sheaf of 
holomorphic differential operators
on $X$.
In this paper, $\nbigd_X$-module
means left $\nbigd_X$-module.
Let $\Hol(X)$ be the category of
holonomic $\nbigd_X$-modules,
and let $D^b_{\hol}(\nbigd_X)$ be the derived
category of cohomologically bounded holonomic
$\nbigd_X$-complexes.
\index{category $\Hol(X)$}
\index{category $D^b_{\hol}(\nbigd_X)$}
Let $\Omega_X^{j}$ denote 
the sheaf of holomorphic $j$-forms.
\index{sheaf $\Omega_X^{j}$}
The invertible sheaf $\Omega_X^{d_X}$
is denoted by $\Omega_X$.
\index{sheaf $\Omega_X$}
The sheaves of $C^{\infty}$-$(p,q)$-forms
are denoted by $\Omega_X^{p,q}$.
\index{sheaf $\Omega_X^{p,q}$}
The dual functor on the derived category
of $\nbigd_X$-modules
is denoted by $\DDD_X$,
i.e.,
$\DDD_X\nbigm^{\bullet}:=
 \nrhom_{\nbigd_X}\bigl(
 \nbigm^{\bullet},
 \nbigd_X\otimes\Omega_X^{\otimes\,-1}
 \bigr)[d_X]$.
\index{dual functor $\DDD_X$}
Recall that
if $\nbigm$ is a holonomic $\nbigd_X$-module,
then $\DDD_X\nbigm$ is a holonomic
$\nbigd_X$-module.
For $\nbigd_X$-modules
$\nbigm_i$ $(i=1,2)$,
the tensor product 
$\nbigm_1\otimes_{\nbigo_X}
 \nbigm_2$
is naturally a $\nbigd_X$-module.
For any tangent vector field $v$,
we have $v(m_1\otimes m_2)
=(vm_1)\otimes m_2+m_1\otimes (vm_2)$.
The $\nbigd_X$-module is
denoted by $\nbigm_1\otimes^D\nbigm_2$.
It is also denoted by
$\nbigm_1\otimes\nbigm_2$
if there is no risk of confusion.
\index{tensor product $\otimes^D$, $\otimes$}

\begin{lem}
\label{lem;09.12.1.1}
Let $\nbigm$ be any holonomic $\nbigd_X$-module.
Let $V$ be any $\nbigd_X$-module,
which is coherent and locally free as an $\nbigo_X$-module.
Its dual is denoted by $V^{\lor}$.
Then, we have a natural isomorphism
\[
 \DDD_X\bigl(\nbigm\otimes^D V\bigr)
\simeq
 (\DDD_X\nbigm)\otimes^D V^{\lor}.
\]
\end{lem}
\pf
We recall Remark 3.4 in \cite{kashiwara_text}.
For any left $\nbigd_X$-module $\nbign$,
we have the left $\nbigd_X$-action
on $\nbigd_X\otimes^D\nbign$.
It is also equipped with
a right $\nbigd_X$-action 
given by the multiplication
$(f\otimes m)\cdot g=fg\otimes m$
for $g\in \nbigd_X$.
The two-sided $(\nbigd_X,\nbigd_X)$-module
is denoted by $\nbign_1$.
Similarly,
we have a left action of $\nbigd_X$
on $\nbigd_X\otimes_{\nbigo_X}\nbign$
(the tensor product $\otimes_{\nbigo_X}$
is taken for the $\nbigo_X$-module structure of
 $\nbigd_X$ given by the right multiplication)
given by the multiplication
$g\cdot (f\otimes m)=gf\otimes m$
for $g\in\nbigd_X$,
and a right $\nbigd_X$-action
given by
$(f\otimes m)\cdot v
=fv\otimes m-f\otimes vm$
for a tangent vector $v$.
The two-sided $(\nbigd_X,\nbigd_X)$-module
is denoted by $\nbign_2$.
We have a naturally defined
$\nbigo_X$-morphism
$\nbign\lrarr\nbign_1$
given by $m\longmapsto 1\otimes m$.
It is naturally extended to a morphism 
of left $\nbigd_X$-modules
$\nbign_2\lrarr\nbign_1$.
Actually, it is an isomorphism
and compatible with the right $\nbigd_X$-action,
as remarked in \cite{kashiwara_text}.

We have two left $\nbigd_X$-actions on
$\nbigd_X\otimes\Omega_X^{\otimes\,-1}$.
The first one is the natural one,
and the second one is induced by 
the right $\nbigd_X$-action.
They induce two $\nbigo_X$-actions.
Let $(\nbigd_X\otimes\Omega_X^{\otimes\,-1})
\otimes_{\nbigo_X}^i\nbign$
denote the tensor product
with respect to the $i$-th one.
Each is equipped with
two left $\nbigd_X$-actions.
From the consideration in the previous paragraph,
we obtain a natural isomorphism
$\iota:
 \nbign\otimes_{\nbigo_X}^1
 (\nbigd_X\otimes\Omega_X^{\otimes\,-1})
\lrarr
 \nbign\otimes_{\nbigo_X}^2
 (\nbigd_X\otimes\Omega_X^{\otimes\,-1})$,
compatible with the $\nbigd_X$-actions.

Let us return to Lemma \ref{lem;09.12.1.1}.
We have the following natural isomorphisms
of $\nbigd_X$-modules:
\begin{multline}
 \DDD_X(\nbigm\otimes^D V)
=
 \nrhom_{\nbigd_X}\bigl(
 \nbigm\otimes^D V,\,
 \nbigd_X
\otimes\Omega_X^{\otimes\,-1}
 \bigr) \\
\simeq
 \nrhom_{\nbigd_X}\Bigl(
 \nbigm,\,V^{\lor}\otimes_{\nbigo_X}^1
 \bigl(
 \nbigd_X
 \otimes\Omega_X^{\otimes\,-1}
 \bigr)
 \Bigr) \\
\simeq
 \nrhom_{\nbigd_X}\Bigl(
 \nbigm,\,V^{\lor}\otimes_{\nbigo_X}^2
 \bigl(
 \nbigd_X\otimes\Omega_X^{\otimes\,-1}
 \bigr)
 \Bigr)
=(\DDD_X\nbigm)\otimes^D V^{\lor}
\end{multline}
Here, the first one is obtained
by using Godement type injective resolution,
and the second one is induced by $\iota$ above.
\hfill\qed

\vspace{.1in}

For any field $R$,
let $R_X$ denote the sheaf on $X$
associated to the constant presheaf valued in $R$.
Let $D^b(R_X)$ (resp. $D^b_c(R_X)$) denote
the derived category of cohomologically 
bounded (resp. bounded constructible) $R_X$-complexes,
and let $\Per(X,R)$ denote the category of
$R$-perverse sheaves.
\index{category $\Per(X,R)$}
\index{category $D^b(R_X)$}
Let $\omega_{X,R}$ denote the dualizing
complex of $R_X$-modules.
\index{complex $\omega_{X,R}$}
It will be denoted by $\omega_{X}$
if there is no risk of confusion.
The dual functor on the derived category
of $R_X$-modules is also denoted by $\DDD_X$,
i.e.,
for an $R_X$-complex $\nbigf^{\bullet}$,
let $\DDD_X\nbigf^{\bullet}:=
 \nrhom_{R_X}\bigl(
 \nbigf^{\bullet},\omega_{X,R}\bigr)$.
\index{dual functor $\DDD_X$}

The de Rham functor is denoted by $\DR_X$,
i.e.,
$\DR_X\nbigm:=
 \Omega_X\otimes_{\nbigd_X}^L\nbigm
=\Omega_X^{\bullet}\otimes_{\nbigo_X}
 \nbigm[d_X]$.
\index{de Rham functor $\DR_X$}
According to \cite{kashiwara_perversity},
it gives a functor of triangulated categories
\[
 \DR_X:D^b_{\hol}(\nbigd_X)
\lrarr D^b_{c}(\cnum_X)
\]
compatible with the $t$-structures,
where the $t$-structure of 
$D^b_{\hol}(\nbigd_X)$ is the natural one,
and the $t$-structure of
$D^b_{c}(\cnum_X)$ is given by
the middle perversity.
In particular, it induces 
an exact functor 
$\DR_X:\Hol(X)\lrarr \Per(X,\cnum)$.
We can identify
$\omega_X=\DR_X\nbigo_X[d_X]$.
It is easy to observe that
$\DR_X\nbigm=0$ implies $\nbigm=0$
for $\nbigm\in\Hol(X)$.
The functor
$\DR_X:\Hol(X)\lrarr\Per(X,\cnum)$
is faithful,
although it is not full in general.

Let $F:X\lrarr Y$ be
a morphism of complex manifolds.
The push-forward for $\cnum_X$-complexes
in the derived category is denoted by $RF_{\ast}$.
(It is also denoted by $F_{\ast}$
if there is
no risk of confusion.)
Its $i$-th perverse cohomology is 
denoted by $F_{\dagger}^i$.
\index{push-forward $F^i_{\dagger}$}
Put 
\[
 \nbigd_{X\rarr Y}:=
 \nbigo_X\otimes_{F^{-1}\nbigo_Y}F^{-1}\nbigd_Y,
\quad\quad
 \nbigd_{Y\larr X}:=
\Omega_X\otimes_{F^{-1}\nbigo_Y} F^{-1}\bigl(
 \nbigd_Y\otimes_{\nbigo_Y}\Omega_Y^{\otimes\,-1}
 \bigr).
\]
The push-forward for $\nbigd_X$-complexes
is denoted by $F_{\dagger}$,
i.e.,
$F_{\dagger}\nbigm
=RF_{\ast}
 \bigl(
 \nbigd_{Y\larr X}\otimes_{\nbigd_X}^{L}
 \nbigm
 \bigr)$.
Its $i$-th cohomology is denoted by
$F_{\dagger}^i$.

Recall that these functors are compatible
on the derived categories.
Let $F:X\lrarr Y$ be a proper morphism
of complex manifolds.
We have natural transformations
\[
 \DR_Y\circ F_{\dagger}
\simeq
 RF_{\ast}\circ \DR_X,
 \quad
\DDD_X\circ \DR_X\simeq
 \DR_X\circ\DDD_X,
\quad
\DDD_Y\circ F_{\dagger}
\simeq
 F_{\dagger}\circ \DDD_X.
\]
In \cite{saito4},
the following diagram is constructed
and it is proved to be commutative
(see Theorem 3.3 of \cite{saito4}):
\begin{equation}
 \label{eq;09.10.3.20}
\begin{CD}
RF_{\ast}\DDD_X\DR_X
@>{\simeq}>>
RF_{\ast}\DR_X\DDD_X
@>{\simeq}>>
\DR_Y F_{\dagger}\DDD_X
 \\
@V{\simeq}VV @. @V{\simeq}VV \\
\DDD_YRF_{\ast}\DR_X
 @>{\simeq}>>
\DDD_Y\DR_Y F_{\dagger}
 @>{\simeq}>>
\DR_Y \DDD_Y F_{\dagger}
\end{CD}
\end{equation}

\subsection{Hypersurfaces}

For any hypersurface $D\subset X$,
let $\nbigo_X(\ast D)$ denote the sheaf of
meromorphic functions whose poles
are contained in $D$.
\index{sheaf $\nbigo_X(\ast D)$}
For $\nbigm\in\Hol(X)$,
we have 
$\nbigm(\ast D),\nbigm(!D)\in\Hol(X)$
given as follows:
\[
\nbigm(\ast D):=
 \nbigm\otimes_{\nbigo_X}\nbigo_X(\ast D),
\quad
\nbigm(!D):=\DDD_X\Bigl(
 \bigl(\DDD_X\nbigm\bigr)(\ast D)
 \Bigr).
\]
\index{sheaf $\nbigm(\ast D)$}
\index{sheaf $\nbigm(\bikkuri D)$}
We have naturally defined morphism
$\nbigm\lrarr\nbigm(\ast D)$.
The morphism
$\DDD_X(\nbigm)\lrarr \DDD_X(\nbigm)(\ast D)$
and the natural transformation $\DDD_X\circ\DDD_X\simeq\id_X$
induce 
$\nbigm(!D)\lrarr\nbigm$.
(See \S3.3 and \S A3.3 of \cite{kashiwara_text}
for $\DDD_X\circ\DDD_X\simeq\id$.)
They are uniquely characterized that 
the restrictions to $X\setminus D$
are the identities.
If $D$ is given as the zero set of a holomorphic
function $f$,
they are denoted by $\nbigm(\ast f)$
and $\nbigm(!f)$, respectively.
\index{sheaf $\nbigm(\ast f)$}
\index{sheaf $\nbigm(\bikkuri f)$}
If we are given two hypersurfaces $D_i$ $(i=1,2)$,
we set
$\nbigm(\star_1D_1)(\star_2D_2):=
 \bigl(\nbigm(\star_1D_1)\bigr)(\star_2D_2)$,
where $\star_i\in\{\ast,!\}$.

\vspace{.1in}

We put 
$\nbigd_{X(\ast D)}:=
 \nbigd_X\otimes\nbigo_X(\ast D)$.
\index{sheaf $\nbigd_{X(\ast D)}$}
A $\nbigd_{X(\ast D)}$-module $\nbigm$
is called holonomic,
if it is holonomic as a $\nbigd_X$-module.
Let $\Hol\bigl(X,\ast D\bigr)$ be the category of
holonomic $\nbigd_{X(\ast D)}$-modules,
which is naturally a full subcategory of $\Hol(X)$.
\index{category $\Hol(X,\ast D)$}
The dual functor on $\Hol\bigl(X,\ast D\bigr)$
is denoted by $\DDD_{X(\ast D)}$,
i.e.,
$\DDD_{X(\ast D)}(\nbigm)
=\DDD_X(\nbigm)(\ast D)$.
\index{dual functor $\DDD_{X(\ast D)}$}

Let $j:X\setminus D\lrarr X$ be the inclusion.
We define a functor
$j^{\ast}:\Hol(X)\lrarr \Hol(X,\ast D)$
by $j^{\ast}(\nbigm)=\nbigm(\ast D)$.
\index{functor $j^{\ast}$}
The natural inclusion $\Hol(X,\ast D)\lrarr \Hol(X)$
is denoted by $j_{\ast}$.
\index{functor $j_{\ast}$}
Another functor $j_!:\Hol(X,\ast D)\lrarr \Hol(X)$
is defined by $j_!(\nbigm):=(j_{\ast}\nbigm)(!D)$.
\index{functor $j_{\bikkuri}$}
The functors $j^{\ast}$, $j_{\ast}$ and $j_!$ are exact.
In this notation,
we have 
$\nbigm(\ast D)=j_{\ast}j^{\ast}\nbigm$
and $\nbigm(!D)=j_!j^{\ast}\nbigm$
for $\nbigm\in\Hol(X)$.

\vspace{.1in}
It is generalized as follows.
Let $H$ be a hypersurface of $X$.
Let $k:X\setminus H\lrarr X$ denote the inclusion.
For $\nbigm\in\Hol(X,\ast D)$,
we define $k^{\ast}\nbigm:=\nbigm(\ast H)$.
We can naturally regard
$\Hol\bigl(X,\ast (D\cup H)\bigr)$
as a full subcategory of
$\Hol(X,\ast D)$.
The natural inclusion is denoted by $k_{\ast}$.
We define another functor
$k_!:\Hol\Bigl(X,\ast (D\cup H)\Bigr)\lrarr\Hol(X,\ast D)$
by $k_!\nbigm=
 j^{\ast}\Bigl(
 \bigl(
 (j\circ k)_{\ast}\nbigm
 \bigr)\bigl(!(D\cup H)\bigr)
 \Bigr)$.

Later (\S\ref{subsection;13.4.25.200}),
we shall consider a successive composition
of the operations.

\subsection{Pre-$K$-holonomic $\nbigd$-modules}

Let $\nbigm$ be any holonomic $\nbigd_X$-module.
Let $K$ be any subfield of $\cnum$.
A pre-$K$-Betti structure of $\nbigm$
is defined to be a $K$-perverse sheaf
$\nbigf$ with an isomorphism
$\lambda:\nbigf\otimes_K\cnum\simeq \DR_X\nbigm$.
\index{pre-$K$-Betti structure}
Such a tuple
$(\nbigm,\nbigf,\lambda)$ is called 
a pre-$K$-holonomic $\nbigd_X$-module.
\index{pre-$K$-holonomic $\nbigd_X$-module}
We will often omit to denote $\lambda$.
A morphism of $K$-holonomic $\nbigd_X$-modules
$(\nbigm_1,\nbigf_1)\lrarr(\nbigm_2,\nbigf_2)$
is defined to be a pair of 
a morphism of $\nbigd_X$-modules
$\nbigm_1\lrarr \nbigm_2$
and a morphism of perverse sheaves
$\nbigf_1\lrarr\nbigf_2$ 
such that the following induced diagram is
commutative:
\[
 \begin{CD}
 \nbigf_1\otimes_K\cnum 
 @>{\simeq}>> \DR_X(\nbigm_1) \\
 @VVV @VVV \\
 \nbigf_2\otimes_K\cnum 
 @>{\simeq}>>\DR_X(\nbigm_2)
 \end{CD}
\]
The category of pre-$K$-holonomic $\nbigd_X$-modules
is denoted by $\Hol^{\pre}(X,K)$.
\index{category $\Hol^{\pre}(X,K)$}

The following lemma is clear.
\begin{lem}
$\Hol^{\pre}(X,K)$ is abelian.
\hfill\qed
\end{lem}

\vspace{.1in}
Let $\nbigf$ be a pre-$K$-Betti structure
of $\nbigm$.
We have induced pre-$K$-Betti structures
$\DDD\nbigf$ and $F_{\dagger}^i\nbigf$
of $\DDD\nbigm$ and $F_{\dagger}^i\nbigm$,
where $F:X\lrarr Y$ be a proper morphism.
We put
$\DDD(\nbigm,\nbigf):=
\bigl(\DDD\nbigm,\DDD\nbigf\bigr)$ and 
$F_{\dagger}^i(\nbigm,\nbigf):=
 \bigl(F_{\dagger}^i\nbigm,
 F_{\dagger}^i\nbigf\bigr)$.

\begin{lem}
The isomorphism
$\DDD F_{\dagger}\nbigm
\simeq
 F_{\dagger}\DDD\nbigm$
is compatible with the induced pre-$K$-Betti
structures.
\end{lem}
\pf
Because (\ref{eq;09.10.3.20}) is commutative,
we have the commutativity of the following
naturally induced diagram:
\[
 \begin{CD}
 \DR\DDD F_{\dagger}\nbigm
 @>{\simeq}>>
 \DDD F_{\dagger}\DR\nbigm
 @>{\simeq}>>
 \DDD F_{\dagger}\nbigf\otimes\cnum\\
 @V{\simeq}VV @V{\simeq}VV @V{\simeq}VV \\
 \DR F_{\dagger}\DDD\nbigm
 @>{\simeq}>>
 F_{\dagger}\DDD\DR\nbigm
 @>{\simeq}>>
 F_{\dagger}\DDD\nbigf\otimes\cnum
 \end{CD}
\]
It means the claim of the lemma.
\hfill\qed

\subsection{Formal completion}
\label{subsection;09.10.29.10}

Let $Y$ be a real analytic manifold.
Let $\nbigc^{\infty}_{Y}$
denote the sheaf of $C^{\infty}$-functions
on $Y$.
\index{sheaf $\nbigc^{\infty}_ZY$}
For any real analytic subset $Z$,
let $\nbigc^{\infty<Z}_Y$ denote the subsheaf
of $\nbigc^{\infty}_Y$ which consists of
the sections $f$ such that
the Taylor series of $f$ at each point $P\in Z$
is $0$.
\index{sheaf $\nbigc^{\infty<Z}_Y$}
We set $\nbigc^{\infty}_{\Zhat}:=
 \nbigc^{\infty}_Y/\nbigc_Y^{\infty<Z}$.
\index{sheaf $\nbigc^{\infty}_{\Zhat}$}
We have other descriptions;
(i)
It is the sheaf of 
Whitney functions of class $C^{\infty}$
on $Z$,
i.e., sections of $\infty$-jets
along $Z$ satisfying the conditions
in Theorem I.2.2 of \cite{malgrange2}.
(ii)
Let $\nbigi_{Z,\infty}$ be the ideal sheaf of 
$\nbigc^{\infty}_Y$ corresponding to $Z$.
Then, $\nbigc^{\infty}_{\Zhat}$
is also isomorphic to
$\varprojlim \nbigc^{\infty}_Y\big/\nbigi_{Z,\infty}^m$.
(See the proof of 
Theorem I.4.1 of \cite{malgrange2}.)
For any $\nbigc^{\infty}_Y$-module $\nbigf$,
let $\nbigf_{|\Zhat}$ denote
$\nbigf\otimes_{\nbigc_Y^{\infty}}
 \nbigc^{\infty}_{\Zhat}$.
Let $Z_i$ $(i=1,2)$ be real analytic subsets
in $Y$.
According to Corollary IV.4.4
with Definition I.5.4 of \cite{malgrange2},
the natural sequence
$0\lrarr \nbigc^{\infty}_{\widehat{Z_1\cup Z_2}}
\lrarr \nbigc^{\infty}_{\Zhat_1}\oplus
 \nbigc^{\infty}_{\Zhat_2}
\lrarr \nbigc^{\infty}_{\widehat{Z_1\cap Z_2}}\lrarr 0$
is exact.

Let $Z_i$ $(i\in \Lambda)$ be real analytic
subsets of $Y$.
For any subset $I\subset \Lambda$,
we put $Z_I:=\bigcap_{i\in I}Z_i$
and $Z(I):=\bigcup_{i\in I}Z_i$.
We fix a total order on $\Lambda$.
For $J\subset K\subset\Lambda$,
we have the restriction
$r_{J,K}:
 \nbigc^{\infty}_{\Zhat_J}
\lrarr
 \nbigc^{\infty}_{\Zhat_K}$.
If $K=J\sqcup\{i\}$,
we put $\kappa(J,K):=\{k\in J\,|\,k<i\}$
and $d_{J,K}:=(-1)^{\kappa(J,K)}r_{J,K}$.
We set 
 $\nbigk^m\bigl(
 \nbigc^{\infty}_{\Zhat(I)}
 \bigr):=
 \bigoplus_{|J|=m+1,\,J\subset I}
 \nbigc^{\infty}_{\Zhat_J}$.
\index{complex $\nbigk^{\bullet}$}
The above morphisms $d_{J,K}$ induce
$d_m:\nbigk^m\bigl(
 \nbigc^{\infty}_{\Zhat(I)}
 \bigr)
 \lrarr
 \nbigk^{m+1}\bigl(
 \nbigc^{\infty}_{\Zhat(I)}
 \bigr)$.
Thus, we obtain a complex 
$\nbigk^{\bullet}\bigl(
 \nbigc^{\infty}_{\Zhat(I)}
 \bigr)$.
By using the exactness 
in the previous paragraph,
it can be proved that
the natural inclusion
$\nbigc^{\infty}_{\Zhat(I)}\lrarr
 \nbigk^0(\nbigc^{\infty}_{\Zhat(I)})$
induces a quasi-isomorphism
$\nbigc^{\infty}_{\Zhat(I)}\simeq
 \nbigk^{\bullet}\bigl(
 \nbigc^{\infty}_{\Zhat(I)}
 \bigr)$.
(See \cite{sabbah4}, for example.)

\vspace{.1in}

Let $X$ be a complex manifold.
For a complex analytic subset $Z$,
we set $\nbigo_{\Zhat}:=
 \varprojlim\nbigo_{X}/\nbigi_Z^m$,
where $\nbigi_Z$ denote the ideal sheaf
of $Z$.
\index{sheaf $\nbigo_{\Zhat}$}
We set $\Omega^{\bullet,\bullet}_{\Zhat}:=
 \Omega^{\bullet,\bullet}_{X|\Zhat}$
which is equipped with 
the differential operators $\del$ and $\delbar$.
\index{sheaf $\Omega^{\bullet,\bullet}_{\Zhat}$}
If $Z$ is smooth,
it is easy to see that
the natural inclusion
$\nbigo_{\Zhat}\lrarr
 \Omega^{0,\bullet}_{\Zhat}$
is a quasi-isomorphism.

Let $D$ be a simple normal crossing hypersurface
with the irreducible decomposition
$D=\bigcup_{i\in\Lambda}D_i$.
By the above procedures,
we obtain the complexes
$\nbigk^{\bullet}\bigl(
 \nbigo_{\Dhat(I)}
 \bigr)$.
It is known that 
the natural inclusion
$\nbigo_{\Dhat(I)}\lrarr
 \nbigk^{0}(\nbigo_{\Dhat(I)})$
induces a quasi-isomorphism
$\nbigo_{\Dhat(I)}\simeq
 \nbigk^{\bullet}\bigl(
 \nbigo_{\Dhat(I)}
 \bigr)$.
(See \cite{har2} and \cite{sabbah4}.)
We also have
$\Omega^{0,\bullet}_{\Dhat(I)}
 \simeq
 \nbigk^{\bullet}\bigl(
 \Omega^{0,\bullet}_{\Dhat(I)}\bigr)$.
Then, we obtain
$\nbigo_{\Dhat(I)}\simeq
 \Omega_{\Dhat(I)}^{0,\bullet}$.

\vspace{.1in}
We recall a useful isomorphism 
due to Z. Mebkhout 
(Lemma 2.2.1.3 of \cite{Mebkhout-positivity}).
\footnote{The author thanks the referee who informed
this result to him.}
\begin{prop}[Z. Mebkhout]
\label{prop;14.1.16.1}
Let $\nbigm$ be any coherent $\nbigd_X$-module.
Let $Z$ be any hypersurface of $X$.
Then, 
$R\nhom_{\nbigd_X}(\nbigm(\ast Z),\nbigo_{\Zhat})=0$
and 
$\nbigm(!Z)\otimes_{\nbigd_X}^L\nbigo_{\Zhat}=0$.
\hfill\qed
\end{prop}
See (3.10) of \cite{kashiwara_text}
to deduce the second vanishing from the first.

%% file: 2.2.tex
Let us recall Beilinson's beautiful construction
of the nearby cycle functor,
the vanishing cycle functor
and the maximal functor,
which is essential for our purpose.
It is particularly convenient for the study of functoriality.
See \cite{beilinson2}
for more details and precisions.
(See also \cite{MacPherson-Vilonen}
 and \cite{Verdier}.)

\subsection{Preliminary}

Let $k$ be any field of characteristic $0$.
Let $A:=k(\!(s)\!)$ and $A^i:=s^ik[\![s]\!]$.
For $a\leq b$,
we put $A^{a,b}:=A^a\big/A^{b}$.
The multiplication of $s$ induces
a nilpotent endomorphism $N_A$ of $A^{a,b}$.
We put $G_m:=\Spec k[t,t^{-1}]$.
We define
$\gbigi^{a,b}:=\nbigo_{G_m}\otimes A^{a,b}$.
\index{sheaf $\gbigi^{a,b}$}
It is equipped with the connection given by
$\nabla \alpha=N_A(\alpha) (dt/t)$
for 
$\alpha\in A^{a,b}$.
We have natural morphisms
$\gbigi^{a,b}\lrarr \gbigi^{c,d}$
for $a\geq c$ and $b\geq d$,
which are compatible with the connections.
We have a natural isomorphism
$\gbigi^{a,a+1}\simeq 
 \gbigi^{0,1}=\nbigo_{G_m}$
given by $s^a\longleftrightarrow 1$.

This construction makes sense
also in the analytic situation.
The multi-valued flat sections are 
formally given by
$\alpha\cdot \exp\bigl(-s\log t \bigr)$
for $\alpha\in A^{a,b}$.

\subsection{Nearby cycle functor
 and maximal functor}

Let $X$ be any complex manifold
with a hypersurface $D$.
Let $f$ be a meromorphic function on $(X,D)$,
i.e., the poles of $f$ are contained in $D$.
We set 
$\gbigi_f^{a,b}:=f^{\ast}\gbigi^{a,b}(\ast D)$,
which are meromorphic flat bundles
on $\bigl(X,f^{-1}(0)\cup D\bigr)$.
\index{sheaf $\gbigi_f^{a,b}$}
Let $j:X-f^{-1}(0)\lrarr X$.
For a holonomic $\nbigd_{X(\ast D)}$-module $\nbigm$,
we obtain the following
holonomic $\nbigd_{X(\ast D)}$-modules:
\[
 \nbigm_f^{a,b}:=\nbigm\otimes \gbigi_f^{a,b}
=j_{\ast}j^{\ast}\bigl(
 \nbigm\otimes\gbigi_f^{a,b}
 \bigr)
\]
We obtain $\nbigd_{X(\ast D)}$-modules
$\Pi_{f!}^{a,b}\nbigm:=
 j_!
 j^{\ast}\nbigm_f^{a,b}$
and
$\Pi_{f\ast}^{a,b}\nbigm:=
 j_{\ast}
 j^{\ast}\nbigm_f^{a,b}$.
\index{sheaf $\Pi_{f\bikkuri}^{a,b}$}
\index{sheaf $\Pi_{f\ast}^{a,b}$}
We define
\[
 \Pi^{a,b}_{f\ast !}(\nbigm):=
\varprojlim_{N\to\infty}
 \Cok\bigl(
 \Pi^{b,N}_!\nbigm
\lrarr
 \Pi^{a,N}_{\ast}\nbigm
 \bigr).
\]
\index{sheaf $\Pi^{a,b}_{f\ast\bikkuri}$}
The following lemma is easy to see.
\begin{lem}
\label{lem;13.4.15.1}
For any point $P\in X$,
there exists a neighbourhood $X_P$ 
and a large integer $N_0$
such that 
the following natural morphisms are isomorphisms
on $X_P$ for any $N\geq N_0$:
\[
  \Cok\bigl(
 \Pi^{b,N+1}_!\nbigm
\lrarr
 \Pi^{a,N+1}_{\ast}\nbigm
 \bigr)
\lrarr 
 \Cok\bigl(
 \Pi^{b,N}_!\nbigm
\lrarr
 \Pi^{a,N}_{\ast}\nbigm
 \bigr)
\]
\end{lem}
\pf
See the proof of Lemma 4.1.1 of \cite{Mochizuki-MTM},
for example.
\hfill\qed

\vspace{.1in}

Beilinson defined the functors
$\psi_f^{(a)}:=\Pi^{a,a}_{f\ast !}$ 
and $\Xi_f^{(a)}:=\Pi^{a,a+1}_{f\ast !}$.
\index{nearby cycle functor $\psi_f^{(a)}$}
\index{maximal functor $\Xi_f^{(a)}$}
In the case $a=0$,
they are denoted by
$\psi_f\nbigm$ and $\Xi_f\nbigm$,
respectively.
The multiplication of $s$ naturally induces
isomorphisms
$\psi_f^{(a)}\nbigm\simeq\psi_f^{(a+1)}\nbigm$
and $\Xi_f^{(a)}\nbigm\simeq\Xi_f^{(a+1)}\nbigm$.
Note that we have natural isomorphisms
$\Pi_{f\star}^{a,a+1}(\nbigm)\simeq
 j_{\star}j^{\ast}\nbigm$
for $\star=\ast,!$
induced by the multiplication of a power of $s$.
They will be implicitly identified.
We have the exact sequences
of holonomic $\nbigd_{X(\ast D)}$-modules:
\[
\begin{CD}
 0 @>>>
 \Pi_{f!}^{a,a+1}\nbigm @>{c_1^{(a)}}>> 
 \Xi_f^{(a)}\nbigm @>{c_2^{(a)}}>>
 \psi_f^{(a)}\nbigm @>>> 0
\end{CD}
\]
\[
\begin{CD}
 0@>>> 
 \psi_f^{(a+1)}\nbigm @>{d_1^{(a)}}>>
 \Xi_f^{(a)}\nbigm @>{d_2^{(a)}}>>
 \Pi_{f\ast}^{a,a+1}\nbigm
 @>>> 0
\end{CD}
\]
The multiplication of $s$
and the endomorphism $c_2^{(a)}\circ d_1^{(a)}$
induce an endomorphism $N^{(a+1)}$
of $\psi_f^{(a+1)}\nbigm$.

Recall the important observation
$\underset{\leftrightarrow}{\lim}
 \Pi_{f!}^{a,b}
 \nbigm
\simeq
 \underset{\leftrightarrow}{\lim}
 \Pi_{f\ast}^{a,b}\nbigm$
due to Beilinson.
(See \cite{beilinson2} for
$\underset{\leftrightarrow}{\lim}$.)
In particular, it implies that
$N^{(a+1)}$ is locally nilpotent.
We also obtain the following isomorphism:
\[
 \Pi^{a,b}_{f\ast !}(\nbigm)
\simeq
 \varinjlim_{N\to\infty}
 \Ker\bigl(
 \Pi^{-N,b}_{f!}\nbigm
\lrarr
 \Pi^{-N,a}_{f\ast}\nbigm
 \bigr)
\]
As in Lemma \ref{lem;13.4.15.1},
$\Ker\bigl(
 \Pi^{-N,b}_{f!}\nbigm
\lrarr
 \Pi^{-N,a}_{f\ast}\nbigm
 \bigr)$
is locally independent of 
the choice of a large $N$.
See \S4.1 of \cite{Mochizuki-MTM}
for an elementary argument.
In particular, we have
the following identifications:
\begin{align}
\label{eq;09.9.21.2}
 \psi^{(a)}_f\nbigm
\simeq
\varinjlim_{N\to\infty}
 \Ker\bigl(
 \Pi^{-N,a}_{f!}\nbigm
\lrarr
 \Pi_{f\ast}^{-N,a}\nbigm
 \bigr),
 \\
\label{eq;13.5.1.2}
 \Xi^{(a)}_f\nbigm\simeq
 \varinjlim_{N\to\infty}
 \Ker\bigl(
 \Pi_{f!}^{-N,a+1}\nbigm
\lrarr
 \Pi_{f\ast}^{-N,a}\nbigm
 \bigr).
\end{align}

\begin{rem}
When we distinguish that we work
on the category of $\nbigd_{X(\ast D)}$-modules,
we will use the symbols
$\psi_f^{(a)}(\nbigm,\ast D)$,
$\Xi_f^{(a)}(\nbigm,\ast D)$,
etc..
\index{nearby cycle functor $\psi_f^{(a)}(\bullet,\ast D)$}
\index{maximal functor $\Xi_f^{(a)}(\bullet,\ast D)$}
\hfill\qed
\end{rem}

\subsection{Vanishing cycle functor
 and gluing}

Let $f$ be as above.
Let $\nbigm_X$ be any holonomic 
$\nbigd_{X(\ast D)}$-module.
We set $\nbigm:=\nbigm_X(\ast f)$.
We have the natural identifications
$\Pi^{a,b}_{f\star}\nbigm_X
=\Pi^{a,b}_{f\star}\nbigm$
for $\star=\ast,!$.
We also have
$\Pi^{a,b}_{f\ast !}\nbigm_X
=\Pi^{a,b}_{f\ast !}\nbigm$.
In particular,
$\psi^{(a)}_{f}\nbigm_X
=\psi^{(a)}_{f}\nbigm$
and 
$\Xi^{(a)}_{f}\nbigm_X
=\Xi^{(a)}_{f}\nbigm$.
We set 
$\nbigm_X^{(a)}:=\nbigm_X\otimes A^{a,a}$.
We have the naturally defined morphisms:
\[
\begin{CD}
 \Pi_{f!}^{a,a+1}\nbigm
 @>{c_{1,X}^{(a)}}>>
 \nbigm^{(a)}_X 
 @>{d_{2,X}^{(a)}}>>
 \Pi_{f\ast}^{a,a+1}\nbigm
\end{CD}
\]
Beilinson defined
the vanishing cycle functor
$\phi_f^{(a)}\nbigm_X$
as the $H^1$-cohomology of 
the following sequence
of holonomic $\nbigd_{X(\ast D)}$-modules:
\[
\begin{CD}
 \Pi^{a,a+1}_{f!}\nbigm
 @>{c_1^{(a)}\oplus c_{1,X}^{(a)}}>>
 \Xi_f^{(a)}\nbigm
 \oplus\nbigm^{(a)}_X 
 @>{d_2^{(a)}-d_{2,X}^{(a)}}>>
 \Pi^{a,a+1}_{f\ast}\nbigm
\end{CD}
\]
\index{vanishing cycle functor $\phi^{(a)}_f$}
The morphisms
$d_1^{(a)}$ and $c_2^{(a)}$ induce
$\can$ and $\var$:
\[
\begin{CD}
 \psi_f^{(a+1)}\nbigm
@>{\can}>>
 \phi_f^{(a)}\nbigm
@>{\var}>>
 \psi_f^{(a)}\nbigm
\end{CD}
\]
By construction,
we have
$\var\circ\can=c_2^{(a)}\circ d_1^{(a)}$.

Conversely,
let $\nbigm_Y$ be a holonomic 
$\nbigd_{X(\ast D)}$-module
whose support is contained in
$Y=f^{-1}(0)$,
with morphisms 
\[
 \psi_f^{(1)}\nbigm
\stackrel{u}{\lrarr}
 \nbigm_Y
\stackrel{v}{\lrarr}
 \psi_f^{(0)}\nbigm,
\quad\quad
 v\circ u=c_2^{(0)}\circ d_1^{(0)}.
\]
Then, 
we obtain a holonomic 
$\nbigd_{X(\ast D)}$-module
$\Glue(\nbigm_Y,u,v)$
as the cohomology of
the complex:
\[
\begin{CD}
 \psi_f^{(1)}\nbigm
 @>{d_1^{(0)}\oplus u}>>
 \Xi_f(\nbigm)\oplus\nbigm_Y
 @>{c_2^{(0)}-v}>>
 \psi_f^{(0)}\nbigm
\end{CD}
\]
Beilinson made an excellent observation
that the above two operations are mutually inverse.
See \cite{beilinson2} for more details.

\subsection{Comparison with
 ordinary definitions}

Let $\psitilde_{f,-1}$ and $\phitilde_{f}$
be the nearby cycle functor
and the vanishing cycle functor
defined in terms of $V$-filtrations,
i.e.,
$\psitilde_{f,-1}(\nbigm)
=\Gr^V_{-1}(\iota_{f\dagger}\nbigm)$
and 
$\phitilde_{f}(\nbigm_X):=
 \Gr^{V}_0(\iota_{f\dagger}\nbigm_X)$,
where $\iota_f:X\lrarr X\times\cnum$
denotes the graph,
and $V$ denotes a $V$-filtration
of $\iota_{f\dagger}\nbigm_X$ along $t$.
For simplicity,
$\psitilde_{f,-1}$ is denoted by
$\psitilde_{f}$ in the following.
\index{nearby cycle functor $\psitilde_{f}$}
\index{vanishing cycle functor $\phitilde_{f}$}

\begin{lem}
We have natural isomorphisms
$\psi_f\simeq\psitilde_{f,}$
and $\phi_f\simeq\phitilde_f$.
\end{lem}
\pf
Recall that $\phitilde_f(\nbigm_X)$
and $\psitilde_f(\nbigm_X)$ are naturally
equipped with the nilpotent endomorphisms $N$,
which are the nilpotent part of 
the multiplication of $-\del_tt$.
We have natural identifications:
\[
 \phitilde_f\bigl(
 \Pi_{f!}^{a,b}\nbigm
 \bigr)
\simeq
 \phitilde_f\bigl(\Pi_{f\ast}^{a,b}\nbigm\bigr)
\simeq
 \psitilde_f\nbigm\otimes A^{a,b}
\]
The natural nilpotent endomorphisms are given by
$N\otimes \id-\id\otimes (s\bullet)$,
which is denoted by $N-s$.
Here, $s\bullet$ denotes the multiplication of
$s$ on $A^{a,b}$.
In the following,
we argue on any compact subset of $X$.

Let us look at the natural morphism
$G^{a,b}:
 \Pi^{a,b}_{f!}\nbigm  
 \lrarr
 \Pi^{a,b}_{f\ast}\nbigm 
$.
The supports of the kernel and the cokernel
are contained in $f^{-1}(0)$.
The morphism
$\phitilde_f(G^{a,b}):
 \phitilde_f\bigl(
 \Pi^{a,b}_{f!}\nbigm
 \bigr)
\lrarr 
  \phitilde_f\bigl(
 \Pi^{a,b}_{f\ast}\nbigm
 \bigr)$ is naturally identified with 
$N-s:
 \psitilde_f\nbigm\otimes A^{a,b}
\lrarr
 \psitilde_f\nbigm\otimes A^{a,b}$.
Hence, 
if $b$ is sufficiently larger than $a$,
$\Cok(G^{a,b})$ is isomorphic to
$\psitilde_f\nbigm\otimes A^{a,a+1}$,
independently of $b$.
Therefore, we obtain 
$\psi^{(a)}_f\nbigm
\simeq 
\psitilde_{f}\nbigm\otimes A^{a,a+1}$.
In particular,
we naturally have 
$\psi_f^{(0)}\nbigm=\psitilde_f\nbigm$.

It follows that
$\Cok\Bigl(
 \Pi_{f!}^{a+1,M}\nbigm
\lrarr
 \Pi_{f\ast}^{a,M}\nbigm
 \Bigr)$
are independent of any sufficiently large $M$,
which should be isomorphic to
$\Xi^{(a)}_f\nbigm$.
We obtain
$\phitilde_f\bigl(
 \Xi_f^{(a)}\nbigm
 \bigr)
\simeq
\Cok\Bigl(
 N-s:
\psi_f\nbigm\otimes A^{a+1,M}\lrarr
 \psi_f\nbigm\otimes A^{a,M}
 \Bigr)$ for any sufficiently large $M$.
Because $\phi_f^{(0)}(\nbigm_X)$
is naturally isomorphic to the cohomology
of the complex
\[
 \phitilde_f\bigl(\Pi_{f!}^{0,1}\nbigm\bigr)
\lrarr
 \phitilde_f\bigl(\Xi_f^{(0)}\nbigm\bigr)
\oplus
 \phitilde_f\bigl(\nbigm_X\bigr)
\lrarr
 \phitilde_f\bigl(\Pi_{f\ast}^{0,1}\nbigm\bigr),
\]
it is easy to obtain
$\phi_f^{(0)}(\nbigm)\simeq\phitilde_f(\nbigm)$
by a direct calculation.
\hfill\qed

\subsection{Compatibility with dual}

In \cite{beilinson2},
the pairing $A\times A\lrarr k= A^{-1}/A^{0}$
is given by
$\bigl\langle
 f(s),\,g(s)
 \bigr\rangle=\Res_{s=0}\bigl(
 f(s)\,g(-s)\,ds\bigr)$.
It induces pairings
$A^{a,b}\otimes
 A^{-b,-a}\lrarr A^{-1}/A^0$.
Then, we obtain flat pairings
$\gbigi^{a,b}\otimes
 \gbigi^{-b,-a}\lrarr \gbigi^{-1,0}$.
We can identify $\gbigi^{a,b}$
with the dual of $\gbigi^{-b,-a}$
by the pairing.

Let $\DDD$ denote the dual functor
on the category of
holonomic $\nbigd_{X(\ast D)}$-modules.
By using the $\nbigd_{X(\ast D)}$-version of 
Lemma \ref{lem;09.12.1.1},
we obtain identifications:
\[
 \DDD\Bigl(
 \Pi_{f\ast}^{a,b}
 \nbigm
 \Bigr)
\simeq
 \Pi_{f!}^{-b,-a}\Bigl(
 \DDD\nbigm
 \Bigr),
\quad\quad
\DDD\Bigl(
 \Pi_{f!}^{a,b}\nbigm\Bigr)
\simeq
 \Pi_{f\ast}^{-b,-a}\Bigl(
 \DDD\nbigm
 \Bigr)
\]
By (\ref{eq;09.9.21.2}) and (\ref{eq;13.5.1.2}),
we obtain the identifications
$\DDD_X\psi^{(a)}_f(\nbigm)
\simeq
 \psi^{(-a)}_f\bigl(\DDD_X\nbigm\bigr)$
and 
$\DDD_X\Xi^{(a)}_f(\nbigm)
\simeq
 \Xi^{(-a-1)}_f\bigl(\DDD_X\nbigm\bigr)$.
We have $\DDD_X(c_1^{(a)})=d_2^{(-a-1)}$,
$\DDD_X(c_2^{(a)})=d_1^{(-a-1)}$
and $\DDD_X(c_{1,X}^{(a)})=d_{2,X}^{(-a-1)}$.
Hence, we obtain
$\DDD_X\phi^{(a)}_f(\nbigm_X)
\simeq
 \phi^{(-a-1)}_f(\DDD_X\nbigm_X)$.
The morphisms
$\DDD_X\psi_f^{(a)}\nbigm
\stackrel{\DDD\var}{\lrarr}
 \DDD_X\phi_f^{(a)}\nbigm_X
\stackrel{\DDD\can}{\lrarr}
 \DDD_X\psi_f^{(a-1)}\nbigm$
are identified with
$\psi_f^{(-a+1)}\nbigm
\stackrel{\can}{\lrarr}
 \phi_f^{(-a)}\nbigm_X
\stackrel{\var}{\lrarr}
  \psi_f^{(-a)}\nbigm$.

The multiplication of $s$ induces
an isomorphism
$\Phi_s:\psi^{(a)}(\nbigm)\simeq
 \psi^{(a+1)}(\nbigm)$, etc.
Under the above identifications,
we have $\DDD\Phi_s=-\Phi_s$.

\begin{rem}
In {\rm\cite{Mochizuki-MTM}},
we use the pairing
$A\times A\lrarr k$
given by
$\langle f(s),g(s)\rangle
=\Res_{s=0}\bigl(f(s)g(-s)ds/s\bigr)$.
It makes an inessential shift of the indexes
in the formulas.
\hfill\qed
\end{rem}

\subsection{Compatibility with push-forward}

Let $F:X\lrarr Y$ be any proper morphism.
Assume that $D=F^{-1}(D_Y)$,
for simplicity.
Let $g$ be any holomorphic function on $Y$.
Let $\nbigm$ be any holonomic 
$\nbigd_{X(\ast D)}$-module.
We set $\gtilde:=F^{\ast}g$.
Let $j_Y:Y-g^{-1}(0)\lrarr Y$
and $j_X:X-\gtilde^{-1}(0)\lrarr X$.
We have natural isomorphisms
$F^i_{\dagger}\bigl(
 \nbigm\otimes \gbigi^{a,b}_{\gtilde}
 \bigr)
\simeq
 F^i_{\dagger}(\nbigm)
 \otimes \gbigi^{a,b}_g$
of $\nbigd_{Y(\ast D_Y)}$-modules.
We naturally have
$(j_{Y\star}j_Y^{\ast})F^i_{\dagger}
\simeq
 F^i_{\dagger}\circ (j_{X\star}j_X^{\ast})$
for $\star=\ast,!$.
Hence, it is easy to obtain the following
identifications:
\[
 F^i_{\dagger}\psi^{(a)}_{\gtilde}\nbigm
=\psi^{(a)}_gF^i_{\dagger}\nbigm,
\quad\quad
 F^i_{\dagger}\Xi^{(a)}_g\nbigm
=\Xi^{(a)}_gF_{\dagger}\nbigm,
\quad\quad
 F^i_{\dagger}\phi^{(a)}_g\nbigm
=\phi^{(a)}_gF^i_{\dagger}\nbigm.
\]
\subsection{Choice of a function}
\label{subsection;09.10.22.5}

Let $f$ and $h$ be meromorphic functions on 
$(X,D)$.
We suppose that $h$ is nowhere vanishing on $X\setminus D$.
We have natural isomorphisms of $\nbigo_X$-modules
$\gbigi_f^{a,b}\simeq
 \gbigi_{hf}^{a,b}
 \simeq
 A^{a,b}\otimes\nbigo_{X(\ast D)}(\ast f)$.
For their flat connections
$\nabla_f$ and $\nabla_{hf}$
and for $\alpha\in A^{a,b}$,
we have the formulas:
\[
 \nabla_f \alpha=\alpha\cdot s\frac{df}{f},
\quad\quad\quad
  \nabla_{hf}\alpha=\alpha\cdot
 s\,\left(
 \frac{df}{f}+\frac{dh}{h}
 \right)
\]
If we have $\log h$ on $X$,
we have a flat isomorphism
$\Phi:\gbigi_f^{a,b}\simeq \gbigi_{hf}^{a,b}$
given by
$\Phi(\alpha)=
 \exp\bigl(-s\log h\bigr)\,
 \alpha$.
It induces isomorphisms:
\begin{equation}
 \label{eq;09.10.22.6}
\Xi_f^{(a)}\simeq \Xi^{(a)}_{hf},
 \quad\quad
\psi_f^{(a)}\simeq\psi^{(a)}_{hf},
 \quad\quad
\phi_f^{(a)}\simeq\phi^{(a)}_{hf}.
\end{equation}
They depend on the choice of 
a branch of $\log h$.

\subsection{$\rnum$-structure of $\gbigi^{a,b}$}

In the analytic case,
the $\rnum$-structure of $A^{a,b}$
is given as follows:
\[
 \cnum\cdot s^j\supset 
 \rnum\cdot(2\pi\sqrt{-1})^j s^j
\]
It gives a $\rnum$-structure
of the fiber of $\gbigi^{a,b}$
over $1\in\cnum^{\ast}$.
We extend it to a flat $\rnum$-structure
of the flat bundle $\gbigi_{|\cnum^{\ast}}$.
Let $u:=2\pi\sqrt{-1}\,s$.
The connection of $\gbigi^{a,b}$ is expressed as
\[
 \nabla (u^{a},\ldots,u^{b-1})
=(u^{a},\ldots,u^{b-1})
 \cdot N \frac{1}{2\pi\sqrt{-1}}
 \frac{dt}{t}
\]
Here, $N$ denotes the constant matrix
such that
$N_{i,i+1}=1$ and $N_{i,j}=0$ otherwise.
Since the monodromy is expressed by $\exp(-N)$,
the $\rnum$-structure is well defined.
More generally,
for any subfield $K\subset\cnum$,
we obtain a $K$-structure of $\gbigi^{a,b}$
in this way.
The pairing
$\langle\cdot,\cdot\rangle:
 \gbigi^{a,b}\otimes\gbigi^{-b,-a}
\lrarr\gbigi^{-1,0}$
is defined over $\rnum$.
Under the identification
$\gbigi^{-1,0}\simeq \gbigi^{0,1}$
by the multiplication of $s$,
the pairing takes values
in $(2\pi\sqrt{-1})^{-1}\rnum$.

\subsection{Comparison with
 the functors for perverse sheaves}

Let $\Loc(\gbigi^{a,b})_{\rnum}$ denote
the $\rnum$-local system associated to $\gbigi^{a,b}$.
The fiber over $1$ is 
 $u^a\rnum[\![u]\!]
 \big/
 u^b\rnum[\![u]\!]$,
and the monodromy along the loop
with the clockwise direction 
is given by the multiplication of $\exp(u)$.
Taking the limit,
we have a $\rnum$-local system
$\Loc(\gbigi)_{\rnum}$,
whose fiber over $1$ is
$\rnum(\!(u)\!)$,
and the monodromy is given by
the multiplication of $\exp(u)$.
We have subsystems
$\Loc(\gbigi^a)_{\rnum}
\subset
 \Loc(\gbigi)_{\rnum}$
whose fiber over $1$ is
$u^a\rnum[\![u]\!]$.
We have
$\Loc(\gbigi^{a,b})_{\rnum}
\simeq
 \Loc(\gbigi^a)_{\rnum}/\Loc(\gbigi^b)_{\rnum}$.
Recall another expression
of these local systems as in \cite{beilinson2}.

Let $A_{\nbigp}:=\rnum(\!(v)\!)$.
We set $t:=v+1$.
The pairing $A_{\nbigp}\times A_{\nbigp}
\lrarr \rnum(-1)$ is given as follows:
\[
 \bigl\langle
 f(t),g(t)
 \bigr\rangle
=\underset{t=1}\Res\Bigl(
 f(t)\,g(t^{-1})
 \frac{dt}{t}
 \Bigr)
\frac{1}{2\pi\sqrt{-1}}
\]
We have a $\rnum$-local system $\gbigi_{\nbigp}$
on $\cnum^{\ast}$
such that 
the fiber over $1$ is $A_{\nbigp}$,
and the monodromy along the loop
with the clockwise direction 
is given by the multiplication of $t=1+v$.
Let us compare
$\gbigi_{\nbigp}$ and $\Loc(\gbigi)_{\rnum}$.
We take an algebra homomorphism
$\Phi:\rnum(\!(u)\!)\longrightarrow
 \rnum(\!(v)\!)$ determined by
$\Phi\bigl(
 \exp(u)
 \bigr)
=1+v$.
We identify the fibers of
$\Loc(\gbigi)_{\rnum}$ 
and $\gbigi_{\nbigp}$
by $\Phi$.
Because it is compatible with
the monodromy,
it induces the identification
$\Loc(\gbigi)_{\rnum}
 \simeq\gbigi_{\nbigp}$.
Note that
$\Phi\bigl(f(-u)\bigr)
=\Phi(f)(t^{-1})$
and 
$\Phi(du)=dt/t$.
Hence the pairing is preserved.

\begin{rem}
Recall that the functors
$\psi$, $\Xi$ and $\phi$ for
perverse sheaves 
are given in terms of
$\gbigi_{\nbigp}$,
according to {\rm\cite{beilinson2}}.
The above comparison
gives the compatibility of
the de Rham functor
$\DR$ with $\phi$, $\psi$ and $\Xi$
in the regular singular case.
\hfill\qed
\end{rem}

%% file: 3.1.tex
We shall introduce
the notion of good holonomic $\nbigd$-modules
on any complex manifold $X$
with a normal crossing hypersurface 
$D=\bigcup_{i\in\Lambda}D_i$.
They are $\nbigd$-modules
locally described as the gluing of
meromorphic flat bundles
on $\bigcap_{j\in J} D_j$ $(J\subset\Lambda)$.
In \S\ref{subsection;09.10.19.5}--\ref{subsection;09.10.23.2},
we study the local case.
We explain the global case 
in \S\ref{subsection;14.1.14.1}.
We explain a kind of quiver description
of good holonomic $\nbigd$-modules 
in the local case
in \S\ref{subsection;14.1.14.2}.

In the local case,
for any good holonomic $\nbigd$-modules,
we have various commutativity of functors
such as $\phi_i^{(a)}\phi_j^{(b)}(\nbigm)
\simeq
 \phi_j^{(b)}\phi_i^{(a)}(\nbigm)$,
for which goodness seems  truly used.

\subsection{$\nbigi$-good
meromorphic flat bundles}
\label{subsection;09.10.19.5}

Let $\Delta^n$ denote a multi-disc in $\cnum^n$,
i.e.,
$\Delta^n:=\{(z_1,\ldots,z_n)\in\cnum^n\,|\,
|z_i|<1 \}$.
We consider the case
$X:=\Delta^n$,
$D_i:=\{z_i=0\}$
and $D:=\bigcup_{i=1}^{\ell}D_i$.
We set $\ellsitabar:=\{1,\ldots,\ell\}$.
For $I\subset\ellsitabar$,
we set $D(I):=\bigcup_{i\in I}D_i$
and $D_I:=\bigcap_{i\in I}D_i$.
\index{divisor $D(I)$}
\index{set $D_I$}
\index{set $\del D_I$}
\index{set $M(X,D)$}
\index{set $H(X)$}
\index{set $\nbigi(I)$}
We put $\del D_I:=D_I\cap D(I^c)$,
where $I^c:=\ellsitabar-I$.
Let $M(X,D)$ be the set of meromorphic
functions on $X$ whose poles are contained in $D$.
Let $H(X)$ be the set of holomorphic functions
on $X$.
We give a review on good meromorphic flat bundles.
See \cite{Mochizuki-DM},
\cite{mochi8}
and \cite{mochi10} for more detailed reviews.

\subsubsection{Good set of irregular values}

Let $f\in M(X,D)$.
Suppose that there exists $\vecm=(m_i)\in\seisuu_{\geq 0}^{\ell}$
such that
(i) $\vecz^{\vecm}f=\prod z_i^{m_i}f$ is holomorphic,
(ii) if $\vecm\neq (0,\ldots,0)$,
we have $(\vecz^{\vecm}f)(O)\neq 0$.
Then, we set $\ord(f):=-\vecm$.
In general, such $\vecm$ does not exist.
For any holomorphic function $f$,
we have $\ord(f)=(0,\ldots,0)$.
If $\ord(g)$ exists for $g\in\nbigo_X(\ast D)$,
then $\ord(g+f)=\ord(g)$ for any 
holomorphic function $f$.
So, the notion $\ord$ is considered
for elements in $M(X,D)/H(X)$.

We use the order $\leq$ on $\seisuu^{\ell}$
given by
$\vecm\leq\vecn$
$\stackrel{\rm def}{\Longleftrightarrow}$
$m_i\leq n_i$ for any $i$.
A finite subset 
$\nbigi\subset M(X,D)/H(X)$
is called good
if the following holds:
\begin{itemize}
\item
 For any $f\in \nbigi$,
 there exists $\ord(f)$.
\item
 For any $f,g\in\nbigi$,
 there exists $\ord(f-g)$,
 and the set
 $\bigl\{
 \ord(f-g)\,\big|\,f,g\in\nbigi
 \bigr\}$ is totally ordered.
\end{itemize}

For any good set of irregular values
$\nbigi\subset M(X,D)/H(X)$
and for any subset $I\subset \ellsitabar$,
let $\nbigi'(I)$ be the set of
the elements $\gminia\in\nbigi$ which are 
regular along $z_i$ $(i\in I)$,
and we put
$\nbigi(I):=\bigl\{
 \gminia_{|D_I}\,\big|\,
 \gminia\in \nbigi'(I)
 \bigr\}$.
It is a good set of irregular values
on $(D_I,\del D_I)$.

\subsubsection{Unramifiedly $\nbigi$-good meromorphic flat bundle}

Let $\nbigi\subset M(X,D)/H(X)$
be a good set of irregular values.
Recall that a meromorphic flat bundle 
$(\nbige,\nabla)$ on $(X,D)$ is called 
unramifiedly $\nbigi$-good
if the following holds:
\begin{itemize}
\item
Let $\nbigi_I$ denote the image of
$\nbigi$ to $M(X,D)/M(X,D(I^c))$.
For any $P\in D_I\setminus\del D_I$,
the formal completion
$(\nbige,\nabla)_{|\Phat}$
is decomposed into
$\bigoplus_{\gminib\in \nbigi_I}
 (\nbigehat_{P,\gminib},\nablahat_{P,\gminib})$
such that
$\nablahat_{P,\gminib}
 -d\gminibtilde\,\id_{\nbigehat_{P,\gminib}}$
are regular singular,
where $\gminibtilde$ are any lifts
of $\gminib$ to $M(X,D)$.
\end{itemize}
In this paper,
we say that a meromorphic flat bundle
$(\nbige,\nabla)$ on $(D_I,\del D_I)$
is unramifiedly $\nbigi$-good
if it is unramifiedly $\nbigi(I)$-good.

\subsubsection{Ramified case}
\label{subsection;14.1.16.30}

For a positive integer $m$,
let $X^{(m)}:=\Delta^n=\{|\zeta_i|<1\}$,
$D_i^{(m)}:=\bigl\{\zeta_i=0\bigr\}$
and 
$D^{(m)}=\bigcup_{i=1}^{\ell} D_i^{(m)}$.
We have a natural ramified covering
$\varphi_m:X^{(m)}\lrarr X$ along $D$
given by
$\varphi_m(\zeta_1,\ldots,\zeta_n)
=(\zeta_1^m,\ldots,\zeta_{\ell}^m,\zeta_{\ell+1},\ldots,\zeta_n)$,
and the induced ramified coverings
$D_I^{(m)}\lrarr D_I$.
Let $\nbigi\subset M(X^{(m)},D^{(m)})/H(X^{(m)})$
be any good set of irregular values
which is preserved by the action of
the Galois group of the ramified covering $X^{(m)}/X$.
In this paper,
a meromorphic flat bundle $\nbige$ on $(D_I,\del D_I)$
is called $\nbigi$-good
if it is the descent of an unramifiedly $\nbigi$-good 
meromorphic flat bundle $\nbige^{(m)}$
on $(D_I^{(m)},\del D_I^{(m)})$.

\subsubsection{Some functors along the divisors}

In this subsection,
we use the following notation
for simplicity of the description.
\begin{notation}
The vanishing cycle functors $\phi^{(a)}_{z_i}$
are denoted by $\phi^{(a)}_i$.
For any $I=(i_1,\ldots,i_m)\in\{1,\ldots,\ell\}^m$
and any $\veca=(a_1,\ldots,a_m)\in \seisuu^{m}$,
 we set
 $\phi^{(\veca)}_I=
 \phi^{(a_1)}_{i_1}\circ\cdots\circ\phi^{(a_m)}_{i_m}$.
 If $\veca=(0,\ldots,0)$, 
 it is often denoted just by $\phi_I$.
We use the symbols
 $\psi^{(\veca)}_I$,
$\Xi^{(\veca)}_{I}$ and $\Pi^{a,b}_{i\star}$
with a similar meaning.
For any holonomic $\nbigd_X$-module $\nbigm$,
we set $\nbigm(\ast i):=\nbigm(\ast D_i)$
and $\nbigm(!i):=\nbigm(!D_i)$.
If we are given a subset $I\subset\ellsitabar$,
we put $\nbigm(!I):=\nbigm\bigl(!D(I)\bigr)$
and $\nbigm(\ast I):=\nbigm\bigl(\ast D(I)\bigr)$.
\hfill\qed
\end{notation}
\index{functor $\psi_i^{(a)}$}
\index{functor $\psi^{(\veca)}_I$}
\index{functor $\phi_i^{(a)}$}
\index{functor $\phi^{(\veca)}_I$}
\index{functor $\Xi_i^{(a)}$}
\index{functor $\Xi^{(\veca)}_I$}
\index{sheaf $\nbigm(\bikkuri I)$}
\index{sheaf $\nbigm(\ast I)$}

\begin{lem}
\label{lem;09.10.19.21}
Let $(\nbige,\nabla)$ be any $\nbigi$-good
meromorphic flat bundle on $(X,D)$.
For $1\leq i,j\leq \ell$ with $i\neq j$,
the natural morphism
$\phi^{(a)}_{i}(\nbige)
\lrarr
\phi^{(a)}_{i}(\nbige)(\ast j)$
is an isomorphism.
\end{lem}
\pf
Because the support of 
$\phi^{(a)}_{i}(\nbige)$
and 
$\phi^{(a)}_{i}(\nbige)(\ast j)$
are contained in $D_i$,
it is enough to prove that
the induced morphism for the formal completions
$\phi^{(a)}_i(\nbige)_{|\Phat}
\lrarr
 \phi^{(a)}_i(\nbige)(\ast j)_{|\Phat}$
is an isomorphism for each $P\in D_i$.
We have only to consider the case 
$P=(0,\ldots,0)$.
We use the notation introduced in \S\ref{subsection;14.1.16.30}.
Take lifts $\gminiatilde$ of $\gminia\in\nbigi$.
We have regular singular meromorphic flat bundles
$(R_{\gminia},\nabla_{\gminia})$ on $(X^{(m)},D^{(m)})$
for $\gminia\in\nbigi$,
and an action of the Galois group $G$ of $\varphi_m$
on $(\nbige',\nabla')=
 \bigoplus_{\gminia\in\nbigi}
 (R_{\gminia},\nabla_{\gminia}+d\gminiatilde)$,
such that
the formal completions of
$(\nbige',\nabla')$
and $\varphi_m^{\ast}(\nbige,\nabla)$
at $(0,\ldots,0)$
are isomorphic in a $G$-equivariant way.
Let $(\nbige'',\nabla'')$ be the meromorphic 
flat bundle on $(X,D)$
obtained as the descent of
$(\nbige',\nabla')$.
The formal completions of
$(\nbige'',\nabla'')$
and $(\nbige,\nabla)$ at 
$P$ are isomorphic.
Then, by using the standard argument
to prove the uniqueness of $V$-filtrations,
the isomorphism
$\nbige''_{|\Phat}\simeq
 \nbige_{|\Phat}$
is compatible with the $V$-filtrations along $z_i$.
Therefore, it is enough to prove the claim for $\nbige''$.

Let $(R,\nabla)$ be a regular singular meromorphic flat bundle
on $(X,D)$.
Let $\gminib\in M(X^{(m)},D^{(m)})$ such that
$\ord(\gminib)$ exists.
We set $L(\gminib):=\nbigo_{X^{(m)}}(\ast D^{(m)})\,e$
with the connection
$\nabla e=e\,d\gminib$.
We obtain a meromorphic flat bundle
$\varphi_{m\ast}(L(\gminib))$
on $(X,D)$.
By the previous consideration,
it is enough to prove the claim for
any direct summand of
the meromorphic flat bundle
$\nbige_1=R\otimes\varphi_{m\ast}L(\gminib)$,
which follows from the claim for
$\nbige_1$.
We may assume that
$\gminib=\prod_{j=1}^{\ell}\zeta_j^{b_j}$
for some $b_j\leq 0$.

Let $V(R)$ denote the $V$-filtration along $z_i$.
For $\vecm\in S:=\{0,1,\ldots,m-1\}^{\ell}$,
let $\veczeta^{\vecm}:=\prod_{k=1}^{\ell}\zeta_k^{m_k}$.
We have
$\varphi_{\ast}L(\gminib)
=\bigoplus_{\vecm\in S}
 \nbigo_X(\ast D)\veczeta^{\vecm}e$.
If $b_i<0$,
the $V$-filtration $V(\nbige_1)$ of $\nbige_1$
is given by
$V_{\alpha}(\nbige_1)=\nbige_1$ for any $\alpha\in\cnum$.
If $b_i=0$,
we have
$V_{\alpha}(\nbige_1)
=\bigoplus
 V_{\alpha+m_i/m}(R)\otimes\nbigo_X\veczeta^{\vecm}e$.
Hence, 
the natural morphism
$\phi_i(\nbige_1)\lrarr \phi_i(\nbige_1)(\ast D_j)$ $(j\neq i)$
is an isomorphism in the both cases.
\hfill\qed

\begin{lem}
\label{lem;09.10.19.30}
If $i\neq j$,
the natural morphism
$\nbige(!i)\lrarr \nbige(!i)(\ast j)$
is an isomorphism.
\end{lem}
\pf
Let $N$ denote the nilpotent part
of the action of $-\del_iz_i$ on $\phi_i(\nbige)$.
We have the following commutative diagram:
\[
 \begin{CD}
 0 @>>> \Ker N @>>>
 \nbige(!i) @>>>
 \nbige @>>>\Cok N@>>> 0\\
 @. @V{a}VV @V{b}VV @V{=}VV @V{c}VV @.\\
 0 @>>> \Ker N(\ast j) @>>>
 \nbige(!i)(\ast j) @>>>
 \nbige @>>>\Cok N(\ast j)@>>> 0
 \end{CD}
\]
By Lemma \ref{lem;09.10.19.21},
we obtain that 
$a$ and $c$ are isomorphisms.
Hence, $b$ is also an isomorphism.
\hfill\qed

\subsection{$\nbigi$-good holonomic 
$\nbigd$-modules}
\label{subsection;09.10.28.3}

We continue to use the notation 
introduced in \S\ref{subsection;09.10.19.5}.
\begin{df}
\label{df;09.10.29.1}
A holonomic $\nbigd_X$-module $\nbigm$
is called $\nbigi$-good on $(X,D)$
if the following holds:
\begin{itemize}
\item
 $\nbigm(\ast D)$ is an $\nbigi$-good meromorphic flat bundle
 on $(X,D)$.
\item
For any 
 $I=(i_1,\ldots,i_m)\in\{1,\ldots,\ell\}^m$,
 $\phi_I(\nbigm)\bigl(\ast I^c\bigr)$
 is the push-forward of
 an $\nbigi$-good meromorphic flat bundle
 on $(D_I,\del D_I)$ 
 by $D_I\lrarr X$.
\hfill\qed
\end{itemize}
\end{df}
\index{good holonomic $\nbigd$-module}
\index{$\nbigi$-good holonomic $\nbigd$-module}

The full subcategory of
$\nbigi$-good holonomic $\nbigd$-modules
is abelian,
and it is closed under extensions.
If $V$ is a good meromorphic flat bundle,
it is a good holonomic $\nbigd_X$-module
in the above sense.
When we do not have to distinguish $\nbigi$,
we will omit to denote it.
We will implicitly use the following obvious lemma.
\begin{lem}
\label{lem;09.10.19.25}
Let $\nbigm$ be a holonomic $\nbigd_X$-module.
Suppose that 
(i) $\nbigm(\ast D)$ is an $\nbigi$-good meromorphic
flat bundle,
(ii) $\phi_i(\nbigm)$ are $\nbigi$-good
 for any $i=1,\ldots,\ell$.
Then, $\nbigm$ is $\nbigi$-good.
\hfill\qed
\end{lem}

\begin{lem}
\label{lem;09.10.19.26}
Let $\nbigm$ be an $\nbigi$-good holonomic $\nbigd$-module
on $(X,D)$.
Then $\DDD_X\nbigm$ is $-\nbigi$-good,
where $-\nbigi=\{-\gminia\,|\,\gminia\in\nbigi\}$.
\end{lem}
\pf
We use an induction on the dimension
of the support of $\nbigm$.
It is easy to check that
$\DDD_X\nbigm(\ast D)$ is a good meromorphic 
flat bundle.
By the inductive assumption,
$\phi^{(a)}_i(\DDD_X\nbigm)\simeq
 \DDD_X\phi^{(-a-1)}_i(\nbigm)$ are also good.
Hence, we obtain that $\nbigm$ is good.
\hfill\qed

\vspace{.1in}

For any good holonomic $\nbigd$-module $\nbigm$,
let $\rho(\nbigm)
 \in\seisuu_{\geq \,0}\times\seisuu_{>0}$ 
denote the pair of $\dim\Supp\nbigm$
and the number of the irreducible components of
$\Supp\nbigm$ with the maximal dimension.
We use the lexicographic order on 
$\seisuu_{\geq \,0}\times\seisuu_{>0}$.
For any good holonomic $\nbigd$-module $\nbigm$,
there exists $J\subset\ellsitabar$
with $\dim\Supp\nbigm=n-|J|$
such that
$\nbigm\bigl(\ast J^c\bigr)\neq 0$.
The kernel $\nbign_1$ and the cokernel
$\nbign_2$ of the natural morphism
$\nbigm\lrarr\nbigm\bigl(\ast J^c\bigr)$
satisfy $\rho\bigl(\nbign_i\bigr)<\rho(\nbigm)$
$(i=1,2)$.

\begin{lem}
\label{lem;09.10.19.20}
Let $\nbigm$ be $\nbigi$-good on $(X,D)$.
Then, $\psi^{(a)}_i(\nbigm)$ are also $\nbigi$-good
for any $i=1,\ldots,\ell$.
\end{lem}
\pf
We use an induction on $\rho(\nbigm)$.
Let $J$ and $\nbign_j$ $(j=1,2)$ be as above.
By the assumption of the induction,
$\psi^{(a)}_i(\nbign_j)$ $(j=1,2)$ are good.
The $\nbigd_X$-module $\nbigm(\ast J^c)$
is the push-forward of an $\nbigi$-good meromorphic
flat bundle $\nbige_J$ on $(D_J,\del D_J)$ by the inclusion
$\iota_J:D_J\lrarr X$.
If $i\in J$, we have
$\psi^{(a)}_i(\nbigm(\ast J^c))=0$.
If $i\not\in J$,
$\psi^{(a)}_i(\nbigm(\ast J^c))$
is isomorphic to
$\iota_{J\dagger}\psi^{(a)}_i(\nbige_J)$.
By computing the formal completion
$\psi^{(a)}_i(\nbige_J)_{|\Phat}$
of $P\in \del D_J$
as in the proof of Lemma \ref{lem;09.10.19.21},
we can prove that
$\psi^{(a)}_i(\nbige_J)_{|\Phat}$ is $\nbigi$-good
on $(D_J,\del D_J)$.
Hence,
we obtain that $\psi^{(a)}_i(\nbigm)$
is also $\nbigi$-good.
\hfill\qed

\subsection{Commutativity of 
the functors along the coordinate functions}
\label{subsection;09.10.23.2}

Let $\nbigm$ be good on $(X,D)$.

\begin{lem}
\label{lem;09.10.19.23}
For any $i\neq j$,
we have natural isomorphisms
$\phi_i\bigl(\nbigm(\ast j)\bigr)
\simeq
 \phi_i(\nbigm)(\ast j)$
and $\phi_i\bigl(\nbigm(!j)\bigr)
\simeq
 \phi_i(\nbigm)(!j)$.
\end{lem}
\pf
The second isomorphism
is obtained as the dual of the first one.
Let us consider the first isomorphism.
We have the following naturally defined morphisms:
\[
 \phi_i\bigl( \nbigm(\ast j) \bigr)
\stackrel{a}{\lrarr}
 \phi_i\bigl(\nbigm(\ast j)\bigr)(\ast j)
\stackrel{b}{\llarr}
 \phi_i\bigl(\nbigm\bigr)(\ast j)
\]
Because the restriction of $b$ to $X-D_j$
is an isomorphism,
it is easy to see that $b$ is an isomorphism.
Let us prove that $a$ is an isomorphism
by using an induction on $\rho(\nbigm)$.
As in the proof of Lemma \ref{lem;09.10.19.20},
the issue can be reduced to the case
where $\nbigm$ is a good meromorphic flat bundle,
which is given in Lemma \ref{lem;09.10.19.21}.
\hfill\qed

\begin{lem}
\label{lem;09.10.19.22}
$\nbigm(\ast j)$ and $\nbigm(!j)$
are also good.
\end{lem}
\pf
Because $\phi_j\bigl(\nbigm(\ast j)\bigr)
\simeq \psi_j(\nbigm)$,
we obtain that 
$\nbigm(\ast j)$ is good
from Lemmas
\ref{lem;09.10.19.25},
\ref{lem;09.10.19.20} and
\ref{lem;09.10.19.23}.
By using Lemma \ref{lem;09.10.19.26},
we obtain that $\nbigm(!j)$
is also good.
\hfill\qed

\vspace{.1in}
We have the following corollary of
Lemma {\rm\ref{lem;09.10.19.22}}.
\begin{cor}
Let $f$ be a meromorphic function on $(X,D)$
whose zeros and poles are contained in $D$.
Take $D^{(1)}\subset D$ such that
the poles of $f$ are contained in $D^{(1)}$.
The holonomic $\nbigd_X$-module
$\Pi^{a,b}_{f\star}(\nbigm,\ast D^{(1)})$
is good on $(X,D)$.
Hence,
$\psi^{(a)}_f(\nbigm,\ast D^{(1)})$,
$\Xi^{(a)}_f(\nbigm,\ast D^{(1)})$
and $\phi^{(a)}_f(\nbigm,\ast D^{(1)})$
are also good on $(X,D)$.
\hfill\qed
\end{cor}

We have the following naturally defined morphisms:
\[
 \nbigm(\ast i)(!j)
\stackrel{a}{\lrarr}
 \nbigm(\ast i)(!j)(\ast i)
\stackrel{b}{\llarr}
 \nbigm(!j)(\ast i)
\]
It is easy to prove that $b$ is an isomorphism
for $i\neq j$.

\begin{lem}
\label{lem;09.10.19.31}
$a$ is also an isomorphism,
by which we can identify
$\nbigm(\ast i)(!j)$
and
$\nbigm(!j)(\ast i)$.
\end{lem}
\pf
By using an induction on $\rho(\nbigm)$,
we can reduce the issue
to the case where
$\nbigm$ is a good meromorphic flat bundle,
which is given in Lemma \ref{lem;09.10.19.30}.
\hfill\qed

\vspace{.1in}
In the following,
we will not distinguish
$\nbigm(\ast i)(!j)$
and $\nbigm(!j)(\ast i)$ for $i\neq j$,
which will be denoted by
$\nbigm(\ast i!j)$.
For $I\sqcup J\subset\ellsitabar$,
we have the natural identification
$\nbigm(!I\!\ast\! J)
\simeq
 \nbigm(\ast J!I)$,
which will be used implicitly.

\begin{lem}
We have the commutativity
$\Xi^{(a)}_i\circ\Xi^{(b)}_j=\Xi^{(b)}_j\circ\Xi^{(a)}_i$,
$\psi^{(a)}_i\circ\psi^{(b)}_j=\psi^{(b)}_j\circ\psi^{(a)}_i$ 
and $\phi^{(a)}_i\circ\phi^{(b)}_j=\phi^{(b)}_j\circ\phi^{(a)}_i$.
Moreover,
the functors
$\Xi^{(a)}_i$, $\psi^{(b)}_j$ and $\phi^{(c)}_k$
are mutually commutative,
where $i$, $j$ and $k$ are mutually distinct.
In the following,
we will not care about the order
of these functors for good holonomic $\nbigd$-modules
on $(X,D)$.
\end{lem}
\pf
We obtain the natural identification
$\Pi_{i\star}^{a,b}\circ
 \Pi_{j\star'}^{c,d}
=\Pi_{j\star'}^{c,d}\circ\Pi_{i\star}^{a,b}$
from Lemma \ref{lem;09.10.19.31}.
Then, the claim of the lemma is clear.
\hfill\qed

\subsection{Globalization}
\label{subsection;14.1.14.1}

Let $X$ be a complex manifold
with a normal crossing hypersurface $D$.

\begin{df}
A holonomic $\nbigd_X$-module $\nbigm$
is called good on $(X,D)$
if the following holds:
\begin{itemize}
\item
Let $P$ be any point of $D$.
Let $(U,z_1,\ldots,z_n)$ be a coordinate
neighbourhood around $P$
such that $D\cap U=\bigcup_{i=1}^{\ell}\{z_i=0\}$.
Then, $\nbigm_{|U}$ is good
in the sense of Definition
{\rm\ref{df;09.10.29.1}}.
\hfill\qed
\end{itemize}
\end{df}
\index{good holonomic $\nbigd$-module}

We obtain the following from the results in
\S\ref{subsection;09.10.28.3}--\S\ref{subsection;09.10.23.2}.
\begin{lem}\mbox{}
Let $\nbigm$ be good on $(X,D)$.
\begin{itemize}
\item
The dual $\DDD_X\nbigm$ is also good on $(X,D)$.
\item
Let $D^{(1)}\subset D$ be the union of 
some irreducible components.
Then, $\nbigm(\ast D^{(1)})$ and $\nbigm(!D^{(1)})$
are also good on $(X,D)$.
\item
Let $D^{(i)}\subset D$ $(i=1,2)$ be the unions of 
some irreducible components
such that $\dim D^{(1)}\cap D^{(2)}<\dim X-1$.
We have a natural isomorphism
$\nbigm(\ast D^{(1)})(!D^{(2)})
\simeq
 \nbigm(! D^{(2)})(\ast D^{(1)})$.
\item
Let $f$ be a meromorphic function on $(X,D)$
which is invertible on $X\setminus D$.
Take $D^{(1)}\subset D$ such that
the poles of $f$ are contained in $D^{(1)}$.
Then,
$\psi^{(a)}_f(\nbigm,\ast D^{(1)})$,
$\Xi^{(a)}_f(\nbigm,\ast D^{(1)})$
and $\phi^{(a)}_f(\nbigm,\ast D^{(1)})$
are also good on $(X,D)$.
\hfill\qed
\end{itemize}
\end{lem}

\subsection{A quiver description in the local case}
\label{subsection;14.1.14.2}

We set $X:=\Delta^n$,
$D_i=\{z_i=0\}$
and $D=\bigcup_{i=1}^{\ell}D_i$.
We use the notation introduced
in \S\ref{subsection;09.10.19.5}.
Let 
$\nbigi\subset M(X^{(m)},D^{(m)})/H(X^{(m)})$
be a good set of irregular values
which is preserved by the action 
of the Galois group of 
the ramified covering $X^{(m)}\lrarr X$.

We consider tuples
of $\nbigi$-good meromorphic flat bundles
$V_I$
on $(D_I,\del D_I)$
$(I\subset\ellsitabar)$,
with a tuple of morphisms
\[
\begin{CD}
 \psi^{(1)}_{i}(V_I)
 @>{g_{I,i}}>>
 V_{Ii}
 @>{f_{I,i}}>>
 \psi^{(0)}_i(V_I)
\end{CD}
\]
for $I\subset\ellsitabar$ and $i\in\ellsitabar\setminus I$.
Here $Ii:=I\cup\{i\}$.
We impose the following conditions:
\begin{itemize}
\item
 $f_{I,i}\circ g_{I,i}$ is equal to
 $\var\circ\can:\psi^{(1)}_i(V_I)\lrarr \psi^{(0)}_i(V_I)$.
\item
 For any $I\sqcup\{i\}\sqcup\{j\}\subset\ellsitabar$,
 we have the commutativity
 $\psi^{(0)}_j(f_{I,i})\circ f_{Ii,j}
 =\psi^{(0)}_i(f_{I,j})\circ f_{Ij,i}$,
 $g_{Ii,j}\circ \psi^{(1)}_j(g_{I,i})
 =g_{Ij,i}\circ \psi^{(1)}_i(g_{I,j})$
 and
 $f_{Ij,i}\circ g_{Ii,j}
=\psi^{(0)}_i(g_{I,j})\circ \psi^{(1)}_j(f_{I,i})$.
\end{itemize}
For such 
$\nbigc^{(a)}
=\bigl((V^{(a)}_I),(f^{(a)}_{I,i},g^{(a)}_{I,i})\bigr)$ $(a=1,2)$,
morphisms
$\nbigc^{(1)}
\lrarr
\nbigc^{(2)}$ 
are defined to be a tuple of morphisms 
$\varphi_I:V_I^{(1)}\lrarr V_I^{(2)}$
of meromorphic flat bundles
such that the following diagram is commutative:
\[
 \begin{CD}
 \psi^{(1)}_{i}(V^{(1)}_I)
 @>{g^{(1)}_{I,i}}>>
 V^{(1)}_{Ii}
 @>{f^{(1)}_{I,i}}>>
 \psi^{(0)}_i(V^{(1)}_I) \\
 @V{\psi^{(1)}_{i}(\varphi_I)}VV 
 @V{\varphi_{Ii}}VV 
 @V{\psi^{(0)}_{i}(\varphi_I)}VV \\
  \psi^{(1)}_{i}(V^{(2)}_I)
 @>{g^{(2)}_{I,i}}>>
 V^{(2)}_{Ii}
 @>{f^{(2)}_{I,i}}>>
 \psi^{(0)}_i(V^{(2)}_I) \\
\end{CD}
\]
Let $C(X,D)$ denote the category of 
such objects and morphisms.
(We do not fix $\nbigi$.)

Let $\nbigm$ be a good holonomic 
$\nbigd$-module on $(X,D)$.
Set 
$V_I(\nbigm):=\phi^{(\veczero)}_I(\nbigm)(\ast \del D_I)$
and $V_{\emptyset}(\nbigm):=\nbigm(\ast D)$,
which are naturally equipped with morphisms
\[
 \begin{CD}
  \psi^{(1)}_{i}(V_I(\nbigm))
 @>{g_{I,i}(\nbigm)}>>
 V_{Ii}(\nbigm)
 @>{f_{I,i}(\nbigm)}>>
 \psi^{(0)}_i(V_I(\nbigm)).
 \end{CD}
\]
Thus, we obtain an object
in $C(X,D)$
denoted by $\Phi(\nbigm)$.
The construction gives a functor
$\nbigm:\Hol^{\good}(X,D)\lrarr 
 C(X,D)$.
\begin{prop}
$\Phi$ is an equivalence of categories.
\end{prop}
\pf
Let us construct a quasi-inverse functor
$\Upsilon:C(X,D)\lrarr \Hol^{\good}(X,D)$.
Let $\iota_I:D_I\lrarr X$ denote the inclusion.
For any $I\subset \ellsitabar$,
we set
$\nbigm_I^{(0)}:=\iota_{I\dagger}V_I$.
For $I\subset\ellsitabar$ with $1\not\in I$,
we define $\nbigm^{(1)}_I$
as the gluing of $V_I$ and $V_{I1}$
by $f_{I,1}$ and $g_{I,1}$,
i.e., 
$\nbigm^{(1)}_I$ is the cohomology of
the complex
\[
\begin{CD}
 \iota_{I1\dagger}\psi_1^{(1)}(V_I)
 @>{d_1^{(1)}+g_{I,1}}>>
 \iota_{I\dagger}\Xi_1^{(1)}(V_I)
 \oplus 
 \iota_{I1\dagger}V_{I1}
 @>{c_2^{(0)}-f_{I,1}}>>
 \iota_{I1\dagger}\psi_1^{(0)}(V_I).
\end{CD}
\]
For $I\sqcup\{i\}\subset\ellsitabar\setminus\{1\}$,
we have naturally induced morphisms
\[
\begin{CD}
 \psi^{(1)}_i(\nbigm^{(1)}_I)
 @>{g^{(1)}_{I,i}}>>
 \nbigm^{(1)}_{I1}
 @>{f^{(1)}_{I,i}}>>
 \psi^{(0)}_i(\nbigm^{(1)}_{I})
\end{CD}
\]
Then, 
(i) $f^{(1)}_{I,i}\circ g^{(1)}_{I,i}$ is equal to
 the canonical morphism,
(ii)
 for any $I\sqcup\{i\}\sqcup\{j\}\subset\ellsitabar\setminus\{1\}$,
 we have the commutativity
 $\psi^{(0)}_j(f^{(1)}_{I,i})\circ f^{(1)}_{Ii,j}
 =\psi^{(0)}_i(f^{(1)}_{I,j})\circ f^{(1)}_{Ij,i}$,
 $g^{(1)}_{Ii,j}\circ \psi^{(0)}_j(g^{(1)}_{I,i})
 =g^{(1)}_{Ij,i}\circ \psi^{(0)}_i(g^{(1)}_{I,j})$,
 and
 $f^{(1)}_{Ij,i}\circ g^{(1)}_{Ii,j}
 =\psi^{(0)}_i(g^{(1)}_{I,j})\circ \psi^{(1)}_j(f^{(1)}_{I,i})$.

Inductively on $m$,
we can introduce 
good holonomic $\nbigd$-modules 
$\nbigm^{(m)}_I$ on $(X,D)$
for $I\subset \ellsitabar\setminus\underline{m}$,
and morphisms for $I\sqcup\{i\}\subset \ellsitabar\setminus\underline{m}$
\begin{equation}
\label{eq;14.1.14.2}
\begin{CD}
 \psi^{(1)}_i(\nbigm^{(m)}_I)
 @>{g^{(m)}_{I,i}}>>
 \nbigm^{(m)}_{I1}
 @>{f^{(m)}_{I,i}}>>
 \psi^{(0)}_i(\nbigm^{(m)}_{I})
\end{CD}
\end{equation}
such that 
 $\psi^{(0)}_j(f^{(m)}_{I,i})\circ f^{(m)}_{Ii,j}
 =\psi^{(0)}_i(f^{(m)}_{I,j})\circ f^{(m)}_{Ij,i}$,
 $g^{(m)}_{Ii,j}\circ \psi^{(0)}_j(g^{(m)}_{I,i})
 =g^{(m)}_{Ij,i}\circ \psi^{(0)}_i(g^{(m)}_{I,j})$,
 and
 $f^{(m)}_{Ij,i}\circ g^{(m)}_{Ii,j}
 =\psi^{(0)}_i(g^{(m)}_{I,j})\circ \psi^{(1)}_j(f^{(m)}_{I,i})$.
Indeed, suppose we are given such holonomic $\nbigd$-modules
for $m-1$,
we define $\nbigm_I^{(m)}$
for $I\subset \ellsitabar\setminus{m}$
as the gluing of 
$\nbigm_I^{(m-1)}$ and $\nbigm^{(m-1)}_{Im}$
by $g^{(m-1)}_{I,m}$ and $f^{(m-1)}_{I,m}$.
By the construction,
we have the induced morphisms
as in (\ref{eq;14.1.14.2})
with the desired property.
After the procedure,
we obtain a good holonomic $\nbigd$-module
$\Upsilon\bigl(
 (V_I\,|\,I\subset\ellsitabar),(
 f_{I,i},g_{I,i}\,|\,I\sqcup\{i\}\subset\ellsitabar)
 \bigr):=
 \nbigm^{(\ell)}$.
Clearly, $\Upsilon$ and $\Phi$ are mutually quasi-inverse.
\hfill\qed

\vspace{.1in}
We can describe some functors on $\Hol^{\good}(X,D)$
in terms of $C(X,D)$.
Let $\nbigc=\bigl((V_I),(g_{I,i},f_{I,i})\bigr)$.
For $i$,
we define 
$\nbigc(\ast D_i)
=\bigl((V_I'),(g'_{I,i},f'_{I,i})\bigr)$
as follows.
We set
$V'_{I}:=V_I$ $(i\not\in I)$
or
$V'_{I}:=\psi^{(0)}_i(V_{I\setminus \{i\}})$ $(i\not\in I)$.
If $j\neq i$,
$g'_{I,j}$ and $f'_{I,j}$
are the naturally induced morphisms,
and $g'_{I,i}$ and $f'_{I,i}$
are given by the canonical morphisms
$\psi^{(1)}_{i}(V_{I})
\stackrel{\can}{\lrarr}
 \psi^{(0)}_{i}(V_{I})
\stackrel{\id}{\lrarr}
 \psi^{(0)}_{i}(V_{I})$.
We define
$\nbigc(!D_i)$
as follows.
We set
$V'_{I}:=V_I$ $(i\not\in I)$
or
$V'_{I}:=\psi^{(1)}_i(V_{I\setminus \{i\}})$ $(i\not\in I)$.
If $j\neq i$,
$g'_{I,j}$ and $f'_{I,j}$
are the naturally induced morphisms,
and $g'_{I,i}$ and $f'_{I,i}$
are given by the canonical morphisms
$\psi^{(1)}_{i}(V_{I})
\stackrel{\id}{\lrarr}
 \psi^{(1)}_{i}(V_{I})
\stackrel{\var}{\lrarr}
 \psi^{(0)}_{i}(V_{I})$.
We have naturally defined morphisms
$\nbigc(!D_i)\lrarr \nbigc\lrarr\nbigc(\ast D_i)$.
It is easy to observe 
$\Phi(\nbigm(\star D_i))
\simeq
 \Phi(\nbigm)(\star D_i)$.

We define 
$\psi_{i}^{(a)}(\nbigc)
=\bigl(
 (V_I'),(g'_{I,i},f'_{I,i})
 \bigr)$
as follows.
If $i\not\in I$,
we set $V_I'=0$.
If $i\in I$,
we set
$V_I':=\psi^{(a)}_i(V_{I\setminus i})$.
The morphisms 
$g'_{I,i}$ and $f'_{I,i}$
are the naturally induced ones.
Then, we have a natural isomorphism
$\Phi\psi^{(a)}_i(\nbigm)
\simeq
 \psi^{(a)}_i\Phi(\nbigm)$.

We define 
$\phi_{i}^{(a)}(\nbigc)
=\bigl(
 (V_I'),(g'_{I,i},f'_{I,i})
 \bigr)$
as follows.
If $i\not\in I$,
we set $V_I'=0$.
If $i\in I$,
we set
$V_I':=V_{I}^{(a)}$.
The morphisms 
$g'_{I,i}$ and $f'_{I,i}$
are the naturally induced ones.
Then, we have a natural isomorphism
$\Phi\phi^{(a)}_i(\nbigm)
\simeq
 \phi^{(a)}_i\Phi(\nbigm)$.

We define 
$\DDD(\nbigc)
=\bigl(
 (V_I'),(g'_{I,i},f'_{I,i})
 \bigr)$
as follows.
We set
$V'_I:=\DDD(V_I^{(-1)})(\ast \del D_I)$.
The morphisms 
$g'_{I,i}$ and $f'_{I,i}$
are the naturally induced ones.
Then, we have a natural isomorphism
$\Phi\DDD(\nbigm)
\simeq
 \DDD\Phi(\nbigm)$.

\subsection{Appendix}
\label{subsection;13.5.5.1}

The category $\Hol^{\good}(X,D)$
of good holonomic $\nbigd$-modules on $(X,D)$
is not abelian.
Indeed, a direct sum of good holonomic $\nbigd$-modules
is not necessarily good.
If we would like to work on an abelian category,
it would be convenient to
restrict ourselves to a smaller category.

We generalize the notion of good system of irregular
values in \S2.4.1 of \cite{mochi7}.
For any point $P\in D$,
we introduce some rings.
To define them,
we introduce a category $C_P$.
Objects in $C_P$ are holomorphic maps
$\varphi:(Z,Q)\lrarr (X,P)$
of smooth complex manifolds
which are coverings with ramification along $D$.
We set $D_Z:=\varphi^{-1}(D)$.
Morphisms 
$F:\bigl((Z,Q),\varphi\bigr)
\lrarr
 \bigl((Z',Q'),\varphi'\bigr)$
are holomorphic maps $F:(Z,Q)\lrarr (Z',Q')$
such that $\varphi'\circ F=\varphi$.
Such morphisms induce the morphisms
$\nbigo_{Z'}(\ast D_{Z'})_{Q'}
\lrarr
\nbigo_{Z}(\ast D_Z)_{Q}$
over $\nbigo_X(\ast D)_Q$.
Let $\nbigotilde_{X}(\ast D)_P$
denote a colimit of
$\nbigo_{Z}(\ast D_Z)_Q$.
Similarly,
let $\nbigotilde_{X,P}$
denote the colimit of
$\nbigo_{Z,Q}$.

We have another more direct description.
Let $\cnum\{z_1,\ldots,z_n\}$
denote the ring of convergent power series.
Let $\cnum\{z_1,\ldots,z_n\}_{z_1\ldots z_{\ell}}$
denote its localization with respect to
$z_1\cdots z_{\ell}$.
For a coordinate system $(z_1,\ldots,z_n)$ such that
$D=\bigcup_{i=1}^{\ell}\{z_i=0\}$,
we have natural isomorphisms
\[
\nbigotilde_{X,P}\simeq
\varinjlim_e
\cnum\bigl\{z_1^{1/e},\ldots,z_{\ell}^{1/e},z_{\ell+1},\ldots,z_n\bigr\},
\]
\[
 \nbigotilde_{X}(\ast D)_P
 \simeq
\varinjlim_e
\cnum\bigl\{z_1^{1/e},\ldots,z_{\ell}^{1/e},z_{\ell+1},\ldots,z_n\bigr\}
 _{z_1^{1/e}\cdots z_{\ell}^{1/e}}.
\]

A finite subset 
$\nbigi\subset
 \nbigotilde_X(\ast D)_P\big/\nbigotilde_{X,P}$
can be regarded as
$\nbigi\subset\nbigo_{Z}(\ast D_Z)_Q
 \big/\nbigo_{Z,Q}$
for some $\bigl((Z,Q),\varphi\bigr)\in C_P$.
It is called a good set of ramified irregular values
if (i) it is a good set of irregular values on $(Z,D_Z)$,
(ii) it is stable under the action of 
the Galois group of $\varphi$.
\index{good set of ramified irregular values}
Note that if $P_1$ is close to $P$,
 we choose $Q_1\in \varphi^{-1}(P_1)$,
 and we obtain a natural map
 $\nbigi_P\lrarr
 \nbigo_{Z}(\ast D_Z)_{Q_1}
 \big/\nbigo_{Z,Q_1}
\lrarr 
 \nbigotilde_{X}(\ast D)_{P_1}
 \big/\nbigotilde_{X,P_1}$.
The image is well defined.

\begin{df}
A good system of ramified irregular values
on $(X,D)$ is a family of good sets of 
ramified irregular values
$\vecnbigi=\bigl\{\nbigi_P\,\big|\,P\in D
 \bigr\}$ satisfying the following condition.
\begin{itemize}
\item
 If $P_1$ is sufficiently close to $P$,
 we impose that the image of
 $\nbigi_P$ in 
the image of $\nbigi_P$
 in 
$\nbigotilde_{X}(\ast D)_{P_1}
 \big/\nbigotilde_{X,P_1}$
 is equal to
 $\nbigi_{P_1}$.
\hfill\qed
\end{itemize}
\end{df}
\index{good system of ramified irregular values}

Let $\vecnbigi=(\nbigi_P\,|\,P\in D)$
be a good system of
ramified irregular values on $(X,D)$.
A holonomic $\nbigd_X$-module $\nbigm$
is called $\nbigi$-good
if for any $P\in D$
there exists a neighborhood $X_P$
such that
$\nbigm_{|X_P}$ is $\nbigi_P$-good.
Then, the category of
$\vecnbigi$-good holonomic $\nbigd$-modules on $(X,D)$
is an abelian full subcategory of
$\Hol(X)$.

%% file: 3.2.tex
\subsection{De Rham complex with
infinite decay}
\label{subsection;10.1.11.2}

For any complex manifold $X$,
let $\Omega_X^{p,q}$
denote the sheaf of $C^{\infty}$-$(p,q)$-forms
on $X$.
\index{sheaf $\Omega_X^{p,q}$}
We set $d_X:=\dim X$.
For any analytic subset $Z\subset X$,
we set $\Omega_{\Zhat}^{p,q}:=
 \Omega_X^{p,q}\otimes_{\nbigc^{\infty}_X}
 \nbigc^{\infty}_{\Zhat}$.
\index{sheaf $\Omega^{p,q}_{\Zhat}$}
For any hypersurface 
$D\subset X$,
we set
$\Omega_{\Zhat}^{p,q}(\ast D):=
 \Omega_{\Zhat}^{p,q}\otimes_{\nbigo_X}
 \nbigo_{X}(\ast D)$.
\index{sheaf $\Omega_{\Zhat}^{p,q}(\ast D)$}
We say that $D_1\cup D_2=D$ 
is a decomposition of $D$
if $D_i\subset X$ $(i=1,2)$ are hypersurfaces
such that $\codim_X(D_1\cap D_2)>1$.
\index{decomposition}
In that situation,
we say that $D_2$ is the complement of $D_1$
in $D$.
\index{complement}
In other words,
the complement of $D_1$ in $D$ is the union
of the irreducible components of $D$
which are not contained in $D_1$.
When we are given a hypersurface 
$D\subset X$ with a decomposition
$D=D_1\cup D_2$,
let $\Omega_X^{p,q}(\ast D_2)^{<D_1}$
denote the kernel of
$\Omega_X^{p,q}(\ast D_2)
\lrarr \Omega_{\Dhat_1}^{p,q}(\ast D_2)$.
\index{sheaf $\Omega_X^{p,q}(\ast D_2)^{<D_1}$}

Let $D_0$ be a normal crossing hypersurface
of $X$ with a decomposition $D_0=D_1\cup D_2$.
For any coherent $\nbigd_X$-module $\nbigm$,
we define 
$\DR^{<D_1\leq D_2}_X\nbigm$
as
\[
 Cone\bigl(
 \DR_{X}\bigl(\nbigm(\ast D_2)\bigr)
\lrarr
 \DR_{\Dhat_1}\bigl(\nbigm(\ast D_2)\bigr)
\bigr)[-1]
\]
in the derived category $D^b(\cnum_X)$.
We have the following natural quasi-isomorphisms:
\[
 \DR^{<D_1\leq D_2}_{X}\nbigm
\simeq
 \Omega_{X}^{d_X,\bullet\,<D_1}(\ast D_2)
 \otimes_{\nbigd_X}^L\nbigm
\simeq
 \Tot
 \Omega_X^{\bullet,\bullet\,<D_1}(\ast D_2)
\otimes_{\nbigo_X}\nbigm[d_X]
\]
\index{functor $\DR^{<D_1\leq D_2}_X$}
Here, $\Tot$ means the total complex
associated to the double complex.
In the following,
we shall often omit to denote $\Tot$.
It is easy to observe that
the natural morphism
$\DR^{<D_1\leq D_2}_X\nbigm
\lrarr
 \DR^{<D_1\leq D_2}_X\bigl(
 \nbigm(\ast D_0)\bigr)$ 
is an isomorphism.
We also have the following natural isomorphisms
in $D^b(\cnum_X)$:
\begin{multline*}
 \DR_X^{<D_1}\bigl(
 \DDD_X\nbigm(\ast D_0)\bigr)
\simeq
 \Omega_{X}^{d_X,\bullet}(\ast D_2)^{<D_1}
\otimes^L_{\nbigd_X}\DDD_X\nbigm(\ast D_0) \\
\simeq 
 \nrhom_{\nbigd_X}\bigl(
 \nbigm,\,\Omega_X^{0,\bullet}(\ast D_2)^{<D_1}
 \bigr)[d_X]
\end{multline*}

The following proposition
is an immediate consequence of
the isomorphism of Mebkhout
recalled in Proposition \ref{prop;14.1.16.1}.
\begin{prop}
\label{prop;14.1.16.2}
If $\bigl(\nbigm(\ast D_2)\bigr) (!D_1)
\simeq\nbigm(\ast D_2)$,
the natural morphism
\[
\DR_X^{<D_1\leq D_2}(\nbigm)
\lrarr
 \DR_X^{\leq D_2}(\nbigm)
\]
is an isomorphism
in $D_c^b(\cnum_X)$.
\hfill\qed
\end{prop}

\subsection{The identification
in the case of good holonomic $\nbigd$-modules}
\label{subsection;09.10.23.3}

Let $X$ be a complex manifold
with a normal crossing hypersurface $D$.
Let $D_0\subset D$ be the union of some irreducible components
with a decomposition $D_0=D_1\cup D_2$.
Let $\nbigm$ be a good holonomic $\nbigd$-module
on $(X,D)$.
The following proposition is a special case of
Proposition \ref{prop;14.1.16.2}.
\footnote{The author thanks the referee
for the simplified proof of the proposition.}
\begin{prop}
If $\nbigm(!D_1)=\nbigm$,
the natural morphism
$\DR_X^{<D_1\leq D_2}\nbigm
\lrarr
 \DR_X^{\leq D_2}\nbigm$
is a quasi-isomorphism.
\hfill\qed
\end{prop}

We obtain the following isomorphisms
in $D^b_c(\cnum_X)$:
\begin{equation}
\label{eq;13.4.17.2}
 \begin{CD}
 \DR_X^{<D_1\leq D_2}(\nbigm)
@<{\simeq}<<
 \DR_X^{<D_1\leq D_2}(\nbigm(!D_1))
@>{\simeq}>>
 \DR_X^{\leq D_2}(\nbigm(!D_1))
 \end{CD}
\end{equation}
We have already seen the right isomorphism.
For the left isomorphism,
we may use
$\Omega^{p,q\,<D_1}_{X}
\simeq
 \Omega^{p,q\,<D_1}_X(\ast D_1)$.
We will identify 
$\DR_X^{<D_1\leq D_2}(\nbigm)$
and
$\DR_X^{\leq D_2}(\nbigm(!D_1))$
by (\ref{eq;13.4.17.2}).

\begin{lem}
\label{lem;13.4.17.1}
If $D_1\subset D_1'\subset D$,
then the following diagram of the natural morphisms
is commutative:
\[
 \begin{CD}
  \DR_X^{<D_1'}\nbigm
  @>{\simeq}>>
 \DR_X\nbigm(!D_1')
 \\
 @VVV @VVV \\
  \DR_X^{<D_1}\nbigm
 @>{\simeq}>>
 \DR_X\nbigm(!D_1)
 \end{CD}
\]
It is also factorized as follows:
\[
 \begin{CD}
  \DR_X^{<D_1'}\nbigm
 @<{\simeq}<<
 \DR^{<D_1}\nbigm(!D_1')
 @>{\simeq}>>
 \DR_X\nbigm(!D_1')
 \\
 @VVV @VVV @VVV \\
  \DR_X^{<D_1}\nbigm
 @<{\simeq}<<
 \DR^{<D_1}\nbigm
 @>{\simeq}>>
 \DR_X\nbigm(!D_1)
 \end{CD}
\]
\end{lem}
\pf
We have the following commutative diagram:
\[
 \begin{CD}
 \DR_X^{<D_1'}\nbigm
 @<{\simeq}<<
 \DR^{<D_1'}\nbigm(!D_1')
 @>{\simeq}>>
 \DR_X\nbigm(!D_1')
 \\
 @VVV @VVV @VVV \\
  \DR_X^{<D_1}\nbigm
 @<{\simeq}<<
 \DR^{<D_1}\nbigm(!D_1)
 @>{\simeq}>>
 \DR_X\nbigm(!D_1)
 \end{CD}
\]
Then, the claim of the lemma is clear.
\hfill\qed

\subsection{Duality}
\label{subsection;13.4.17.30}

We continue to use the notation in
\S\ref{subsection;09.10.23.3}.
For simplicity, we assume $D=D_0$.
We have a morphism of complexes
\begin{equation}
 \label{eq;13.4.17.10}
 \Tot
 \Bigl(
 \Tot\Omega^{\bullet,\bullet\,<D_2}(\ast D_1)[d_X]
 \otimes
 \Tot\Omega^{0,\bullet\,<D_1}(\ast D_2)[d_X]
 \Bigr)
\lrarr
 \Tot\Omega^{\bullet,\bullet}[2d_X]
\end{equation}
by
$\xi\otimes\eta
\!\longmapsto\!
(-1)^{pd_X}
 \xi\wedge\eta$,
where
$\xi$ and $\eta$
are local sections of
$\bigl(
 \Tot\Omega^{\bullet,\bullet\,<D_2}(\ast D_1)
 \bigr)^{p+d_X}$
and 
$\bigl(
 \Tot\Omega^{0,\bullet\,<D_1}(\ast D_2)
 \bigr)^{q+d_X}$
respectively.
Let $\nbigi^{\bullet}_1$ be a $\nbigd_X$-injective resolution
of $\Tot\Omega^{0,\bullet\,<D_1}(\ast D_2)[d_X]$,
and let $\nbigi^{\bullet}_2$ be a $\cnum_X$-injective resolution
of $\Tot\Omega^{\bullet,\bullet}[2d_X]$.
Then, the morphism is extended to a $\cnum_X$-homomorphism
$\DR^{\leq D_1<D_2}_X(\nbigi_1^{\bullet})
\lrarr\nbigi_2^{\bullet}$.

For any coherent $\nbigd_X$-module $\nbigm$,
we have the following natural morphism:
\begin{equation}
\label{eq;13.4.17.20}
 \DR_X^{<D_1\leq D_2}
 (\DDD_X\nbigm)
\lrarr
 \DDD_X\DR_X^{<D_2\leq D_1}(\nbigm).
\end{equation}
Indeed, 
$\DR_X^{<D_1\leq D_2}\DDD_X\nbigm$
is represented by
$\nhom_{\nbigd_X}\bigl(
 \nbigm,\nbigi^{\bullet}_1
 \bigr)$.
Hence, we have the desired morphism
given as follows:
\begin{multline*}
 \nhom_{\nbigd_X}\bigl(
 \nbigm,\nbigi^{\bullet}_1
 \bigr)
\lrarr
 \nhom_{\cnum_X}\bigl(
 \DR_X^{<D_2\leq D_1}\nbigm,
 \DR_X^{<D_2\leq D_1}\nbigi_1^{\bullet}
 \bigr) \\
\lrarr
 \nhom_{\cnum_X}\bigl(
  \DR_X^{<D_2\leq D_1}\nbigm,
 \nbigi_2^{\bullet}
 \bigr)
\end{multline*}

\begin{thm}
\label{thm;09.10.3.50}
Let $V$ be a good meromorphic flat bundle on $(X,D)$.
The following diagram is commutative:
\begin{equation}
\label{eq;09.10.23.11}
\begin{CD}
 \DR^{<D_1\leq D_2}(V^{\lor})
 @>{G_1}>>
 \DDD_X \DR^{<D_2\leq D_1}(V)\\
 @V{\simeq}VV @A{\simeq}AA \\
 \DR V^{\lor}(!D_1)
 @>{G_2}>{\simeq}>
 \DDD_X\DR_X\bigl(V(!D_2)\bigr)
\end{CD}
\end{equation}
Here,
 $G_1$ is induced by {\rm(\ref{eq;13.4.17.20})}
and 
 $\DR_X^{<D_1\leq D_2}(\DDD_X V)
\simeq
 \DR_X^{<D_1\leq D_2}(V^{\lor})$.
The vertical isomorphisms are given
by {\rm(\ref{eq;13.4.17.2})},
and 
$G_2$ is induced by the natural isomorphism
of $\nbigd$-modules
$V^{\lor}(!D_1)\simeq
 \DDD_X\bigl(V(!D_2)\bigr)$.
(See \S{\rm\ref{subsection;09.10.23.2}}.)
In particular,
$G_1$ is also an isomorphism.
\end{thm}
\pf
We have the commutativity of the following natural morphisms:
{\footnotesize
\[
 \begin{CD}
 \DR^{<D_1\leq D_2}_X(V^{\lor})
 @>{\simeq}>>
  \DR_X^{<D_1\leq D_2}(\DDD_X V)
 @>>>
 \DDD_X\DR_X^{<D_2\leq D_1}(V)\\
 @A{\simeq}AA @V{\simeq}VV @V{\simeq}VV
 \\
  \DR_X^{<D_1\leq D_2}(V^{\lor}(!D_1))
 @>{\simeq}>>
  \DR_X^{<D_1\leq D_2}\bigl(\DDD_X (V(!D_2))\bigr)
 @>>>
 \DDD_X\DR_X^{<D_2\leq D_1}(V(!D_2))\\
 @VVV @VVV @VVV \\
 \DR_X(V^{\lor}(!D_1))
 @>{\simeq}>>
  \DR_X\bigl(\DDD_X (V(!D_2))\bigr)
 @>>>
 \DDD_X\DR_X(V(!D_2))
 \end{CD}
\]
}
Then, the claim of the theorem is clear.
\hfill\qed

\subsection{Functoriality for birational morphisms}

Let $X$ be a complex manifold,
and let $D$ be a normal crossing hypersurface 
with a decomposition $D=D_1\cup D_2$.
Let $D_3$ be a hypersurface of $X$.
Let $\varphi:X'\lrarr X$ be a proper birational
morphism such that
(i) $D'=\varphi^{-1}\bigl(D\cup D_3\bigr)$
is normal crossing,
(ii) $X'\setminus D'\simeq 
 X\setminus\bigl(D\cup D_3\bigr)$.
We put $D_1':=\varphi^{-1}(D_1)$.
Let $D_2'$ be the complement of $D_1'$
in $D'$.

Let $\nbigm'$ be any coherent $\nbigd_{X'}$-module
having a good filtration in the neighbourhood of
fibers of $\varphi$.
We have the following natural morphism:
\begin{equation}
 \label{eq;13.4.17.30}
 \DR_X^{<D_1\leq D_2}
 \varphi_{\dagger}\nbigm'
\lrarr
 R\varphi_{\ast}\DR^{<D_1'\leq D_2'}_{X'}\nbigm'
\end{equation}
Indeed,
we have the following:
\begin{multline}
 \DR^{<D_1\leq D_2}_X
 \varphi_{\dagger}\nbigm'
\simeq
 R\varphi_{\ast}\Bigl(
 \Omega_{X'}\otimes_{\varphi^{-1}\nbigo_{X}}
 \varphi^{-1}(\Omega_{X}^{0,\bullet<D_1}(\ast D_2))
 \otimes^L_{\nbigd_{X'}}\nbigm'
 \Bigr) 
 \\
\lrarr
 R\varphi_{\ast}\Bigl(
 \bigl(
 \Omega_{X'}\otimes
 \Omega_{X'}^{0,\bullet<D_1'}(\ast D_2')
\bigr)
 \otimes^L_{\nbigd_{X'}}\nbigm'
 \Bigr)
\simeq
 R\varphi_{\ast}\bigl(
 \DR_{X'}^{<D_1'\leq D_2'}
 (\nbigm')
 \bigr)
\end{multline}

Let $V$ be a good meromorphic flat bundle on $(X,D)$,
and we set 
$V':=\varphi^{\ast}V\otimes\nbigo_{X'}(\ast D')$.
We have a natural isomorphism
$\bigl(V(\ast D_3)\bigr)(!D_1)
 \simeq
 \varphi_{\dagger}\bigl(V'(!D_1')\bigr)$.
Hence,
we have a morphism of $\nbigd_X$-modules
$V(!D_1)\lrarr \varphi_{\dagger}\bigl(
 V'(!D_1')\bigr)$.
We obtain the following morphism
from (\ref{eq;13.4.17.30})
and $V\lrarr \varphi_{\dagger}V'$:
\begin{equation}
\label{eq;09.10.3.151}
\DR_X^{<D_1\leq D_2}(V)\lrarr 
R\varphi_{\ast}\DR_{X'}^{<D_1'\leq D_2'}(V')
\end{equation}
It is equal to the one induced by
$\varphi^{-1}
 \bigl(
 \Omega_{X}^{\bullet,\bullet\,<D_1}(\ast D_2)
\otimes V
 \bigr)
\lrarr
 \Omega_{X'}^{\bullet,\bullet\,<D_1'}(\ast D_2')
\otimes V'$.
Note that we have natural isomorphisms
\begin{multline}
 (\Omega_{X'}\otimes V')
\otimes^L_{\nbigd_{X'}}
 \bigl(\nbigo_{X'}\otimes_{\varphi^{-1}\nbigo_X}
 \varphi^{-1}(\nbigd_X\otimes\Omega_X^{-1})
 \bigr)
\simeq \\
 (\Omega_{X'}\otimes V')
\otimes^L_{\nbigd_{X'}(\ast D')}
 \bigl(\nbigo_{X'}(\ast D')\otimes_{\varphi^{-1}\nbigo_X}
 \varphi^{-1}(\nbigd_X\otimes\Omega_X^{-1})
 \bigr) 
\simeq \\
 (\Omega_{X'}\otimes V')\otimes^L_{\nbigd_{X'}(\ast D')}
 \bigl(\nbigd_{X'}(\ast D')
 \otimes_{\varphi^{-1}\nbigo_X}
 \varphi^{-1}\Omega_X\bigr)
\simeq
 V'
\end{multline}

By considering the dual with $V^{\lor}$
(see Theorem \ref{thm;09.10.3.50}),
we also obtain the following morphism:
\begin{equation}
\label{eq;09.10.3.220}
 R\varphi_{\ast}\DR_{X'}^{<D_2'\leq D_1'}(V')
\lrarr
 \DR_X^{<D_2\leq D_1}(V)
\end{equation}

\begin{thm}
\label{thm;09.10.3.55}
We have the following commutative diagram:
\begin{equation}
\label{eq;09.10.3.231}
 \begin{CD}
 \DR_X^{<D_1\leq D_2}V
 @>>>
 R\varphi_{\ast}
 \DR_{X'}^{<D_1'\leq D_2'}V'\\
 @V{\simeq}VV @V{\simeq}VV \\
 \DR_XV(!D_1)
 @>>>
 R\varphi_{\ast}\DR_{X'}V'(!D_1')
 \end{CD}
\end{equation}
Here, the vertical isomorphisms
are given in {\rm(\ref{eq;13.4.17.2})},
the upper horizontal arrow is
{\rm (\ref{eq;09.10.3.151})},
and the lower horizontal arrow is
induced by the morphism of 
$\nbigd_X$-modules
$V(!D_1)\lrarr
 \varphi_{\dagger}\bigl(
 V'(!D_1')
 \bigr)$.

Similarly,
we have the following commutative diagram:
\begin{equation}
\label{eq;09.10.3.230}
 \begin{CD}
 R\varphi_{\ast}
 \DR_{X'}^{<D_2'\leq D_1'}V'
@>>>
 \DR^{<D_2\leq D_1}_XV  \\
 @V{\simeq}VV @V{\simeq}VV \\
 R\varphi_{\ast}\DR_{X'}V'(!D_2')
@>>>
 \DR_XV(!D_2)
 \end{CD}
\end{equation}
Here, the vertical isomorphisms
are given in {\rm(\ref{eq;13.4.17.2})},
the upper horizontal arrow is
{\rm (\ref{eq;09.10.3.220})},
and the lower horizontal arrow is
induced by the natural morphism of
$\nbigd_X$-modules
$\varphi_{\dagger}\bigl(
 V'(!D_2')\bigr)
\lrarr
 V(!D_2)$.
\end{thm}
\pf
We have the following commutative diagram:
{\footnotesize
\[
 \begin{CD}
 \DR^{<D_1\leq D_2}_X(V)
 @>>>
 \DR^{<D_1\leq D_2}_X(\varphi_{\dagger}V')
 @>>>
 R\varphi_{\ast}\DR_X^{<D_1'\leq D_2'}V'
 \\
 @A{\simeq}AA @AAA @A{\simeq}AA \\
 \DR^{<D_1\leq D_2}_X(V(!D_1))
 @>>>
 \DR^{<D_1\leq D_2}_X(\varphi_{\dagger}V'(!D_1'))
 @>>>
 R\varphi_{\ast}\DR_X^{<D_1'\leq D_2'}V'(!D_1')
 \\
 @VVV @VVV @VVV \\
 \DR_XV(!D_1)
 @>>>
 \DR_X\varphi_{\dagger}V'(!D_1')
 @>>>
 R\varphi_{\ast}\DR_XV'(!D_1')
 \end{CD}
\]
}
Then, we obtain the commutativity of
(\ref{eq;09.10.3.231}).

\vspace{.1in}

Let us consider the commutativity of
(\ref{eq;09.10.3.230}).
Recall the commutativity of (\ref{eq;09.10.3.20}).
We have the following commutative diagram
for $\nbign\lrarr\varphi_{\dagger}\nbign'$,
where $\nbign$ (resp. $\nbign'$)
is a coherent $\nbigd_X$-module
(resp. $\nbigd_{X'}$-module):
\[
 \begin{array}{cccccccccccc}
 R\varphi_{\ast}\DR_X\DDD\nbign'
 & \simeq &
  \DR \varphi_{\dagger}\DDD\nbign' 
 & \simeq &
  \DR \DDD\varphi_{\dagger}\nbign' 
  & \lrarr &
  \DR \DDD\nbign \\
 \darr & & & 
 & \darr & & \darr \\
 R\varphi_{\ast}\DDD \DR_X\nbign'
 & \simeq &
  \DDD R\varphi_{\ast}\DR\nbign'
 & \simeq &
 \DDD\DR\varphi_{\dagger}\nbign'
 & \lrarr &
 \DDD \DR\nbign
 \end{array}
\]
The vertical arrows are also isomorphisms.
Hence,
the lower horizontal arrow in 
(\ref{eq;09.10.3.230})
is obtained as the dual of 
$\DR_X V^{\lor}(D_1)\lrarr
 R\varphi_{\ast}\DR_{X'}V^{\prime\lor}(!D_1')$
in $D_c^b(\cnum_X)$.
Then, the commutativity of (\ref{eq;09.10.3.230})
follows from the commutativity of
(\ref{eq;09.10.3.231}).
Thus, the proof of 
Theorem \ref{thm;09.10.3.55}
is finished.
\hfill\qed

%% file: 4.1.tex
We shall introduce the sheaves of
holomorphic functions of various types.
We give some statements mainly on flatness.
The proof will be given later.

\subsection{Preliminary}
\label{subsection;14.1.17.30}

Let $X$ be an $n$-dimensional complex manifold
with a simply normal crossing hypersurface
$D$ with the irreducible decomposition
$\bigcup_{i\in\Lambda} D_i$.
In this paper,
the real blow up $\pi:\Xtilde(D)\lrarr X$
means the fiber product of $\Xtilde(D_i)$
over $X$.
\index{real blow up $\Xtilde(D)$}
For any subset $I\subset\Lambda$,
we set $D_I:=\bigcap_{i\in I}D_i$
and $D(I):=\bigcup_{i\in I}D_i$.
Formally, $D_{\emptyset}:=X$.
For $J\subset I^c:=\Lambda\setminus I$,
we put
$D_I(J):=D_I\cap D(J)$.
In particular,
$\del D_I:=D_I(I^c)$.

\subsection{Holomorphic functions
with moderate growth or rapid decay}

Recall that holomorphic functions
on an open subset $U\subset \Xtilde(D)$
are defined to be 
$C^{\infty}$-functions on $U$
whose restriction to $U\setminus\pi^{-1}(D)$
are holomorphic.
A holomorphic function $f$ on $U$
is called of rapid decay
if the following holds:
\index{rapid decay}
\begin{itemize}
\item
Let $P$ be any point of $\pi^{-1}(D)\cap U$.
We take a holomorphic coordinate system
$(z_1,\ldots,z_n)$ around $\pi(P)$
such that $D=\bigcup_{i=1}^{\ell}\{z_i=0\}$.
Then,
we have $f=O\Bigl(\prod_{i=1}^{\ell}|z_i|^{N}\Bigr)$
for any $N$ around $P$.
\end{itemize}
In this paper,
the sheaf of holomorphic functions on $\Xtilde(D)$
is denoted by $\nbigo_{\Xtilde(D)}$.
The sheaf of holomorphic functions with rapid decay
is denoted by $\nbiga^{\rapid}_{\Xtilde(D)}$.
\index{sheaf $\nbigo_{\Xtilde(D)}$}
\index{sheaf $\nbiga^{\rapid}(\Xtilde(D))$}

Let $U$ be any open subset in $\Xtilde(D)$.
A holomorphic function $f$ on $U\setminus \pi^{-1}(D)$
is called of moderate growth
if the following holds:
\index{moderate growth}
\begin{itemize}
\item
Let $P$ be any point of
 $\pi^{-1}(P)\cap U\neq\emptyset$.
We take a holomorphic 
coordinate system $(z_1,\ldots,z_n)$ around $\pi(P)$
such that $D=\bigcup_{i=1}^{\ell}\{z_i=0\}$.
Then, we have $f=O\Bigl(\prod_{i=1}^{\ell}|z_i|^{-N}\Bigr)$
for some $N$ around $P$.
\end{itemize}
In this paper,
the sheaf of holomorphic functions with moderate growth
is denoted by $\nbiga^{\moderate}_{\Xtilde(D)}$.
\index{sheaf $\nbiga^{\moderate}_{\Xtilde(D)}$}
We shall prove the following
(Proposition \ref{prop;09.10.25.32},
Theorem \ref{thm;12.9.18.10}).
\begin{thm}
\label{thm;14.1.22.1}
The sheaves
$\nbigo_{\Xtilde(D)}$,
$\nbiga^{\rapid}_{\Xtilde(D)}$
and $\nbiga^{\moderate}_{\Xtilde(D)}$
are flat over $\pi^{-1}(\nbigo_X)$.
\end{thm}

\subsection{Partially rapid decay functions on completions}

Suppose that $Z$ is $\pi^{-1}\bigl(D_I(J)\bigr)$
for some $I\sqcup J\subset\Lambda$.
Let $\nbigi_Z\subset\nbigo_{\Xtilde(D)}$
be the ideal sheaf of $Z$,
and put $\nbigo_{\Zhat}:=
 \varprojlim\nbigo_X\big/\nbigi_Z^m$.
\index{sheaf $\nbigo_{\Zhat}$}
For a given $\nbigo_{\Xtilde(D)}$-module $\nbigf$,
we set $\nbigf_{|\Zhat}:=
 \nbigf\otimes_{\nbigo_{\Xtilde(D)}}
 \nbigo_{\Zhat}$.
\index{sheaf $\nbigf_{|\Zhat}$}
According to a generalized Borel-Ritt theorem
due to Majima and Sabbah
(\cite{majima}, 
Proposition II.1.1.16 of \cite{sabbah4}),
the natural morphism
$\nbigo_{\widehat{\pi^{-1}(D_I)}}\lrarr
 \nbigo_{\widehat{\pi^{-1}(D_I(J))}}$ 
is surjective.
The kernel is denoted by 
$\nbigo_{\widehat{\pi^{-1}(D_I)}}^{<D(J)}$.
\index{sheaf $\nbigo_{\widehat{\pi^{-1}(D_I)}}^{<D(J)}$}
If $D_I=X$ and $D(J)=D$,
it is equal to $\nbiga^{\rapid}_{\Xtilde(D)}$.
We shall prove the following theorem.
(See Proposition \ref{prop;09.10.25.32}
for a refined claim.)

\begin{prop}
The sheaves
$\nbigo_{\widehat{\pi^{-1}(D_I)}}^{<D(J)}$
and 
$\nbigo_{\widehat{\pi^{-1}D_I(J)}}$
are flat over $\pi^{-1}(\nbigo_X)$.
\end{prop}

\subsection{Holomorphic functions of Nilsson type}

\subsubsection{Preliminary}

We set
$\Nil(z):=\bigoplus_{\alpha\in\cnum}
 z^{\alpha}\cnum[\log z]$.
\index{set $\Nil(z)$}
For $(\alpha,k)\in\cnum\times\seisuu_{\geq \,0}$,
we put $\varphi_{\alpha,k}(z):=z^{\alpha}(\log z)^k
 \in\Nil(z)$.
\index{function $\varphi_{\alpha,k}$}
Let $T$ be any finite subset 
contained in $\bigl\{
 \alpha\in\cnum\,\big|\,0\leq\Re(\alpha)<1
 \bigr\}$.
For simplicity, we assume $0\in T$.
Let $N$ be a non-negative integer.
We set
\[
 \Nil_{T,N}(z):=
 \Bigl\{
 \sum a_{\alpha,j,k}\,
 \varphi_{\alpha+j,k}(z)
 \in\Nil(z)
 \,\Big|\,
 a_{\alpha,j,k}\in\cnum,\,\,
 j\geq -N,\,\,k\leq N,\,\,
 \alpha\in T
 \Bigr\}.
\]
\index{set $\Nil_{T,N}(z)$}
Note that
$\Nil_{T,N}(z)$ is a finitely 
generated free $\cnum[z]$-module.
For $T\subset T'$ and $N\leq N'$,
we have a natural inclusion
$\Nil_{T,N}(z)\subset
 \Nil_{T',N'}(z)$.
We have 
$\Nil(z)=\varinjlim\Nil_{T,N}(z)$.

Let $\cnumtilde_z$ be the real blow up
of $\cnum_z$ along $0$.
Let $\iota$ be the inclusion
$\iota:\cnum_z^{\ast}\lrarr\cnumtilde_z$.
We have the subsheaves of
$\iota_{\ast}\nbigo_{\cnum^{\ast}}$
on $\cnumtilde$
corresponding to 
$\Nil(z)$ and $\Nil_{T,N}(z)$.
The sheaves are also denoted by
$\Nil(z)$ and $\Nil_{T,N}(z)$.

For $\ell\geq 1$,
put
$\Nil(z_1,\ldots,z_{\ell}):=
 \Nil(z_1)\otimes_{\cnum}\cdots\otimes_{\cnum}
 \Nil(z_{\ell})$
and
$\Nil_{T,N}(z_1,\ldots,z_{\ell}):=
 \Nil_{T,N}(z_1)\otimes_{\cnum}
 \cdots\otimes_{\cnum}
 \Nil_{T,N}(z_{\ell})$.
\index{set $\Nil_{T,N}(z_1,\ldots,z_{\ell})$}
\index{set $\Nil(z_1,\ldots,z_{\ell})$}
We naturally regard
$\Nil(z_1,\ldots,z_{\ell})$
as a subsheaf of
$\iota_{\ast}\nbigo_{\cnum^n-D}$
on the real blow up $\cnumtilde(D)$,
where
$D=\bigcup_{i=1}^{\ell}\{z_i=0\}$
and $\iota:\cnum^n-D\lrarr\cnumtilde^n(D)$.
For $(\vecalpha,\veck)\in
 \cnum^{\ell}\times\seisuu_{\geq 0}^{\ell}$,
we put 
$\varphi_{\vecalpha,\veck}(z_1,\ldots,z_n):=
 \prod_{i=1}^{\ell}
 \varphi_{\alpha_i,k_i}(z_i)$,
which are regarded as multi-valued flat sections
of $\Nil(z_1,\ldots,z_{\ell})$.
\index{function $\varphi_{\vecalpha,\veck}$}

\subsubsection{Holomorphic functions of Nilsson type}

Let $X$ be an $n$-dimensional complex manifold
with a simply normal crossing hypersurface $D$.
Let $D=D^{(1)}\cup D^{(2)}$ be a decomposition.
We shall introduce a sheaf
$\nbiga^{<D^{(1)}\leq D^{(2)}}_{\Xtilde(D)}$
on $\Xtilde(D)$.
First, let us consider the case $X=\Delta^n$,
$D=\bigcup_{i=1}^{\ell}\{z_i=0\}$.
Let $\ellsitabar=I_1\sqcup I_2$
be determined by
$D^{(j)}=\bigcup_{i\in I_j}\{z_i=0\}$
for $j=1,2$.
Let $\jtilde$ denote the inclusion
$X-D\lrarr\Xtilde(D)$.
Let $\nbiga^{<D^{(1)}\leq D^{(2)}}_{\Xtilde(D)}$
be the image of the naturally defined morphisms:
\index{sheaf $\nbiga^{<D^{(1)}\leq D^{(2)}}_{\Xtilde(D)}$}
\[
 \nbigo_{\Xtilde(D)}^{<D^{(1)}}
 \otimes\Nil(z_i\,|\,i\in I_2)
\lrarr
 \jtilde_{\ast}\nbigo_{X-D}.
\]
We can observe that
they are independent of the choice of a coordinate system
$(z_1,\ldots,z_n)$.
Hence, we obtain globally defined sheaf
$\nbiga^{<D^{(1)}\leq D^{(2)}}_{\Xtilde(D)}$
on $\Xtilde(D)$.
It is also denoted by
$\nbiga^{\nil\,<D^{(1)}}_{\Xtilde(D)}$.
We shall prove the following.
(See Theorem \ref{thm;09.12.4.5}
and Corollary \ref{cor;14.1.19.3}
for refined claims.)
\begin{thm}
$\nbiga^{<D^{(1)}\leq D^{(2)}}$
is flat over $\pi^{-1}\nbigo_X$.
We also have
$R\pi_{\ast}\nbiga^{\nil}_{\Xtilde(D)}
\simeq
 \nbigo_X(\ast D)$.
\end{thm}

\begin{rem}
This type of sheaves are useful
when we study the de Rham complex of
$V(!D^{(1)}\ast D^{(2)})$
for a good meromorphic flat bundle
on $(X,D)$.
Compared with functions with moderate growth,
we may consider 
functions with rapid decay along some direction
and of Nilsson type along other direction.
\hfill\qed
\end{rem}

\subsection{Real blow up along holomorphic functions}

\subsubsection{Category of complex manifolds over $\cnum^{\ell}$}

It is convenient to consider the category $\Cat_{\ell}$
of complex manifolds over $\cnum^{\ell}$
given as follows.
\index{category $\Cat_{\ell}$}
An object of $\Cat_{\ell}$ is 
a morphism $f:X\lrarr \cnum^{\ell}$ of complex manifolds.
Morphisms
$\varphi:(X_1,f_1)\lrarr (X_2,f_2)$ in $\Cat_{\ell}$
are morphisms of complex manifolds $\varphi:X_1\lrarr X_2$
such that $f_1=f_2\circ\varphi$.
We say that $\varphi$ has {\em some property}
when the underlying $\varphi$ has the property.
For example,
we say that $\varphi:(X_1,f_1)\lrarr (X_2,f_2)$ is a closed immersion
when $\varphi:X_1\lrarr X_2$ is a closed immersion.
For a given object $(X,f)$ in $\Cat_{\ell}$,
we set $D_X:=f^{-1}(D_0)$,
where $D_0:=\bigcup_{i=1}^{\ell}\{z_i=0\}$.
Let $\cnumtilde$ denote the real blow up of
$\cnum$ along $z=0$.
We have $\cnumtilde^{\ell}(D_0)=\cnumtilde^{\ell}$.

For any object $(X,f)$ in $\Cat_{\ell}$,
we have the naturally defined map
$\Gamma_f:X\lrarr X\times\cnum^{\ell}$
given by $\Gamma_f(x)=(x,f(x))$.
A morphism $\varphi:(X_1,f_1)\lrarr (X_2,f_2)$
induces maps
$X_1\times\cnum^{\ell}\lrarr X_2\times\cnum^{\ell}$
and
$X_1\times\cnumtilde^{\ell}\lrarr X_2\times\cnumtilde^{\ell}$,
which are denoted by
$\varphi_1$ and $\varphitilde_1$,
respectively.

\subsubsection{Real blow up along functions}
\label{subsection;14.1.19.10}

Let $(X,f)$ be an object in $\Cat_{\ell}$.
Let $j:X\times(\cnum^{\ast})^{\ell}\lrarr
 X\times\cnumtilde^{\ell}$
denote the inclusion.
Let $\Xtilde(f)$ denote the topological space
obtained as the closure of
$j\bigl(\Gamma_f(X\setminus D_X)\bigr)$
in $X\times\cnumtilde^{\ell}$,
which is called the real blow up of $X$ along $f$
\cite{sabbah_lecture_Stokes}.
\index{real blow up $\Xtilde(f)$}
The projection $\Xtilde(f)\lrarr X$
is denoted by $\pi_{f}$.
The inclusion
$\Xtilde(f)\lrarr X\times\cnumtilde^{\ell}$
is denoted by $\Gammatilde_f$.
If there is no risk of confusion,
we shall omit to denote the subscript $f$
to simplify the notation.
If $f$ is submersive, 
$\Xtilde(f)$ is naturally diffeomorphic to
$\Xtilde(D_X)$.
A morphism
$\varphi:(X_1,f_1)\lrarr (X_2,f_2)$ 
in $\Cat_{\ell}$ naturally induces a continuous map
$\varphitilde:\Xtilde_1(f_1)\lrarr \Xtilde_2(f_2)$.

\subsubsection{Moderate growth and rapid decay}

Let $(X,f)\in\Cat_{\ell}$.
Let $U$ be any open subset of $\Xtilde(f)$.
A holomorphic function $s$ on
$U\setminus\pi_f^{-1}(D_X)$
is called of moderate growth
if we have
$|s|=O\bigl(\prod|f_i|^{-N}\bigr)$
for some $N$ locally around any point of 
$U\cap\pi^{-1}(D_X)$.
A holomorphic function $s$ on
$U\setminus\pi_f^{-1}(D_X)$
is called of rapid decay
if we have
$|s|=O\bigl(\prod|f_i|^{N}\bigr)$
for any $N$ locally around any point of 
$U\cap\pi^{-1}(D_X)$.
\index{holomorphic functions with moderate growth}
\index{holomorphic functions with rapid decay}
The sheaf of holomorphic functions 
with moderate growth (resp. rapid decay)
is denoted by
$\nbiga^{\moderate}_{\Xtilde(f)}$
(resp. $\nbiga^{\rapid}_{\Xtilde(f)}$).
\index{sheaf $\nbiga^{\moderate}_{\Xtilde(f)}$}
\index{sheaf $\nbiga^{\rapid}_{\Xtilde(f)}$}
We shall prove the following theorem.
(See Theorems \ref{thm;13.4.20.1}, \ref{thm;14.1.19.1},
and Theorems \ref{thm;13.4.19.42}, \ref{thm;13.4.19.210}
for refined claims.)
\begin{thm}
\mbox{{}}
\begin{itemize}
\item
The sheaves
$\nbiga^{\moderate}_{\Xtilde(f)}$
and 
$\nbiga^{\rapid}_{\Xtilde(f)}$
are flat over $\pi_f^{-1}(\nbigo_X)$.
\item
Let $\Gammatilde_f:\Xtilde(f)\lrarr X\times\cnumtilde^{\ell}$
denote the inclusion.
Then, we naturally have
\[
 \Gammatilde_{f\ast}\nbiga^{\rapid}_{\Xtilde(f)}
\simeq
 \pi^{-1}\nbigo_{\Gamma_f(X)}
 \otimes_{\pi^{-1}\nbigo_{X\times\cnum^{\ell}}}
 \nbiga^{\rapid}_{X\times\cnumtilde^{\ell}},
\]
\[
 \Gammatilde_{f\ast}\nbiga^{\moderate}_{\Xtilde(f)}
\simeq
 \pi^{-1}\nbigo_{\Gamma_f(X)}
 \otimes_{\pi^{-1}\nbigo_{X\times\cnum^{\ell}}}
 \nbiga^{\moderate}_{X\times\cnumtilde^{\ell}}.
\]
\item
Let $\rho_0:\Xtilde(D_X)\lrarr \Xtilde(f)$
denote the naturally induced map.
Then, we naturally have
\[
 R\rho_{0\ast}\nbiga^{\rapid}_{\Xtilde(D_X)}
\simeq
 \nbiga^{\rapid}_{\Xtilde(f)},
\quad\quad
 R\rho_{0\ast}\nbiga^{\moderate}_{\Xtilde(D_X)}
\simeq
 \nbiga^{\moderate}_{\Xtilde(f)}.
\]
\item
Let $\varphi:(Y,g)\lrarr (X,f)$
be a projective morphism in $\Cat_{\ell}$.
Let $M$ be a coherent $\nbigo_Y$-module.
Then, we have the following natural isomorphism:
\[
 \nbiga^{\moderate}_{\Xtilde(f)}
 \otimes_{\pi_f^{-1}\nbigo_X}
 \pi_f^{-1}R\varphi_{\ast}M
\lrarr
 R\varphitilde_{\ast}\bigl(
 \nbiga^{\moderate}_{\Ytilde(g)}
 \otimes_{\pi_g^{-1}\nbigo_Y}
 \pi_g^{-1}M
 \bigr)
\]
\end{itemize}
\end{thm}

%% file: 4.2.tex
\subsection{Preliminary}
\label{subsection;14.1.19.2}

Let $X$ be any $n$-dimensional complex manifold
with a simply normal crossing hypersurface
$D$ with the irreducible decomposition
$\bigcup_{i\in\Lambda} D_i$.
We use the notation in \S\ref{subsection;14.1.17.30}.
Let $D^{\circ}$ be a (possibly empty) hypersurface of $X$
such that 
(i) $D\cup D^{\circ}$ is simply normal crossing,
(ii) $\dim D\cap D^{\circ}<n-1$.
For $J\subset\Lambda$,
we set $D(\Jbar):=D(J)\cup D^{\circ}$.
For $I\sqcup J\subset\Lambda$,
we put $D_I(\Jbar):=D_I\cap D(\Jbar)$.

Let $\Omega^{0,q}_{\Xtilde(D)}$
denote the sheaf of 
$C^{\infty}$-logarithmic $(0,q)$-forms
on $\Xtilde(D)$,
i.e.,
a section of $\Omega^{0,q}_{\Xtilde(D)}$
is locally described as a linear combination of
\[
 f\cdot
 d\zbar_{i_1}/\zbar_{i_1}
 \cdots d\zbar_{i_m}/\zbar_{i_m}\cdot
 \,d\zbar_{j_1}\cdots d\zbar_{j_k}
\quad
(1\leq i_p\leq \ell,\,\,
\ell+1\leq j_q\leq n,\,\,
f\in\nbigc^{\infty}_{\Xtilde(D)})
\]
in terms of 
a local holomorphic coordinate system $(z_1,\ldots,z_n)$
such that $D$ is locally described as
$\bigcup_{i=1}^{\ell}\{z_i=0\}$.
\index{sheaf $\Omega^{p,q}_{\Xtilde(D)}$}
We have the naturally defined operator
$\delbar:\Omega_{\Xtilde(D)}^{0,q}
\lrarr \Omega_{\Xtilde(D)}^{0,q+1}$.
The complex
$\Omega^{0,\bullet}_{\Xtilde(D)}$
is called the Dolbeault complex of $\Xtilde(D)$.
We put 
$\Omega^{0,\bullet}_{\Zhat}:=
 \Omega_{\Xtilde(D)|\Zhat}^{0,\bullet}$
for any real analytic subset $Z\subset\Xtilde(D)$.
\index{sheaf $\Omega^{0,\bullet}_{\Zhat}$}

For a given $C^{\infty}$-manifold $Y$
and a real analytic subset $W\subset X$,
let $\nbigc^{\infty\,<W}
 _{\widehat{\pi^{-1}(D_I)}\times Y}$
denote the sheaf
$\nbigc^{\infty\,<\pi^{-1}(W)\times Y}
 _{\widehat{\pi^{-1}(D_I)}\times Y}$
on $\Xtilde(D)\times Y$,
for simplicity of the description.
\index{sheaf $\nbigc^{\infty\,<W}
 _{\widehat{\pi^{-1}(D_I)}\times Y}$}
We also put 
$\Omega^{0,\bullet<W}_{
 \widehat{\pi^{-1}(D_I)}\times Y}:=
 \Omega^{0,\bullet}_{\Xtilde(D)}
 \otimes_{\nbigc^{\infty}_{\Xtilde(D)}}
 \nbigc^{\infty<W}_{\widehat{\pi^{-1}(D_I)}\times Y}$
on $\Xtilde(D)\times Y$.
\index{sheaf $\Omega^{0,\bullet<W}_{
 \widehat{\pi^{-1}(D_I)}\times Y}$}

Let $q_I$ denote the projection
$\pi^{-1}(D_I)\lrarr \Dtilde_I(\del D_I)$.
If we are given a holomorphic coordinate system
$(z_1,\ldots,z_n)$ as above,
then
$\nbigo_{\widehat{\pi^{-1}(D_I)}}^{<D(J)}
=q_I^{-1}\nbigo_{\Dtilde_I(\del D_I)}^{<D_I(J)}
 [\![z_i\,|\,i\in I]\!]$.
By a natural diffeomorphism
$\pi^{-1}(D_I)\simeq
 \Dtilde_I(\del D_I)\times (S^1)^{|I|}$,
we can locally identify
$\nbigc^{\infty\,<D(\Jbar)}
 _{\widehat{\pi^{-1}(D_I)}}
=\nbigc^{\infty\,<D_I(\Jbar)}
 _{\Dtilde_I(\del D_I)\times (S^1)^{|I|}}
 [\![z_i\,|\,i\in I]\!]$.

\vspace{.1in}
For $I\subset J$,
put $\nbigt(m,I,J):=\bigl\{
 K\subset J\,\big|\,
 I\subset K,\,|K|=|I|+m+1\bigr\}$
for $m\geq 0$.
We set
$\nbigk^{m}\Bigl(
 \nbigo_{\widehat{\pi^{-1}(D_I(J))}}
 \Bigr)
:=\bigoplus_{K\in\nbigt(m,I,J)}
 \nbigo_{\widehat{\pi^{-1}(D_K)}}$.
We obtain a complex
$\nbigk^{\bullet}\Bigl(
 \nbigo_{\widehat{\pi^{-1}(D_I(J))}}
 \Bigr)$
as in \S\ref{subsection;09.10.29.10}.
Similarly,
we obtain a complex
$\nbigk^{\bullet}\bigl(
 \Omega^{0,\bullet\,<D^{\circ}}
 _{\widehat{\pi^{-1}(D_I(J))}\times Y}\bigr)$.
See \S I.5 of \cite{malgrange2}
and \S II.1.1 of \cite{sabbah4}
for the following.

\begin{lem}
\label{lem;09.10.24.10}
Let $\nbigb$ be
$\nbigo_{\widehat{\pi^{-1}(D_I(J))}}$
or $\Omega^{0,\bullet\,<D^{\circ}}
 _{\widehat{\pi^{-1}(D_I(J))}\times Y}$.
The natural inclusion
$\nbigb\lrarr \nbigk^{0}(\nbigb)$
induces a quasi-isomorphism
$\nbigb\lrarr\nbigk^{\bullet}(\nbigb)$.
\hfill\qed
\end{lem}

\subsection{Dolbeault resolution}

In this subsection,
we suppose $D^{\circ}=\emptyset$.
\begin{prop}[\cite{majima}, \cite{sabbah4}]
\label{prop;09.10.26.1}
$\Omega^{0,\bullet}_{\widehat{\pi^{-1}(D_I(J))}}$
and $\Omega^{0,\bullet\,<D(J)}_{\widehat{\pi^{-1}(D_I)}}$
are $c$-soft resolutions of
$\nbigo_{\widehat{\pi^{-1}(D_I(J))}}$
and $\nbigo_{\widehat{\pi^{-1}(D_I)}}^{<D(J)}$
respectively,
where $J\subset I^c$.
\end{prop}
\pf
We give only an outline.
In each case,
it is easy to compute
the $0$-th cohomology of
the Dolbeault complexes.
It is enough to prove the vanishing
of the higher cohomology.
We may assume
$X=\Delta^n$,
$D_i=\{z_i=0\}$
and $D=\bigcup_{i=1}^{\ell}D_i$.
First, let us look at
$\Omega^{0,\bullet}_{\Xtilde(D)}$.
For $1\leq j\leq n$,
let $\nbigp^0_{\leq j}$
be the sheaf of $C^{\infty}$-functions
on $\Xtilde(D)$ which are 
$\delbar_i$-holomorphic for $i>j$.
We set 
$X_j:=\Delta^j=\{(z_1,\ldots,z_j)\}$
and $D_{j,\ell}:=
 \bigcup_{i\leq \min\{j,\ell\}}\{z_i=0\}$.
Let $q_{\leq j}$ be the projection 
$\Xtilde(D)\lrarr \Xtilde_j(D_{j,\ell})$.
Let $\nbigp^1_{\leq j}$ be the sheaf of
$C^{\infty}$-sections of
$q_{\leq j}^{-1}
 \Omega^{0,1}_{\Xtilde_j(D_{j,\ell})}$,
which are $\delbar_i$-holomorphic for $i>j$.
We set
$\nbigp^{\bullet}_{\leq j}:=
 \bigwedge^{\bullet}\nbigp_{\leq j}^1$
over $\nbigp^0_{\leq j}$.
We have the naturally defined operator
 $\delbar:
 \nbigp^{\bullet}_{\leq j}\lrarr
 \nbigp^{\bullet+1}_{\leq j}$.

Because
$\nbigp^{\bullet}_{\leq 0}=
 \nbigo_{\Xtilde(D)}$
and 
$\nbigp^{\bullet}_{\leq n}
=\Omega^{0,\bullet}_{\Xtilde(D)}$,
it is enough to prove that
the natural inclusions
$\nbigp^{\bullet}_{\leq j}
\lrarr
 \nbigp^{\bullet}_{\leq j+1}$
are quasi-isomorphisms
for the vanishing of the higher cohomology of
$\Omega^{0,\bullet}_{\Xtilde(D)}$.
Let $\nbigq^0_{\leq j}=\nbigp^0_{\leq j+1}$.
Let $\nbigq^1_{\leq j}$ be the sheaf of
$q_{\leq j}^{-1}
 \Omega^{0,1}_{\Xtilde_j(D_{j,\ell})}$
which are $\delbar_i$-holomorphic
for $i>j+1$.
We take the exterior product
$\nbigq^{\bullet}_{\leq j}=
 \bigwedge^{\bullet}\nbigq_{\leq j}^1$
over $\nbigq_{\leq j}^0$.
We have the naturally defined operator
$\delbar_{j+1}:
 \nbigq^{\bullet}_{\leq j}
\lrarr
 \nbigq^{\bullet}_{\leq j}
\wedge
 d\zbar_{j+1}\big/\zbar_{j+1}$
$(j-1\leq \ell)$
or 
$\delbar_{j+1}:
 \nbigq^{\bullet}_{\leq j}
\lrarr
 \nbigq^{\bullet}_{\leq j}
\wedge
 d\zbar_{j+1}$
$(j\geq \ell)$.
We clearly have
$\Ker\delbar_{j+1}=\nbigp_{\leq j}^{\bullet}$.
Let us prove $\Cok\delbar_{j+1}=0$.
In the case $j\geq \ell$,
it can be proved by the argument
for the standard Dolbeault's lemma.
Let us consider the case $j<\ell$.

\begin{lem}
\label{lem;09.10.22.11}
$\delbar_{j+1}:
 \nbigq^{\bullet}_{\leq j|\widehat{\pi^{-1}(D_{j+1})}}
\lrarr
 \nbigq^{\bullet}_{\leq j|\widehat{\pi^{-1}(D_{j+1})}}
 \wedge
 d\zbar_{j+1}\big/\zbar_{j+1}$
is an epimorphism.
\end{lem}
\pf
We use the polar coordinate system
$z_{j+1}=r_{j+1}\,e^{\sqrt{-1}\theta_{j+1}}$.
The action of $\delbar_{j+1}$ 
is expressed as follows:
\[
 \delbar_{j+1}\Bigl(
 \sum_n f_n(\theta_{j+1})\,z_{j+1}^n
 \Bigr)
=\sum_n
 \Bigl(
 \frac{\sqrt{-1}}{2}\del_{\theta_{j+1}}
 \Bigr)f_n(\theta_{j+1})\,z_{j+1}^n
 \cdot d\zbar_{j+1}\big/\zbar_{j+1}
\]
Then, it is easy to prove the claim of
Lemma \ref{lem;09.10.22.11}.
\hfill\qed

\vspace{.1in}

Put $D':=\bigcup_{i=1,i\neq j+1}^{\ell}\{z_i=0\}$,
and let us consider the real blow up
$\pi':\Xtilde(D')\lrarr X$.
We have a naturally induced morphism
$q_{\leq j}':\Xtilde(D')\lrarr \Xtilde_j(D_{j,\ell})$.
Let $\nbigs^1_{\leq j,X}$ be
the sheaf of sections of
$(q'_{\leq j})^{-1}
 \Omega^{0,1}_{\Xtilde_j(D_{j,\ell})}$
on $\Xtilde(D')$,
which are $\delbar_i$-holomorphic for $i>j+1$.
Let $\nbigs^0_{\leq j,X}$ be
the sheaf of $C^{\infty}$-functions
on $\Xtilde(D')$,
which are $\delbar_{i}$-holomorphic
for $i>j+1$.
We set
$\nbigs^{\bullet}_{\leq j}:=
 \bigwedge^{\bullet}\nbigs^{1}_{\leq j}$.
It is easy to prove the vanishing of
the cokernel of 
$\delbar_{j+1}:
 \nbigs^{\bullet}_{\leq j}\lrarr
 \nbigs^{\bullet}_{\leq j}\wedge d\zbar_{j+1}$
by using the argument
for standard Dolbeault's lemma.

Let $P\in \pi^{-1}(D)$.
Let $U$ be a small neighbourhood around $P$.
We will shrink it in the following argument.
According to Lemma \ref{lem;09.10.22.11},
for any section $\varphi$ of
$\nbigq^{\bullet}_{\leq j}
 \wedge d\zbar_{j+1}/\zbar_{j+1}$ on $U$,
we can take a local section
$\psi$ of $\nbigq^{\bullet}_{\leq j}$
such that
\[
\bigl(
 \varphi-\delbar_j\psi
\bigr)_{|\widehat{\pi^{-1}(D_j)}\cap U}=0.
\]
We put 
$\lambda:=\varphi-\delbar_j\psi$.
We take a cut function $\rho$ around $P$,
i.e.,
$\rho$ is constantly $1$ around $P$
and constantly $0$ near the boundary of $U$.
We can regard $\rho\,\lambda$
as a section of 
$\nbigs^{\bullet}_{\leq j}\wedge d\zbar_{j+1}$.
Then, we can find a section
$\kappa$  of $\nbigs^{\bullet}_{\leq j}$ 
around $\pi_j(P)$ such that
$\delbar_{j+1}\kappa=\rho\lambda$,
where $\pi_j$ denotes the natural projection
$\Xtilde(D)\lrarr \Xtilde(D')$.
We obtain
$\varphi=\delbar_j(\psi+\kappa)$
around $P$.
Thus, we obtain the vanishing
of the cokernel of
$\delbar_{j+1}:
 \nbigq^{\bullet}_{\leq j}
\lrarr
 \nbigq^{\bullet}_{\leq j}
\wedge
 d\zbar_{j+1}\big/\zbar_{j+1}$,
and hence the vanishing of
the higher cohomology of
$\Omega^{0,\bullet}_{\Xtilde(D)}$.

\vspace{.1in}

Because 
$\pi^{-1}(D_I)=
 \Dtilde_I(\del D_I)
 \times (S^1)^{|I|}$,
we can reduce
the vanishing of the higher cohomology
of $\Omega^{0,\bullet}_{\widehat{\pi^{-1}(D_I)}}$
to the vanishing of
$\Omega^{0,\bullet}_{\Dtilde_I(\del D_I)}$
by a formal calculation as in 
Lemma \ref{lem;09.10.22.11}.
By using the resolution in
Lemma \ref{lem;09.10.24.10},
we obtain the vanishing of 
the higher cohomology of 
$\Omega^{0,\bullet}_{\widehat{\pi^{-1}(D(I))}}$.
We have the following diagram
of exact sequences:
\[
 \begin{CD}
 0 @>>>\nbigo_{\Xtilde(D)}^{<D(I)}
 @>>> \nbigo_{\Xtilde(D)}
 @>>> \nbigo_{\widehat{\pi^{-1}(D(I))}}
 @>>> 0 \\
 @. @VVV @VVV @VVV @.\\
 0 @>>> 
 \Omega^{0,\bullet\,<D(I)}_{\Xtilde(D)}
 @>>>
 \Omega^{0,\bullet}_{\Xtilde(D)}
 @>>>
 \Omega^{0,\bullet}_{\widehat{\pi^{-1}(D(I))}}
 @>>> 0
 \end{CD}
\]
Then, we obtain the vanishing
of the higher cohomology of
$\Omega^{0,\bullet\,<D(I)}_{\Xtilde(D)}$.
By a formal calculation as in 
Lemma \ref{lem;09.10.22.11},
we obtain the vanishing of
the higher cohomology of
$\Omega_{\widehat{\pi^{-1}(D_I(J))}}
 ^{0,\bullet}$
and
$\Omega_{\widehat{\pi^{-1}(D_I)}}
 ^{0,\bullet\,<D(J)}$.
\hfill\qed

\subsection{Flatness}

In this subsection,
$D^{\circ}$ is not necessarily empty.

\begin{prop}
\label{prop;09.10.25.32}
Let $I\sqcup J\subset\Lambda$.
The sheaves
$\nbigc^{\infty\,<D(\Jbar)}_{\widehat{\pi^{-1}(D_I)}}$,
$\nbigc^{\infty\,<D^{\circ}}_{\widehat{\pi^{-1}(D_I(J))}}$,
$\nbigo^{<D(J)}_{\widehat{\pi^{-1}(D_I)}}$
and 
$\nbigo_{\widehat{\pi^{-1}(D_I(J))}}$
are flat over $\pi^{-1}\nbigo_X$.
In particular,
the sheaves
$\nbigo_{\Xtilde(D)}$ and
$\nbiga^{\rapid}_{\Xtilde(D)}$
are flat over $\pi^{-1}\nbigo_X$.
\end{prop}
\pf
Let us recall a general result.
For a real analytic manifold $Y$,
let $\nbigo^{\real}_Y$ denote the sheaf of
real analytic functions on $Y$.
If $Y$ is the product of a complex manifold $Y_1$
and a real analytic manifold $Y_2$,
let $\nbigo^{Y_1-\hol}_{Y}$ denote the sheaf of
real analytic functions which are holomorphic
in the $Y_1$-direction.
The extension
$\nbigo^{Y_1-\hol}_{Y}\subset
 \nbigo^{\real}_{Y}$ is faithfully flat.

\begin{lem}
\label{lem;09.10.25.30}
Let $W_1\subset W_2\subset Y$ be 
real analytic subsets.
Then,
$\nbigc^{\infty\,<W_i}_Y$
and $\nbigc^{\infty\,<W_1}_Y\big/
 \nbigc^{\infty\,<W_2}_Y$
are flat over $\nbigo_{Y}^{\real}$.
\end{lem}
\pf
The sheaf $\nbigc^{\infty}_Y$
is faithfully flat over $\nbigo^{\real}_Y$
(Corollary 1.12 of \cite{malgrange2}).
Theorem VI.1.2 of \cite{malgrange2} implies
$\gminia\,\nbigc^{\infty\,<W_1}_Y\cap
 \nbigc^{\infty\,<W_2}_Y
=\gminia\,\nbigc^{\infty\,<W_2}_Y$
for any real analytic subsets
$W_1\subset W_2\subset Y$
and for any ideal sheaf $\gminia$ 
of $\nbigo_Y^{\real}$.
By using the argument in the proof of 
Proposition III.4.7 in \cite{malgrange2},
we can prove the following:
\begin{itemize}
\item
Let $A$ be a ring.
Let $M$ be an $A$-flat module.
Let $N$ be an $A$-submodule of $M$.
If $\gminia M\cap N=\gminia N$
for any ideal $\gminia$ of $A$,
then $N$ and $M/N$ are also $A$-flat.
\end{itemize}
We immediately obtain the claim of
Lemma \ref{lem;09.10.25.30}
from these results.
\hfill\qed

\vspace{.1in}

Let $Z_0$ be a complex manifold
with a normal crossing hypersurface $D_0$.
Let $Z_1$ be a real analytic manifold.
We put $Z:=Z_0\times Z_1$
and $D:=D_0\times Z_1$.
Let $G$ denote the composite of 
the maps
$Z\lrarr Z_0\lrarr Z_0\times\cnum^n$,
where the latter is induced by
the inclusion
$\{(0,\ldots,0)\}\subset\cnum^n$.
Let $(t_1,\ldots,t_n)$ be the standard
holomorphic coordinate system of $\cnum^n$.
\begin{lem}
\label{lem;09.10.25.31}
$\nbigc^{\infty\,<D}_{Z}
 [\![t_1,\ldots,t_n]\!]$
is flat over 
$G^{-1}\nbigo_{Z_0\times \cnum^n}$.
\end{lem}
\pf
Let $\iota_1$ denote the inclusion
$Z\lrarr Z_2:=Z\times\real^n$
induced by $\{(0,\ldots,0)\}\lrarr\real^n$.
We put $D_2:=D\times\real^n$.
We regard that $(t_1,\ldots,t_n)$
is a real coordinate system of
$\real^n\subset\cnum^n$.
We have the natural identification
$\nbigc^{\infty\,<D}_Z[\![t_1,\ldots,t_n]\!]
=\nbigc^{\infty\,<D_2}_{Z_2}\big/
 \nbigc^{\infty\,<D_2\cup Z}_{Z_2}$.
According to Lemma \ref{lem;09.10.25.30},
it is flat over $\iota_1^{-1}\nbigo_{Z_2}^{\real}$.
Let $G_1$ be the composite of
$Z\lrarr Z_0\lrarr Z_0\times\real^n$.
We have a natural isomorphism
$G_1^{-1}\nbigo_{Z_0\times\real^n}^{Z_0-\hol}
\simeq
 G^{-1}\nbigo_{Z_0\times\cnum^n}$.
Since the extension
$G_1^{-1}\nbigo_{Z_0\times\real^n}^{Z_0-\hol}
 \subset
 \nbigo_{Z_2}^{\real}$ is faithfully flat,
we obtain the claim of Lemma \ref{lem;09.10.25.31}.
\hfill\qed

\vspace{.1in}

Let us return to the proof of 
Proposition \ref{prop;09.10.25.32}.
We may assume that
$X=\Delta^n$, $D_i=\{z_i=0\}$,
$D=\bigcup_{i=1}^{\ell}D_i$
and $D^{\circ}=\bigcup_{i=\ell+1}^mD_i$.
For $I\subset \ellsitabar$,
let $\pi_I:\Xtilde(D(I))\lrarr X$ be the real blow up.
We have the natural identification
$\pi_I^{-1}(D_I)=D_I\times (S^1)^{|I|}$
and $\pi_I^{-1}\bigl(D_I(\Ibar^c)\bigr)=
 D_I(\Ibar^c)\times (S^1)^{|I|}$.
From Lemma \ref{lem;09.10.25.31},
we obtain that
$\nbigc^{\infty\,<D(\Ibar^c)}_{
 \widehat{\pi_I^{-1}(D_I)}}
=\nbigc^{\infty\,<D_I(\Ibar^c)}_{
 \pi_I^{-1}(D_I)}[\![z_i\,|\,i\in I]\!]$
is flat over $\pi_I^{-1}\nbigo_X$.

\begin{lem}
\label{lem;09.10.25.35}
$\nbigc^{\infty\,<D(\Ibar^c)}
 _{\widehat{\pi^{-1}(D_I)}}$
is flat over $\pi^{-1}\nbigo_X$.
(Note that $\pi:\Xtilde(D)\lrarr X$.)
\end{lem}
\pf
The claim is clear outside of
$\pi^{-1}(\del D_I)$.
Let $P$ be any point of $\del D_I$.
Let $\gminia$ be any finitely generated ideal
of $\nbigo_{X,P}$.
We take a free resolution
$\nbigq_{\bullet}$ of $\gminia$,
i.e.,
$\cdots \rarr \nbigq_{1}\rarr \nbigq_{0}
\lrarr\gminia$.
We obtain a $\pi^{-1}\nbigo_X$-free resolution
$\pi^{-1}\nbigq_{\bullet}$ of
$\pi^{-1}\gminia$.
We set $\nbigqtilde_j=\nbigq_j$ for $j\geq 0$
and $\nbigqtilde_{-1}:=\gminia$
for simplicity of the description.
It is enough to prove that
$\pi^{-1}\nbigqtilde_{\bullet}
 \otimes 
 \nbigc^{\infty\,<D(\Ibar^c)}
 _{\widehat{\pi^{-1}(D_I)}}$ is exact.
Let $\rho:\Xtilde(D)\lrarr\Xtilde(D(I))$
be the naturally induced map.
Note 
\[
 \rho_{\ast}\bigl(
 \pi^{-1}\nbigqtilde_{\bullet}
 \otimes\nbigc^{\infty\,<D(\Ibar^c)}_{
 \widehat{\pi^{-1}(D_I)}}
 \bigr)
=\pi_I^{-1}(\nbigqtilde_{\bullet})\otimes
 \rho_{\ast}\bigl(
 \nbigc^{\infty\,<D(\Ibar^c)}_{
 \widehat{\pi^{-1}(D_I)}}
 \bigr)
=\pi_I^{-1}\bigl(\nbigqtilde_{\bullet}\bigr)
 \otimes
 \nbigc^{\infty\,<D(\Ibar^c)}
 _{\widehat{\pi_I^{-1}(D_I)}}.
\]
The first equality is the projection formula.
As for the second one,
it is enough to observe
that the natural morphism
$\nbigc^{\infty <D(\Ibar^c)}_{\pi_I^{-1}(D_I)}
\lrarr
\rho_{\ast}\nbigc^{\infty <D(\Ibar^c)}_{\pi^{-1}(D_I)}$
is an isomorphism.
It is clearly injective.
Let $f$ be a section of
$\rho_{\ast}\nbigc^{<\infty D(\Ibar^c)}_{\pi^{-1}(D_I)}$.
The restriction
$g:=f_{|\pi_I^{-1}(D_I\setminus D(\Ibar^c))}$
gives a $C^{\infty}$-function on
$\pi_I^{-1}(D_I\setminus \del D_I)$.
For any differential operator $R$ 
on $\pi_I^{-1}(D_I)$,
$R(g)(P)$ goes to $0$ 
when $P$ goes to a point in 
$\pi_I^{-1}(\del D_I)$.
Hence, $g$ gives a section of
$\nbigc^{\infty <D(\Ibar^c)}_{\pi_I^{-1}(D_I)}$
which is mapped to $f$.
Let $Q\in \pi^{-1}(P)$.
Take any cycle $\varphi$ of
$ \pi^{-1}\nbigqtilde_{i}\otimes
 \nbigc^{\infty\,<D(\Ibar^c)}
 _{\widehat{\pi^{-1}(D_I)}}$ at $Q$.
By using a cut function around $Q$,
we can regard it as a global cycle of
$\pi^{-1}\nbigqtilde_{i}\otimes
 \nbigc^{\infty\,<D(\Ibar^c)}
 _{\widehat{\pi^{-1}(D_I)}}$
whose support is a small neighbourhood of $Q$.
Then, it can be regarded 
as a cycle of
$\pi_I^{-1}(\nbigqtilde_{i})\otimes
 \nbigc^{\infty\,<D(\Ibar^c)}
 _{\widehat{\pi_I^{-1}(D_I)}}$
around $\rho(Q)$.
Because 
$\nbigc^{\infty\,<D(\Ibar^c)}
 _{\widehat{\pi_I^{-1}(D_I)}}$
is flat over $\pi_I^{-1}\nbigo_X$,
we obtain that $\varphi$ is a boundary
in the complex
$\pi_I^{-1}(\nbigqtilde_{\bullet})\otimes
 \nbigc^{\infty\,<D(\Ibar^c)}
 _{\widehat{\pi_I^{-1}(D_I)}}$.
Then, it is easy to deduce that
$\varphi$ is a boundary in the complex
$\pi^{-1}(\nbigqtilde_{\bullet})\otimes
 \nbigc^{\infty\,<D(\Ibar^c)}
 _{\widehat{\pi^{-1}(D_I)}}$.
Thus, the proof of Lemma 
\ref{lem;09.10.25.35} is finished.
\hfill\qed

\vspace{.1in}
Let us prove that
$\nbigc^{\infty\,<D(\Jbar)}_{\widehat{\pi^{-1}(D_I)}}$
is flat over $\pi^{-1}\nbigo_X$,
where $I\sqcup J\subset\ellsitabar$.
We put
\[
 \nbigs(I,J,m):=\bigl\{
 K\subset \ellsitabar-J\,\big|\,
 \,I\subset K,\,|K|=m
 \bigr\}. 
\]
Put 
$\nbigg_{I,\ell+1}:=
 \nbigc^{\infty\,<D(\Jbar)}
_{\widehat{\pi^{-1}(D_I)}}$,
and descending inductively
we set
\[
 \nbigg_{I,m}:=
 \Ker\Bigl(
 \nbigg_{I,m+1}\lrarr
 \bigoplus_{K\in\nbigs(I,J,m)}
 \nbigc^{\infty\,<D(\Kbar^c)}_{\widehat{\pi^{-1}(D_K)}}
 \Bigr).
\]
We have
$\nbigg_{I,|I|+1}=
 \nbigc^{\infty\,<D(\Ibar^c)}
 _{\widehat{\pi^{-1}(D_I)}}$,
which is flat over $\pi^{-1}\nbigo_X$.
By an induction,
we obtain that
$\nbigg_{I,m}$ are flat over $\pi^{-1}\nbigo_X$.
Hence, we obtain that
$\nbigc^{\infty\,<D(\Jbar)}_{\widehat{\pi^{-1}(D_I)}}$
is flat over $\pi^{-1}\nbigo_X$.
By using the resolution of
$\nbigc^{\infty\,<D^{\circ}}_{\widehat{\pi^{-1}(D_I(J))}}$
in Lemma \ref{lem;09.10.24.10},
we obtain that
$\nbigc^{\infty\,<D^{\circ}}
 _{\widehat{\pi^{-1}(D_I(J))}}$
is flat over $\pi^{-1}\nbigo_X$.
As a result,
we obtain that 
$\Omega^{0,\bullet\,<D(\Jbar)}
 _{\widehat{\pi^{-1}(D_I)}}$
and 
$\Omega^{0,\bullet\,<D^{\circ}}
 _{\widehat{\pi^{-1}(D_I(J))}}$
are flat over $\pi^{-1}\nbigo_X$,
where $J\subset I^c$.
In particular,
$\Omega^{0,\bullet\,<D(J)}
 _{\widehat{\pi^{-1}(D_I)}}$
and 
$\Omega^{0,\bullet}
 _{\widehat{\pi^{-1}(D_I(J))}}$
are flat over $\pi^{-1}\nbigo_X$.
Then, we obtain the $\pi^{-1}\nbigo_X$-flatness of
$\nbigo^{<D(J)}_{\widehat{\pi^{-1}(D_I)}}$
and $\nbigo_{\widehat{\pi^{-1}(D_I(J))}}$
by using Proposition \ref{prop;09.10.26.1}.
Thus, the proof of Proposition
\ref{prop;09.10.25.32} is finished.
\hfill\qed

%% file: 4.3.tex
\subsection{$C^{\infty}$-functions
of Nilsson type}

Let $X$, $D$ and $D^{\circ}$
be as in \S\ref{subsection;14.1.19.2}.
We put $D^{(3)}:=D^{(1)}\cup D^{\circ}$.
We shall introduce a sheaf
$\nbigc^{\infty<D^{(3)}\leq D^{(2)}}_{\Xtilde(D)}$
on $\Xtilde(D)$.
First, let us consider the case $X=\Delta^n$,
$D=\bigcup_{i=1}^{\ell}\{z_i=0\}$
and
$D^{\circ}=\bigcup_{i=\ell+1}^m\{z_i=0\}$.
Let $\ellsitabar=I_1\sqcup I_2$
be determined by
$D^{(j)}=\bigcup_{i\in I_j}\{z_i=0\}$
for $j=1,2$.
Let $\jtilde$ denote the inclusion
$X-D\lrarr\Xtilde(D)$.
\index{sheaf $\nbiga^{<D^{(1)}\leq D^{(2)}}_{\Xtilde(D)}$}
Let $\nbigc^{\infty\,<D^{(3)}\leq D^{(2)}}_{\Xtilde(D)}$
be the image of the naturally defined morphisms:
\[
 \nbigc_{\Xtilde(D)}^{\infty\,<D^{(3)}}
 \otimes\Nil(z_i\,|\,i\in I_2)
\lrarr
 \jtilde_{\ast}\nbigc^{\infty\,<D^{\circ}}_{X-D}.
\]
\index{sheaf $\nbigc^{\infty\,<D^{(3)}\leq D^{(2)}}_{\Xtilde(D)}$}
We can observe that
they are independent of the choice of a coordinate system
$(z_1,\ldots,z_n)$.
Hence, we obtain a globally defined sheaf
$\nbigc^{\infty<D^{(3)}\leq D^{(2)}}_{\Xtilde(D)}$
on $\Xtilde(D)$.
It is also denoted by
$\nbigc^{\infty\,\nil\,<D^{(3)}}_{\Xtilde(D)}$.
Put
$\Omega^{0,\bullet\,<D^{(3)}\leq D^{(2)}}
 _{\Xtilde(D)}:=
 \Omega^{0,\bullet}_{\Xtilde(D)}
\otimes_{\nbigc^{\infty}_{\Xtilde(D)}}
 \nbigc^{\infty<D^{(3)}\leq D^{(2)}}_{\Xtilde(D)}$.
\index{sheaf $\nbiga^{\nil\,<D^{(1)}}_{\Xtilde(D)}$}
\index{sheaf $\nbigc^{\infty\,\nil\,<D^{(3)}}_{\Xtilde(D)}$}
We will prove the following theorem
in \S\ref{subsection;09.12.4.13}.
(More refined claims will be proved.)

\begin{thm}
\mbox{{}}\label{thm;09.12.4.5}
\begin{itemize}
\item
 $\Omega^{0,\bullet\,<D^{(1)}\leq D^{(2)}}_{\Xtilde(D)}$
is naturally a c-soft resolution 
 of $\nbiga^{<D^{(1)}\leq D^{(2)}}_{\Xtilde(D)}$
in the case $D^{\circ}=\emptyset$.
\item
 The sheaves 
$\nbiga^{<D^{(1)}\leq D^{(2)}}_{\Xtilde(D)}$
and $\Omega^{0,\bullet\,<D^{(3)}\leq D^{(2)}}
 _{\Xtilde(D)}$
are flat over $\pi^{-1}\nbigo_X$.
\end{itemize}
\end{thm}

Let $D^{(i)}=\bigcup_{j\in \Lambda_i} D^{(i)}_j$
$(i=1,2)$ be the irreducible decomposition.
Fix $k\in \Lambda_1\sqcup \Lambda_2$.
We put
\[
 E^{(i)}:=\bigcup_{j\in\Lambda_i\setminus \{k\}}
 D^{(i)}_j
\quad
 (i=1,2).
\]
We put 
$E:=E^{(1)}\cup E^{(2)}$
and $E^{(3)}:=D^{(3)}$.
We have the naturally defined projection
$\rho:\Xtilde(D)\lrarr \Xtilde(E)$.
We will prove the following theorem
in \S\ref{subsection;09.12.4.40}.
\begin{thm}
\label{thm;09.12.4.30}
If $k\in \Lambda_1$,
the following naturally defined morphism
is an isomorphism:
\[
 \Omega_{\Xtilde(E)}^{0,\bullet\,<E^{(3)}\leq E^{(2)}}
\lrarr
 \rho_{\ast}
 \Omega_{\Xtilde(D)}^{0,\bullet\,<D^{(3)}\leq D^{(2)}}
\]
If $k\in\Lambda_2$,
the following naturally defined morphism
is a quasi-isomorphism:
\[
 \Omega_{\Xtilde(E)}^{0,\bullet\,<E^{(3)}\leq E^{(2)}}
 (\ast D^{(2)}_k)
\lrarr
 \rho_{\ast}
 \Omega_{\Xtilde(D)}^{0,\bullet\,<D^{(3)}\leq D^{(2)}}
\]
\end{thm}

\begin{cor}
\label{cor;14.1.19.3}
The natural morphism
\[
 \Omega_X^{0,\bullet\,<D^{(1)}}(\ast D^{(2)})
\lrarr
 \pi_{\ast}
 \Omega^{0,\bullet\,<D^{(1)}\leq D^{(2)}}
 _{\Xtilde(D)} 
\]
is a quasi-isomorphism.
In particular,
$R\pi_{\ast}\nbiga_{\Xtilde(D)}^{\nil}
\simeq
 \nbigo_X(\ast D)$.
\end{cor}

For the proof of the theorems,
we may assume $X=\Delta^n$
and $D=\bigcup_{i=1}^{\ell}\{z_i=0\}$
and $D^{\circ}=\bigcup_{i=\ell+1}^m\{z_i=0\}$,
where $1\leq \ell\leq m\leq n$.
We set $D_i:=\{z_i=0\}$ for $i=1,\ldots,m$.
We use the notation in \S\ref{subsection;14.1.17.30}.
For a subset $J\subset \ellsitabar$,
we set $\Jbar:=J\sqcup\bigl(\mbar\setminus\ellsitabar\bigr)$.

\subsection{Refinements}
\label{subsection;09.12.4.2}

For any locally closed 
real analytic subset $Z\subset \Xtilde(D)$,
we implicitly regard $\nbigo_{\Zhat}$
as a sheaf on $\Xtilde(D)$ in a natural way.
For any $I\sqcup J\subset\ellsitabar$,
let $\nbiga^{\nil\,<D(J)}_{\widehat{\pi^{-1}(D_I)}}$
denote the image of the following naturally 
defined morphism:
\[
 \nbigo_{\widehat{\pi^{-1}(D_I)}}
 ^{<D(J)}
 \otimes_{\cnum[z_1,\ldots,z_{\ell}]}
 \Nil(z_1,\ldots,z_{\ell})
\lrarr
 \nbigo_{\widehat{
 \pi^{-1}(D_I\setminus \del D_I)}}
 \otimes_{\cnum[z_i|i\in I]}
 \Nil(z_i\,|\,i\in I)
\]
\index{sheaf $\nbiga^{\nil\,<D(J)}_{\widehat{\pi^{-1}(D_I)}}$}
In the case $I=\emptyset$,
it is 
$\nbiga^{\nil<D(J)}_{\Xtilde(D)}$.
For $I\sqcup J\subset\ellsitabar$,
let $\nbiga^{\nil}_{\widehat{\pi^{-1}(D_I(J))}}$
denote the image of the following naturally
defined morphism:
\[
 \nbigo_{\widehat{\pi^{-1}(D_I(J))}}
 \otimes_{\cnum[z_1,\ldots,z_{\ell}]}
 \Nil(z_1,\ldots,z_{\ell})
\lrarr
\bigoplus_{j\in J}
 \nbigo_{\widehat{
 \pi^{-1}(D_{Ij}\setminus \del D_{Ij})}}
 \otimes_{\cnum[z_i|i\in Ij]}
 \Nil(z_i\,|\,i\in Ij)
\]
\index{sheaf $\nbiga^{\nil}_{\widehat{\pi^{-1}(D_I(J))}}$}
Here, $Ij:=I\sqcup\{j\}$.
In particular,
$\nbiga^{\nil}_{\widehat{\pi^{-1}(D(J))}}$
is the image of the following morphism:
\[
 \nbigo_{\widehat{\pi^{-1}(D(J))}}
 \otimes_{\cnum[z_1,\ldots,z_{\ell}]}
 \Nil(z_1,\ldots,z_{\ell})
\lrarr
 \bigoplus_{j\in J}
 \nbigo_{\widehat{\pi^{-1}(D_j\setminus\del D_j)}}
 \otimes_{\cnum[z_j]}
 \Nil\bigl(z_j\bigr)
\]
\index{sheaf $\nbiga^{\nil}_{\widehat{\pi^{-1}(D(J))}}$}
Let 
$\nbiga^{\nil\,<D(J)}_{\widehat{\pi^{-1}(D_I)},T,N}$
and $\nbiga^{\nil}_{\widehat{\pi^{-1}(D_I(J))},T,N}$
be the sheaves obtained from
$\Nil_{T,N}(z_1,\ldots,z_{\ell})$
instead of $\Nil(z_1,\ldots,z_{\ell})$.
\index{sheaf $\nbiga^{\nil\,<D(J)}_{\widehat{\pi^{-1}(D_I)},T,N}$}
\index{sheaf $\nbiga^{\nil}_{\widehat{\pi^{-1}(D_I(J))},T,N}$}
For $T\subset T'$ and $N\leq N'$,
we have natural inclusions
$\nbiga^{\nil\,<D(J)}_{\widehat{\pi^{-1}(D_I)},T,N}
\subset
  \nbiga^{\nil\,<D(J)}_{\widehat{\pi^{-1}(D_I)},T',N'}$
and 
$\nbiga^{\nil}_{\widehat{\pi^{-1}(D_I(J))},T,N}
\subset
 \nbiga^{\nil}_{\widehat{\pi^{-1}(D_I(J))},T',N'}$.
We have the following natural isomorphisms:
\begin{equation}
\label{eq;09.10.30.15}
 \nbiga^{\nil\,<D(J)}_{\widehat{\pi^{-1}(D_I)}}
\simeq\varinjlim
 \nbiga^{\nil\,<D(J)}_{\widehat{\pi^{-1}(D_I)},T,N}
\quad\quad
 \nbiga^{\nil}_{\widehat{\pi^{-1}(D_I(J))}}
\simeq
 \varinjlim
 \nbiga^{\nil}_{\widehat{\pi^{-1}(D_I(J))},T,N}
\end{equation}

Let $q_I:\pi^{-1}(D_I)\lrarr \Dtilde_I(\del D_I)$
denote the projection.
Let $\pi_I:\Dtilde_I(\del D_I)\lrarr D_I$
be the real blow up.
Then, we have 
\begin{equation}
\label{eq;09.10.30.1}
 \nbiga^{\nil\,<D(J)}_{\widehat{\pi^{-1}(D_I)},T,N}
=q_I^{-1}\nbiga^{\nil\,<D_I(J)}
 _{\Dtilde_I(\del D_I),T,N}
 [\![z_i\,|\,i\in I]\!]\otimes_{\cnum[z_i|i\in I]}
 \Nil_{T,N}(z_i\,|\,i\in I)
\end{equation}
\begin{equation}
 \label{eq;09.10.29.10}
 \nbiga^{\nil}_{\widehat{\pi^{-1}(D_I(J))},T,N}
=q_I^{-1}
 \nbiga^{\nil}_{\widehat{\pi_I^{-1}(D_I(J))},T,N}
 [\![z_i\,|\,i\in I]\!]\otimes_{\cnum[z_i|i\in I]}
 \Nil_{T,N}(z_i\,|\,i\in I)
\end{equation}

\subsection{Specialization}
\label{subsection;09.10.25.1}

Let us construct a morphism
$\nbiga^{\nil}_{\widehat{\pi^{-1}(D_I)}}\lrarr
 \nbiga^{\nil}_{\widehat{\pi^{-1}(D_I(J))}}$
for any $I\sqcup J\subset\ellsitabar$.
First, let us construct
$\nbiga^{\nil}_{\Xtilde(D)}\lrarr
 \nbiga^{\nil}_{\widehat{\pi^{-1}(D)}}$
in the case $D=D_1$.
Let $\Phi$ denote the natural morphism
$\Phi:\nbigo_{\Xtilde(D)}\otimes\Nil(z_1)
\lrarr
 \jtilde_{\ast}\nbigo_{X-D}$,
where $\jtilde:X-D\lrarr \Xtilde(D)$.

\begin{lem}
\label{lem;09.10.16.50}
Assume that $D=D_1$.
Let $\nbigs\subset\cnum$ be a finite subset
such that the induced map
$\nbigs\lrarr\cnum/\seisuu$ is injective.
Assume that we are given
 $f=\sum_{\alpha\in \nbigs}\sum_{j=0}^M
 f_{\alpha,j}\otimes
 \varphi_{\alpha,j}(z_1)
 \in \nbigo_{\Xtilde(D)}\otimes \Nil(z_1)$
such that $\Phi(f)\in\nbigo^{<D}_{\Xtilde(D)}$.
Then, we have
$f_{\alpha,j}\in\nbigo^{<D}_{\Xtilde(D)}$.
In particular,
we have the well defined map
$\nbiga^{\nil}_{\Xtilde(D)}\lrarr
 \nbiga^{\nil}_{\piinverseDhat}$
in the case $D=\{z_1=0\}$.
\end{lem}
\pf
Let us consider the growth order of
$f_{\alpha,j}\,z_1^{\alpha}(\log z_1)^j$.
For the polar coordinate system
$z_1=re^{\sqrt{-1}\theta}$,
we have
$z_1^{\alpha}
=\exp\bigl( 
 \beta\log r
 -\gamma\theta
+\sqrt{-1}(\gamma\log r+\beta\theta)
\bigr)$,
where
$\beta=\Re\alpha$ and
$\gamma=\Image\alpha$.
Let $V$ be the set of 
$(\alpha,j)\in \nbigs\times\seisuu_{\geq 0}$
such that 
$f_{\alpha,j}$ is not contained in
$\nbigo_{\Xtilde(D)}^{<D}$.
We will derive a contradiction
by assuming $V\neq\emptyset$.
For each $(\alpha,j)\in V$,
there exists a unique integer $m(\alpha,j)$
such that 
(i) 
$h_{\alpha,j}:=
z_1^{-m(\alpha,j)}f_{\alpha,j}
 \in \nbigo_{\Xtilde(D)}$,
(ii) $h_{\alpha,j|\pi^{-1}(D)}$
is not constantly $0$.
We set
\[
 \kappa:=\max_{(\alpha,j)\in V}
 \bigl\{
 \Re\alpha+m(\alpha,j)
 \bigr\},
\quad
 S:=\bigl\{
 (\alpha,j)\in V\,\big|\,
 \Re\alpha+m(\alpha,j)=\kappa
 \bigr\}.
\]
For $(\alpha_1,j_1),(\alpha_2,j_2)\in S$,
we have
$\Re\alpha_1=\Re\alpha_2$
and $m(\alpha_1,j_1)=m(\alpha_2,j_2)$.
We also have 
$\Image\alpha_1\neq\Image\alpha_2$
if $\alpha_1\neq\alpha_2$.
We obtain the following estimate for some $\epsilon>0$:
\begin{multline}
 \label{eq;09.10.16.30}
 \sum_{(\alpha,j)\in V}
 h_{\alpha,j|\pi^{-1}(D)}\,
 z_1^{\alpha+m(\alpha,j)}
 (\log z_1)^j
= \\
 r^{\kappa}\Bigl(
 \sum_{(\alpha,j)\in V}
 h_{\alpha,j|\pi^{-1}(D)}\,
 e^{-\Image\alpha\theta
+\sqrt{-1}(\Image\alpha\log r+\Re\alpha\theta)}
(\log z_1)^j
 \Bigr) 
=O(r^{\kappa+\epsilon})
\end{multline}
Let us deduce that
$h_{\alpha,j|\pi^{-1}(D)}$
are constantly $0$ from (\ref{eq;09.10.16.30}).
Assume the contrary.
Let $Q\in \pi^{-1}(D)$
at which $h_{\alpha,j}(Q)\neq 0$
for one of $(\alpha,j)\in V$.
We may assume $\theta(Q)=0$.
We obtain the following from (\ref{eq;09.10.16.30}):
\begin{equation}
 \label{eq;09.12.3.1}
\sum_{(\alpha,j)\in V} 
 h_{\alpha,j}(Q)\cdot
 e^{\sqrt{-1}\Image\alpha\log r}
 (\log r)^j=O(r^{\epsilon})
\end{equation}
But, for any $\delta>0$,
we can take $0<r<\delta$
such that the amplitudes of the complex numbers
\[
 (-1)^jh_{\alpha,j}(Q)\,e^{\sqrt{-1}\Image\alpha\log r}
\quad
 (\alpha,j)\in V
\]
are sufficiently close,
which contradicts with (\ref{eq;09.12.3.1}).
Hence, $h_{\alpha,j}$ $(\alpha,j)\in V$
are constantly $0$.
Thus, we obtain Lemma \ref{lem;09.10.16.50}.
\hfill\qed

\vspace{.1in}

Let us return to the general case.
We take $\nbigs\subset\cnum$ 
such that the induced map
$\nbigs\lrarr\cnum/\seisuu$ is bijective.
Let $q_i:(\nbigs\times\seisuu)^{\ell}
 \lrarr \nbigs\times\seisuu$ be the projection
onto the $i$-th component,
and $\pi_i:(\nbigs\times\seisuu)^{\ell}
\lrarr (\nbigs\times\seisuu)^{\ell-1}$
be the projection forgetting the $i$-th component.
For a given 
\[
 \sum_{(\vecalpha,\veck)\in
 \nbigs^{\ell}\times\seisuu_{\geq\,0}^{\ell}}
 A_{\vecalpha,\veck}\otimes
 \varphi_{\vecalpha,\veck}
 \in \nbigo_{\Xtilde(D)}\otimes 
 \Nil(z_1,\ldots,z_{\ell}),
\]
we set
$\lefttop{i}F_{\beta,j}:=
 \sum_{q_i(\vecalpha,\veck)=(\beta,j)}
 A_{\vecalpha,\veck}\cdot
 \varphi_{\pi_i(\vecalpha,\veck)} (z_j|\,j\!\neq\! i)$.
Put $i^c:=\ellsitabar-\{i\}$.
If $\sum A_{\vecalpha,\veck}\cdot
 \varphi_{\vecalpha,\veck}\in
 \nbigo_{\Xtilde(D)\setminus
 \pi^{-1}(D(i^c))}^{<D_i}$,
we obtain
$\lefttop{i}F_{\beta,j|
 \widehat{\pi^{-1}(D_i\setminus\del D_i)}}=0$
by applying Lemma \ref{lem;09.10.16.50} to
$\sum \lefttop{i}F_{\beta,j}\cdot
 \varphi_{\beta,j}(z_i)$.
It implies that 
the morphism
\[
 \nbigo_{\Xtilde(D)}\otimes\Nil(z_1,\ldots,z_{\ell})
\lrarr
 \nbigo_{\widehat{\pi^{-1}(D_i)}}
 \otimes \Nil(z_1,\ldots,z_{\ell})
\lrarr
\nbiga^{\nil}_{\widehat{\pi^{-1}(D_i)}}
\]
factors through $\nbiga^{\nil}_{\Xtilde(D)}$.
Hence, we have a well defined morphism
$\nbiga^{\nil}_{\Xtilde(D)}
\lrarr
 \nbiga^{\nil}_{\widehat{\pi^{-1}(D_i)}}$.
By construction,
it is an epimorphism.
We also obtain that the following morphism
factors through $\nbiga^{\nil}_{\Xtilde(D)}$:
\[
 \nbigo_{\Xtilde(D)}\otimes\Nil(z_1,\ldots,z_{\ell})
\lrarr
 \nbigo_{\widehat{\pi^{-1}(D(I))}}
 \otimes \Nil(z_1,\ldots,z_{\ell})
\lrarr
 \nbiga^{\nil}_{\widehat{\pi^{-1}(D(I))}}
\subset
\bigoplus_{i\in I}
\nbiga^{\nil}_{\widehat{\pi^{-1}(D_i)}}
\]
Hence, we obtain the well defined morphism
$\nbiga^{\nil}_{\Xtilde(D)}
\lrarr
 \nbiga^{\nil}_{\widehat{\pi^{-1}(D(I))}}$.
We also obtain 
$\nbiga^{\nil}_{\Xtilde(D),T,N}
\lrarr
 \nbiga^{\nil}_{\widehat{\pi^{-1}(D(I))},T,N}$.
They are surjective by construction.
By using (\ref{eq;09.10.30.15}),
(\ref{eq;09.10.30.1})
and (\ref{eq;09.10.29.10}),
we also obtain epimorphisms
$\nbiga^{\nil}_{\widehat{\pi^{-1}(D_I)}}
\lrarr
 \nbiga^{\nil}_{\widehat{\pi^{-1}(D_I(J))}}$
and 
$\nbiga^{\nil}_{\widehat{\pi^{-1}(D_I)},T,N}
\lrarr
 \nbiga^{\nil}_{\widehat{\pi^{-1}(D_I(J))},T,N}$.

\begin{lem}
\label{lem;09.10.30.2}
We have the following:
\[
 \nbiga^{\nil<D(J)}_{\widehat{\pi^{-1}(D_I)}}
=\Ker\Bigl(
 \nbiga^{\nil}_{\widehat{\pi^{-1}(D_I)}}
\lrarr
 \nbiga^{\nil}_{\widehat{\pi^{-1}(D_I(J))}}
 \Bigr)
\]
\[
\nbiga^{\nil<D(J)}_{\widehat{\pi^{-1}(D_I)},T,N}
=\Ker\Bigl(
 \nbiga^{\nil}_{\widehat{\pi^{-1}(D_I)},T,N}
\lrarr
 \nbiga^{\nil}_{\widehat{\pi^{-1}(D_I(J))},T,N}
 \Bigr)
\]
\end{lem}
\pf
The implication $\subset$ is clear.
Let us prove the converse.
First, we consider the case $I=\emptyset$.
Let $f=\sum A_{\vecalpha,\veck}\,
 \varphi_{\vecalpha,\veck}$
be any section of
$\Ker\Bigl(
 \nbiga^{\nil}_{\Xtilde(D)}
\lrarr
 \nbiga^{\nil}_{\widehat{\pi^{-1}(D(J))}}
 \Bigr)$.
Let us prove the following equality
on $\widehat{\pi^{-1}(D_K-\del D_K)}$
for any subset $K\subset \ellsitabar$
such that $K\cap J\neq\emptyset$:
\begin{equation}
 \label{eq;09.10.16.100}
 \sum_{q_K(\vecalpha,\veck)=(\vecbeta,\vecj)}
 A_{\vecalpha,\veck|\widehat{\pi^{-1}(D_K)}}
 \prod_{i\not\in K}\varphi_{\alpha_i,k_i}(z_i)=0
\end{equation}
We use an induction on $|K|$.
In the case $|K|=1$,
it follows from the assumption.
Let $K=K'\sqcup\{j\}$.
Assume that we have already known
(\ref{eq;09.10.16.100}) for $K'$.
By using Lemma \ref{lem;09.10.16.50},
we obtain the claim for $K$.
As a special case of (\ref{eq;09.10.16.100}),
we have
$A_{\vecalpha,\veck|
 \widehat{\pi^{-1}(D_{\ellsitabar})}}=0$.

Note that the expression of $f$ is not unique.
We would like to replace $A_{\vecalpha,\veck}$
such that the following holds:
\begin{description}
\item[P(m)]
 $A_{\vecalpha,\veck|\widehat{\pi^{-1}(D_K)}}=0$
 if $|K|\geq m$ and $K\cap J\neq\emptyset$.
\end{description}
We use a descending induction on $m$.
In the case $m=\ell$, it holds
as was already proved.
Assume that $P(m+1)$ holds.
Take $K\subset\ellsitabar$ such that
$|K|=m$ and $K\cap J\neq\emptyset$.
We have
\[
 A_{\vecalpha,\veck|\widehat{\pi^{-1}(D_K)}}
 \prod_{i\not\in K}\varphi_{\alpha_i,k_i}(z_i)
 \in 
 \nbigo^{<D(K^c)}_{\widehat{\pi^{-1}(D_K)}}.
\]
By a generalized Borel-Ritt theorem
due to Majima and Sabbah,
we can take $G_{\vecalpha,\veck}
 \in \nbigo^{<D(K^c)}_{\Xtilde(D)}$
satisfying
$G_{\vecalpha,\veck|\widehat{\pi^{-1}(D_K)}}
=A_{\vecalpha,\veck|\widehat{\pi^{-1}(D_K)}}
 \prod_{i\not\in K}
 \varphi_{\alpha_i,k_i}(z_i)$.
By (\ref{eq;09.10.16.100}),
the following holds:
\[
 \sum_{q_K(\vecalpha,\veck)=(\vecbeta,\vecj)}
 G_{\vecalpha,\veck|\widehat{\pi^{-1}(D_K)}}=0
\]
We have the following equality:
\begin{multline*}
 f=
 \sum_{\vecalpha,\veck}
 \Bigl(
 A_{\vecalpha,\veck}
-\frac{G_{\vecalpha,\veck}}
 {\prod_{i\not\in K}\varphi_{\alpha_i,k_i}(z_i)}
 \Bigr)
 \varphi_{\vecalpha,\veck}(z_1,\ldots,z_{\ell})
 \\
+\sum_{\vecbeta,\vecj}
 \Bigl(
 \sum_{q_K(\vecalpha,\veck)=(\vecbeta,\vecj)}
 G_{\vecalpha,\veck}
 \Bigr)
 \varphi_{\vecbeta,\vecj}(z_i|i\in K)
\end{multline*}
Note that 
 $\sum_{q_K(\vecalpha,\veck)=(\vecbeta,\vecj)}
 G_{\vecalpha,\veck}$ is $0$
 on $\widehat{\pi^{-1}(D_{K})}\cup
 \widehat{\pi^{-1}(D(K^c))}$.
 In particular, it is $0$
 on $\bigcup_{|K_1|=m}
 \widehat{\pi^{-1}(D_{K_1})}$.
By construction,
 $A_{\vecalpha,\veck}
-G_{\vecalpha,\veck}
 \prod_{i\not\in K}\varphi_{\alpha_i,k_i}(z_i)^{-1}$
 vanishes on $\widehat{\pi^{-1}(D_K)}$.
 Moreover, if 
 $A_{\vecalpha,\veck|\widehat{\pi^{-1}(D_L)}}=0$
 for some $|L|=m$ with $L\cap J\neq\emptyset$,
$A_{\vecalpha,\veck}
-G_{\vecalpha,\veck}
 \prod_{i\not\in K}\varphi_{\alpha_i,k_i}(z_i)^{-1}$
also vanishes on $\widehat{\pi^{-1}(D_L)}$.
Hence, by applying the above procedure
to each $K$ 
satisfying $|K|=m$
and $K\cap J\neq\emptyset$,
we can arrive at $P(m)$.
The status $P(0)$ means
$f=\sum A_{\vecalpha,\veck}\,
 \varphi_{\vecalpha,\veck}$
with
$A_{\vecalpha,\veck}\in
 \nbigo^{<D(J)}_{\Xtilde(D)}$,
which implies that
$f\in \nbiga^{\nil\,<D(J)}_{\Xtilde(D)}$.
Thus, we are done in the case $I=\emptyset$.
We can reduce
the general case to the case $I=\emptyset$
by using (\ref{eq;09.10.30.15}),
(\ref{eq;09.10.30.1}) and
(\ref{eq;09.10.29.10}).
Thus, the proof of Lemma \ref{lem;09.10.30.2}
is finished.
\hfill\qed

\subsection{A resolution}

For $I\subset J$,
put $\nbigt(m,I,J):=\bigl\{
 K\subset J\,\big|\,
 I\subset K,\,|K|=|I|+m+1\bigr\}$
for $m\geq 0$.
We set 
\[
\nbigk^{m}\Bigl(
 \nbiga^{\nil}_{\widehat{\pi^{-1}(D_I(J))}}
 \Bigr)
:=\bigoplus_{K\in\nbigt(m,I,J)}
 \nbiga^{\nil}_{\widehat{\pi^{-1}(D_K)}}. 
\]
We obtain a complex
$\nbigk^{\bullet}\Bigl(
 \nbiga^{\nil}_{\widehat{\pi^{-1}(D_I(J))}}
 \Bigr)$
as in \S\ref{subsection;09.10.29.10}.
\begin{lem}
\label{lem;09.10.25.3}
The $0$-th cohomology of 
$\nbigk^{\bullet}\Bigl(
 \nbiga^{\nil}_{\widehat{\pi^{-1}(D_I(J))}}
 \Bigr)$ is 
$\nbiga^{\nil}_{\widehat{\pi^{-1}(D_I(J))}}$,
and the higher cohomology sheaves are $0$.
A similar claim holds for
$\nbiga^{\nil}_{\widehat{\pi^{-1}(D_I(J))},T,N}$.
\end{lem}
\pf
It is enough to consider the issue for
$\nbigk^{\bullet}\Bigl(
 \nbiga^{\nil}_{\widehat{\pi^{-1}(D_I(J)),T,N}}
 \Bigr)$.
First, let us consider the case $I=\emptyset$.
We use an induction on $|J|$
and the dimension of $X$.
The cases $|J|=1$ or $\dim X=1$ are clear.
Let $J=J_0\sqcup\{j\}$.
Assume that the claim holds for $J_0$.
We set
$\nbigl^m_{T,N}:=
 \bigoplus_{|K|=m+1,j\in K\subset J}
 \nbiga^{\nil}_{\widehat{\pi^{-1}(D_K)},T,N}$.
We have the exact sequence:
\[
 0\lrarr \nbigl^{\bullet}_{T,N}\lrarr
 \nbigk^{\bullet}\Bigl(
 \nbiga^{\nil}_{\widehat{\pi^{-1}(D(J))},T,N}
 \Bigr)
 \lrarr
 \nbigk^{\bullet}\Bigl(
 \nbiga^{\nil}_{\widehat{\pi^{-1}(D(J_0))},T,N}
 \Bigr)
 \lrarr 0
\]
Let $q_j:\pi^{-1}(D_j)\lrarr \Dtilde_j(\del D_j)$
and $\pi_j:\Dtilde_j(\del D_j)\lrarr D_j$
be the projections.
We have a natural isomorphism:
{\small
\[
 \nbigl^{\bullet}_{T,N}
 \simeq
 \Cone\Bigl(
 \nbiga^{\nil}_{\widehat{\pi^{-1}(D_j)},T,N}
\lrarr
 q_j^{-1}\nbigk^{\bullet}\Bigl(
 \nbiga^{\nil}_{
 \widehat{\pi_j^{-1}(D_j\cap D(J_0))},T,N}
 \Bigr)[\![z_j]\!]
 \otimes_{\cnum[z_j]}\Nil_{T,N}(z_j)
\Bigr)[-1]
\]
}
By the inductive assumption,
we obtain the vanishing of the higher cohomology
sheaves of $\nbigl^{\bullet}_{T,N}$
and $\nbigk^{\bullet}\Bigl(
 \nbiga^{\nil}_{
 \widehat{\pi^{-1}(D(J_0))},T,N}
 \Bigr)$.
Hence, we obtain the vanishing of
the higher cohomology
of $\nbigk^{\bullet}\Bigl(
 \nbiga^{\nil}_{\widehat{\pi^{-1}(D(J))},T,N}
 \Bigr)$.
The calculation of the $0$-th cohomology is easy.
The general case can be easily reduced to the case
$I=\emptyset$
by (\ref{eq;09.10.30.15}),
(\ref{eq;09.10.30.1}) and (\ref{eq;09.10.29.10}).
\hfill\qed

\subsection{The $C^{\infty}$-version}
\label{subsection;09.12.4.3}

Let $Y$ be a $C^{\infty}$-manifold.
For $I\sqcup J\subset \ellsitabar$,
let $\nbigc^{\infty\,\nil\,<D(\Jbar)}_{
 \widehat{\pi^{-1}(D_I)}\times Y}$
denote the image of the following morphism:
\[
 \nbigc^{\infty\,<D(\Jbar)}
 _{\widehat{\pi^{-1}(D_I)}\times Y}
 \otimes_{\cnum[z_i|i\in J^c]}
 \Nil(z_i\,|\,i\in J^c)
\lrarr
 \nbigc^{\infty\,<D(\Jbar)}_{\widehat{
 \pi^{-1}(D_I\setminus \del D_I)}\times Y}
 \otimes_{\cnum[z_i|i\in I]}
 \Nil(z_i\,|\,i\in I)
\]
Let 
$\nbigc^{\infty\,\nil\,<D^{\circ}}_{
 \widehat{\pi^{-1}(D_I(J))}\times Y}$
be the image of the following morphism:
\begin{multline*}
 \nbigc^{\infty\,<D^{\circ}}_{
 \widehat{\pi^{-1}(D_I(J))}\times Y}
 \otimes_{\cnum[z_1,\ldots,z_{\ell}]}
 \Nil(z_1,\ldots,z_{\ell})
\lrarr \\
 \bigoplus_{j\in J}
 \nbigc^{\infty\,<D^{\circ}}_{
 \widehat{\pi^{-1}(D_{Ij}-\del D_{Ij})}\times Y}
 \otimes_{\cnum[z_i|i\in Ij]}
 \Nil\bigl(z_i\,\big|\,i\in Ij\bigr)
\end{multline*}
In particular,
$\nbigc^{\infty\,\nil\,<D^{\circ}}_{
 \widehat{\pi^{-1}(D(J))}\times Y}$
is the image of the following morphism:
\[
  \nbigc^{\infty\,<D^{\circ}}_{
 \widehat{\pi^{-1}(D(J))}\times Y}
 \otimes_{\cnum[z_1,\ldots,z_{\ell}]}
 \Nil(z_1,\ldots,z_{\ell})
\lrarr
 \bigoplus_{j\in J}
 \nbigc^{\infty\,<D^{\circ}}_{
 \widehat{\pi^{-1}(D_{j}-\del D_{j})}\times Y}
 \otimes_{\cnum[z_j]}
 \Nil\bigl(z_j\bigr)
\]
Similarly,
$\nbigc^{\infty\,\nil\,<D(\Jbar)}
 _{\widehat{\pi^{-1}(D_I)}\times Y,T,N}$
and 
$\nbigc^{\infty\,\nil\,<D^{\circ}}
 _{\widehat{\pi^{-1}(D_I(J))}\times Y,T,N}$
denote the sheaves obtained from
$\Nil_{T,N}(z_1,\ldots,z_{\ell})$
instead of $\Nil(z_1,\ldots,z_{\ell})$.
We have 
\begin{equation}
 \label{eq;09.10.30.10}
\nbigc^{\infty\,\nil\,<D(\Jbar)}_{
 \widehat{\pi^{-1}(D_I)}\times Y,T,N}
=
 \nbigc^{\infty\,\nil\,<D_I(\Jbar)}_{
 \Dtilde_I(\del D_I)\times
 (S^1)^{|I|}\times Y,T,N}
 [\![z_i|i\in I]\!]
 \otimes_{\cnum[z_i|i\in I]}
 \Nil_{T,N}(z_i|i\in I)
\end{equation}
\begin{equation}
 \label{eq;09.10.30.11}
 \nbigc^{\infty\,\nil\,<D^{\circ}}_{
 \widehat{\pi^{-1}(D_I(J))}\times Y,T,N}
=
 \nbigc^{\infty\,\nil\,<D^{\circ}\cap D_I}_{
 \widehat{\pi_I^{-1}(D_I(J))}
 \times
 (S^1)^{|I|}\times Y,T,N}
 [\![z_i|i\in I]\!]
 \otimes_{\cnum[z_i|i\in I]}
 \Nil_{T,N}(z_i|i\in I)
\end{equation}

By the argument in \S\ref{subsection;09.10.25.1},
we obtain the well defined surjective morphisms:
\begin{equation}
 \label{eq;09.10.30.20}
 \nbigc^{\infty\,\nil\,<D^{\circ}}
 _{\widehat{\pi^{-1}(D_I)}\times Y}
\lrarr
 \nbigc^{\infty\,\nil\,<D^{\circ}}
 _{\widehat{\pi^{-1}(D_I(J))}\times Y},
 \quad\quad
 \nbigc^{\infty\,\nil\,<D^{\circ}}_{
 \widehat{\pi^{-1}(D_I)}\times Y,T,N}
\lrarr
 \nbigc^{\infty\,\nil\,<D^{\circ}}_{
 \widehat{\pi^{-1}(D_I(J))}\times Y,T,N}
\end{equation}
By the argument in the proof of 
Lemma \ref{lem;09.10.30.2},
we can prove that the kernels of the morphisms
in (\ref{eq;09.10.30.20})
are 
$\nbigc^{\infty\,\nil\,<D(\Jbar)}_{
 \widehat{\pi^{-1}(D_I)}\times Y}$
and 
$\nbigc^{\infty\,\nil\,<D(\Jbar)}_{
 \widehat{\pi^{-1}(D_I)}\times Y,T,N}$,
respectively.

We set
$\nbigk^{m}\Bigl(
 \nbigc^{\infty\,\nil\,<D^{\circ}}_{
 \widehat{\pi^{-1}(D_I(J))}\times Y}
 \Bigr)
:=\bigoplus_{K\in\nbigt(m,I,J)}
 \nbigc^{\infty\,\nil\,<D^{\circ}}
 _{\widehat{\pi^{-1}(D_K)}\times Y}$.
We obtain a complex
$\nbigk^{\bullet}\Bigl(
 \nbigc^{\infty\,\nil\,<D^{\circ}}_{
 \widehat{\pi^{-1}(D_I(J))}\times Y}
 \Bigr)$.
It is easy to see that
the $0$-th cohomology is
$\nbigc^{\infty\,\nil\,<D^{\circ}}_{
 \widehat{\pi^{-1}(D_I(J))}\times Y}$.
By using an argument in the proof of
Lemma \ref{lem;09.10.25.3},
we can prove the vanishing of
the higher cohomology.
Similar claims hold for
$\nbigk^{\bullet}\Bigl(
 \nbigc^{\infty\,\nil\,<D^{\circ}}_{
 \widehat{\pi^{-1}(D_I(J))}\times Y,T,N}
 \Bigr)$.

\subsection{Proof of Theorem \ref{thm;09.12.4.5}}
\label{subsection;09.12.4.13}

In this subsection,
we do not consider $D^{\circ}$.
We put
$\Omega^{0,\bullet\,\nil\,<D(J)}
 _{\widehat{\pi^{-1}(D_I)}}:=
 \Omega^{0,\bullet}_{\Xtilde(D)}
\otimes_{\nbigc^{\infty}_{\Xtilde(D)}}
 \nbigc^{\infty\nil<D(J)}_{\widehat{\pi^{-1}(D_I)}}$
and 
$\Omega^{0,\bullet\,\nil}_{\widehat{\pi^{-1}(D_I(J))}}
 :=
 \Omega^{0,\bullet}_{\Xtilde(D)}
 \otimes_{\nbigc^{\infty}_{\Xtilde(D)}}
 \nbigc^{\infty\,\nil}_{\widehat{\pi^{-1}(D_I(J))}}$.
We use the symbols
$\Omega^{0,\bullet\,\nil\,<D(J)}
 _{\widehat{\pi^{-1}(D_I)},T,N}$
and 
$\Omega^{0,\bullet\,\nil}
 _{\widehat{\pi^{-1}(D_I(J))},T,N}$
with a similar meaning.
The following proposition implies
the first claim of Theorem {\rm\ref{thm;09.12.4.5}}.

\begin{prop}
\label{prop;09.12.4.4}
The complexes
$\Omega^{0,\bullet\,\nil\,<D(J)}
 _{\widehat{\pi^{-1}(D_I)}}$
and 
$\Omega^{0,\bullet\,\nil}
 _{\widehat{\pi^{-1}(D_I(J))}}$
are c-soft resolutions of
the sheaves
$\nbiga^{\nil\,<D(J)}_{\widehat{\pi^{-1}(D_I)}}$
and 
$\nbiga^{\nil}_{\widehat{\pi^{-1}(D_I(J))}}$
respectively.
Similar claims hold for
$\nbiga^{\nil\,<D(J)}
 _{\widehat{\pi^{-1}(D_I)},T,N}$
and 
$\nbiga^{\nil}
 _{\widehat{\pi^{-1}(D_I(J))},T,N}$.
\end{prop}
\pf
We use an induction on $\dim X$.
In the case $\dim X=0$,
the claim is trivial.
Let us prove the claim for
$\widehat{\pi^{-1}(D_I)}$.
For $I\neq\emptyset$,
let $q_I:\pi^{-1}(D_I)\lrarr \Dtilde_I(\del D_I)$
denote the naturally induced morphism.
We put
$\Nil_{T,N}(I):=\Nil_{T,N}(z_i|i\in I)$.
By using the inductive assumption
and a formal calculation as in Lemma
\ref{lem;09.10.22.11},
we can prove that the following morphisms
are quasi-isomorphisms:
\begin{multline}
 q_I^{-1}
 \nbiga^{\nil\,<D_I(J)}_{\Dtilde_I(\del D_I),T,N}
 [\![z_i|i\in I]\!]
\otimes\Nil_{T,N}(I)
\lrarr \\
 q_I^{-1}
 \Omega^{0,\bullet\,\nil\,<D_I(J)}
 _{\Dtilde_I(\del D_I),T,N}[\![z_i|i\in I]\!]
 \otimes\Nil_{T,N}(I)
\lrarr
 \Omega^{0,\bullet\,\nil\,<D(J)}
 _{\widehat{\pi^{-1}(D_I)},T,N}
\end{multline}
It implies the claim for
$\nbiga^{\nil\,<D(J)}_{\widehat{\pi^{-1}(D_I)},T,N}$.
We obtain the claim for
$\nbiga^{\nil\,<D(J)}_{\widehat{\pi^{-1}(D_I)}}$
from (\ref{eq;09.10.30.15}).
For any subset $I\subset \ellsitabar$
($I$ can be $\emptyset$),
by using the resolutions
$\nbigk^{\bullet}\bigl(
 \nbiga^{\nil}_{\widehat{\pi^{-1}(D_I(J))}}
 \bigr)$
and 
$\nbigk^{\bullet}\bigl(
 \Omega^{0,\bullet\,\nil}_{\widehat{\pi^{-1}(D_I(J))}}
 \bigr)$,
we can reduce the claim for
$\nbiga^{\nil}_{\widehat{\pi^{-1}(D_I(J))}}$
to the claims for 
$\nbiga^{\nil}_{\widehat{\pi^{-1}(D_K)}}$
($I\subsetneq K$).
The claim for
$\nbiga^{\nil}_{\widehat{\pi^{-1}(D_I(J))},T,N}$
can be obtained in a similar way.
By using the exact sequences
\begin{eqnarray*}
0\lrarr 
 \Omega_{\Xtilde(D)}^{0,\bullet\,<D}
 \lrarr
 \Omega_{\Xtilde(D)}^{0,\bullet\,\nil}
 \lrarr
 \Omega_{\widehat{\pi^{-1}(D)}}^{0,\bullet\,\nil}
\lrarr 0,
\\
0\lrarr\nbigo_{\Xtilde(D)}^{<D}
\lrarr \nbiga^{\nil}_{\Xtilde(D)}
\lrarr \nbiga^{\nil}_{\widehat{\pi^{-1}(D)}}
\lrarr 0,
\end{eqnarray*}
we obtain the claim for 
$\nbiga^{\nil}_{\Xtilde(D)}$.
By using the exact sequences
\begin{eqnarray*}
 0\lrarr
 \Omega_{\Xtilde(D)}^{0,\bullet\,<D(J)}
 \lrarr
 \Omega_{\Xtilde(D)}^{0,\bullet\,\nil}
 \lrarr
 \Omega_{\widehat{\pi^{-1}(D(J))}}^{0,\bullet\,\nil}
\lrarr 0,
\\
0\lrarr\nbiga_{\Xtilde(D)}^{\nil\,<D(J)}
\lrarr \nbiga^{\nil}_{\Xtilde(D)}
\lrarr \nbiga^{\nil}_{\widehat{\pi^{-1}(D(J))}}
\lrarr 0,
\end{eqnarray*}
we obtain the claim for
$\nbiga^{\nil\,<D(J)}_{\Xtilde(D)}$.
The claims for
$\nbiga^{\nil}_{\Xtilde(D),T,N}$
and
$\nbiga^{\nil\,<D(J)}_{\Xtilde(D),T,N}$
can be obtained similarly.
\hfill\qed

\vspace{.1in}

The following proposition implies 
the second claim of Theorem \ref{thm;09.12.4.5}.
\begin{prop}
\label{prop;09.10.25.40}
$\nbigc^{\infty\,\nil\,<D(\Jbar)}_{
 \widehat{\pi^{-1}(D_I)}}$,
$\nbigc^{\infty\,\nil\,<D^{\circ}}_{
 \widehat{\pi^{-1}(D_I(J))}}$,
$\nbiga^{\nil\,<D(J)}_{\widehat{\pi^{-1}(D_I)}}$
and 
$\nbiga^{\nil}_{\widehat{\pi^{-1}(D_I(J))}}$
are flat over $\pi^{-1}\nbigo_X$.
Similar claims hold for
$\nbigc^{\infty\,\nil\,<D(\Jbar)}
 _{\widehat{\pi^{-1}(D_I)},T,N}$,
$\nbigc^{\infty\,\nil\,<D^{\circ}}
 _{\widehat{\pi^{-1}(D_I(J))},T,N}$,
$\nbiga^{\nil\,<D(J)}
 _{\widehat{\pi^{-1}(D_I)},T,N}$
and 
$\nbiga^{\nil}
 _{\widehat{\pi^{-1}(D_I(J))},T,N}$
are also flat over $\pi^{-1}\nbigo_X$.
\end{prop}
\pf
We have
$\nbigc^{\infty\nil\,<D(I^c)}
 _{\widehat{\pi^{-1}(D_I)}}
=\nbigc^{\infty<D(I^c)}_{\widehat{\pi^{-1}(D_I)}}
\otimes_{\cnum[z_i|i\in I]}\Nil(z_i|\,i\in I)$,
which is flat over
$\pi^{-1}\nbigo_X$,
according to 
Lemma \ref{lem;09.10.25.35}.
Then, we can prove Proposition 
\ref{prop;09.10.25.40}
by the arguments in the last part of
the proof of Proposition \ref{prop;09.10.25.32}.
\hfill\qed

\subsection{Proof of Theorem \ref{thm;09.12.4.30}}
\label{subsection;09.12.4.40}

The first claim of Theorem \ref{thm;09.12.4.30}
is obvious.
We give a preliminary for the second claim.
Put $X':=\cnum\times X$,
$X'_0:=\{0\}\times X$
and $D':=(\cnum\times D)\cup (\{0\}\times X)$.
Let $J\subset \ellsitabar$.
Put $D'(\Jbar):=\cnum\times D(\Jbar)$.
Let $\pi_0:\Xtilde'(D')\lrarr X'$ and
$\pi_1:\cnum\times\Xtilde(D)\lrarr \cnum\times X$
be the real blow up.
We have a natural diffeomorphism
$\pi_0^{-1}(X_0')\simeq S^1\times\Xtilde(D)$.
Let $\rho_0:\Xtilde'(D')\lrarr 
 \cnum\times\Xtilde(D)$
be the naturally induced map.
We use the coordinate system $z=r\,e^{\sqrt{-1}\theta}$
of $\cnum$.
We have a natural inclusion:
\begin{equation}
\label{eq;10.1.11.1}
 \nbigc^{\infty\,\nil\,<D'(\Jbar)}
 _{\widehat{\pi_1^{-1}(X_0')}}(\ast X_0')
\lrarr
 \rho_{0\ast}\Bigl(
 \nbigc^{\infty\,\nil\,<D'(\Jbar)}
 _{\widehat{\pi_0^{-1}(X_0')}}\Bigr)
\end{equation}
The operator
$\zbar\delbar_z$
induces
endomorphisms of
$\nbigc^{\infty\nil<D'(\Jbar)}
 _{\widehat{\pi_1^{-1}(X_0')}}(\ast X_0')$
and 
$\rho_{0\ast}\Bigl(
 \nbigc^{\infty\nil<D'(\Jbar)}
 _{\widehat{\pi_0^{-1}(X_0')}}\Bigr)$,
which are denoted by 
$F_1$ and $F_2$, respectively.

\begin{lem}
\label{lem;09.10.25.5}
The cokernel of $F_i$ $(i=1,2)$ are $0$,
and {\rm(\ref{eq;10.1.11.1})} induces
an isomorphism
$\Ker F_1\simeq \Ker F_2$.
\end{lem}
\pf
It is easy to obtain the vanishing of
$\Cok F_1$
by a formal calculation.
Let us prove the other claims.
We take $\nbigs\subset\cnum$ such that
(i) the induced map
$\nbigs\lrarr\cnum/\seisuu$ is bijective,
(ii) $0\in \nbigs$.
According to the decomposition
$\Nil(z)=\bigoplus_{\alpha\in \nbigs}
 z^{\alpha}\cnum[z,z^{-1}]\,[\log z]$,
we have the decomposition
$\nbigc^{\infty\,\nil\,<D'(\Jbar)}_{
 \widehat{\pi_0^{-1}(X_0')}}
=\bigoplus_{\alpha\in\nbigs}
 \nbigc^{\infty\,\nil\,<D'(\Jbar)}_{
 \widehat{\pi_0^{-1}(X_0')},
 \alpha}$.
Let $U\subset \Xtilde(D)$ be an open subset.
Let $f$ be a section of
$\nbigc^{\infty\,\nil\,<D'(\Jbar)}
 _{\widehat{\pi_0^{-1}(X_0')},\alpha}$
on $S^1\times U
\subset \pi_0^{-1}(X_0')$
expressed as follows:
\[
f=
 \sum_{\vecbeta,\veck}
 \sum_{n,j}
 f_{\vecbeta,\veck,n,j}\,
 \varphi_{\vecbeta,\veck}
 \,e^{-\sqrt{-1}\theta\alpha}
 z^{\alpha+n}\bigl( \log|z|^2\bigr)^j
\quad\quad
 \bigl(
 f_{\vecbeta,\veck,n,j}\in
 \nbigc^{\infty\,<D(\Jbar)}_{S^1\times\Xtilde(D)}
 \bigr)
\]
We have the following equality:
\begin{multline}
 \label{eq;09.9.15.2}
\zbar\delbar_z f
=\sum_{\vecbeta,\veck}
 \sum_{n,j}
 \Bigl(
 \frac{\sqrt{-1}}{2}\del_{\theta}
+\frac{\alpha}{2}
 \Bigr)f_{\vecbeta,\veck,n,j}\, 
 \varphi_{\vecbeta,\veck}\,
 e^{-\sqrt{-1}\theta\alpha}
 z^{\alpha+n}(\log|z|^2)^j
 \\
+\sum_{\vecbeta,\veck}
 \sum_{n,j}
 f_{\vecbeta,\veck,n,j}
 \varphi_{\vecbeta,\veck}
 \,e^{-\sqrt{-1}\theta\alpha}
 z^{\alpha+n}
 j(\log|z|^2)^{j-1}
\end{multline}
For any section $g$ of
$\nbigc^{\infty\,<D(\Jbar)}_{S^1\times\Xtilde(D)}$
on $S^1\times U$,
we can solve the equation
\[
 \del_{\theta}G
-\sqrt{-1}\alpha\, G=g
\quad\quad
(\alpha\neq 0)
\]
in $\nbigc^{\infty\,\nil\,<D(\Jbar)}_{S^1\times\Xtilde(D)}$.
We remark
$\int_0^{2\pi}
 e^{-\sqrt{-1}\alpha\theta}g(\theta)\,d\theta=0$.
It is easy to obtain $\Cok(\zbar\del_{\zbar})=0$
and $\Ker(\zbar\del_{\zbar})=0$
in the part $\alpha\neq 0$
by using (\ref{eq;09.9.15.2}).
Let us consider the part $\alpha=0$.
We use the filtration
with respect to the order of $\log|z|^2$.
If we take Gr with respect to this filtration,
the second term in (\ref{eq;09.9.15.2})
with $\alpha=0$
disappears.
We obtain $\nbigh^0\Gr_j=\nbigh^1\Gr_j$
for each $j$,
and they are represented by 
constants with respect to $\theta$.
Then, the second term in (\ref{eq;09.9.15.2})
induces $\nbigh^0\Gr_j\simeq \nbigh^1\Gr_{j-1}$
for $j\geq 1$.
Hence, we obtain the vanishing
of the cokernel of $\zbar\del_{\zbar}$,
and the kernel is $\nbigh^0\Gr_0$.
Then, the remaining claims of  Lemma \ref{lem;09.10.25.5}
are clear. 
\hfill\qed

\vspace{.1in}

We have the following morphism
of exact sequences:
{\small
\[
 \begin{CD}
 0 @>>>
 \Omega_{\cnum\times\Xtilde(D)}
 ^{0,\bullet\,<D'(\Jbar)\cup X_0'}
 @>>>
  \Omega_{\cnum\times\Xtilde(D)}
 ^{0,\bullet<D'(\Jbar)}(\ast X_0')
 @>>>
 \Omega_{\widehat{\pi_1^{-1}(X_0')}}^{0,\bullet\,
 <D'(\Jbar)}
 @>>> 0\\
 @. @V{=}VV @VVV @V{\simeq}VV @.\\
 0 @>>>
 \rho_{0\ast}\Omega_{\Xtilde'(D')}
 ^{0,\bullet\,<D'(\Jbar)\cup X_0'}
 @>>>
 \rho_{0\ast}\Omega_{\Xtilde(D')}
 ^{0,\bullet\,<D'(\Jbar)}
 @>>>
 \rho_{0\ast}
 \Omega_{\widehat{\pi_0^{-1}(X_0')}}^{0,\bullet\,
 <D'(\Jbar)}
 @>>> 0\\
 \end{CD}
\]
}
The left vertical arrow is an isomorphism.
According to Lemma \ref{lem;09.10.25.5},
the right vertical arrow is a quasi-isomorphism.
Thus, the central vertical arrow is also a quasi-isomorphism,
which is the second claim of Theorem \ref{thm;09.12.4.30}.
\hfill\qed

%% file: 4.4.tex
\subsection{Preliminary}

We shall freely use the notation in \S\ref{subsection;14.1.19.10}.
Let $(t_1,\ldots,t_{\ell})$
denote the standard coordinate system of 
$\cnum^{\ell}$.
We set $D_0:=\bigcup_{i=1}^{\ell}\{t_i=0\}$.
We have 
$\cnumtilde^{\ell}(D_0)=\cnumtilde^{\ell}$.
Let $X$ be any complex manifold.
The projection
$X\times\cnumtilde^{\ell}\lrarr X\times\cnum^{\ell}$
is denoted by $\pi$.
We put $H_X:=X\times D_0$.
\index{set $H_X$}

For any closed complex submanifold $Y\subset X$,
we have naturally defined morphisms:
\begin{equation}
\label{eq;13.4.19.1}
 \pi^{-1}\nbigo_{Y\times\cnum^{\ell}}
 \otimes^L_{\pi^{-1}\nbigo_{X\times\cnum^{\ell}}}
 \nbiga^{\moderate}_{X\times\cnumtilde^{\ell}}
\lrarr
 \itilde_{\ast}\nbiga^{\moderate}_{Y\times\cnumtilde^{\ell}}
\end{equation}
\begin{equation}
\label{eq;13.4.19.2}
 \pi^{-1}\nbigo_{Y\times\cnum^{\ell}}
 \otimes^L_{\pi^{-1}\nbigo_{X\times\cnum^{\ell}}}
 \nbiga^{\rapid}_{X\times\cnumtilde^{\ell}}
\lrarr
 \itilde_{\ast}\nbiga^{\rapid}_{Y\times\cnumtilde^{\ell}}
\end{equation}
Here, $\itilde:Y\times\cnumtilde^{\ell}\lrarr X\times\cnumtilde^{\ell}$
denotes the map induced by the inclusion $Y\subset X$.

\begin{lem}
\label{lem;13.4.19.3}
The morphisms {\rm(\ref{eq;13.4.19.1})}
and {\rm(\ref{eq;13.4.19.2})} are isomorphisms.
\end{lem}
\pf
Let us prove the claim for (\ref{eq;13.4.19.1}).
The other case can be proved similarly.
It is enough to argue it locally around each point of 
$H_X$.
It is easy to reduce the case 
$X=\Delta^n=\bigl\{(z_1,\ldots,z_n)\,\big|\,
 |z_i|<1
 \bigr\}$
and $Y=\{z_1=0\}$.
Let $F$ be the endomorphism of
$\nbiga^{\moderate}_{X\times\cnumtilde^{\ell}}$
given by $F(x)=z_1x$.
The complex 
$\nbiga^{\moderate}_{X\times\cnumtilde^{\ell}}
 \stackrel{F}{\lrarr}
 \nbiga^{\moderate}_{X\times\cnumtilde^{\ell}}$
expresses 
$\pi^{-1}\nbigo_{Y\times\cnum^{\ell}}
 \otimes^L_{\pi^{-1}\nbigo_{X\times\cnum^{\ell}}}
 \nbiga^{\moderate}_{X\times\cnumtilde^{\ell}}$.
Clearly, $F$ is injective.
It is enough to prove that
the induced map
$\rho:\Cok(F)\lrarr 
\nbiga^{\moderate}_{Y\times\cnum^{\ell}}$
is an isomorphism.
It is clearly surjective.
Let $f$ be any section of 
$\nbiga^{\moderate}_{X\times\cnumtilde^{\ell}}$
on $\nbigu\subset X\times\cnumtilde^{\ell}$
such that $\rho(f)=0$.
Then, $z_1^{-1}f$ naturally gives
a holomorphic function on
$\nbigu\setminus \pi^{-1}(H_X)$.
Let us prove that $z_1^{-1}f$ is of moderate growth.
We may assume that $\nbigu$ is the product
of a multi-sector
\[
\begin{array}{l}
 S_t=\bigl\{
 (t_1,\ldots,t_{\ell})
 \,\big|\,
 |\arg(t_i)-\theta_{0i}|\leq\delta_{0i},\,0<|t_i|<r_{0i}\,\,
 (i=1,\ldots,\ell)
 \bigr\}
 \\
 (\theta_{0i}\in\real,\,\,\delta_{0i}>0,\,\,r_{0i}>0)
\end{array}
\]
in $(\cnum^{\ast})^{\ell}$,
and multi-discs
$U_1=\{|z_1|\leq r_{1}\}$
and $U=\{(z_2,\ldots,z_n)\,|\,|z_i|\leq r_2\}$.
We put $U_1':=\{r_1/2\leq|z_1|\leq r_1\}$.
On $U_1'\times U\times S_t$,
we have 
$|z_1^{-1}f|\leq C\prod_{i=1}^{\ell}|t_i|^{-N}$.
By using the maximum principle,
we obtain the estimate of $z_1^{-1}f$
on $U_1\times U\times S_t$.
\hfill\qed

\subsection{The push-forward of coherent $\nbigo_X$-modules}

For any $\pi^{-1}\nbigo_{X\times\cnum^{\ell}}$-module $\nbigm$,
we canonically have 
a standard $\pi^{-1}\nbigo_{X\times\cnum^{\ell}}$-flat resolution
$\nbign_{\bullet}(\nbigm)$ of $\nbigm$
given as follows.
For any open subset $U\subset X\times\cnumtilde^{\ell}$,
let $\nbign_U$ be the free 
$\pi^{-1}(\nbigo_{X\times\cnum^{\ell}})_{|U}$-module
generated by $\nbigm(U)$,
and let $\nbign_U'$ denote its $0$-extension on 
$X\times\cnumtilde^{\ell}$.
It is naturally equipped with a morphism
$a_U:\nbign_U'\lrarr \nbigm$.
We put $\nbign_0(\nbigm):=\bigoplus_{U}\nbign_U'$,
and then $a:=\bigoplus_U a_U$ gives a surjection
$\nbign_0(\nbigm)\lrarr\nbigm$.
By applying the same procedure to
$\Ker a$,
we obtain a flat $\pi^{-1}\nbigo_{X\times\cnum^{\ell}}$-module
$\nbign_1(\nbigm)$ with a surjection
$\nbign_1(\nbigm)\lrarr \Ker a$.
By the standard inductive procedure,
we obtain the flat resolution.
In particular, we obtain a canonical flat resolution
$\nbign_{\bullet}(\nbiga^{\moderate}_{X\times\cnumtilde^{\ell}})$.

Let $\varphi:(Y,g)\lrarr (X,f)$ be a morphism in $\Cat_{\ell}$.
We have a canonical morphism
$\varphitilde_1^{-1}
 \nbign_{\bullet}(\nbiga^{\moderate}_{X\times\cnumtilde^{\ell}})
\lrarr
 \nbign_{\bullet}(\nbiga^{\moderate}_{Y\times\cnumtilde^{\ell}})$.
Hence,
for any $\nbigo_{Y}$-sheaf $M$,
we obtain the following morphism:
\[
 \varphitilde_1^{-1}
 \nbign_{\bullet}(\nbiga^{\moderate}_{X\times\cnumtilde^{\ell}})
 \otimes_{\varphitilde_1^{-1}\pi^{-1}\nbigo_{X\times\cnum^{\ell}}}
 \pi^{-1}(\Gamma_{g\ast}M)
\lrarr
 \nbign_{\bullet}(\nbiga^{\moderate}_{Y\times\cnumtilde^{\ell}})
\otimes_{\pi^{-1}\nbigo_{Y\times\cnum^{\ell}}}
 \pi^{-1}(\Gamma_{g\ast}M).
\]
It induces the following morphism:
\begin{equation}
\label{eq;13.4.19.40}
 \nbiga^{\moderate}_{X\times\cnumtilde^{\ell}}
 \otimes^L_{\pi^{-1}\nbigo_{X\times\cnum^{\ell}}}
 \pi^{-1}(\Gamma_{f\ast}R\varphi_{!}M)
\lrarr
 R\varphitilde_{1!}\bigl(
 \nbiga^{\moderate}_{Y\times\cnumtilde^{\ell}}
 \otimes^L_{\pi^{-1}\nbigo_{Y\times\cnum^{\ell}}}
 \pi^{-1}\Gamma_{g\ast}M
 \bigr)
\end{equation}
Similarly, 
we have the following natural morphism:
\begin{equation}
\label{eq;13.4.19.41}
 \nbiga^{\rapid}_{X\times\cnumtilde^{\ell}}
 \otimes^L_{\pi^{-1}\nbigo_{X\times\cnum^{\ell}}}
 \pi^{-1}(\Gamma_{f\ast}R\varphi_{!}M)
\lrarr
  R\varphitilde_{1!}\bigl(
  \nbiga^{\rapid}_{Y\times\cnumtilde^{\ell}}
 \otimes_{\pi^{-1}\nbigo_{Y\times\cnum^{\ell}}}^L
 \pi^{-1}(\Gamma_{g\ast}M)
 \bigr)
\end{equation}
\begin{rem}
Because $\nbiga^{\rapid}_{X\times\cnumtilde^{\ell}}$ is flat over
$\pi^{-1}\nbigo_{X\times\cnum^{\ell}}$
(Proposition {\rm\ref{prop;09.10.25.32}}),
we may replace $\otimes^L$
in {\rm(\ref{eq;13.4.19.41})} with $\otimes$.
Later, we shall prove that
$\nbiga^{\moderate}_{X\times\cnumtilde^{\ell}}$
is also flat over $\pi^{-1}\nbigo_{X\times\cnum^{\ell}}$
(Theorem {\rm\ref{thm;12.9.18.10}}).
\hfill\qed
\end{rem}

\begin{thm}
\label{thm;13.4.19.42}
Suppose that $M$ is $\nbigo_Y$-coherent
and that $\varphi$ is projective.
Then,
the morphisms
{\rm(\ref{eq;13.4.19.40})}
and {\rm(\ref{eq;13.4.19.41})}
are isomorphisms.
\end{thm}
\pf
We shall give details for (\ref{eq;13.4.19.40}).
Because the other case can be argued in a similar way,
we give only an indication in the last.
It is enough to consider the cases
(i) $\varphi$ is a closed immersion,
(ii) $\varphi$ is the projection
$Y=\proj^n\times X\lrarr X$.

\subsubsection{The case (i)}

The following natural morphisms
are isomorphisms:
\begin{multline}
\pi^{-1}\bigl(\Gamma_{f\ast}\varphi_{\ast}M\bigr)
\otimes^L_{\pi^{-1}\nbigo_{X\times\cnum^{\ell}}}
\nbiga^{\moderate}_{X\times\cnumtilde^{\ell}}
 \\
\simeq 
\pi^{-1}\bigl(
 \varphi_{1\ast}\Gamma_{g\ast}M
 \bigr)
\otimes^L_{\pi^{-1}\varphi_{1\ast}\nbigo_{Y\times\cnum^{\ell}}}
\Bigl(
\pi^{-1}\varphi_{1\ast}\nbigo_{Y\times\cnum^{\ell}}
\otimes^L_{\pi^{-1}\nbigo_{X\times\cnum^{\ell}}}
\nbiga^{\moderate}_{X\times\cnumtilde^{\ell}}
\Bigr)
 \\
\simeq
\varphitilde_{1\ast}
\Bigl(
\pi^{-1}(\Gamma_{g\ast}M)
 \otimes^L_{\pi^{-1}\nbigo_{Y\times\cnum^{\ell}}}
\nbiga^{\moderate}_{Y\times\cnumtilde^{\ell}}
\Bigr) 
\end{multline}
Here, we have used Lemma \ref{lem;13.4.19.3}.
Thus, we are done in the case (i).

\subsubsection{The case (ii)}

Let us consider the case where
$\varphi$ is the projection
$Y=\proj^n\times X\lrarr X$.
Let $L$ be a line bundle on $\proj^n$.
Its pull back to $Y\times\cnum^{\ell}=\proj^n\times X\times\cnum^{\ell}$
is denoted by $L_Y$.
\begin{lem}
\label{lem;13.4.19.100}
Let $q>0$.
If $H^q\bigl(\proj^n,L\bigr)=0$,
we have
\[
 R^q\varphitilde_{1\ast}\bigl(
 \pi^{-1}\Gamma_{g\ast}L_Y
 \otimes_{\pi^{-1}\nbigo_{Y\times\cnum^{\ell}}}
 \nbiga^{\moderate}_{Y\times\cnumtilde^{\ell}}
 \bigr)=0.
\]
\end{lem}
\pf
We have the natural decomposition 
$\delbar_{Y\times\cnum^{\ell}}=
 \delbar_{\proj^n}+\delbar_{X}+\delbar_{\cnum^{\ell}}$
into the differentials
of the $\proj^n$-direction,
the $X$-direction
and the $\cnum^{\ell}$-direction.
Let $\nbigb_{Y\times\cnumtilde^{\ell}}$ 
be the sheaf of $C^{\infty}$-functions
$\kappa$ on $Y\times\cnumtilde^{\ell}$
satisfying $(\delbar_X+\delbar_{\cnum^{\ell}})\kappa=0$
and the following condition locally:
\begin{description}
\item[(Moderate)]
For any differential operator $\nbigr$ on $\proj^n$,
there exists $N>0$ such that 
$\nbigr(\kappa)=O\Bigl(\prod_{i=1}^{\ell}|t_i|^{-N}\Bigr)$.
\end{description}
We naturally have
$\nbiga^{\moderate}_{Y\times\cnumtilde^{\ell}}
 \subset\nbigb_{Y\times\cnumtilde^{\ell}}$.
We set
$\nbigb^{0,\bullet}_{Y\times\cnumtilde^{\ell}}:=
 \nbigb_{Y\times\cnumtilde^{\ell}}\otimes
 \pi^{-1}(\Omega_{Y/X}^{0,\bullet})$.
The naturally defined morphism
$\nbiga^{\moderate}_{Y\times\cnumtilde^{\ell}}
\lrarr
 \nbigb^{0,\bullet}_{Y\times\cnumtilde^{\ell}}$
is a quasi isomorphism,
which can be proved by a standard argument
for Dolbeault's lemma.
Hence, we obtain the following
$\varphitilde_1$-soft resolution of
$\pi^{-1}(L_Y)\otimes_{\pi^{-1}\nbigo_{Y\times\cnum^{\ell}}}
 \nbiga^{\moderate}_{Y\times\cnumtilde^{\ell}}$:
\[
 \pi^{-1}(L_Y)\otimes_{\pi^{-1}\nbigo_{Y\times\cnum^{\ell}}}
 \nbiga^{\moderate}_{Y\times\cnumtilde^{\ell}}
\lrarr
  \pi^{-1}(L_Y)\otimes_{\pi^{-1}\nbigo_{Y\times\cnum^{\ell}}}
 \nbigb^{0,\bullet}_{Y\times\cnumtilde^{\ell}}
\]

We take a hermitian metric $h_L$ of $L$.
We fix a K\"ahler metric $g_{\proj^n}$ of $\proj^n$.
Let $\delbar_L^{\ast}$ denote the formal adjoint of
$\delbar_L:
 C^{\infty}(L\otimes\Omega_{\proj^n}^{0,\bullet})
\lrarr
 C^{\infty}(L\otimes\Omega_{\proj^n}^{0,\bullet+1})$.
Let $\Delta_L^{0,\bullet}$ denote the Laplacian
on $\Gamma(\proj^n,L\otimes
 \Omega_{\proj^n}^{0,\bullet})$
associated to $h_L$ and $g_{\proj^n}$.
Let $G^{0,\bullet}$ be the Green operator.
By the assumption $H^q(\proj^n,L)=0$ for $q>0$,
we have 
$\Delta^{0,q}\circ G^{0,q}
=G^{0,q}\circ\Delta^{0,q}=\id$
if $q>0$.
We have 
$[G^{0,\bullet},\delbar_L]=
 [G^{0,\bullet},\delbar^{\ast}_L]=0$.
In particular, if $\delbar_L\tau=0$
for 
$\tau\in \Gamma(\proj^n,L\otimes\Omega^{0,q})$
$(q>0)$,
we have
$\delbar_L\delbar_L^{\ast}G(\tau)=\tau$.
Recall the following standard results
for elliptic operators:
\begin{itemize}
\item
$G^{0,q}$ are integral operators.
\item
For any non-negative integer $m$,
there exists $C_m>0$ such that
$\|G^{0,q}(\tau)\|_{L_{m+2}^2}
\leq
 C_m\,\| \tau\|_{L_m^2}$
for any 
$\tau\in\Gamma(\proj^n,L\otimes\Omega^{0,q})$,
where $\|\cdot\|_{L_m^2}$ denotes
the Sobolev norm.
\end{itemize}

Let $P\in\pi^{-1}(H_X)$.
Let $\nbigu_P$ be an open neighbourhood
of $P$ in $X\times\cnumtilde^{\ell}$.
Put $\nbigu_P^{\circ}:=
 \nbigu_P\setminus\pi^{-1}(H_X)$.
We have 
$\varphitilde_1^{-1}(\nbigu_P)
=\proj^n\times\nbigu_P$.
Let $\tau\in\Gamma(\proj^n\times\nbigu_P,
 \pi^{-1}L_Y
 \otimes_{\pi^{-1}\nbigo_{Y\times\cnum^{\ell}}}
 \nbigb^{0,q})$.
We obtain a $C^{\infty}$-function $G(\tau)$ on
$\proj^n\times\nbigu_P^{\circ}$,
and we have
$\delbar_{z_i} G(\tau)=0$
and $\del_{z_i}G(\tau)=G(\del_{z_i}\tau)$
for any local coordinate system
$(z_1,\ldots,z_n)$ on $X\times\cnum^{\ell}$.
Then, by the estimate of the Green operator,
we obtain that
$G(\tau)\in
 \Gamma(\proj^n\times\nbigu_P,
 \pi^{-1}L_Y\otimes_{\pi^{-1}\nbigo_Y}
 \nbigb^{0,q})$.
Moreover, if $\delbar_L\tau=0$ and $q>0$,
we have
$\delbar_L\bigl(
 \delbar_L^{\ast}G(\tau)
 \bigr)=\tau$.
Thus, we obtain Lemma \ref{lem;13.4.19.100}.
\hfill\qed

\begin{lem}
\label{lem;13.4.19.50}
We have
$\varphitilde_{1\ast}
 \nbiga^{\moderate}_{Y\times\cnumtilde^{\ell}}
\simeq
 \nbiga^{\moderate}_{X\times\cnumtilde^{\ell}}$,
i.e.,
the morphism {\rm(\ref{eq;13.4.19.40})} is an isomorphism
for $\nbigo_Y$.
\end{lem}
\pf
Let $P\in \pi^{-1}(H_X)$.
Let $\nbigu_P$ be a small neighbourhood 
of $P$ in $X\times\cnumtilde^{\ell}$.
Let $\kappa\in\Gamma(\proj^n\times\nbigu_P,
 \nbiga^{\moderate}_{Y\times\cnumtilde^{\ell}})$.
Take any point $Q$ of $\proj^n$.
We consider the inclusion
$\iota_Q:\nbigu_P\simeq\nbigu_P\times\{Q\}\lrarr
 \proj^n\times\nbigu_P$.
We have 
$\mu:=\iota_Q^{-1}(\kappa)\in
 \Gamma(\nbigu_P,
 \nbiga^{\moderate}_{X\times\cnumtilde^{\ell}})$.
It is easy to deduce that
$\kappa=\varphitilde(\mu)$.
Then, we obtain Lemma \ref{lem;13.4.19.50}.
\hfill\qed

\begin{lem}
\label{lem;13.4.19.51}
Let $L$ be a line bundle on $\proj^n$.
Then, {\rm(\ref{eq;13.4.19.40})} is 
an isomorphism for $L_Y$.
\end{lem}
\pf
We use an induction on $n$.
In the case $n=0$, the claim is trivial.
Assume that we have already obtained
the claim in the case $n-1$.
Let $L=\nbigo_{\proj^n}(m)$.
If $m=0$, the claim follows from
Lemma \ref{lem;13.4.19.50}.
We fix a hyperplane 
$\proj^{n-1}_{\infty}\subset\proj^n$.
If $m>0$, we can reduce the claim
to the case $m-1$,
by using the exact sequence
$0\lrarr
 \nbigo_{\proj^n}(m-1)
\lrarr
 \nbigo_{\proj^n}(m)
\lrarr
 \nbigo_{\proj_{\infty}^{n-1}}(m)\lrarr 0$.
If $m<0$,
we can reduce the claim to the case
$m+1$,
by using the exact sequence
$0\lrarr
 \nbigo_{\proj^n}(m)
\lrarr
 \nbigo_{\proj^n}(m+1)
\lrarr
 \nbigo_{\proj_{\infty}^{n-1}}(m+1)\lrarr 0$.
\hfill\qed

\vspace{.1in}

Let us finish the proof in the case (ii).
It is enough to prove that
(\ref{eq;13.4.19.40})
is an isomorphism around 
any point of $X\times\cnum^{\ell}$,
which we shall implicitly use.
We may assume to have a resolution
\[
\Bigl(
\cdots\lrarr
\nbigq_{p}
\lrarr
\nbigq_{p-1}
\lrarr\cdots
\lrarr\nbigq_1\lrarr\nbigq_0
\Bigr)\simeq M,
\]
such that
$\nbigq_p$ are of the form
$\bigoplus_{i=1}^{N_p}(L_{p,i})_Y$,
where $L_{p,i}$ are line bundles on $\proj^n$.
By Lemma \ref{lem;13.4.19.51},
the morphisms (\ref{eq;13.4.19.40})
for $\nbigq_p$ are isomorphisms.
Hence, (\ref{eq;13.4.19.40})
for $M$ is also an isomorphism.
Thus, the proof for (\ref{eq;13.4.19.40}) 
is finished.

\vspace{.1in}

Let us give an indication
to prove that (\ref{eq;13.4.19.41}) is an isomorphism.
We can argue the case (i) in the same way.
In the case (ii),
we replace the condition 
{\bf(Moderate)} in the proof of Lemma \ref{lem;13.4.19.100}
with the following:
\begin{description}
\item[(Rapid)]
Let $\nbigr$ be any differential operators on $\proj^n$.
Then, $\nbigr(\kappa)=O\Bigl(\prod|t_i|^{N}\Bigr)$
for any $N$.
\end{description}
Then, we can prove that
(\ref{eq;13.4.19.41}) is an isomorphism
in the case (ii).
Thus, the proof of Theorem \ref{thm;13.4.19.42}
is finished.
\hfill\qed

%% file: 4.5.tex
\subsection{Statements}

\begin{thm}
\label{thm;13.4.20.1}
Let $(X,f)$ be an object in $\Cat_{\ell}$.
\begin{itemize}
\item
 $\Tor_i^{\pi^{-1}\nbigo_{X\times\cnum^{\ell}}}
  \bigl(
 \nbiga^{\moderate}_{X\times\cnumtilde^{\ell}},\,
 \pi^{-1}\nbigo_{\Gamma_f(X)}
 \bigr)=0$
 for $i\neq 0$.
Namely, 
\[
  \nbiga^{\moderate}_{X\times\cnumtilde^{\ell}}
 \otimes^L_{\pi^{-1}\nbigo_{X\times\cnum^{\ell}}}
 \pi^{-1}\nbigo_{\Gamma_f(X)}
\simeq
 \nbiga^{\moderate}_{X\times\cnumtilde^{\ell}}
 \otimes_{\pi^{-1}\nbigo_{X\times\cnum^{\ell}}}
 \pi^{-1}\nbigo_{\Gamma_f(X)}.
\]
\item
Let $\varphi:(Y,g)\lrarr (X,f)$ be a projective birational morphism
such that (i) $D_Y$ is normal crossing,
(ii) $Y\setminus D_Y\simeq X\setminus D_X$.
For the naturally induced map
$\rho:\Ytilde(D_Y)\lrarr X\times\cnumtilde^{\ell}$,
we have
\begin{equation}
\label{eq;13.4.20.50}
 R\rho_{\ast}\nbiga^{\moderate}_{\Ytilde(D_Y)}
\simeq
 \nbiga^{\moderate}_{X\times\cnumtilde^{\ell}}
 \otimes_{\pi^{-1}\nbigo_{X\times\cnum^{\ell}}}
 \pi^{-1}\nbigo_{\Gamma_f(X)}
\end{equation}
\begin{equation}
\label{eq;13.4.20.51}
 R\rho_{\ast}\nbiga^{\rapid}_{\Ytilde(D_Y)}
\simeq
 \nbiga^{\rapid}_{X\times\cnumtilde^{\ell}}
 \otimes_{\pi^{-1}\nbigo_{X\times\cnum^{\ell}}}
 \pi^{-1}\nbigo_{\Gamma_f(X)}
\end{equation}
\item
 The support of 
  $\nbiga^{\moderate}_{X\times\cnumtilde^{\ell}}
 \otimes_{\pi^{-1}\nbigo_{X\times\cnum^{\ell}}}
 \pi^{-1}\nbigo_{\Gamma_f(X)}$
and 
  $\nbiga^{\rapid}_{X\times\cnumtilde^{\ell}}
 \otimes_{\pi^{-1}\nbigo_{X\times\cnum^{\ell}}}
 \pi^{-1}\nbigo_{\Gamma_f(X)}$
 are $\Xtilde(f)$.
\end{itemize}
\end{thm}

\begin{rem}
Note that
$\nbiga^{\rapid}_{X\times\cnumtilde^{\ell}}$
is flat over $\pi^{-1}\nbigo_{X\times\cnum^{\ell}}$,
according to Proposition {\rm\ref{prop;09.10.25.32}}.
The first claim of the theorem
is a special case of the flatness of
$\nbiga^{\moderate}_{X\times\cnumtilde^{\ell}}$
over $\pi^{-1}\nbigo_{X\times\cnum^{\ell}}$
(Theorem {\rm\ref{thm;12.9.18.10}}).
\hfill\qed
\end{rem}

Let us state some consequences.
We have the sheaves of algebras $\nbiga^{\moderate}_{X,f}$
and $\nbiga^{\rapid}_{X,f}$
on $\Xtilde(f)$ determined by the following conditions:
\index{sheaf $\nbiga^{\rapid}_{X,f}$}
\index{sheaf $\nbiga^{\moderate}_{X,f}$}
\begin{eqnarray*}
 \Gammatilde_{f\ast}\nbiga^{\moderate}_{X,f}=
 \pi^{-1}(\nbigo_{\Gamma_f(X)})
 \otimes_{\pi^{-1}\nbigo_{X\times\cnum^{\ell}}}
 \nbiga^{\moderate}_{X\times\cnumtilde^{\ell}}
 \\
  \Gammatilde_{f\ast}\nbiga^{\rapid}_{X,f}=
 \pi^{-1}(\nbigo_{\Gamma_f(X)})
 \otimes_{\pi^{-1}\nbigo_{X\times\cnum^{\ell}}}
 \nbiga^{\rapid}_{X\times\cnumtilde^{\ell}}
\end{eqnarray*}

\begin{thm}
\label{thm;14.1.19.1}
Let $(X,f)\in\Cat_{\ell}$.
\begin{itemize}
\item
For the inclusion
$j:X\setminus D_X\lrarr \Xtilde(f)$,
the natural morphism
$\nbiga^{\moderate}_{X,f}\lrarr
 j_{\ast}\nbigo_{X\setminus D_X}$
is a monomorphism.
The image is 
$\nbiga^{\moderate}_{\Xtilde(f)}$.
\item
The natural morphism
 $\nbiga^{\rapid}_{X,f}\lrarr
 j_{\ast}\nbigo_{X\setminus D_X}$
is a monomorphism.
The image is
$\nbiga^{\rapid}_{\Xtilde(f)}$.

\item
In particular,
if $f$ is submersive,
then we naturally have
$\nbiga^{\moderate}_{X,f}\simeq
 \nbiga^{\moderate}_{\Xtilde(D_X)}$
and 
$\nbiga^{\rapid}_{X,f}\simeq
 \nbiga^{\rapid}_{\Xtilde(D_X)}$.
\hfill\qed
\end{itemize}
\end{thm}
\pf
It follows from the descriptions (\ref{eq;13.4.20.50})
and (\ref{eq;13.4.20.51}).
\hfill\qed

\vspace{.1in}

Theorem \ref{thm;13.4.19.42}
can be reformulated 
in terms of 
$\nbiga^{\moderate}_{\Xtilde(f)}$
and 
$\nbiga^{\rapid}_{\Xtilde(f)}$.
\begin{thm}
\label{thm;13.4.19.210}
Let $\varphi:(Y,g)\lrarr (X,f)$ be a projective morphism in $\Cat_{\ell}$.
Let $M$ be any coherent $\nbigo_Y$-module.
Then, the following natural morphisms are isomorphisms:
\begin{equation}
 \label{eq;14.1.14.11}
 \nbiga^{\moderate}_{\Xtilde(f)}\otimes^L_{\pi^{-1}\nbigo_X}
 \pi^{-1}R\varphi_{\ast}M
\simeq
 R\varphitilde_{\ast}\Bigl(
  \nbiga^{\moderate}_{\Ytilde(g)}\otimes^L_{\pi^{-1}\nbigo_Y}
 \pi^{-1}M
 \Bigr)
\end{equation}
\begin{equation}
 \label{eq;14.1.14.12}
 \nbiga^{\rapid}_{\Xtilde(f)}\otimes^L_{\pi^{-1}\nbigo_X}
 \pi^{-1}R\varphi_{\ast}M
\simeq
 R\varphitilde_{\ast}\Bigl(
  \nbiga^{\rapid}_{\Ytilde(g)}\otimes^L_{\pi^{-1}\nbigo_Y}
 \pi^{-1}M
 \Bigr)
\end{equation}
\hfill\qed
\end{thm}
After the flatness results in 
Proposition \ref{prop;09.10.25.32}
and Theorem \ref{thm;12.9.18.10} below,
we may replace $\otimes^L$
with $\otimes$ 
in (\ref{eq;14.1.14.11}) and (\ref{eq;14.1.14.12}).

\subsection{Proof of Theorem \ref{thm;13.4.20.1}}

Let us begin with the simplest case.

\begin{lem}
Suppose that $f$ is submersive.
For the naturally induced closed immersion
$\rho:\Xtilde(D_X)\lrarr  X\times\cnumtilde^{\ell}$,
the following natural morphisms are isomorphisms:
\begin{equation}
 \label{eq;13.4.20.10}
 \pi^{-1}\nbigo_{\Gamma_f(X)}
 \otimes^L_{\pi^{-1}\nbigo_{X\times\cnum^{\ell}}}
 \nbiga^{\moderate}_{X\times\cnumtilde^{\ell}}
\lrarr
 \rho_{\ast}\nbiga^{\moderate}_{\Xtilde(D_X)}
\end{equation}
\begin{equation}
 \label{eq;13.4.20.11}
  \pi^{-1}\nbigo_{\Gamma_f(X)}
 \otimes^L_{\pi^{-1}\nbigo_{X\times\cnum^{\ell}}}
 \nbiga^{\rapid}_{X\times\cnumtilde^{\ell}}
\lrarr
 \rho_{\ast}\nbiga^{\rapid}_{\Xtilde(D_X)}
\end{equation}
\end{lem}
\pf
It is enough to argue it locally around any point of $H_X$.
We may assume
$X=\{(z_1,\ldots,z_n)\}$ and $f=(z_1,\ldots,z_{\ell})$.
Let $G:X\times\cnum^{\ell}\lrarr \cnum^n\times\cnum^{\ell}$
be given by
\[
 G(z_1,\ldots,z_n,t_1,\ldots,t_{\ell})=
 (z_1-t_1,z_2-t_2,\ldots,z_{\ell}-t_{\ell},z_{\ell+1},\ldots,z_n,t_1,\ldots,t_{\ell}).
\]
Then, 
$G\circ \Gamma_{f}(z_1,\ldots,z_n)=
 (0,\ldots,0,z_{\ell+1},\ldots,z_n,z_1,\ldots,z_{\ell})$.
By using $G$,
it is easy to prove that the morphisms
(\ref{eq;13.4.20.10})
and 
(\ref{eq;13.4.20.11})
are isomorphisms.
\hfill\qed

\vspace{.1in}
Let us consider the case where $D_X$ is normal crossing.
We have a naturally defined map
$X\setminus D_X\lrarr X\times(\cnum^{\ast})^{\ell}$
as the graph.
Let us observe that it is extended to
$\rho_1:\Xtilde(D_X)\lrarr X\times\cnumtilde^{\ell}$.
Let $f_i$ be the composite of
$f:X\lrarr\cnum^{\ell}$
and the projection $\cnum^{\ell}\lrarr\cnum$
onto the $i$-th component.
It induced a map
$g_i:X\setminus D_X\lrarr \cnum^{\ast}$.
It is enough to observe that it is extended to
a map $\Xtilde(D_X)\lrarr\cnumtilde$.
Let $P$ be any point of $D_X$.
Because $f_i^{-1}(0)$ is contained in 
the normal crossing hypersurface $D_X$,
we can take a holomorphic coordinate neighbourhood
$(X_P;z_1,\ldots,z_n)$ around $P$
such that
$D_X=\bigcup_{i=1}^{\ell}\{z_i=0\}$
and $f=\prod_{i=1}^{\ell}z_i^{m_i}$,
where $m_i\geq 0$.
Let $z_i=r_ie^{\sqrt{-1}\theta_i}$.
Because the map
$\Xtilde(D_X)\lrarr \cnum^{\ast}$
is described as
$(r_1,e^{\sqrt{-1}\theta_1},\ldots,
 r_{\ell}e^{\sqrt{-1}\theta{\ell}},
z_{\ell+1},\ldots,z_n)
\longmapsto
 \prod r_i^{m_i}e^{\sqrt{-1}m\theta_i}$,
we obtain that
$g_{i|X_P\setminus D_X}$ 
is extended to
$\Xtilde_P(D_X\cap X_P)
\lrarr \cnumtilde$.
Then, the claim follows.

We have the naturally defined morphism:
\begin{equation}
 \label{eq;13.4.19.20}
\nbiga^{\moderate}_{X\times\cnumtilde^{\ell}}
 \otimes_{\pi^{-1}\nbigo_{X\times\cnum^{\ell}}}
 \pi^{-1}\nbigo_{\Gamma_f(X)}
 \lrarr
 \rho_{1\ast}\nbiga^{\moderate}_{\Xtilde(D_X)}
\end{equation}

\begin{prop}
\label{prop;13.4.19.30}
Suppose that $D_X:=f^{-1}(D_0)$ is normal crossing.
 The morphism {\rm(\ref{eq;13.4.19.20})}
 is an isomorphism.
 Moreover, we have the following isomorphisms:
\[
 R\rho_{1\ast}\nbiga^{\moderate}_{\Xtilde(D_X)}
\simeq
 \rho_{1\ast}\nbiga^{\moderate}_{\Xtilde(D_X)}
\]
\[
 \nbiga^{\moderate}_{X\times\cnumtilde^{\ell}}
 \otimes^L_{\pi^{-1}\nbigo_{X\times\cnum^{\ell}}}
 \pi^{-1}\nbigo_{\Gamma_f(X)}
\simeq
 \nbiga^{\moderate}_{X\times\cnumtilde^{\ell}}
 \otimes_{\pi^{-1}\nbigo_{X\times\cnum^{\ell}}}
  \pi^{-1}\nbigo_{\Gamma_f(X)}
\]
\end{prop}
\pf
In the proof,
we omit to denote $\pi^{-1}$.
We have the maps
$\Gammatilde^{(1)}_{f}:
 \Xtilde(D_X)\lrarr \Xtilde(D_X)\times\cnumtilde^{\ell}$
and
$\Gammatilde^{(2)}_{f}:
 \Xtilde(D_X)\lrarr \Xtilde(D_X)\times\cnum^{\ell}$
induced by $f$.
We have the projections:
\[
 \nu_1:
 \Xtilde(D_X)\times\cnumtilde^{\ell}
\lrarr
 X\times\cnumtilde^{\ell},
\quad\quad
 \nu_2:
 \Xtilde(D_X)\times\cnumtilde^{\ell}
 \lrarr
 \Xtilde(D_X)\times\cnum^{\ell}.
\]
We set
$D_X':=D_X\times\cnum^{\ell}$.
According to \S II.1.1 of \cite{sabbah4},
we have the following isomorphisms:
\[
 R\nu_{1\ast}\nbiga^{\moderate}_{\Xtilde(D_X)\times\cnumtilde^{\ell}}
\simeq
 \nbiga^{\moderate}_{X\times\cnumtilde^{\ell}}(\ast D_X'),
\quad\quad
 R\nu_{2\ast}\nbiga^{\moderate}_{\Xtilde(D_X)\times\cnumtilde^{\ell}}
\simeq
 \nbiga^{\moderate}_{\Xtilde(D_X)\times\cnum^{\ell}}(\ast H_X)
\]
Hence, we have the following natural isomorphisms:
\begin{multline}
 \label{eq;13.4.19.21}
  R\nu_{1\ast}\bigl(
 \nbiga^{\moderate}_{\Xtilde(D_X)\times\cnumtilde^{\ell}}
\otimes^L_{\nbigo_{X\times\cnum^{\ell}}}
 \nbigo_{\Gamma_f(X)}
 \bigr)
\simeq
 \nbiga^{\moderate}_{X\times\cnumtilde^{\ell}}(\ast D_X')
\otimes^L_{\nbigo_{X\times\cnum^{\ell}}}
 \nbigo_{\Gamma_f(X)}
 \\
\simeq
  \nbiga^{\moderate}_{X\times\cnumtilde^{\ell}}
\otimes^L_{\nbigo_{X\times\cnum^{\ell}}}
 \nbigo_{\Gamma_f(X)}(\ast D_{X'})
\simeq
  \nbiga^{\moderate}_{X\times\cnumtilde^{\ell}}
\otimes^L_{\nbigo_{X\times\cnum^{\ell}}}
 \nbigo_{\Gamma_f(X)}(\ast H_X)
 \\
\simeq
 \nbiga^{\moderate}_{X\times\cnumtilde^{\ell}}(\ast H_X)
\otimes^L_{\nbigo_{X\times\cnum^{\ell}}}
 \nbigo_{\Gamma_f(X)}
\simeq
  \nbiga^{\moderate}_{X\times\cnumtilde^{\ell}}
\otimes^L_{\nbigo_{X\times\cnum^{\ell}}}
 \nbigo_{\Gamma_f(X)}
\end{multline}
We also have the following:
\begin{multline}
 \label{eq;13.4.19.22}
 R\nu_{2\ast}\bigl(
 \nbiga^{\moderate}_{\Xtilde(D_X)\times\cnumtilde^{\ell}}
 \otimes^L_{\nbigo_{X\times\cnum^{\ell}}}
 \nbigo_{\Gamma_f(X)}
 \bigr)
\simeq
 \nbiga^{\moderate}_{\Xtilde(D_X)\times\cnum^{\ell}}(\ast H_X)
 \otimes^L_{\nbigo_{X\times\cnum^{\ell}}}
 \nbigo_{\Gamma_f(X)}
 \\
\simeq
 \nbiga^{\moderate}_{\Xtilde(D_X)\times\cnum^{\ell}}
 \otimes^L_{\nbigo_{X\times\cnum^{\ell}}}
 \nbigo_{\Gamma_f(X)}
\end{multline}

\begin{lem}
$\Gammatilde^{(2)}_f$ is a closed embedding,
and that we have
\begin{equation}
\label{eq;13.4.20.20}
 \nbiga^{\moderate}_{\Xtilde(D_X)\times\cnum^{\ell}}
 \otimes^L_{\nbigo_{X\times\cnum^{\ell}}}
 \nbigo_{\Gamma_f(X)}
\simeq
\nbiga^{\moderate}_{\Xtilde(D_X)\times\cnum^{\ell}}
 \otimes_{\nbigo_{X\times\cnum^{\ell}}}
 \nbigo_{\Gamma_f(X)}
\simeq
 \Gammatilde_{f\ast}^{(2)}
 \nbiga^{\moderate}_{\Xtilde(D_X)}.
\end{equation}
\end{lem}
\pf
For the expression
$f=(f_1,\ldots,f_{\ell})$,
we define 
$G':X\times\cnum^{\ell}\lrarr X\times\cnum^{\ell}$
by
$G'(P,t_1,\ldots,t_{\ell}):=
(P,t_1-f_1(P),\ldots,t_{\ell}-f_{\ell}(P))$.
We have $G'\circ\Gamma_f(P)=(P,0,\ldots,0)$.
Then, we can prove (\ref{eq;13.4.20.20})
by an induction on $\ell$.
\hfill\qed

\begin{lem}
\label{lem;13.4.19.23}
The support of 
$\Tor_{\ast}^{\pi^{-1}\nbigo_{X\times\cnum^{\ell}}}
 \bigl(
 \nbiga^{\moderate}_{\Xtilde(D_X)\times\cnumtilde^{\ell}},\!
 \pi^{-1}\nbigo_{\Gamma_f(X)}
 \bigr)$
is contained in 
$\Gammatilde^{(1)}_f(\Xtilde(D_X))$.
\end{lem}
\pf
Let $U$ denote an $\ell$-dimensional vector space with
a basis $e_1,\ldots,e_{\ell}$.
We set $C^{k-\ell}:=\bigwedge^kU\otimes \nbigo_{X\times\cnum^{\ell}}$.
Let $\del:C^m\lrarr C^{m+1}$ be given by
$\del\alpha=\sum (t_i-f_i)e_i\wedge \alpha$.
Then, we obtain a complex of $\nbigo_{X\times\cnum^{\ell}}$-modules
$C^{\bullet}$,
and it gives a free resolution of 
$\nbigo_{X\times\cnum^{\ell}}$-module
$\nbigo_{\Gamma_f(X)}$.
If $Q\in\Xtilde(D_X)\times\cnumtilde^{\ell}$ is not 
contained in $\Gammatilde_f^{(1)}(X)$,
then one of $t_i-f_i$ are invertible in
$\nbiga^{\moderate}_{\Xtilde(D_X)\times\cnumtilde^{\ell}}$
around $Q$.
Hence, the complex
$\nbiga^{\moderate}_{\Xtilde(D_X)\times\cnumtilde^{\ell}}
 \otimes C^{\bullet}$ is acyclic around $Q$.
It implies the claim of Lemma \ref{lem;13.4.19.23}.
\hfill\qed

\vspace{.1in}
Note that $\nu_2$ induces a homeomorphism
$\Gammatilde^{(1)}_f(X)\simeq
 \Gammatilde^{(2)}_f(X)$.
By Lemma \ref{lem;13.4.19.23},
we obtain that
\[
 R^p\nu_{2\ast}
 \Tor_j^{\pi^{-1}\nbigo_{X\times\cnum^{\ell}}}
 \bigl(
 \nbiga^{\moderate}_{\Xtilde(D)\times\cnumtilde^{\ell}},\,
 \pi^{-1}\nbigo_{\Gamma_f(X)}
 \bigr)=0
\]
 for $p\neq 0$.
By applying the argument of the spectral sequence
with (\ref{eq;13.4.20.20}) to (\ref{eq;13.4.19.22}),
we obtain that 
\[
 \Tor_j^{\pi^{-1}\nbigo_{X\times\cnum^{\ell}}}
 \bigl(
 \nbiga^{\moderate}_{\Xtilde(D)\times\cnumtilde^{\ell}},\,
 \pi^{-1}\nbigo_{\Gamma_f(X)}
 \bigr)=0
\]
for $j\neq 0$,
i.e.,
$\nbiga^{\moderate}_{\Xtilde(D)\times\cnumtilde^{\ell}}
 \otimes_{\pi^{-1}\nbigo_{X\times\cnum^{\ell}}}^L
 \pi^{-1}\nbigo_{\Gamma_f(X)}
\simeq
\nbiga^{\moderate}_{\Xtilde(D)\times\cnumtilde^{\ell}}
 \otimes_{\pi^{-1}\nbigo_{X\times\cnum^{\ell}}}
 \pi^{-1}\nbigo_{\Gamma_f(X)}$
on $\Xtilde(D_X)\times\cnumtilde^{\ell}$.
We also obtain an isomorphism
of sheaves on $\Xtilde(D_X)\simeq
 \Gammatilde_f^{(i)}(X)$:
\[
 \nbiga^{\moderate}_{\Xtilde(D)\times\cnumtilde^{\ell}}
 \otimes_{\nbigo_{X\times\cnum^{\ell}}}
 \nbigo_{\Gamma_f(X)}
\simeq
 \nbiga^{\moderate}_{\Xtilde(D_X)}
\]
From (\ref{eq;13.4.19.21}),
we obtain 
\begin{equation}
 \label{eq;14.1.14.10}
 R\rho_{1\ast}\nbiga^{\moderate}_{\Xtilde(D_X)}
\simeq
 R\nu_{1\ast}\bigl(\nbiga^{\moderate}_{\Xtilde(D_X)\times\cnumtilde^{\ell}}
 \otimes_{\nbigo_{X\times\cnum^{\ell}}}
 \nbigo_{\Gamma_f(X)}\bigr)
\simeq
 \nbiga^{\moderate}_{X\times\cnumtilde^{\ell}}
 \otimes^L_{\nbigo_{X\times\cnum^{\ell}}}
 \nbigo_{\Gamma_f(X)}.
\end{equation}
Note 
 $R^p\nu_{1\ast}\bigl(\nbiga^{\moderate}_{\Xtilde(D_X)\times\cnumtilde^{\ell}}
 \otimes_{\nbigo_{X\times\cnum^{\ell}}}
 \nbigo_{\Gamma_f(X)}\bigr)=0$
unless $p\geq 0$,
and 
the $p$-th cohomology sheaf of 
$\nbiga^{\moderate}_{X\times\cnumtilde^{\ell}}
 \otimes^L_{\nbigo_{X\times\cnum^{\ell}}}
 \nbigo_{\Gamma_f(X)}$
is $0$ unless $p\leq 0$.
Hence, (\ref{eq;14.1.14.10})
implies the claims of 
Proposition \ref{prop;13.4.19.30}.
\hfill\qed

\begin{prop}
\label{prop;13.4.20.21}
Suppose that 
$D_X$ is normal crossing.
Then, the natural map
\[
 \nbiga^{\rapid}_{X\times\cnumtilde^{\ell}}
 \otimes_{\pi^{-1}\nbigo_{X\times\cnum^{\ell}}}
 \pi^{-1}\nbigo_{\Gamma_f(X)}
\simeq
 \rho_{1\ast}\nbiga^{\rapid}_{\Xtilde(D_X)}
\]
is an isomorphism.
Moreover, we have
$R\rho_{1\ast}\nbiga^{\rapid}_{\Xtilde(D_X)}
\simeq
 \rho_{1\ast}\nbiga^{\rapid}_{\Xtilde(D_X)}$.
\end{prop}
\pf
It is proved by the arguments
in the proof of Proposition \ref{prop;13.4.19.30}.
We omit to denote $\pi^{-1}$.
We have the following isomorphisms:
\[
 R\nu_{1\ast}\nbiga^{<H_X\leq D_X'}_{\Xtilde(D_X)\times\cnumtilde^{\ell}}
\simeq
 \nbiga^{<H_X}_{X\times\cnumtilde^{\ell}}(\ast D_X'),
\quad\quad
 R\nu_{2\ast}\nbiga^{\leq H_X<D_X'}_{\Xtilde(D_X)\times\cnumtilde^{\ell}}
\simeq
 \nbiga^{<D_X'}_{\Xtilde(D_X)\times\cnum^{\ell}}(\ast H_X)
\]
Hence, we have the following natural isomorphisms:
\begin{multline}
\label{eq;13.4.20.40}
  R\nu_{1\ast}\bigl(
 \nbiga^{<H_X\leq D_X'}_{\Xtilde(D_X)\times\cnumtilde^{\ell}}
\otimes_{\nbigo_{X\times\cnum^{\ell}}}
 \nbigo_{\Gamma_f(X)}
 \bigr)
\simeq
 \nbiga^{<H_X}_{X\times\cnumtilde^{\ell}}(\ast D_X')
\otimes_{\nbigo_{X\times\cnum^{\ell}}}
 \nbigo_{\Gamma_f(X)}
 \\
\simeq
  \nbiga^{<H_X}_{X\times\cnumtilde^{\ell}}
\otimes_{\nbigo_{X\times\cnum^{\ell}}}
 \nbigo_{\Gamma_f(X)}
\end{multline}
\begin{multline}
 R\nu_{2\ast}\bigl(
 \nbiga^{\leq H_X<D_X'}_{\Xtilde(D_X)\times\cnumtilde^{\ell}}
 \otimes_{\nbigo_{X\times\cnum^{\ell}}}
 \nbigo_{\Gamma_f(X)}
 \bigr)
\simeq
 \nbiga^{<D_X'}_{\Xtilde(D_X)\times\cnum^{\ell}}(\ast H_X)
 \otimes_{\nbigo_{X\times\cnum^{\ell}}}
 \nbigo_{\Gamma_f(X)}
 \\
\simeq
 \nbiga^{<D_X'}_{\Xtilde(D_X)\times\cnum^{\ell}}
 \otimes_{\nbigo_{X\times\cnum^{\ell}}}
 \nbigo_{\Gamma_f(X)}
\simeq
 \Gammatilde^{(2)}_{f\ast}
 \nbiga^{<D_X}_{\Xtilde(D_X)}
\end{multline}
Let us consider the following morphisms:
\begin{equation}
 \label{eq;13.4.20.30}
 \nbiga^{<H_X\leq D_X'}_{\Xtilde(D)\times\cnumtilde^{\ell}}
 \otimes_{\nbigo_{X\times\cnum^{\ell}}}\nbigo_{\Gamma_f}
 \llarr
 \nbiga^{<(H_X\cup D_X')}_{\Xtilde(D)\times\cnumtilde^{\ell}}
 \otimes_{\nbigo_{X\times\cnum^{\ell}}}\nbigo_{\Gamma_f}
 \lrarr
  \nbiga^{<D_X'\leq H_X}_{\Xtilde(D)\times\cnumtilde^{\ell}}
 \otimes_{\nbigo_{X\times\cnum^{\ell}}}\nbigo_{\Gamma_f}
\end{equation}
Because $t_i-f_i$ are invertible
on $\nbiga^{<H_X}_{\widehat{\pi^{-1}(D_X')}}$
and 
$\nbiga^{<D'_X}_{\widehat{\pi^{-1}(H_X)}}$,
we have
\[
\nbiga^{<H_X}_{\widehat{\pi^{-1}(D_X')}}
\otimes_{\nbigo_{X\times\cnum^{\ell}}}\nbigo_{\Gamma_f(X)}
=0,
\quad\quad
\nbiga^{<D'_X}_{\widehat{\pi^{-1}(H_X)}}
\otimes_{\nbigo_{X\times\cnum^{\ell}}}\nbigo_{\Gamma_f(X)}
=0.
\]
Hence,
the morphisms in {\rm(\ref{eq;13.4.20.30})}
are isomorphisms.
By the argument in the proof of
Lemma \ref{lem;13.4.19.23},
we obtain that
the support of the sheaves in (\ref{eq;13.4.20.30})
are contained in $\Gammatilde^{(1)}_f(X)$.
Because $\nu_2$ gives a homeomorphism
$\Gammatilde_f^{(1)}(\Xtilde(D))
\simeq
 \Gammatilde_f^{(2)}(\Xtilde(D))$,
we identify 
$\nbiga^{<H_X\leq D_X'}_{\Xtilde(D)\times\cnumtilde^{\ell}}
 \otimes_{\nbigo_{X\times\cnum^{\ell}}}
 \nbigo_{\Gamma_f(X)}$
with 
$\nbiga^{<D_X}_{\Xtilde(D_X)}$
as sheaves on $\Xtilde(D_X)$.
Then, the claim of Proposition \ref{prop;13.4.20.21}
follows from (\ref{eq;13.4.20.40}).
\hfill\qed

\vspace{.1in}

Let us finish the proof of Theorem \ref{thm;13.4.20.1}.
Let $(X,f)$ be any object in $\Cat_{\ell}$.
We take any projective birational morphism
$\varphi:(Y,g)\lrarr (X,f)$
such that (i) $D_Y$ is normal crossing,
(ii) $Y\setminus D_Y\simeq X\setminus D_X$.
We set $D_Y':=D_Y\times\cnum^{\ell}$
and $D_X':=D_X\times\cnum^{\ell}$.
We have
$R\varphi_{\ast}\nbigo_{Y}\bigl(\ast D_Y\bigr)
\simeq
 \nbigo_{X}(\ast D_X)$.
By using Theorem \ref{thm;13.4.19.42},
we obtain 
\[
 R\varphitilde_{1\ast}\Bigl(\!
 \nbiga^{\moderate}_{Y\times\cnumtilde^{\ell}}
 \otimes^L_{\pi^{-1}\nbigo_{Y\times\cnum^{\ell}}}
 \pi^{-1}\Gamma_{g\ast}(\nbigo_Y(\ast D_Y))
 \!
 \Bigr)\!
\simeq\!
 \nbiga^{\moderate}_{X\times\cnumtilde^{\ell}}
 \otimes^L_{\pi^{-1}\nbigo_{X\times\cnum^{\ell}}}
 \pi^{-1}\Gamma_{f\ast}\bigl(
 \nbigo_X(\ast D_X)
 \bigr).
\]
By using Proposition \ref{prop;13.4.19.30},
we obtain
\begin{multline}
  R\varphitilde_{1\ast}\Bigl(
 \nbiga^{\moderate}_{Y\times\cnumtilde^{\ell}}
 \otimes^L_{\pi^{-1}\nbigo_{Y\times\cnum^{\ell}}}
 \pi^{-1}\Gamma_{g\ast}(\nbigo_Y(\ast D_Y))
 \Bigr)
\simeq
 \\
 R\varphitilde_{1\ast}\Bigl(
 \nbiga^{\moderate}_{Y\times\cnumtilde^{\ell}}
 \otimes^L_{\pi^{-1}\nbigo_{Y\times\cnum^{\ell}}}
 \pi^{-1}\nbigo_{\Gamma_g(Y)}
 \Bigr)
\simeq
 R\varphitilde_{1\ast}\Bigl(
 \nbiga^{\moderate}_{Y\times\cnumtilde^{\ell}}
 \otimes_{\pi^{-1}\nbigo_{Y\times\cnum^{\ell}}}
 \pi^{-1}\nbigo_{\Gamma_g(Y)}
 \Bigr).
\end{multline}
We also have
\[
 \nbiga^{\moderate}_{X\times\cnumtilde^{\ell}}
 \otimes^L_{\pi^{-1}\nbigo_{X\times\cnum^{\ell}}}
 \pi^{-1}\Gamma_{f\ast}\bigl(
 \nbigo_X(\ast D_X)
 \bigr)
\simeq
 \nbiga^{\moderate}_{X\times\cnumtilde^{\ell}}
 \otimes^L_{\pi^{-1}\nbigo_{X\times\cnum^{\ell}}}
 \pi^{-1}\nbigo_{\Gamma_f(X)}.
\]
We obtain
$R\varphitilde_{1\ast}\Bigl(
 \nbiga^{\moderate}_{Y\times\cnumtilde^{\ell}}
 \otimes_{\pi^{-1}\nbigo_{Y\times\cnum^{\ell}}}
 \pi^{-1}\nbigo_{\Gamma_g(Y)}
 \Bigr)
\simeq
 \nbiga^{\moderate}_{X\times\cnumtilde^{\ell}}
 \otimes^L_{\pi^{-1}\nbigo_{X\times\cnum^{\ell}}}
 \pi^{-1}\nbigo_{\Gamma_f(X)}$.
It implies that
the claims for $\nbiga^{\moderate}$
in Theorem \ref{thm;13.4.20.1}.
The claims for $\nbiga^{\rapid}$
can be proved similarly.
\hfill\qed

\subsection{Complement
for the sheaf of Nilsson type functions (Appendix)}

Let us consider an analogue for the sheaves of Nilsson type functions.
We restrict ourselves to the case $\ell=1$.
Let $\nbiga^{\nil}_{X\times\cnumtilde}$
denote the sheaf of holomorphic functions of Nilsson type
on $X\times\cnumtilde$.
\begin{lem}
For any complex manifold $i:(Y,g)\subset (X,f)$ in $\Cat_{1}$,
the naturally defined morphism
\[
 \nbiga^{\nil}_{X\times\cnumtilde}
\otimes^L_{\pi^{-1}\nbigo_{X\times\cnum}}
 \pi^{-1}\nbigo_{Y\times\cnum}
\lrarr
 \itilde_{\ast}
 \nbiga^{\nil}_{Y\times\cnumtilde}
\]
is an isomorphism.
\end{lem}
\pf
As in Lemma \ref{lem;13.4.19.3},
we have an isomorphism
$\nbiga^{\rapid}_{X\times\cnumtilde}
\otimes^L_{\nbigo_{X\times\cnum}}
 \nbigo_{Y\times\cnum}
\simeq
 \nbiga^{\rapid}_{Y\times\cnumtilde}$.
We can check
$\nbiga^{\nil}_{\widehat{\pi^{-1}(H_X)}}
\otimes^L_{\nbigo_{\widehat{H}_X}}
 \nbigo_{\widehat{H_Y}}
\simeq
 \nbiga^{\nil}_{\widehat{\pi^{-1}(H_Y)}}$
directly.
Then, the claim of the lemma follows.
\hfill\qed

\vspace{.1in}

Let $\varphi:(Y,g)\lrarr (X,f)$ be a morphism in $\Cat_{1}$.
For any $\nbigo_Y$-coherent sheaf $M$,
we have the following naturally defined morphism
\begin{equation}
 \label{eq;13.4.20.100}
 \nbiga^{\nil}_{X\times\cnumtilde}
\otimes^L_{\pi^{-1}\nbigo_{X\times\cnum}}
 \pi^{-1}\bigl(
 \Gamma_{f\ast}R\varphi_{\ast}M
 \bigr)
\lrarr
 R\varphitilde_{1\ast}\Bigl(
 \nbiga^{\nil}_{Y\times\cnumtilde}
 \otimes^L_{\pi^{-1}\nbigo_{Y\times\cnum}}
 \pi^{-1}\Gamma_{g\ast}M
 \Bigr).
\end{equation}

\begin{prop}
Suppose that $M$ is $\nbigo_X$-coherent,
and that $\varphi$ is projective.
Then, the morphism {\rm(\ref{eq;13.4.20.100})}
is an isomorphism.
\end{prop}
\pf
By Theorem \ref{thm;13.4.19.42},
we have an isomorphism
\[
\nbiga^{\rapid}_{X\times\cnumtilde}
\otimes^L_{\pi^{-1}\nbigo_{X\times\cnum}}
 \pi^{-1}\bigl(
 \Gamma_{f\ast}R\varphi_{\ast}M
 \bigr)
\simeq
 R\varphitilde_{1\ast}\Bigl(
 \nbiga^{\rapid}_{Y\times\cnumtilde}
 \otimes^L_{\pi^{-1}\nbigo_{Y\times\cnum}}
 \pi^{-1}\Gamma_{g\ast}M
 \Bigr). 
\]
We also have the following formal isomorphism:
\[
\nbiga^{\nil}_{\widehat{H}_X}
\otimes^L_{\pi^{-1}\nbigo_{X\times\cnum}}
 \pi^{-1}\bigl(
 \Gamma_{f\ast}R\varphi_{\ast}M
 \bigr)
\simeq
 R\varphitilde_{1\ast}\Bigl(
 \nbiga^{\nil}_{\Hhat_Y}
 \otimes^L_{\pi^{-1}\nbigo_{Y\times\cnum}}
 \pi^{-1}\Gamma_{g\ast}M
 \Bigr). 
\]
Then, the claim of the proposition follows.
\hfill\qed

\begin{thm}
\label{thm;13.4.20.150}
Let $(X,f)$ be an object in $\Cat_1$.
Let $\varphi:(Y,g)\lrarr (X,f)$
be a projective birational morphism such that 
(i) $D_Y$ is normal crossing,
(ii) $Y\setminus D_Y\simeq X\setminus D_X$.
For the naturally induced map
$\rho:\Ytilde(D_Y)\lrarr X\times\cnumtilde$,
we have
\[
 R\rho_{\ast}\nbiga^{\nil}_{\Ytilde(D_Y)}
\simeq
 \nbiga^{\nil}_{X\times\cnumtilde}
\otimes_{\pi^{-1}\nbigo_{X\times\cnum}}
 \pi^{-1}\nbigo_{\Gamma_f(X)}.
\]
\end{thm}
\pf
As in the proof of Theorem \ref{thm;13.4.20.1},
it is enough to consider the case where
$\varphi=\id$.
We use the notation in the proof of 
Proposition \ref{prop;13.4.19.30}.
We have the isomorphism
$R\nu_{1\ast}
 \nbiga^{\nil}_{\Xtilde(D_X)\times\cnumtilde}
\simeq
 \nbiga^{\nil}_{X\times\cnumtilde}(\ast D_X')$.
Hence, we have the following natural isomorphism
\[
 R\nu_{1\ast}\bigl(
 \nbiga^{\nil}_{\Xtilde(D_X)\times\cnumtilde}
 \otimes^L_{\nbigo_{X\times\cnum}}
 \nbigo_{\Gamma_f(X)}
 \bigr)
\simeq
 \nbiga^{\nil}_{X\times\cnumtilde}
 \otimes_{\nbigo_{X\times\cnum}}
 \nbigo_{\Gamma_f(X)}.
\]
We have the naturally defined morphism
$\nbiga^{\nil}_{\Xtilde(D_X)\times\cnumtilde}
 \otimes_{\nbigo_{X\times\cnum}}
 \nbigo_{\Gamma_f(X)}
\lrarr
 \Gammatilde^{(1)}_{f\ast}
 \nbiga^{\nil}_{\Xtilde(D_X)}$.
It is enough to prove that the following induced morphism
is an isomorphism:
\begin{equation}
 \label{eq;13.4.20.120}
 R\nu_{1\ast}\Bigl(
 \nbiga^{\nil}_{\Xtilde(D_X)\times\cnumtilde}
 \otimes_{\nbigo_{X\times\cnum}}
 \nbigo_{\Gamma_f(X)}
 \Bigr)
\lrarr
 R\nu_{1\ast}\Bigl(
 \Gammatilde^{(1)}_{f\ast}
 \nbiga^{\nil}_{\Xtilde(D_X)}
 \Bigr).
\end{equation}
We have already known that
the following is an isomorphism,
by Proposition \ref{prop;13.4.20.21}:
\[
 R\nu_{1\ast}\Bigl(
 \nbiga^{\rapid}_{\Xtilde(D_X)\times\cnumtilde}
 \otimes_{\nbigo_{X\times\cnum}}
 \nbigo_{\Gamma_f(X)}
 \Bigr)
\lrarr
 R\nu_{1\ast}\Bigl(
 \Gammatilde^{(1)}_{f\ast}
 \nbiga^{\rapid}_{\Xtilde(D_X)}
 \Bigr).
\]
Let $D_X=\bigcup_{i\in\Lambda}D_{i}$ be
the irreducible decomposition.
For any $I\subset \Lambda$,
we set $D_{I0}:=\bigcap_{i\in I}(D_{i}\times\{0\})$.
To prove that (\ref{eq;13.4.20.120})
is an isomorphism,
it is enough to prove that
the following natural morphisms are isomorphisms:
\begin{equation}
 \label{eq;13.4.20.110}
 R\nu_{1\ast}
 \nbiga^{<\del D_{I0}}_{\widehat{\pi^{-1}(D_{I0})}}
 \otimes_{\nbigo_{X\times\cnum}}
 \nbigo_{\Gamma_f(X)}
\lrarr
 R\nu_{1\ast}\Gammatilde^{(1)}_{f\ast}
 \nbiga^{<\del D_I}_{\widehat{\pi_1^{-1}(D_I)}}
\end{equation}
It is enough to consider the issue locally 
around any point of $D_X\times\{0\}$.
We may assume that
$X=\Delta^n$,
$D_X=\bigcup_{i=1}^{\ell}\{z_i=0\}$
and $f=\prod_{i=1}^{\ell}z_i^{m_i}$.

\begin{lem}
\label{lem;13.4.25.100}
We may assume that
$\gcd(m_i\,|\,i\in I)=1$.
\end{lem}
\pf
Let $p:=\gcd(m_i\,|\,i\in I)$.
We set $X':=\Delta^n$
and $D':=\bigcup_{i=1}^{\ell}\{w_i=0\}$.
We define $D_I':=\bigcap_{i\in I}\{w_i=0\}$.
On $X'$,
we set $g:=\prod_{i\not\in I}z_i^{m_i}
 \times \prod_{i\in I}z_i^{m_i/p}$.
We define
$\psi:X\lrarr X'$ by
$z_i\longmapsto z_i^p$ $(i\in I)$
and 
$z_i\longmapsto z_i$ $(i\not\in I)$.
We have $f=g\circ\psi$.
The map $\psi$ gives
$D_I\simeq D_I'$
and 
$\Dtilde_I(\del D_I)\simeq \Dtilde'(\del D_I')$.

Let $\Gammatilde^{(1)}_g:
 \Xtilde'(D')\lrarr 
 \Xtilde'(D')\times\cnumtilde$ 
and 
$\nu'_1:\Xtilde'(D')\times\cnumtilde\lrarr 
 X'\times\cnumtilde$
be given similarly to
$\Gammatilde^{(1)}_f$
and $\nu_1$.
We have the following natural commutative diagram
of the sheaves on $\Dtilde_I(\del D_I)$:
\[
 \begin{CD}
 R\nu_{1\ast}
 \nbiga^{<\del D_{I0}}_{\widehat{\pi^{-1}(D_{I0})}}
 \otimes_{\nbigo_{X\times\cnum}}
 \nbigo_{\Gamma_f(X)}
@>>>
 R\nu_{1\ast}\Gammatilde^{(1)}_{f\ast}
 \nbiga^{<\del D_I}_{\widehat{\pi_1^{-1}(D_I)}} \\
 @A{\simeq}AA @A{\simeq}AA \\
 R\nu'_{1\ast}
 \nbiga^{<\del D'_{I0}}_{\widehat{\pi^{-1}(D'_{I0})}}
 \otimes_{\nbigo_{X'\times\cnum}}
 \nbigo_{\Gamma_g(X')}
@>>>
 R\nu'_{1\ast}\Gammatilde^{(1)}_{g\ast}
 \nbiga^{<\del D'_I}_{\widehat{\pi_1^{-1}(D'_I)}}
 \end{CD}
\]
It is easy to check that
the vertical arrows are isomorphisms.
Then, we obtain the claim of Lemma \ref{lem;13.4.25.100}.
\hfill\qed

\vspace{.1in}

Let $\pi_1:\Xtilde(D_X)\lrarr X$,
$\pi_2:\Xtilde(f)\lrarr X$
and $\pi:\Xtilde(D_X)\times\cnumtilde\lrarr X\times \cnum$
be the projections.
We have 
\[
 \pi_1^{-1}(D_I)
\simeq \Dtilde_I(\del D_I)\times(S^1)^{|I|},\,\,\,
 \pi^{-1}(D_{I0})
\simeq
 \Dtilde_{I}(\del D_I)\times(S^1)^{|I|+1},\,\,\,
 \pi_2^{-1}(D_I)
\simeq
 D_I\times S^1.
\]
We decompose the map
$\nu_{1|\pi^{-1}(D_{I0})}:
 \pi^{-1}(D_{I0})
\lrarr
 \pi_2^{-1}(D_I)$
into
\[
 \Dtilde_I(\del D_I)\times (S^1)^{|I|+1}
\stackrel{\mu_1}{\lrarr}
 \Dtilde_I(\del D_I)\times S^1
\stackrel{\mu_2}{\lrarr}
 D_I\times S^1
\]
To prove that 
(\ref{eq;13.4.20.110}) are isomorphisms,
it is enough to prove that
\begin{equation}
\label{eq;13.4.20.130}
 R\mu_{1\ast}
 \nbiga^{<\del D_{I0}}_{\widehat{\pi^{-1}(D_{I0})}}
 \otimes_{\nbigo_{X\times\cnum}}
 \nbigo_{\Gamma_f}
\lrarr
 R\mu_{1\ast}\Gammatilde^{(1)}_{f\ast}
 \nbiga^{<\del D_I}_{\widehat{\pi_1^{-1}(D_I)}}
\end{equation}
is an isomorphism.

We have the following expression:
\[
 \nbiga^{<\del D_{I0}}_{\widehat{\pi^{-1}(D_{I0})}}
\simeq 
 \varinjlim_{T,N}
\left(
 \nbiga^{\del D_I}_{\Dtilde_I(\del D_I),T,N}
 [\![t,z_i\,|\,i\in I]\!]
 \otimes_{\cnum[t,z_i|i\in i]}
 \Nil(t,z_i|i\in I)
\right)
\]
By the argument in Lemma \ref{lem;09.10.25.5},
or by a direct computation of the cohomology of the sheaves
on the fiber of $\mu_1$,
we obtain
\[
 R\mu_{1\ast}
 \nbiga^{<\del D_{I0}}_{\widehat{\pi^{-1}(D_{I0})}}
\simeq
\varinjlim_{T,N}
\left(
 \nbiga^{<\del D_I}_{\Dtilde_I(\del D_I),T,N}
 [\![t,z_i\,|\,i\in I]\!]
\otimes_{\cnum[t]}
 \Nil(t)
\right)
\]
Hence,
we obtain the following natural isomorphism:
\begin{equation}
 \label{eq;13.4.25.20}
 R\mu_{1\ast}
 \nbiga^{<\del D_{I0}}_{\widehat{\pi^{-1}(D_{I0})}}
 \otimes_{\nbigo_{X\times\cnum}}\nbigo_{\Gamma_f}
\simeq
 \varinjlim_{T,N}
 \left(
 \nbiga^{<\del D_I}_{\Dtilde_I(\del D_I),T,N}
 [\![z_i\,|\,i\in I]\!]
\otimes_{\cnum[t]}
 \Nil(t)
\right)
\end{equation}
Here, $t$ acts as $f$ on 
$\nbiga^{<\del D_I}_{\Dtilde_I(\del D_I),T,N}
 [\![z_i\,|\,i\in I]\!]$.
We have the following expression:
\[
  \nbiga^{<\del D_I}_{\widehat{\pi_1^{-1}(D_I)}}
\simeq
 \varinjlim_{T,N}
 \left(
 \nbiga^{<\del D_I}_{\Dtilde_I(\del D_I),T,N}
 [\![z_i\,|\,i\in I]\!]
 \otimes_{\cnum[z_i|i\in I]}
 \Nil(z_i|i\in I)
\right)
\]
We take $T_0\subset \cnum$
such that $T_0\lrarr \cnum/\seisuu$ is bijective.
We have the decomposition
\[
 \Nil(z_i|i\in I)
=\bigoplus_{\vecalpha\in T_0^I}
 \vecz^{\vecalpha}
 \cnum[z_i,\log z_i|i\in I].
\]
We have the corresponding decomposition:
\begin{multline}
\label{eq;13.4.25.50}
  \nbiga^{<\del D_I}_{\Dtilde_I(\del D_I),T,N}
 [\![z_i\,|\,i\in I]\!]
 \otimes_{\cnum[z_i|i\in I]}
 \Nil(z_i|i\in I)
= \\
 \bigoplus_{\vecalpha\in T_0^I}
  \nbiga^{<\del D_I}_{\Dtilde_I(\del D_I),T,N}
 [\![z_i\,|\,i\in I]\!]
 \vecz^{\vecalpha}
 \otimes\cnum[\log z_i|i\in I]
\end{multline}
Recall $f=\prod_{i=1}^{\ell}z_i^{m_i}$
with $\gcd(m_i\,|\,i\in I)=1$.
Under the assumption,
the map
$\cnum/\seisuu
\lrarr
 (\cnum/\seisuu)^I$
given by
$\beta\longmapsto
 (\beta m_i\,|\,i\in I)$
is injective.
We have the subsheaf
\begin{equation}
\label{eq;13.4.25.51}
 \bigoplus_{\beta\in T_0}
  \nbiga^{<\del D_I}_{\Dtilde_I(\del D_I),T,N}
 [\![z_i\,|\,i\in I]\!]
 \prod_{i\in I}z_i^{\beta m_i}
 \otimes\cnum[\log z_i|i\in I]
\end{equation}
Let $\nbigq$ be the quotient of (\ref{eq;13.4.25.50})
by (\ref{eq;13.4.25.51}).
Note that the fibers of 
$\mu_1\circ\Gammatilde_f^{(1)}$
are connected.
By a direct computation of the sheaves
on the fibers of 
$\mu_1\circ\Gammatilde_f^{(1)}$,
we obtain the push-forward of $\nbigq$
by $\mu_1\circ\Gammatilde_f^{(1)}$
is $0$.
Moreover,
we obtain that the push-forward of (\ref{eq;13.4.25.51})
is naturally isomorphic to
\begin{equation}
\label{eq;13.4.25.52}
 \bigoplus_{\beta\in T_0}
  \nbiga^{<\del D_I}_{\Dtilde_I(\del D_I),T,N}
 [\![z_i\,|\,i\in I]\!]
 \prod_{i\in I}z_i^{\beta m_i}
 \Bigl(
 \log\bigl(\prod_{i=1}^{\ell}z_i^{m_i}\bigr)
 \Bigr)
\end{equation}
Hence, the push-forward of
$\nbiga^{<\del D_I}_{\widehat{\pi_1^{-1}(D_I)}}$
by $\mu_1\circ\Gammatilde_f^{(1)}$ is 
isomorphic to the limit of (\ref{eq;13.4.25.52}).
Together with (\ref{eq;13.4.25.20})
we obtain Theorem \ref{thm;13.4.20.150}.
\hfill\qed

\vspace{.1in}
For any object $(X,f)$ in $\Cat_1$,
we have the sheaves $\nbiga^{\nil}_{X,f}$
on $\Xtilde(f)$
determined by the condition
$\Gammatilde_{f\ast}\nbiga^{\nil}_{X,f}
=\pi^{-1}\nbigo_{\Gamma_f(X)}
 \otimes_{\pi^{-1}\nbigo_{X\times\cnum}}
 \nbiga^{\nil}_{X\times\cnumtilde}$.
\index{sheaf $\nbiga^{\nil}_{X,f}$}
For a morphism
$\varphi:(X_1,f_1)\lrarr (X_2,f_2)$
in $\Cat_1$,
we naturally have
$\varphitilde^{-1}\nbiga^{\nil}_{X_2,f_2}\lrarr
 \nbiga^{\nil}_{X_1,f_1}$.
We obtain the following propositions
as in the case of
$\nbiga^{\nil}$ and $\nbiga^{\moderate}$.

\begin{prop}
For the inclusion $j:X\setminus D_X\lrarr \Xtilde(D_X)$,
the natural morphism 
$\nbiga^{\nil}_{X,f}\lrarr
 j_{\ast}\nbigo_{X\setminus D_X}$
is a monomorphism.
\hfill\qed
\end{prop}

\begin{prop}
\label{prop;13.4.20.202}
Let $\varphi:(Y,g)\lrarr (X,f)$ be a projective morphism in $\Cat_{1}$.
Let $M$ be any coherent $\nbigo_Y$-module.
Then, the following natural morphism is an isomorphism:
\[
 \nbiga^{\nil}_{X,f}\otimes_{\pi^{-1}\nbigo_X}
 \pi^{-1}R\varphi_{\ast}M
\simeq
 R\varphitilde_{\ast}\Bigl(
  \nbiga^{\nil}_{Y,g}\otimes_{\pi^{-1}\nbigo_Y}
 \pi^{-1}M
 \Bigr)
\]
\hfill\qed
\end{prop}

%% file: 4.6.tex
\subsection{Statement}

Let $(X,f)$ be any object in $\Cat_{\ell}$. 
Let $j:X\setminus D_X\lrarr \Xtilde(f)$ denote 
the natural inclusion.
For any $\nbigo_X$-module $M$,
we set 
\[
 \pi^{\ast}_{\moderate}M:=
 \nbiga^{\moderate}_{\Xtilde(f)}
 \otimes_{\pi^{-1}\nbigo_X}\pi^{-1}M. 
\]
\index{functor $\pi^{\ast}_{\moderate}$}
It is also denoted by
$\pi^{\ast}_{f\moderate}M$,
when we would like to emphasize the dependence on $f$.
We shall prove the following theorem.

\begin{thm}
\label{thm;12.9.18.10}
$\nbiga^{\moderate}_{\Xtilde(f)}$ is flat over 
$\pi^{-1}\nbigo_X$,
i.e.,
$\pi^{\ast}_{\moderate}M\simeq
 \nbiga^{\moderate}_{\Xtilde(f)}
 \otimes^L_{\pi^{-1}\nbigo_X}\pi^{-1}M$
for any coherent $\nbigo_X$-module $M$.
Moreover, 
the natural morphism
$\pi^{\ast}_{\moderate}(M)
\lrarr
 j_{\ast}(M_{|X\setminus D_X})$
is injective.
\end{thm}

\begin{cor}
$\nbiga^{\moderate}_{\Xtilde(f)}$
is faithfully flat over
$\pi^{-1}\nbigo_{X}(\ast D_X)$.
\hfill\qed
\end{cor}

We define 
$\pi_{\rapid}^{\ast}M:=
 \nbiga^{\rapid}_{\Xtilde(f)}
\otimes_{\pi^{-1}\nbigo_X}
 \pi_f^{-1}M$.
We can prove the following
by a similar argument.
\begin{prop}
The natural morphism
$\pi^{\ast}_{\rapid}(M)
\lrarr
 j_{\ast}(M_{|X\setminus D_X})$
is injective.
\hfill\qed
\end{prop}

By Theorem {\rm\ref{thm;09.12.4.5}},
$\nbiga_{\Xtilde(f)}^{\rapid}$ is flat over 
$\pi^{-1}\nbigo_X$.
So, we have the following.
\begin{prop}
$\nbiga^{\rapid}_{\Xtilde(f)}$
is faithfully flat over $\pi^{-1}\nbigo_{X}(\ast D_X)$.
\hfill\qed
\end{prop}

\subsection{Induction}

We consider the following conditions for 
any coherent $\nbigo_X$-module $M$:
\begin{description}
\item[($\nbigp1$)]
$\pi^{-1}M\otimes^L_{\pi^{-1}\nbigo_X}\nbiga^{\moderate}_{\Xtilde(f)}
\simeq
\pi^{-1}M\otimes_{\pi^{-1}\nbigo_X}\nbiga^{\moderate}_{\Xtilde(f)}$.
\item[($\nbigp2$)]
$\pi^{\ast}_{\moderate}(M)\lrarr
 j_{\ast}(M_{|X\setminus D_X})$
is injective.
\end{description}
Let $\nbigp(X)$ denote the class of
coherent $\nbigo_X$-modules 
satisfying the conditions $(\nbigp1)$ and $(\nbigp2)$.
It is our purpose to prove
that any coherent $\nbigo_X$-modules
are members of $\nbigp(X)$.
We shall implicitly use that the conditions are local.

\vspace{.1in}
We shall prove the following claim 
by using an induction on $k$:
\begin{description}
\item[($Q_k$)]
Let $(X,f)$ be any object in $\Cat_{\ell}$.
Let $M$ be any coherent $\nbigo_X$-module
such that $\dim\Supp M\leq k$.
Then, $M$ is a member of
$\nbigp(X)$.
\end{description}

\subsection{Preliminary}
The following lemma is easy to prove.
\begin{lem}
Let $0\lrarr M_1\lrarr M_2\lrarr M_3\lrarr 0$
be an exact sequence of coherent $\nbigo_X$-modules.
\begin{itemize}
\item
If $M_2$ and $M_3$ are members of $\nbigp(X)$,
then $M_1$ is also a member of $\nbigp(X)$.
\item
If $M_1$ and $M_3$ are members of $\nbigp(X)$,
then $M_2$ is also a member of $\nbigp(X)$.
\hfill\qed
\end{itemize}
\end{lem}
The following direct corollary will be used
implicitly.
\begin{cor}
Let $\rho:M_1\lrarr M_2$ be any morphism of coherent
$\nbigo_X$-modules
such that 
$\Cok(\rho),\Ker(\rho)\in\nbigp(X)$.
If $M_2$ is contained in $\nbigp(X)$,
then $M_1$ is also contained in $\nbigp(X)$.
\hfill\qed
\end{cor}

\begin{lem}
\label{lem;12.9.17.1}
Let $Z$ be any complex submanifold of $X$
with the inclusion $i_Z:Z\lrarr X$.
Let $M_Z$ be any locally free $\nbigo_Z$-module.
Then,
we have $i_{Z\ast}M_Z\in\nbigp(X)$.
\end{lem}
\pf
It follows from Theorem \ref{thm;13.4.20.1}
and Theorem \ref{thm;13.4.19.210}.
\hfill\qed

\subsection{Functoriality for the push-forward}

Let $\varphi:(X',f')\lrarr (X,f)$ be a morphism in $\Cat_{\ell}$
such that $\varphi:X'\lrarr X$
is projective and birational.
We do not assume that
$X'\setminus D_{X'}$ is isomorphic to $X\setminus D_X$.
Let $D''$ be the exceptional divisor of $\varphi$.
Let $M$ be a coherent $\nbigo_{X'}$-module
such that $M\in\nbigp(X')$.
Assume that $\dim\Supp M= k$
and $\dim\varphi(\Supp M\cap D'')<k$.

\begin{lem}
\label{lem;12.9.17.5}
Assume that $Q_{k-1}$ holds.
Then, we obtain
$\varphi_{\ast}(M')\in\nbigp(X)$.
\end{lem}
\pf
According to Theorem \ref{thm;13.4.19.210},
we have the following isomorphism:
\begin{equation}
 \label{eq;12.9.17.3}
 R\varphitilde_{\ast}\bigl(
 \nbiga^{\moderate}_{\Xtilde'(f')}
\otimes^L_{\pi^{-1}\nbigo_{X'}}
 \pi^{-1}M
 \bigr)
\simeq
 \nbiga^{\moderate}_{\Xtilde(f)}
\otimes^L_{\pi^{-1}\nbigo_X}
 \pi^{-1}R\varphi_{\ast}M
\end{equation}
If $i>0$,
we have
$R^i\varphi_{\ast}M
 \in\nbigp(X)$,
because
$\dim\Supp R^i\varphi_{\ast}M<k$.
By using the degeneration of the spectral sequence,
we obtain
\[
 H^i\bigl(
 \nbiga^{\moderate}_{\Xtilde(f)}
 \otimes^L_{\pi^{-1}\nbigo_X}
 \pi^{-1}R\varphi_{\ast}M
 \bigr)
\simeq
 \left\{
 \begin{array}{ll}
 \Tor_{-i}^{\pi^{-1}\nbigo_X}\bigl(
 \nbiga^{\moderate}_{\Xtilde(f)},
 \pi^{-1}\varphi_{\ast}M
 \bigr)
 & (i<0)\\
 \mbox{{}}\\
 \pi_{\moderate}^{\ast}
 R^i\varphi_{\ast}M
 & (i\geq 0)
 \end{array}
 \right.
\]
By (\ref{eq;12.9.17.3})
and the isomorphism
$\nbiga^{\moderate}_{\Xtilde'(f')}
\otimes^L_{\pi^{-1}\nbigo_{X'}}M
\simeq 
\nbiga^{\moderate}_{\Xtilde'(f')}
\otimes_{\pi^{-1}\nbigo_{X'}}M$,
we have $H^i=0$ for $i<0$.
Hence, we obtain that 
$\varphi_{\ast}M$ satisfies $(\nbigp1)$.
Because 
$\pi_{\moderate}^{\ast}\varphi_{\ast}M
\simeq
 \varphitilde_{\ast}(\pi_{\moderate}^{\ast}M)$,
$(\nbigp2)$ for $\varphi_{\ast}M$
follows from $(\nbigp2)$ for $M$.
\hfill\qed

\vspace{.1in}
We have a direct consequence.
Let $(X',f')\lrarr (X,f)$ be a morphism in $\Cat_{\ell}$
such that 
$\varphi:X'\lrarr X$ is a projective birational morphism.
We do not assume that
$X'\setminus D_{X'}$ is isomorphic to $X\setminus D_X$.
Let $Z'\subset X'$ be a $k$-dimensional 
irreducible complex submanifold.
We assume that $Z'$ is not contained in the exceptional divisor
of $\varphi$,
in particular, $Z'$ is birational to $\varphi(Z')$.
We obtain the following lemma 
from Lemma \ref{lem;12.9.17.1}
and Lemma \ref{lem;12.9.17.5}.

\begin{cor}
Let $M_{Z'}$ be any locally free $\nbigo_{Z'}$-module.
Suppose $Q_{k-1}$.
Then, we have
$\varphi_{\ast}(i_{Z'\ast}M_{Z'})\in\nbigp(X)$.
\hfill\qed
\end{cor}

\subsection{Coherent sheaves on submanifolds}

Let $Z$ be any $k$-dimensional irreducible submanifold of $X$
with the inclusion $i_Z:Z\lrarr X$.

\begin{lem}
\label{lem;12.9.18.50}
Let $M$ be any coherent $\nbigo_X$-module
such that $\Supp(M)\subset Z$.
Assume that $Q_{k-1}$ holds.
Then, we have $M\in\nbigp(X)$.
\end{lem}
\pf
It is enough to consider locally around each point 
$P$ of $X$.
We shall shrink $X$ around $P$ without mention.

First, let us consider the case where
$M=i_{Z\ast}M_Z$.
We may assume that $M_Z$ is a torsion-free $\nbigo_Z$-module.
We can find a projective birational morphism
$\varphi:(X',f')\lrarr (X,f)$ in $\Cat_{\ell}$
such that
(i) the strict transform $Z'$ of $Z$ is a complex
 submanifold of $X'$,
(ii) there exists a locally free $\nbigo_{Z'}$-module $M'$
 with a morphism $\psi:\varphi^{\ast}M\lrarr M'$
 such that $\psi_{|X'\setminus D''}$ is an isomorphism.
We obtain a morphism
$\psi_1:M\lrarr \varphi_{\ast}M'$,
which is an isomorphism on $Z\setminus \varphi(D'')$.
By $Q_{k-1}$,
$\Ker\psi_1$ and $\Cok\psi_1$ are contained in
$\nbigp(X)$.
Then, we obtain
$\iota_{Z\ast}M\in\nbigp(X)$.

In the general case,
we have a finite increasing filtration 
$F=\bigl\{F_i(M)\,\big|\,i=0,\ldots,N\bigr\}$ of $M$
by $\nbigo_X$-modules
such that each $F_i(M)/F_{i-1}(M)$ comes from
an $\nbigo_Z$-module.
Then, the claim of the lemma is reduced to the result
in the previous paragraph.
\hfill\qed

\subsection{End of the proof of Theorem
\ref{thm;12.9.18.10}}

Let $Z$ be any $k$-dimensional irreducible reduced analytic subset
of $X$ such that $Z\not\subset D_X$.
\begin{lem}
\label{lem;12.9.17.12}
Let $M$ be any coherent $\nbigo_X$-module
such that $\Supp(M)\subset Z$.
Assume that $Q_{k-1}$ holds.
Then, we have $M\in\nbigp(X)$.
\end{lem}
\pf
It is enough to consider the issue
locally around any point $P$ of $X$.
Hence, we shall shrink $X$ around $P$
without mention.
Let $Z_1$ denote the union of
the singular points of $Z$ and $D_X\cap Z$.
There exists a projective birational morphism
$\varphi_P:(X',f')\lrarr (X,f)$ in $\Cat_{\ell}$
with the following property:
\begin{itemize}
\item
 The induced morphism
 $X\setminus D''\lrarr X\setminus(Z_1\cup D)$
 is an isomorphism.
\item
 The strict transform $Z'$ of $Z$ is a complex
 submanifold of $X'$.
\end{itemize}
We have $M\lrarr \varphi_{\ast}\varphi^{\ast}M$,
which is an isomorphism outside the singular locus of $Z$.
Hence, we obtain $M\in\nbigp(X)$
by Lemma \ref{lem;12.9.17.5}
and Lemma \ref{lem;12.9.18.50}.
\hfill\qed

\vspace{.1in}

Let $M$ be any coherent $\nbigo_X$-module
such that $\dim\Supp(M)\leq k$.
If we have a decomposition
$\Supp(M)=Z_1\cup Z_2$
such that $Z_1\cap Z_2\subsetneq Z_i$,
then we have an exact sequence
$0\lrarr M_1\lrarr M\lrarr M_2\lrarr 0$
of coherent $\nbigo_X$-modules,
such that
$\Supp(M_i)\subset Z_i$.
Hence, by an easy induction,
we obtain $M\in\nbigp(X)$
from Lemma \ref{lem;12.9.17.12}.
Thus, our induction can proceed,
and the proof of Theorem \ref{thm;12.9.18.10}
is finished.
\hfill\qed

%% file: 4.7.tex
\subsection{Rapid decay and moderate growth}

Let $(X,f)$ be any object in $\Cat_{\ell}$.
We put $\nbigd^{\moderate}_{\Xtilde(f)}:=
\pi^{-1}(\nbigd_X)\otimes_{\pi^{-1}\nbigo_X}
 \nbiga^{\moderate}_{\Xtilde(f)}$.
For any $\nbigd_X$-module $\nbigm$,
we set 
\[
\pi_{f\moderate}^{\ast}(\nbigm):=
 \pi^{-1}\nbigm\otimes_{\pi^{-1}\nbigo_X}
 \nbiga^{\moderate}_{\Xtilde(f)},
\quad\quad
\pi_{f\rapid}^{\ast}(\nbigm):=
 \pi^{-1}\nbigm\otimes_{\pi^{-1}\nbigo_X}
 \nbiga^{\rapid}_{\Xtilde(f)}.
\]
They are naturally $\nbigd^{\moderate}_{\Xtilde(f)}$-modules.
\index{functor $\pi_{f\moderate}^{\ast}$}
\index{functor $\pi_{f\rapid}^{\ast}$}

Let $\varphi:(X,f)\lrarr (Y,g)$ be any morphism in $\Cat_{\ell}$.
For any $\nbigd^{\moderate}_{\Xtilde(f)}$-module $\nbigmtilde$,
we put
\begin{equation}
 \label{eq;13.4.20.200}
 \varphitilde_{\dagger}(\nbigmtilde):=
 R\varphitilde_{!}\bigl(
 \pi^{-1}(\nbigd_{Y\larr X})
 \otimes^L_{\pi^{-1}\nbigd_X}
 \nbigmtilde
 \bigr).
\end{equation}
\index{functor $\varphitilde_{\dagger}$}

Let $\nbigm$ be any $\nbigd_X$-module.
We have the following naturally defined morphism:
\[
 \varphitilde^{-1}\nbiga^{\moderate}_{\Ytilde(g)}
 \otimes_{\varphitilde^{-1}\pi^{-1}\nbigo_Y}
 \pi^{-1}(\nbigm)
\lrarr
 \pi_{f\moderate}^{\ast}(\nbigm)
\]
It induces the following morphism
in the derived category of $\nbigd_{Y,g}$-modules:
\[
 R\varphitilde_{!}\Bigl(
 \pi^{-1}(\nbigd_{Y\larr X}\otimes_{\nbigd_X}^L\nbigm)
 \otimes_{\varphitilde^{-1}\pi^{-1}\nbigo_Y}
 \varphitilde^{-1}\nbiga^{\moderate}_{\Ytilde(g)}
 \Bigr)
\lrarr
 \varphitilde_{\dagger}
 \pi_{f\moderate}^{\ast}(\nbigm)
\]
We also have the following isomorphisms:
\begin{multline*}
 R\varphitilde_{!}\Bigl(
 \pi^{-1}(\nbigd_{Y\larr X}\otimes_{\nbigd_X}^L\nbigm) 
 \otimes_{\varphitilde^{-1}\pi^{-1}\nbigo_Y}
 \varphitilde^{-1}\nbiga^{\moderate}_{\Ytilde(g)}
 \Bigr)
\simeq \\
 R\varphitilde_{!}\Bigl(
 \pi^{-1}(\nbigd_{Y\larr X}\otimes_{\nbigd_X}^L\nbigm)
 \Bigr)
 \otimes_{\pi^{-1}\nbigo_Y}
\nbiga^{\moderate}_{\Ytilde(g)}
\simeq \\
 \pi^{-1}R\varphi_{!}(\nbigd_{Y\larr X}
 \otimes_{\nbigd_X}^L\nbigm)
 \otimes_{\pi^{-1}\nbigo_Y}\nbiga^{\moderate}_{\Ytilde(g)}
\simeq
 \pi_{g\moderate}^{\ast}\varphi_{\dagger}\nbigm
\end{multline*}
Hence, we obtain the following morphism
in the derived category of 
$\nbigd_{Y,g}$-modules:
\begin{equation}
\label{eq;10.12.29.21}
 \pi_{g\moderate}^{\ast}\varphi_{\dagger}\nbigm 
\lrarr
 \varphitilde_{\dagger}\pi_{f\moderate}^{\ast}(\nbigm)
\end{equation}
Similarly, we obtain the following morphism:
\begin{equation}
\label{eq;12.9.19.100}
 \pi_{g\rapid}^{\ast}\varphi_{\dagger}\nbigm 
\lrarr
 \varphitilde_{\dagger}\pi_{f\rapid}^{\ast}(\nbigm)
\end{equation}

\begin{prop}
Assume that $\varphi$ is projective,
and that $M$ has a good filtration
in the neighbourhood of fibers of $\varphi$.
Then, the morphisms {\rm(\ref{eq;10.12.29.21})}
and {\rm(\ref{eq;12.9.19.100})}
are isomorphisms.
\end{prop}
\pf
By considering a resolution,
it is enough to consider the case
$\nbigm=
 M\otimes_{\nbigo_X}\nbigd_X\otimes\Omega_X^{-1}$,
and $M$ is an $\nbigo_X$-coherent sheaf.
Then, the claim is reduced to Theorem \ref{thm;13.4.19.210}.
\hfill\qed

\vspace{.1in}

Let $(X,f)$ be an object in $\Cat_{\ell}$
such that $D_X$ is normal crossing.
We set
$\nbigd^{\moderate}_{\Xtilde(D_X)}:=
 \nbiga^{\moderate}_{\Xtilde(D_X)}
 \otimes_{\pi^{-1}\nbigo_X}\pi^{-1}\nbigd_X$.
Let $\pi_1:\Xtilde(D_X)\lrarr X$ be the projection.
For any $\nbigd_X$-module $\nbigm$,
we define
\[
 \pi_{1\moderate}^{\ast}\nbigm:=
 \nbiga^{\moderate}_{\Xtilde(D_X)}
 \otimes_{\pi^{-1}\nbigo_X}
 \nbigm,
\quad
 \pi_{1\rapid}^{\ast}\nbigm:=
 \nbiga^{\rapid}_{\Xtilde(D_X)}
 \otimes_{\pi^{-1}\nbigo_X}
 \nbigm.
\]
We have the naturally defined proper map
$\rho:\Xtilde(D_X)\lrarr \Xtilde(f)$.
We obtain the following proposition
from Theorem \ref{thm;13.4.20.1}.
\begin{prop}
\label{prop;13.4.20.220}
We have the following natural isomorphisms
for any coherent $\nbigd_X$-module $\nbigm$:
\[
 R\rho_{\ast}
 \pi^{\ast}_{1\moderate}\nbigm
\simeq
\pi^{\ast}_{f\moderate}\nbigm,
\quad\quad
 R\rho_{\ast}
 \pi^{\ast}_{1\rapid}\nbigm
\simeq
\pi^{\ast}_{f\rapid}\nbigm.
\]
\hfill\qed
\end{prop}

\subsection{Compatibility with the de Rham functor}

For any $\nbigd_X$-module $\nbigm$,
we put
\[
\DR^{\moderate}_{X,f}(\nbigm)
:=\pi^{-1}(\DR_X\nbigm)\otimes
 _{\pi^{-1}\nbigo_X}\nbiga^{\moderate}_{\Xtilde(f)}
\simeq
 \pi^{-1}(\Omega_X)
 \otimes_{\pi^{-1}\nbigd_X}^L
 \pi_{f\moderate}^{\ast}(\nbigm),
\]
\index{functor $\DR^{\moderate}_{X,f}$}
\[
\DR^{\rapid}_{X,f}(\nbigm)
:=\pi^{-1}(\DR_X\nbigm)\otimes
 _{\pi^{-1}\nbigo_X}\nbiga_{\Xtilde(f)}^{\rapid}
\simeq
 \pi^{-1}(\Omega_X)
 \otimes_{\pi^{-1}\nbigd_X}^L
 \pi_{f\rapid}^{\ast}(\nbigm).
\]
\index{functor $\DR^{\rapid}_{X,f}$}

\begin{cor}
\label{cor;13.4.20.230}
Suppose that  $\nbigm$ has a good filtration
in the neighbourhood of fibers of $\varphi$.
Assume that $\varphi$ is projective.
Then, we have natural isomorphisms:
\[
 R\varphitilde_{!}
 \DR^{\moderate}_{X,f}(\nbigm)
\simeq
 \DR^{\moderate}_{Y,g}
 \varphi_{\dagger}(\nbigm),
\quad\quad
 R\varphitilde_{!}
 \DR^{\rapid}_{X,f}(\nbigm)
\simeq
 \DR^{\rapid}_{Y,g}
 \varphi_{\dagger}(\nbigm)
\]
\end{cor}
\pf
From 
$\varphitilde_{\dagger}\pi_{\moderate}^{\ast}\nbigm
\simeq
 \pi_{\moderate}^{\ast}\varphi_{\dagger}\nbigm$,
we obtain the following isomorphisms:
\begin{multline}
R\varphitilde_{!}
 \DR^{\moderate}_{X,f}\nbigm
\simeq
R\varphitilde_{!}\bigl(
 \pi^{-1}
 \Omega_X\otimes_{\pi^{-1}\nbigd_X}^L
 \pi^{\ast}_{f\moderate}\nbigm
 \bigr)
\simeq
\pi^{-1}\Omega_Y
 \otimes_{\pi^{-1}\nbigd_Y}^L
 \varphitilde_{\dagger}
 \pi_{f\moderate}^{\ast}\nbigm \\
\simeq
 \pi^{-1}\Omega_Y
 \otimes_{\pi^{-1}\nbigd_Y}^L
 \bigl(
 \pi_{g\moderate}^{\ast}\varphi_{\dagger}\nbigm
 \bigr)
\simeq
 \DR^{\moderate}_{Y,g}\varphi_{\dagger}\nbigm
\end{multline}
Thus, we obtain the first isomorphism.
We obtain the second one similarly.
\hfill\qed

\vspace{.1in}

Let $(X,f)$ be an object in $\Cat_{\ell}$
such that $D_X$ is normal crossing.
We consider the real blow up
$\pi_1:\Xtilde(D_X)\lrarr X$.
We define 
$\DR^{\moderate}_{\Xtilde(D_X)}(\nbigm)$
and 
$\DR^{\rapid}_{\Xtilde(D_X)}(\nbigm)$
as follows:
\[
 \DR^{\moderate}_{\Xtilde(D_X)}(\nbigm)
:=\pi^{-1}\Omega\otimes^L_{\pi^{-1}\nbigd_X} 
 \pi_{1\moderate}^{\ast}(\nbigm)
\]
\[
 \DR^{\rapid}_{\Xtilde(D_X)}(\nbigm)
:=\pi^{-1}\Omega\otimes^L_{\pi^{-1}\nbigd_X} 
 \pi_{1\rapid}^{\ast}(\nbigm)
\]
\index{functor $\DR^{\moderate}_{\Xtilde(D_X)}$}
\index{functor $\DR^{\rapid}_{\Xtilde(D_X)}$}
We have the naturally defined proper map
$\rho:\Xtilde(D_X)\lrarr \Xtilde(f)$.

\begin{prop}
\label{prop;13.4.20.231}
The following natural morphisms
are isomorphisms:
\[
 R\rho_{\ast}
 \DR^{\moderate}_{\Xtilde(D_X)}(\nbigm)
\simeq
 \DR^{\moderate}_{X,f}(\nbigm),
\quad
 R\rho_{\ast}
 \DR^{\rapid}_{\Xtilde(D_X)}(\nbigm)
\simeq
 \DR^{\rapid}_{X,f}(\nbigm).
\]
\end{prop}
\pf
It immediately follows from
Proposition \ref{prop;13.4.20.220}.
\hfill\qed

\vspace{.1in}

We obtain the following corollary
from Corollary \ref{cor;13.4.20.230}
and Proposition \ref{prop;13.4.20.231}.
\begin{cor}
\label{cor;13.4.20.400}
Let $\varphi:X\lrarr Y$ be any projective
morphism of complex manifolds.
Let $D_Y$ be a normal crossing hypersurface of $Y$
such that $D_X:=\varphi^{-1}(D_Y)$ is normal crossing.
Let $\varphitilde:\Xtilde(D_X)\lrarr\Ytilde(D_Y)$
be the induced map.
Then, 
for any coherent $\nbigd_X$-module
having a good filtration in the neighbourhood of
fibers of $\varphi$,
we have the following natural isomorphisms:
\begin{equation}
 \label{eq;13.4.20.300}
 R\varphitilde_{!}
 \DR^{\moderate}_{\Xtilde(D_X)}(\nbigm)
\simeq
 \DR^{\moderate}_{\Ytilde(D_Y)}
 \varphi_{\dagger}\nbigm
\end{equation}
\begin{equation}
  R\varphitilde_{!}
 \DR^{\rapid}_{\Xtilde(D_X)}(\nbigm)
\simeq
 \DR^{\rapid}_{\Ytilde(D_Y)}
 \varphi_{\dagger}\nbigm.
\end{equation}
\hfill\qed
\end{cor}

\begin{rem}
G. Morando informed the author
that the isomorphism {\rm(\ref{eq;13.4.20.300})}
and its generalizations can be deduced from
some results
in {\rm\cite{Kashiwara-Schapira-Ind-sheaves}}.
While the author hopes that the generalization would 
make the subject more transparent,
he also hopes that
our direct method would be also significant 
for our understanding.
\hfill\qed
\end{rem}

\subsection{Nilsson type (Appendix)}

We have variants in the case of Nilsson type.
Let $(X,f)$ be an object in $\Cat_1$.
We set $\nbigd^{\nil}_{\Xtilde(f)}:=
 \nbiga^{\nil}_{\Xtilde(f)}\otimes_{\pi^{-1}\nbigo_X}
 \pi^{-1}\nbigd_X$.
For any $\nbigd_X$-module $\nbigm$,
we set 
$\pi_{\nil}^{\ast}(\nbigm):=
 \pi^{-1}\nbigm\otimes_{\pi^{-1}\nbigo_X}
 \nbiga^{\nil}_{\Xtilde(f)}$.
They are naturally $\nbigd^{\nil}_{\Xtilde(f)}$-modules.

Let $\varphi:(X,f)\lrarr (Y,g)$ be a morphism in $\Cat_1$.
For any $\nbigd^{\nil}_{\Xtilde(f)}$-module $\nbigmtilde$,
we define 
$\varphitilde_{\dagger}(\nbigmtilde)$
by the formula (\ref{eq;13.4.20.200}).
We also define 
$\DR^{\nil}_{X,f}(\nbigm):=
 \pi^{-1}\Omega_X\otimes^L_{\pi^{-1}\nbigd_X}
 \pi_{f\nil}^{\ast}\nbigm$.
We obtain the following
from Proposition \ref{prop;13.4.20.202}.

\begin{prop}
Suppose that $\varphi$ is projective
and that $\nbigm$ has a good filtration
in the neighbourhood of fibers of $\varphi$.
Then, the natural morphism
\begin{equation}
 \label{eq;13.4.20.201}
 \pi_{\nil}^{\ast}\varphi_{\dagger}\nbigm 
\lrarr
 \varphitilde_{\dagger}\pi_{\nil}^{\ast}(\nbigm)
\end{equation}
is an isomorphism.
In particular,
a natural morphism
$R\varphitilde_{!}
 \DR^{\nil}_{X,f}(\nbigm)
\simeq
 \DR^{\nil}_{Y,g}
 \varphi_{\dagger}\nbigm$
is an isomorphism.
\hfill\qed
\end{prop}

Let $(X,f)$ be an object in $\Cat_1$
such that $D_X$ is normal crossing.
We consider the real blow up
$\pi_1:\Xtilde(D_X)\lrarr X$.
We define 
$\DR^{\nil}_{\Xtilde(D_X)}(\nbigm):=
 \pi_1^{-1}\Omega
 \otimes_{\pi_1^{-1}\nbigo_X}
 \pi_{1\nil}^{\ast}\nbigm$
for any $\nbigd_X$-module $\nbigm$.
We obtain the following proposition from
Theorem \ref{thm;13.4.20.150}.
\begin{prop}
Let $\rho:\Xtilde(D_X)\lrarr \Xtilde(f)$ 
be the natural map.
We have a natural isomorphism
\[
 R\rho_{\ast}
 \pi^{\ast}_{1\nil}(\nbigm)
\simeq
 \pi^{\ast}_{\nil}(\nbigm).
\]
In particular,
we obtain an isomorphism
 $R\rho_{\ast}\DR^{\nil}_{\Xtilde(D)}(\nbigm)
\simeq
 \DR^{\nil}_{X,f}(\nbigm)$.
\hfill\qed
\end{prop}

\begin{cor}
\label{cor;13.4.20.301}
Let $\varphi:X\lrarr Y$ be any projective morphism
of complex manifolds.
Let $D_Y$ be a smooth hypersurface of $Y$
such that $\varphi^{-1}(D_Y)$ is normal crossing.
Let $\varphitilde:\Xtilde(D_X)\lrarr \Ytilde(D_Y)$
be the induced map.
Then, for any coherent $\nbigd_X$-module $\nbigm$
having a good filtration in the neighbourhood of
fibers of $\varphi$,
we have the natural isomorphism
$R\varphitilde_{!}\DR^{\nil}_{\Xtilde(D_X)}(\nbigm)
\simeq
 \DR^{\nil}_{\Ytilde(D_Y)}(\nbigm)$.
\hfill\qed
\end{cor}

%% file: 5.1.tex
\subsection{De Rham complex
and a description by dual}
\label{subsection;09.10.26.10}

Let $X$ be a complex manifold
and $D$ be a normal crossing hypersurface
with a decomposition $D=D_1\cup D_2$.
(Note that $D_i$ are not necessarily irreducible.
 See \S\ref{subsection;10.1.11.2}.)
We set $d_X:=\dim X$.
Let $\pi:\Xtilde(D)\lrarr X$ be the real blow up.
Let $\Omega_X^{\bullet}$ denote the sheaf of holomorphic
$1$-forms on $X$.
We put
\[
\Omega^{\bullet\,<D_1\leq D_2}_{\Xtilde(D)}:=
 \nbiga^{<D_1\leq D_2}_{\Xtilde(D)}
 \otimes_{\pi^{-1}\nbigo_X}
 \pi^{-1}\Omega_X^{\bullet},
\]
\[
\Omega^{\bullet,\bullet\,<D_1\leq D_2}
 _{\Xtilde(D)}
:=
 \Omega^{0,\bullet\,<D_1\leq D_2}_{\Xtilde(D)}
 \otimes_{\pi^{-1}\nbigo_X}
 \pi^{-1}\Omega_X^{\bullet}.
\]
\index{sheaf $\Omega^{\bullet\,<D_1\leq D_2}_{\Xtilde(D)}$}
\index{sheaf $\Omega^{\bullet,\bullet\,<D_1\leq D_2}
 _{\Xtilde(D)}$}
For any holonomic $\nbigd$-module
$\nbigm$ on $X$,
we define
\begin{multline*}
 \DR^{<D_1\leq D_2}_{\Xtilde(D)}(\nbigm):=
 \nbiga^{<D_1\leq D_2}_{\Xtilde(D)}
 \otimes_{\pi^{-1}\nbigo_X}
 \pi^{-1}\DR_X(\nbigm)
 \\
\simeq
 \Omega^{\bullet\,<D_1\leq D_2}_{\Xtilde(D)}
 [d_X]
\otimes_{\pi^{-1}\nbigo_X}
 \pi^{-1}\nbigm
\simeq
\Tot\Bigl(
 \Omega^{\bullet,\bullet\,<D_1\leq D_2}
 _{\Xtilde(D)}
 \otimes_{\pi^{-1}\nbigo_X}
 \pi^{-1}\nbigm
\Bigr)[d_X].
\end{multline*}
\index{functor $\DR_{\Xtilde(D)}^{<D_1\leq D_2}$}
Note 
$\DR^{<D_1\leq D_2}_{\Xtilde(D)}\bigl(
 \nbigm\bigr)
\simeq
 \DR^{<D_1\leq D_2}_{\Xtilde(D)}\bigl(
 \nbigm(\ast D)
 \bigr)$
because
$\Omega_{\Xtilde(D)}^{\bullet\,<D_1\leq D_2}
 (\ast D)=
 \Omega_{\Xtilde(D)}^{\bullet\,<D_1\leq D_2}$.

We have a natural isomorphism
$R\pi_{\ast}\DR_{\Xtilde(D)}^{<D_1\leq D_2}\bigl(
 \nbigm \bigr)
\simeq
 \DR_X^{<D_1\leq D_2}\nbigm$
induced as follows,
by Theorem \ref{thm;09.12.4.5}:
\begin{multline}
\label{eq;09.10.4.1}
 R\pi_{\ast}\Tot\Bigl(
 \Omega_{\Xtilde(D)}^{\bullet,\bullet\,<D_1\leq D_2}
 \otimes_{\pi^{-1}\nbigo_X}
 \pi^{-1}\nbigm
 \Bigr)[d_X]
\simeq
\Tot\Bigl(
 R\pi_{\ast}
 \Omega_{\Xtilde(D)}^{\bullet,\bullet\,
 <D_1\leq D_2} 
\otimes_{\nbigo_X}\nbigm
\Bigr)
[d_X]
 \\
\simeq 
 \Tot\Bigl(
 \Omega_{X}^{\bullet,\bullet\,<D_1}(\ast D_2)
\otimes_{\nbigo_X}\nbigm
 \Bigr)[d_X]
\end{multline}

\begin{lem}
\label{lem;09.10.3.40}
We have a natural isomorphism
\[
 \nrhom_{\pi^{-1}\nbigd_X}
 \bigl(
 \pi^{-1}\nbigm,\nbiga_{\Xtilde(D)}^{<D_1\leq D_2}
 \bigr)[d_X]
\simeq
 \DR^{<D_1\leq D_2}_{\Xtilde(D)}
 \bigl(\DDD\nbigm\bigr).
\]
\end{lem}
\pf
Since $\nbigm$ is $\nbigd_X$-coherent,
we have the following isomorphisms:
\begin{multline}
\label{eq;09.10.2.1}
 \nrhom_{\pi^{-1}\nbigd_X}\bigl(
 \pi^{-1}\nbigm,\,
 \nbiga^{<D_1\,\leq D_2}_{\Xtilde(D)}
 \bigr)[d_X] \\
\simeq
 \nrhom_{\pi^{-1}\nbigd_X}\bigl(
 \pi^{-1}\nbigm,\,\pi^{-1}\nbigd_X
 \bigr)
 \otimes_{\pi^{-1}\nbigd_X}^L
 \nbiga_{\Xtilde(D)}^{<D_1\,\leq D_2}[d_X]
 \\
=\pi^{-1}\Bigl(
 \Omega_X\otimes_{\nbigo_X}
 \DDD\nbigm
 \Bigr)
\otimes^L_{\pi^{-1}\nbigd_X}
 \nbiga^{<D_1\leq D_2}_{\Xtilde(D)}
 \\
\simeq
 \Bigl(
 \pi^{-1}\Omega_X\otimes_{\pi^{-1}\nbigo_X}
 \nbiga^{<D_1\leq D_2}_{\Xtilde(D)}
\Bigr)
\otimes^L_{\pi^{-1}\nbigd_X}
 \pi^{-1}
 \DDD\nbigm
\end{multline}
Because
$\nbiga^{<D_1\leq D_2}_{\Xtilde(D)}$
is flat over $\pi^{-1}\nbigo_X$
(Theorem \ref{thm;09.12.4.5}),
$\pi^{-1}\nbigd_X\otimes_{\pi^{-1}\nbigo_X}
 \nbiga^{<D_1\leq D_2}_{\Xtilde(D)}$
is flat over $\pi^{-1}\nbigd_X$.
Therefore,
\[
 \nbiga_{\Xtilde(D)}^{<D_1\leq D_2}
\simeq
 \pi^{-1}\Bigl(
 \nbigd_X\otimes_{\nbigo_X}\Theta_X^{-\bullet}
 \Bigr)
\otimes_{\pi^{-1}\nbigo_X}
 \nbiga^{<D_1\leq D_2}_{\Xtilde(D)}
\]
is a $\pi^{-1}\nbigd_X$-flat resolution.
Hence, (\ref{eq;09.10.2.1}) is quasi-isomorphic to
the following:
\begin{multline}
 \left(
 \pi^{-1}\Bigl(
 \Omega_X^{\bullet}\otimes \nbigd_X
 \Bigr)
 \otimes_{\pi^{-1}\nbigo_X}
 \nbiga_{\Xtilde(D)}^{<D_1\leq D_2}
 \right)
 \otimes_{\pi^{-1}\nbigd_X}
 \pi^{-1}\DDD\nbigm[d_X]
 \\
\simeq
 \Omega^{\bullet\,<D_1\leq D_2}_{\Xtilde(D)}
 \otimes_{\pi^{-1}\nbigo_X}
 \pi^{-1}\DDD\nbigm[d_X]
\end{multline}
Thus, we obtain the desired isomorphism.
\hfill\qed

\vspace{.1in}

According to Lemma \ref{lem;09.10.3.40},
we have a natural isomorphism
\begin{multline}
 \label{eq;09.10.3.120}
 \DR_{\Xtilde(D)}^{<D_1\leq D_2}\bigl(\nbigm\bigr)
\simeq
 \nrhom_{\pi^{-1}\nbigd_X}
 \bigl(
 \pi^{-1}\DDD\nbigm,\,
 \nbiga_{\Xtilde(D)}^{<D_1\leq D_2}
 \bigr)[d_X] 
\\ 
\simeq 
  \nrhom_{\pi^{-1}\nbigd_X}
 \bigl(
 \pi^{-1}\DDD\bigl(
 \nbigm(\ast D)\bigr),\,
 \nbiga_{\Xtilde(D)}^{<D_1\leq D_2}
 \bigr)[d_X]
\end{multline}
We will implicitly identify them
in the following argument.

\subsection{A combinatorial description
in the case of good meromorphic flat bundles}
\label{subsection;09.10.28.1}

Let $X$ be a complex manifold
with a normal crossing hypersurface $D$.
Let $\pi:\Xtilde(D)\lrarr X$ be the real blow up.
Let $V$ be a good meromorphic flat bundle
on $(X,D)$.
We have the local system on $X-D$
associated to $V_{|X-D}$.
Its prolongment over $\Xtilde(D)$
is denoted by $\nbigl$.
If $V$ is unramifiedly good,
for any $P\in \pi^{-1}(D)$,
we have the Stokes filtration $\nbigf^P$
of the stalk $\nbigl_{P}$ indexed by
the set of the irregular values
$\Irr(V,\pi(P))\subset
 \nbigo_{X}(\ast D)_{\pi(P)}\big/\nbigo_{X,\pi(P)}$
with the order $\leq_P$.
The system of filtrations
$\bigl\{\nbigf^P\,\big|\,P\in\pi^{-1}(D)\bigr\}$
satisfies some compatibility condition.
See \cite{mochi7}, \cite{mochi8} 
or \S3 of \cite{mochi10} 
for more details.

Let $D=D_1\cup D_2$ be a decomposition.
Let us describe 
$\DR^{<D_1\leq D_2}_{\Xtilde(D)}(V)$
in terms of the Stokes filtrations.
If $V$ is unramifiedly good,
for $P\in \Xtilde(D)$,
let $\nbigl^{<D_1\leq D_2}_P$
be the union of the subspaces
$\nbigf^P_{\gminia}(\nbigl_P)
 \subset\nbigl_P$
such that
(i) $\gminia\leq_P 0$,
(ii) the poles of $\gminia$
contain the germ of $D_1$ at $\pi(P)$.
\index{sheaf $\nbigl^{<D_1\leq D_2}$}
If $V$ is not unramifiedly good,
we take a ramified covering
$\varphi:(X',D')\lrarr (X,D)$
such that
$V'=\varphi^{\ast}V$ is unramifiedly good.
We obtain the local system
$\nbigl'$ 
and a sheaf $\nbigl^{\prime< D_1'\leq D_2'}$
on $\Xtilde'(D')$ associated to $V'$
with the Stokes structure.
By taking the descent,
we obtain a subsheaf 
$\nbigl^{< D_1\leq D_2}
\subset\nbigl$.

\begin{lem}
The family 
$\bigl\{
 \nbigl_{P}^{<D_1\leq D_2}
 \bigr\}$ gives a constructible sheaf
$\nbigl^{<D_1\leq D_2}$ on $\Xtilde(D)$.
\end{lem}
\pf
It is enough to consider the case
$X=\Delta^n$ and 
$D=\bigcup_{i=1}^{\ell}\{z_i=0\}$.
We may also assume that $V$
is unramifiedly good.
By using a decomposition around $P$
as in Theorem 4.1 of \cite{mochi10},
it is easy to observe that
it is enough to consider 
the case $V=\nbigo_X(\ast D)$
with a flat connection
$\nabla e=e\,d\gminia$,
where $\gminia=\prod_{i=1}^{m}z_i^{-m_i}$
$(m_i>0)$
for some $1\leq m\leq \ell$.
We have a decomposition
$\ellsitabar=I_1\sqcup I_2$
such that $D_j=\bigcup_{i\in I_j}\{z_i=0\}$.
For $P\in\Xtilde(D)$,
we set $I_j(P):=\bigl\{
 i\in I_j\,\big|\,z_i\bigl(\pi(P)\bigr)=0
 \bigr\}$.
We set $F_{\gminia}:=-|\gminia|^{-1}\Re\gminia$.
We put 
$R_0:=\bigcup_{i=1}^{m}\{z_i=0\}$
and 
$R_1:=\bigcup_{i=m+1}^{\ell}\{z_i=0\}
 \setminus R_0$.
\begin{itemize}
\item
For $P\in X-D$,
we have $\nbigl^{<D_1\leq D_2}_P\neq 0$.
\item
For $P\in \pi^{-1}(R_1)$,
we have $\nbigl^{<D_1\leq D_2}_P\neq 0$
if and only if 
$I_1(P)=\emptyset$.
\item
For $P\in \pi^{-1}(R_0)$,
we have 
$\nbigl^{<D_1\leq D_2}_P\neq 0$
if and only if 
(i) $F_{\gminia}(P)<0$,
(ii) $I_1(P)\subset \mbar$.
\end{itemize}
Then, the claim of the lemma is clear.
\hfill\qed

\vspace{.1in}
We recall the following proposition.
(See \cite{majima} and \cite{sabbah4}.
See also \cite{hien1}.)

\begin{prop}
\label{prop;09.10.26.2}
The natural inclusion
$\nbigl^{<D_1\leq D_2}[d_X]
\lrarr
\DR^{<D_1\leq D_2}_{\Xtilde(D)}(V)$
is a quasi-isomorphism.
\end{prop}
\pf
We give a preparation from elementary analysis
on multi-sectors.
We set $Y:=\Delta_z\times\Delta^n_{\vecw}$
and $D_Y=\{z=0\}\cup\bigcup_{i=1}^{\ell}\{w_i=0\}$.
Let $\pi:\Ytilde(D_Y)\lrarr Y$ be the real blow up.
For $m>0$ and 
$\vecm=(m_1,\ldots,m_k)\in \seisuu_{>0}^k$
$(0\leq k\leq \ell)$,
we put $\gminia=z^{-m}\prod_{i=1}^kw_i^{-m_i}$.
We put $F_{\gminia}=-|\gminia^{-1}|\Re(\gminia)$,
which naturally gives a $C^{\infty}$-function
on $\Ytilde(D_Y)$.
Take a point $P\in \pi^{-1}(O)\subset \Ytilde(D_Y)$.
Let $S=S_z\times S_{\vecw}$ be a small
multi-sector in $Y-D_Y$
such that $P$ is contained in the interior
part of the closure of $\Sbar$ in $\Ytilde(D_Y)$.
\begin{itemize}
\item
If $F_{\gminia}(P)<0$
(resp. $F_{\gminia}(P)>0$),
we assume that 
$F_{\gminia}<0$
(resp. $F_{\gminia}>0$)  on $\Sbar$.
\item
If $F_{\gminia}(P)=0$,
we assume that
$F_{\gminia}$ is monotonous
with respect to $\theta$,
where $z=re^{\sqrt{-1}\theta}$ 
is the polar coordinate system.
Let $\theta_{i}$ $(i=1,2)$ be the arguments
of the edges of $S_z$,
i.e.,
$S_z=\bigl\{
 (r,\theta)\,\big|\,\theta_1\leq \theta\leq \theta_2,\,
 0<r\leq r_0 \bigr\}$.
Let $\theta_{+}$ be one of $\theta_i$
such that
$F_{\gminia}>0$ on
$\{re^{\sqrt{-1}\theta_+}\}
 \times \Sbar_{\vecw}$.
\end{itemize}

Let $f$ be a holomorphic function on $S$
of moderate growth with respect to
$z$ and $\vecw$.
We set
\begin{equation}
 \label{eq;09.10.26.3}
 \Phi(f)(z,\vecw):=
 \int_{\gamma(z,\vecw)}
 \exp\bigl(
 -\gminia(z,\vecw)+\gminia(\zeta,\vecw)
 \bigr)\,
 f(\zeta,\vecw)\,d\zeta.
\end{equation}
Here, $\gamma(z,\vecw)$ is a path 
contained in $S_z\times\{\vecw\}$
taken as follows.
\begin{description}
\item[Case $F_{\gminia}(P)<0$]
We fix a point $z_0\in S_z$,
and $\gamma(z,\vecw)$ is a path 
from $z_0$ to $z$.
\item[Case $F_{\gminia}(P)>0$]
Let $\gamma(z,\vecw)$ be the segment
from $0$ to $z$.
\item[Case $F_{\gminia}(P)=0$]
Let $\theta_+$ be as above.
For the polar coordinate system $z=re^{\sqrt{-1}\theta}$,
let $\gamma(z,\vecw)$ be the union of the ray
$\{\rho e^{\sqrt{-1}\theta_+}\,|\,
 0\leq \rho\leq r
 \}$
and the arc
connecting $r\,e^{\sqrt{-1}\theta_+}$
and $z$.
\end{description}
\begin{lem}
\label{lem;09.10.26.5}
For each $N>0$, there exists 
$C_N>0$ such that 
$\bigl|\Phi(f)(z,\vecw) \bigr|
\leq
 C_N\cdot C\, 
 |z|^N\prod_{i=1}^{\ell}|w_i|^{N_i}$
if $|f(z,\vecw)|\leq 
 C\,|z|^N\prod_{i=1}^{\ell}|w_i|^{N_i}$.
\end{lem}
\pf
We give only an outline.
Let us consider the case $F_{\gminia}(P)<0$.
Let $z_0=r_0e^{\sqrt{-1}\theta_0}$
and $z=r e^{\sqrt{-1}\theta}$.
We may assume that the path $\gamma$
is the union of 
(i) the arc $\gamma_1$ connecting $z_0$ and 
$z_1=r_0e^{\sqrt{-1}\theta}$,
(ii) the segment $\gamma_2$ connecting
$z_1$ and $z$.
The segment $\gamma_2$
is divided into 
$\gamma_{2,1}=\gamma_1\cap
 \bigl\{|\zeta|>3|z|/2\bigr\}$
and 
$\gamma_{2,2}=\gamma_1\cap
 \bigl\{
 |\zeta|\leq 3|z|/2
 \bigr\}$.
The contributions of $\gamma_1$
and $\gamma_{2,1}$ are dominated by
$\bigl|\exp\bigl(-\gminia(z,\vecw)\bigr)\bigr|
 \prod_{i=k+1}^{\ell}|w_i|^{N_i}$.
The function $\Re\gminia$ is monotone on
$\gamma_{2,2}$.
We also have $|f(\zeta,\vecw)|\leq 
 C'\,|z^N|\prod_{i=1}^{\ell}|w_i|^{N_i}$
on $\gamma_{2,2}$.
Hence, the contribution of
$\gamma_{2,2}$ is dominated by
$|z|^{N}\prod_{i=1}^{\ell}|w_i|^{N_i}$.
Let us consider the case $F_{\gminia}(P)\geq 0$.
On $\gamma$,
we have $\bigl|f(\zeta,\vecw)\bigr|
\leq
 C'\,|z^N|\prod_{i=1}^{\ell}|w_i|^{N_i}$,
and $\Re(\gminia)$ is monotone.
Hence, it is easy to obtain the desired estimate.
\hfill\qed

\vspace{.1in}

Let us return to the proof of 
Proposition \ref{prop;09.10.26.2}.
It is enough to consider the case
$X=\Delta^n$ and 
$D=\bigcup_{i=1}^{\ell}\{z_i=0\}$.
We may assume that $V$ is unramifiedly good.
Let $P\in \pi^{-1}(0,\ldots,0)$.
By using the local decomposition around $P$
as in Theorem 4.1 of \cite{mochi10},
we can reduce the issue to the case
$V=\bigoplus_{i=1}^M
 \nbigo_X(\ast D)\,e_i$
with a flat connection
\[
\nabla \vece=\vece\,
\Bigl(d\gminia+
 \sum_{i=1}^{\ell}\bigl(\alpha_i\,I_M+ N_i\bigr)
 \frac{dz_i}{z_i}
 \Bigr),
\]
where 
$I_M$ denotes the identity matrix,
$N_i$ $(i=1,\ldots,\ell)$
are mutually commuting nilpotent matrices,
$\alpha_i$ are complex numbers,
and we put
$\vece:=(e_1,\ldots,e_n)$ and
$\gminia:=\prod_{i=1}^{m}z_i^{-m_i}$.
Then, it is easy to observe that
$\nbigl^{<D_1\leq D_2}$
is naturally isomorphic to
the $0$-th cohomology of
$\DR^{<D_1\leq D_2}_{\Xtilde(D)}(V)[-d_X]$.
Hence, it is enough to show the vanishing
of the higher cohomology of
$\DR^{<D_1\leq D_2}_{\Xtilde(D)}(V)[-d_X]$.
It is enough to consider the case
$\rank V=1$,
and we put $v=e_1$.

First, let us consider the case $D_1=D$.
For a subset $J\subset\{1,\ldots,n\}$,
we set $dz_J=dz_{j_1}\wedge\cdots\wedge dz_{j_k}$.
For a section $\omega$ of
$\Omega^{\bullet\,<D}_{\Xtilde(D)}$,
we have the unique decomposition
$\omega=\sum \omega_J\,dz_J$,
where $\omega_J\in\nbiga^{<D}_{\Xtilde(D)}$.
Let $S_{i}$ $(i=1,\ldots,\ell)$
be a small sector in $\Delta_{z_i}^{\ast}$,
and let $U$ be a small neighbourhood
of $(0,\ldots,0)$ in 
$\prod_{i=\ell+1}^n\Delta_{z_i}$,
such that the closure $\Sbar$ of
$S:=\prod S_{i}\times U$ in $\Xtilde(D)$
is a neighbourhood of $P$.
In the following,
we will shrink $S$ without mention.
It is easy to observe that 
it is enough to consider the case
$\alpha_i=0$ $(i=1,\ldots,\ell)$.

Take $h=1,\ldots,n$.
Assume $\nabla (\omega\, v)=0$
for some section $\omega$
of $\Omega^{\bullet\,<D}_{\Xtilde(D)}$ on $S$
such that 
$\omega_{J}=0$ unless
$J\subset \{1,\ldots,h\}$.
We have
$d\bigl(\exp(\gminia)\omega\bigr)=0$.
For the expression
$\exp(\gminia)\,\omega
=\sum_{h\not\in J} f_J\,dz_h\,dz_J
+\sum_{h\not\in J}f_J\,dz_J$,
we set
$\tau(\vecz)=\sum_{h\not\in J}
 \exp(-\gminia)
 \,\Bigl(
 \int_{\gamma(\vecz)}f_J\, dz_h
\Bigr)\,dz_J$,
where $\gamma(\vecz)$ is 
a path taken as follows:
\begin{itemize}
\item
If $h\leq m$,
the condition is similar to
that for the path in (\ref{eq;09.10.26.3}).
\item
If $m<h$,
$\gamma$ is a path connecting
$(z_1,\ldots,z_{h-1},0,z_{h+1},\ldots,z_n)$
and $(z_1,\ldots,z_n)$.
\end{itemize}
By using Lemma \ref{lem;09.10.26.5},
we obtain that 
$\tau\in
 \Omega^{\bullet\,<D}_{\Xtilde(D)}\otimes V$.
By a formal computation,
we can show that
$\omega\,v-\nabla(\tau\,v)$
does not contain
$dz_{j}$ for $j\geq h$.
Hence, we can show the vanishing
of the higher cohomology of
$\Omega^{\bullet\,<D}_{\Xtilde(D)}\otimes V$
by an induction.

We have the decomposition
$I_1\sqcup I_2=\ellsitabar$ such that
$D_j=\bigcup_{i\in I_j}\{z_i=0\}$.
Let us consider
$\Omega^{\bullet<D(J^c)\leq D(J)}
 _{\widehat{\pi^{-1}(D_J)}}
 \otimes V$
for any subset $J\subset I_2$,
where $J^c:=\ellsitabar\setminus J$.
If $\mbar\cap J\neq\emptyset$,
it is easy to show that
$\Omega^{\bullet<D(J^c)\leq D(J)}
 _{\widehat{\pi^{-1}(D_J)}}
 \otimes V$ is acyclic 
by a formal computation.
Assume $\mbar\cap J=\emptyset$.
Let $V_J=\nbigo_{D_J}(\ast\del D_J)\, v_J$
be equipped with the flat connection
$\nabla v_J=v_J\cdot d\gminia_{|D_J}$ on $D_J$.
Let $q_J$ be the projection
$\pi^{-1}(D_J)\lrarr \Dtilde_J(\del D_J)$.
Then, it is easy to obtain a natural quasi-isomorphism
$q_J^{-1}\bigl(
 \Omega^{\bullet\,<\del D_J}
 _{\Dtilde_J(\del D_J)}\otimes V_J
 \bigr)
\simeq
\Omega^{\bullet\,<D(J^c)\leq D(J)} 
 _{\widehat{\pi^{-1}(D_J)}}
 \otimes V$
by a formal computation.
Hence, we obtain the vanishing of
the higher cohomology of
$\Omega^{\bullet\,<D(J^c)\leq D(J)}
 _{\widehat{\pi^{-1}(D_J)}}
 \otimes V$.

We put $h:=|I_2|$.
Let $\nbigg_h^{\bullet}$
denote the kernel of the surjection
$\Omega_{\Xtilde(D)}^{\bullet\,<D_1\leq D_2}
 \otimes V
\lrarr
 \Omega^{\bullet\,<D_1\leq D_2}
 _{\widehat{\pi^{-1}(D_{I_2})}}
 \otimes V$.
Inductively, let $\nbigg_{k}^{\bullet}$ 
be the kernel of the following surjection:
\[
 \nbigg_{k+1}^{\bullet}
\lrarr
 \bigoplus_{\substack{
 J\subset I_2\\
 |J|=k }}
 \Omega^{\bullet\,<D(J^c)\leq D(J)}
 _{\widehat{\pi^{-1}(D_{J})}}
 \otimes V
\]
Because
$\nbigg_{1}^{\bullet}
=\Omega_{\Xtilde(D)}^{<D}\otimes V$,
we obtain the vanishing of the higher cohomology
by an induction on $k$.
Thus, the proof of Proposition \ref{prop;09.10.26.2}
is finished.
\hfill\qed

\vspace{.1in}
Similarly, we also obtain the following.
(See also \cite{sabbah_lecture_Stokes}.)
\begin{prop}
The natural inclusion
$\nbigl^{\leq D}[d_X]
\lrarr
 \DR^{\moderate}_{\Xtilde(D)}(V)$
is an isomorphism
in $D^b_c(\cnum_{\Xtilde(D)})$.
\hfill\qed
\end{prop}

\subsection{Isomorphisms}

Let $X$ and $D$ be as in the beginning of
\S\ref{subsection;09.10.26.10}.
Let $H$ be hypersurfaces of $X$
contained in $D_1$.
We have the naturally defined projection
$\rho:\Xtilde(D)\lrarr\Xtilde(H)$.

\begin{lem}
\label{lem;13.4.27.1}
For any good meromorphic flat bundle $V$
on $(X,D)$,
the following natural morphisms
are isomorphisms.
\begin{multline}
 R\rho_{\ast}\DR^{<D_1\leq D_2}_{\Xtilde(D)}(V)
 \stackrel{a_1}{\llarr}
 \DR^{<D_1}_{\Xtilde(H)}(V)
 \stackrel{a_2}{\llarr}
 \DR^{<D_1}_{\Xtilde(H)}\bigl(V(!D_1)\bigr)
 \\
 \stackrel{a_3}{\lrarr}
  \DR^{<H}_{\Xtilde(H)}\bigl(V(!D_1)\bigr)
\end{multline}
\end{lem}
\pf
The claim for $a_1$
follows from Theorem \ref{thm;09.12.4.30}.
The claim for $a_2$ is clear.
Let us look at $a_3$.
We use an induction on $\dim X$
and the number of the irreducible components of
$D_1\setminus H$.
We may assume $X=\Delta^n$
and $D=\bigcup_{i=1}^{\ell}\{z_i=0\}$.
We set $L_i:=\{z_i=0\}$.
We may assume
$D_1=\bigcup_{i=1}^{\ell_1}L_i$,
$H=\bigcup_{i=1}^{m_1}L_i$ and
$D_2=\bigcup_{i=\ell_1+1}^{\ell_1+m_1}L_i$.
We set 
$D_3:=\bigcup_{i=2}^{\ell_1}\{z_i=0\}$.
We set
$X':=L_1$ and $D_2':=D_2\cap X'$. 
We set $D_3':=X'\cap D_3$
and $H':=X'\cap\bigcup_{i=2}^{m_1}L_i$.
Let $\iota:X'\lrarr X$ denote the inclusion.
There exist good meromorphic flat bundles
$V_3'$ and $V_3''$
with the following exact sequence:
\[
 0\lrarr
\iota_{\dagger}V_3'(!D_3')
\lrarr
 V(!D_1)
\stackrel{c}{\lrarr}
 V(!D_3)
\lrarr
\iota_{\dagger}V_3''(!D_3')
\lrarr 0
\]
Let $\nbigk$ denote the image of $c$.
We have the following:
{\small
\[
\begin{array}{ccccccc}
 0\lrarr
 &
 \DR^{<D_3}_{\Xtilde(H)}(\iota_{\dagger}V'_3(!D_3'))
 &
\!\lrarr\!
 &
 \DR^{<D_3}_{\Xtilde(H)}(V(!D_1))
 &
 \!\lrarr\!
 &
 \DR^{<D_3}_{\Xtilde(H)}(\nbigk)
 &
 \lrarr 0\\
  & \darr & & \darr & & \darr & \\
0\lrarr
 &
 \DR^{}_{\Xtilde(H)}(\iota_{\dagger}V'_3(!D_3'))
 &
\!\lrarr\!
 &
 \DR^{}_{\Xtilde(H)}(V(!D_1))
 &
\!\lrarr\!
 &
 \DR^{}_{\Xtilde(H)}(\nbigk)
 &
 \lrarr 0
\end{array}
\]
}
{\small
\[
\begin{array}{ccccccc}
 0\lrarr \!&
\DR^{<D_3}_{\Xtilde(H)}(\nbigk)
 &\!\lrarr\!&
 \DR^{<D_3}_{\Xtilde(H)}(V(!D_3))
 &\!\lrarr\!&
  \DR^{<D_3}_{\Xtilde(H)}(\iota_{\dagger}V''_3(!D_3'))
&\!\lrarr 0
 \\
 & \darr & & \darr & & \darr & \\
 0\lrarr \!&
\DR^{}_{\Xtilde(H)}(\nbigk)
 &\!\lrarr\!&
 \DR^{}_{\Xtilde(H)}(V(!D_3))
 &\!\lrarr\!&
  \DR^{}_{\Xtilde(H)}(\iota_{\dagger}V''_3(!D_3'))
&\!\lrarr 0
\end{array}
\]
}
By using the inductive assumption,
we obtain that
\[
 \DR^{<D_3}_{\Xtilde(H)}
 (V(!D_1))
\lrarr
 \DR_{\Xtilde(H)}
 (V(!D_1))
\]
is a quasi-isomorphism.
Because we have
$\DR^{<D_3}_{\Xtilde(H)}
 (V(!D_1))
\simeq
 \DR^{<D_3}_{\Xtilde(H)}
 (V(!L_1))$
and 
$\DR^{<D_1}_{\Xtilde(H)}
 (V(!D_1))
\simeq
 \DR^{<D_1}_{\Xtilde(H)}
 (V(!L_1))$,
it is enough to prove the natural morphism
\begin{equation}
 \label{eq;14.1.16.21}
 \DR^{<D_1}_{\Xtilde(H)}
 (V(!L_1))
\lrarr 
 \DR^{<D_3}_{\Xtilde(H)}
 (V(!L_1))
\end{equation}
is a quasi-isomorphism.

Let $I\subset\{1,\ldots,\ell\}=:\ellsitabar$
be any subset with $1\in I$.
Let $\pi_H:\Xtilde(H)\lrarr X$ denote the projection.
We set $L_{I}:=\bigcap_{i\in I}L_i$
and $\del L_{I}:=L_{I}\cap
 \bigcup_{j\in\ellsitabar\setminus I}L_j$.
\begin{lem}
\label{lem;14.1.16.20}
$\DR^{< \del L_I}_{\widehat{\pi_I^{-1}(L_I)}}(V(!L_1))=0$.
\end{lem}
\pf
By using the pull back and the push-forward
with respect to a ramified covering,
we may assume that $V$ is unramifiedly good.
Let $\nbigi\subset M(X,D)/H(X)$ 
denote the set of irregular values of $V$.
We set $L(I^c):=
 \bigcup_{j\in\ellsitabar\setminus I}L_j$.
Let $\nbigi_I$ denote the image of
$\nbigi$ in $M(X,D)/M(X,L(I^c))$.
For each element of $[\gminia]\in \nbigi_I$,
we fix a representative $\gminia$ in $M(X,D)$.
There exist meromorphic 
$\nbigo_{\Lhat_I}(\ast\del D)$-subbundles
$\Vhat_{[\gminia]}$ of $V_{|\Lhat_I}$
stable by the connection
and a decomposition
$V_{|\Lhat_I}=
 \bigoplus_{[\gminia]\in\nbigi_I}
 \Vhat_{[\gminia]}$
compatible with the connection,
such that
$\nablahat_{\gminia}^{\reg}
:=\nablahat_{\gminia}-d\gminia\id_{\Vhat_{[\gminia]}}$
are regular along $L_i$ $(i\in I)$,
where $\nablahat_{\gminia}$ denotes
the induced connection on $\Vhat_{[\gminia]}$.

Let $j\in I$.
Suppose $\ord_{z_j}\gminia<0$.
We consider the Deligne-Malgrange filtration
$\nbigp_{\ast}$ on $\Vhat_{[\gminia]}$.
(See \cite{Mochizuki-DM} for a survey.)
We have
$(\del_j\gminia)^{-1}\nablahat^{\reg}_{\gminia,\del_j}
 \nbigp_{\vecb}\Vhat_{[\gminia]}\subset
 \nbigp_{\vecb}\Vhat_{[\gminia]}$
for any $\vecb\in\real^{\ell}$.
Hence we obtain that
$\nablahat_{\gminia,\del_j}$ is invertible on
$\nbigc^{\infty<\del L_I}_{\widehat{\pi_H^{-1}(L_I)}}
 \otimes
 \Vhat_{[\gminia]}$.
Suppose moreover that
$j\neq 1$ and that $\ord_{z_1}(\gminia)=0$.
Let $\leq$ denote the total order on $\cnum$
defined by the lexicographic order 
on $(\Re(\alpha),\Image(\alpha))\in\real\times\real$.
We have the $V$-filtration
$\nbigptilde$ of 
$\Vhat_{[\gminia]}$
along $z_1$
indexed by $(\cnum,\leq)$ 
such that 
(i) $z_1\nablahat_{\gminia,\del_1}$ preserves
 the filtration $\nbigptilde$
(ii) the endomorphisms of
 $\Gr^{\nbigptilde}_{\beta}(\Vhat_{[\gminia]})$
 induced by 
 $-\nablahat_{\gminia,\del_1}z_1-\beta$
 are nilpotent for any $\beta$.
The induced morphisms
$\nablahat_{\gminia,\del_1}:
 \Gr^{\nbigptilde}_{\beta}(\Vhat_{[\gminia]})
\lrarr
 \Gr^{\nbigptilde}_{\beta+1}(\Vhat_{[\gminia]})$
are isomorphisms unless $\beta=-1$.
We can observe that
the filtration $\nbigptilde$ is preserved by
$\nablahat_{[\gminia],\del_j}$
and the multiplication of $\del_j\gminia$.
Hence, 
$\nablahat_{\gminia,\del_j}$ is invertible on
$\nbigc^{\infty<\del L_I}_{\widehat{\pi_H^{-1}(L_I)}}
 \otimes
 \nbigp_a\Vhat_{[\gminia]}$
and 
$\nbigc^{\infty<\del L_I}_{\widehat{\pi_H^{-1}(L_I)}}
 \otimes
 \Gr^{\nbigp}_a\Vhat_{[\gminia]}$.

Suppose $\ord_{z_j}\gminia=0$ for any $j\in I$,
i.e., $[\gminia]=[0]$.
For the Deligne-Malgrange filtration $\nbigp_{\ast}$
of $\Vhat_{[0]}$,
we have
$\nabla_{[0],\del_1}\bigl(
 \nbigp_{\vecb}\Vhat_{[0]} 
 (\ast \del L_I)
 \bigr)
\subset
 \nbigp_{\vecb+(1,0,\ldots,0)}
 \Vhat_{[0]}(\ast \del L_I)$.
For the $V$-filtration $\nbigptilde$ along $z_1$,
we obtain that if $\beta<-1$,
the morphism
$\nablahat_{0,\del_1}:
\nbigc_{\widehat{\pi_H^{-1}(L_I)}}^{\infty<\del L_I}
 \otimes
 \nbigptilde_{\beta}(\Vhat_{[0]})
\lrarr
\nbigc_{\widehat{\pi_H^{-1}(L_I)}}^{\infty<\del L_I}
 \otimes
\nbigptilde_{\beta+1}(\Vhat_{[0]})$
is an isomorphism.

We have the decomposition 
$V(!L_1)_{|\Lhat_I}
=\bigoplus_{[\gminia]}
 \widehat{V(!L_1)}_{[\gminia]}$,
compatible with the decomposition
of $V_{|\Lhat_I}$.
If $\ord_{z_1}\gminia<0$,
we have
$\widehat{V(!L_1)}_{[\gminia]}
=\Vhat_{[\gminia]}$.
The action of $\nablahat_{\gminia,\del_1}$
on $\nbigc^{\infty<\del L_I}_{\widehat{\pi_H^{-1}(L_I)}}
 \otimes\widehat{V(!L_1)}_{[\gminia]}$
is invertible.
If $\ord_{z_1}\gminia=0$,
for the $V$-filtration $\nbigptilde$ along $z_1$,
we have
$\nbigptilde_{\beta}\bigl(
 \widehat{V(!L_1)}_{[\gminia]}\bigr)
=\nbigptilde_{\beta}(\Vhat_{[\gminia]})$
for $\beta<0$,
and that 
$\nablahat_{\gminia,\del_1}:
 \Gr^{\nbigptilde}_{\beta}
 (\widehat{V(!L_1)_{[\gminia]}})
\lrarr
\Gr^{\nbigptilde}_{\beta+1}
(\widehat{V(!L_1)_{[\gminia]}})$
are isomorphisms 
for $\beta\geq -1$.
If $[\gminia]\neq [0]$,
take $j\in I$ such that $\ord_{z_j}\gminia<0$,
and then the action of
$\nablahat_{\gminia,\del_j}$
on $\nbigc^{\infty<\del L_I}_{\widehat{\pi_H^{-1}(L_I)}}
 \otimes\widehat{V(!L_1)}_{[\gminia]}$
is invertible.
If $[\gminia]=[0]$,
the action of $\nablahat_{0,\del_1}$
on $\nbigc^{\infty<\del L_I}_{\widehat{\pi_H^{-1}(L_I)}}
 \otimes\widehat{V(!L_1)}_{[\gminia]}$
is invertible.
Then, the claim of Lemma \ref{lem;14.1.16.20}
follows.
\hfill\qed

\vspace{.1in}

Then,
by an easy inductive argument,
we obtain that (\ref{eq;14.1.16.21})
is a quasi-isomorphism,
and the proof of Lemma \ref{lem;13.4.27.1}
is finished.
\hfill\qed

\vspace{.1in}
Suppose 
that we are given a holomorphic function
$G:X\lrarr\cnum^{\ell}$
such that $G^{-1}(D_0)=H$,
where $D_0=\bigcup_{i=1}^{\ell}\{z_i=0\}$.
\begin{lem}
\label{lem;13.4.27.3}
For the naturally defined map
$\rho_1:\Xtilde(D)\lrarr \Xtilde(G)$,
we obtain the following natural isomorphism:
\begin{equation}
\label{eq;13.4.27.2}
 R\rho_{1\ast}\DR^{<D_1\leq D_2}_{\Xtilde(D)}(V)
\simeq
 \DR^{\rapid}_{X,G}\bigl(V(!D_1)\bigr).
\end{equation}
\end{lem}
\pf
It follows from Lemma \ref{lem;13.4.27.1}
and Proposition \ref{prop;13.4.20.231}.
\hfill\qed

\vspace{.1in}
Let $\varphi:X'\lrarr X$ be a projective birational morphism
such that 
(i) $D':=\varphi^{-1}(D)$ is normal crossing,
(ii) $X'\setminus D'\simeq X\setminus D$.
We put $D_1':=\varphi^{-1}(D_1)$
and $H_1':=\varphi^{-1}(H_1)$.
Let $D_2'$ be the complement of $D_1'$
in $D'$.
We set 
$G':=G\circ\varphi$.
We put $V':=\varphi^{\ast}V$.
We have the following natural commutative diagram:
\[
 \begin{CD}
 \Xtilde'(D')  @>{\varphitilde_1}>>
 \Xtilde(D)\\
 @V{\rho_1'}VV @V{\rho_1}VV \\
 \Xtilde'(G') @>{\varphitilde}>>
 \Xtilde(G)
 \end{CD}
\]
We set $\rho_2:=\varphitilde\circ\rho_1'$.
Correspondingly,
we have the following commutative diagram
of isomorphisms
by the construction:
\begin{equation}
 \label{eq;13.4.27.10}
 \begin{CD}
 R\rho_{2\ast}\DR^{<D_1'\leq D_2'}_{\Xtilde'(D')}(V')
 @>>>
 R\rho_{1\ast}\DR^{<D_1\leq D_2}_{\Xtilde(D)}(V)
 \\
 @V{\simeq}VV @V{\simeq}VV \\
 R\varphitilde_{\ast}\DR^{\rapid}_{X',G'}
 \bigl(V'(!D_1')\bigr)
 @>{\simeq}>>
 \DR^{\rapid}_{X,G}\bigl(V(!D_1)\bigr)
 \end{CD}
\end{equation}
The lower horizontal arrow is an isomorphism
according to Corollary \ref{cor;13.4.20.400}.

%% file: 5.2.tex
\subsection{Duality morphisms}

Let $X$, $D$ and $\nbigm$ be as in 
\S\ref{subsection;09.10.26.10}.
We have the following natural morphism
given in a way parallel to that of (\ref{eq;13.4.17.20}):
\begin{equation}
\label{eq;09.10.25.120}
 \DR_{\Xtilde(D)}^{<D_1\leq D_2}\bigl(
 \DDD\nbigm\bigr)
\lrarr
 \DDD \DR^{<D_2\leq D_1}_{\Xtilde(D)}
 (\nbigm)
\end{equation}
Namely,
we take a $\pi^{-1}(\nbigd_X)$-injective resolution
$\nbigitilde_1^{\bullet}$
of $\Omega^{0,\bullet<D_1\leq D_2}_{\Xtilde(D)}[d_X]$,
and a $\cnum_{\Xtilde(D)}$-injective resolution 
$\nbigitilde_2^{\bullet}$ of
$\Tot\Omega^{\bullet,\bullet<D}_{\Xtilde(D)}[2d_X]$
with a morphism
$\DR^{\leq D_1<D_2}_{\Xtilde(D)}\nbigitilde_1^{\bullet}
\lrarr
 \nbigitilde_2^{\bullet}$
extending a natural morphism
\[
 \DR^{\leq D_1<D_2}_{\Xtilde(D)}
 \bigl(
 \Omega^{0,\bullet<D_1\leq D_2}_{\Xtilde(D)}[d_X]
 \bigr)
\lrarr
 \Tot\Omega^{\bullet,\bullet<D}_{\Xtilde(D)}[2d_X].
\]
Then, (\ref{eq;09.10.25.120}) is given as the composite of
the following morphisms:
\begin{multline}
 \nhom_{\pi^{-1}(\nbigd_X)}
 \bigl(
 \pi^{-1}\nbigm,\nbigitilde_1^{\bullet}
 \bigr)
\lrarr
 \nhom_{\cnum_{\Xtilde(D)}}
 \bigl(
 \DR^{<D_2\leq D_1}_{\Xtilde(D)}\nbigm,
 \DR^{<D_2\leq D_1}_{\Xtilde(D)}\nbigitilde_1
 \bigr)
 \\
\lrarr
  \nhom_{\cnum_{\Xtilde(D)}}
 \bigl(
 \DR^{<D_2\leq D_1}_{\Xtilde(D)}\nbigm,
 \nbigitilde_2
 \bigr)
\end{multline}

\begin{prop}
\label{prop;09.10.27.1}
The following diagram is commutative:
\begin{equation}
 \label{eq;09.10.3.50}
 \begin{CD}
 R\pi_{\ast}\DR_{\Xtilde(D)}
 ^{<D_1\leq D_2}\bigl(
 \DDD\nbigm\bigr)
 @>>>
 R\pi_{\ast}\DDD
 \DR_{\Xtilde(D)}^{<D_2\leq D_1}(\nbigm)\\
 @V{\simeq}VV @V{\simeq}VV \\
 \DR_X^{<D_1\leq D_2}(\DDD\nbigm)
 @>>>
 \DDD\DR_X^{<D_2\leq D_1}(\nbigm)
 \end{CD}
\end{equation}
Here, 
the upper horizontal arrow is induced by
{\rm(\ref{eq;09.10.25.120})},
the lower horizontal arrow is
given as in {\rm(\ref{eq;13.4.17.20})},
the left vertical arrow is given in
{\rm(\ref{eq;09.10.4.1})},
and the right vertical arrow is given by
$R\pi_{\ast}\DDD 
 \DR^{<D_2\leq D_1}_{\Xtilde(D)}\nbigm
\simeq
 \DDD R\pi_{\ast}\DR^{<D_2\leq D_1}_{\Xtilde(D)}
 \nbigm
\simeq
 \DDD \DR^{<D_2\leq D_1}_X(\nbigm)$.
\end{prop}
\pf
We have a morphism
$R\pi_{\ast}\DR_{\Xtilde(D)}^{<D_1\leq D_2}
 (\DDD\nbigm)
\lrarr
 \DR_X^{<D_1\leq D_2}(\DDD \nbigm) $
given as follows,
by Lemma \ref{lem;09.10.3.40}:
{\small
\begin{multline}
\label{eq;09.10.4.2}
 R\pi_{\ast}
 \nrhom_{\pi^{-1}\nbigd_X}\bigl(
 \pi^{-1}\nbigm,\,
 \Omega_{\Xtilde(D)}^{0,\bullet\,<D_1\leq D_2}
 \bigr)[d_X]
\simeq
 \nrhom_{\nbigd_X}\bigl(
 \nbigm,\,
 R\pi_{\ast}\Omega_{\Xtilde(D)}^{0,\bullet\,<D_1\leq D_2}
 \bigr)[d_X] \\
\simeq
  \nrhom_{\nbigd_X}\bigl(
 \nbigm,\,
  \Omega_X^{0,\bullet}(\ast D_2)^{<D_1}
 \bigr)[d_X]
\end{multline}
}
It is equal to the morphism
obtained as in {\rm(\ref{eq;09.10.4.1})}.
Then, the claim of the proposition can be checked easily.
\hfill\qed

\subsection{The case of good meromorphic
flat bundles}

Let us consider the case where
$\nbigm$ is a good meromorphic flat bundle $V$
on $(X,D)$.

\begin{thm}
\label{thm;09.10.26.11}
The duality morphism
$\DR_{\Xtilde(D)}^{<D_1\leq D_2}\DDD V
\lrarr 
 \DDD\DR^{<D_2\leq D_1}_{\Xtilde(D)}V$
is an isomorphism.
\end{thm}
\pf
We begin with elementary preparations.
Let $\real^2=S_0\cup S_1\cup S_2$
be a decomposition given as follows:
\[
 S_0:=\bigl\{
 (x,y)\,\big|\,y\geq 0
 \bigr\}
\quad
 S_1:=\bigl\{
 (x,y)\,\big|\,y\leq 0,\,x\leq 0
 \bigr\}
\quad
 S_2:=\big\{
 (x,y)\,\big|\,y\leq 0,\,x\geq 0
 \}
\]
We put 
$X_1:=(\real\times S_1)\cup (\real_{\geq 0}\times S_0)$
and
$X_2:=(\real\times S_2)\cup (\real_{\leq 0}\times S_0)$.
The following lemma is easy to see.
\begin{lem}
\label{lem;09.10.26.10}
$X_i\subset\real^3$ $(i=1,2)$ are 
closed $C^0$-submanifolds with boundaries.
We have
$X_1\cup X_2=\real^3$
and $X_1\cap X_2=\del X_i$.
\hfill\qed
\end{lem}

We put $\nbigj:=\openopen{-1}{1}$,
$\nbigj_+:=\closedopen{0}{1}$,
$\nbigj_-:=\openclosed{-1}{0}$,
and $\nbigi_i:=\closedopen{0}{1}$ $(i=1,2,3)$.
We have a homeomorphism
$\del(\nbigi_1\times\nbigi_2\times\nbigi_3)
\simeq \real^2$,
and we can identify
the decomposition
\[
\del(\nbigi_1\times\nbigi_2\times\nbigi_3)=
\bigl(
 \del\nbigi_1\times\nbigi_2\times\nbigi_3
\bigr)
\cup
\bigl(
 \nbigi_1\times\del\nbigi_2\times\nbigi_3
\bigr)
\cup
\bigl(
\nbigi_1\times\nbigi_2\times\del\nbigi_3
\bigr)
\]
with $\real^2=S_0\cup S_1\cup S_2$.
We put
\[
 X_1':=
 \Bigl(
 \nbigj\times\nbigi_1\times\del\nbigi_2\times\nbigi_3
\Bigr)
\cup
 \Bigl(
 \nbigj_+\times\del\nbigi_1\times\nbigi_2\times\nbigi_3
 \Bigr)
\]
\[
 X_2':=\Bigl(
 \nbigj\times\nbigi_1\times\nbigi_2\times\del\nbigi_3
 \Bigr)
\cup
 \Bigl(
 \nbigj_-\times\del\nbigi_1\times\nbigi_2\times\nbigi_3
 \Bigr)
\]
They are closed subsets of
$\nbigj\times \del(\nbigi_1\times\nbigi_2\times\nbigi_3)$.
We obtain the following lemma
from Lemma \ref{lem;09.10.26.10}.
\begin{lem}
\label{lem;09.10.26.21}
$X_i'\subset 
\nbigj\times \del(\nbigi_1\times\nbigi_2\times\nbigi_3)$
are $C^0$-submanifolds with boundaries.
We have $X_1'\cup X_2'=
 \nbigj\times
 \del\bigl(\nbigi_1\times\nbigi_2\times\nbigi_3\bigr)$
and $X_1'\cap X_2'=\del X_i'$.
\hfill\qed
\end{lem}

We recall some elementary facts
on constructible sheaves.
Let $Y$ be an oriented $\ell$-dimensional 
$C^0$-manifold with the boundary $\del Y$.
For a closed $C^0$-submanifold $W\subset \del Y$
with boundary such that $\dim W=\ell-1$,
let $j_{W}$ denote the inclusion $Y-W\lrarr Y$.
We have the following natural isomorphisms:
\[
 \nrhom_{\cnum_Y}\bigl(j_{W!}\cnum_{Y-W},K \bigr)
\simeq 
 Rj_{W\ast}\nrhom_{\cnum_{Y-W}}
 (\cnum_{Y-W},Rj_{W}^!K)
\simeq
 Rj_{W\ast}j_W^{\ast}K
\]
The dualizing complex $\omega_Y$ of $Y$ is given by
$j_{\del Y!}\cnum_{Y-\del Y}[\ell]$.

\begin{lem}
\label{lem;09.10.26.30}
Let $Y_i\subset\del Y$ be
closed $C^0$-submanifolds
with boundaries such that
$Y_1\cup Y_2=Y$ and
$Y_1\cap Y_2=\del Y_i$.
Then, we have
$\DDD j_{Y_1!}\cnum_{Y-Y_1}
\simeq
 j_{Y_2!}\cnum_{Y-Y_2}$.
\end{lem}
\pf
The left hand side is naturally isomorphic to
$j_{Y_1\ast}j_{Y_1}^{\ast}
 \omega_Y
\simeq
 j_{Y_1\ast}
 j_{0!}
 \cnum_{Y-\del Y}[\ell]$,
where $j_0$ denotes the inclusion
$Y-\del Y\lrarr Y-Y_1$.
Then, we can check the claim directly.
\hfill\qed

\vspace{.1in}

Let us return to the proof of Theorem
\ref{thm;09.10.26.11}.
It is enough to consider the case
$X=\Delta^n$ and
$D=\bigcup_{i=1}^{\ell}\{z_i=0\}$.
As in the proof of Proposition \ref{prop;09.10.26.2},
we can reduce the issue to the case where
$V=\nbigo_X(\ast D)\,v$ with a meromorphic
flat connection $\nabla v=v\,d\gminia$,
where $\gminia=\prod_{i=1}^mz_i^{-m_i}$
($m_i>0$).
We put $F_{\gminia}:=-|\gminia|^{-1}\Re\gminia$.
We have the decomposition
$I_1\sqcup I_2=\ellsitabar$ such that
$D_j=\bigcup_{i\in I_j}\{z_i=0\}$ $(j=1,2)$.
We set
$I_j(>m):=\bigl\{
 i\in I_j\,\big|\,i>m \bigr\}$.
We also put
$D(>m):=\bigcup_{i=m+1}^{\ell}\{z_i=0\}$
and $D(\leq m):=\bigcup_{i=1}^m\{z_i=0\}$.
We consider the closed subsets $W_i\subset\pi^{-1}(D)$
$(i=1,2)$ given as follows:
\[
 W_1:=\pi^{-1}\Bigl(
 D_1\cap D(>m)
 \Bigr)
\cup
 \Bigl[
 \pi^{-1}\bigl(D(\leq m)\bigr)
\cap\{F_{\gminia}\geq 0\}
 \Bigr]
\]
\[
 W_2:=\pi^{-1}\Bigl(
 D_2\cap D(>m)
 \Bigr)
\cup
 \Bigl[
 \pi^{-1}\bigl(D(\leq m)\bigr)
\cap\{F_{\gminia}\leq 0\}
 \Bigr]
\]
\begin{lem}
\label{lem;09.10.26.20}
$W_i\subset \pi^{-1}(D)$ are 
closed $C^0$-submanifolds
with boundaries,
and we have
$W_1\cup W_2=\pi^{-1}(D)$
and $W_1\cap W_2=\del W_i$.
\end{lem}
\pf
It is easy to observe that it is enough to consider 
the case $n=\ell$.
We have the natural identification
$\Xtilde(D)\simeq
 (S^1)^{\ell}\times\real_{\geq\,0}^{\ell}$.
By the decomposition
$\ellsitabar=\mbar\sqcup I_1(>m)\sqcup I_2(>m)$,
we identify
$\real_{\geq 0}^{\ell}
=\real_{\geq 0}^m\times
 \real_{\geq 0}^{I_1(>m)}\times
 \real_{\geq 0}^{I_2(>m)}$.
We argue the case $I_j(>m)\neq\emptyset$ $(j=1,2)$.
The other cases are easier.
We fix homeomorphisms
\[
\real^m_{\geq 0}\simeq
 \nbigi_1\times\real^{m-1},
\quad 
\real^{I_1(>m)}_{\geq 0}\simeq
 \nbigi_2\times\real^{|I_1(>m)|-1},
\quad
\real^{I_2(>m)} _{\geq 0}\simeq
 \nbigi_3\times\real^{|I_2(>m)|-1}.
\]
We put $N:=m+|I_1(>m)|+|I_2(>m)|-3$.
Let $H_{\pm}$ be the subsets of $(S^1)^{\ell}$
given as follows:
\[
 H_+:=\Bigl\{
 \cos\Bigl(\sum m_i\theta_i\Bigr)\geq 0
 \Bigr\}
\quad\quad
 H_-:=\Bigl\{
 \cos\Bigl(\sum m_i\theta_i\Bigr)\leq 0
 \Bigr\}
\]
Then, 
$\pi^{-1}(D)$ is identified with
$(S^1)^{\ell}\times
 \del(\nbigi_1\times\nbigi_2\times\nbigi_3)
 \times\real^{N}$, under which we have
\[
 W_1\simeq \Bigl(
 \bigl(
 (S^1)^{\ell}\times
 \nbigi_1\times\del\nbigi_2\times\nbigi_3
 \bigr)
\cup
 \bigl(
 H_-\times\del\nbigi_1
 \times\nbigi_2\times\nbigi_3
 \bigr)
 \Bigr)
 \times\real^N
\]
\[
 W_2=\Bigl(\bigl(
 (S^1)^{\ell}\times
 \nbigi_1\times\nbigi_2\times
 \del\nbigi_3\bigr)
\cup
 \bigl(
 H_+\times\del \nbigi_1\times
 \nbigi_2\times\nbigi_3\bigr)
 \Bigr)
 \times\real^{N}
\]
For a point $Q\in H_+\cap H_-$,
we can take a neighbourhood $U_Q$
such that
$U\simeq \nbigj\times\real^{\ell-1}$
under which
$H_{\pm}\cap U_Q=
 \nbigj_{\pm}\times\real^{\ell-1}$.
Then, we obtain Lemma \ref{lem;09.10.26.20}
from Lemma \ref{lem;09.10.26.21}.
\hfill\qed

\vspace{.1in}

Let $j_{W_i}$ be the inclusion
$\Xtilde(D)\setminus W_i\lrarr\Xtilde(D)$.
Let $\nbigl$ and $\nbigl^{\lor}$
be the local systems on $\Xtilde(D)$
associated to $V$ and $V^{\lor}$,
respectively.
According to the description of
$\nbigl^{<D_1\leq D_2}$
and
$\nbigl^{\lor<D_2\leq D_1}$,
we have the following natural
isomorphisms:
\[
 \nbigl^{<D_1\leq D_2}\simeq
 j_{W_1!}\bigl(
 \nbigl_{\Xtilde(D)\setminus W_1}
 \bigr)
\quad\quad\quad
 \nbigl^{\lor<D_2\leq D_1}\simeq
 j_{W_2!}\bigl(
 \nbigl^{\lor}_{\Xtilde(D)\setminus W_2}
 \bigr)
\]
Lemma \ref{lem;09.10.26.30} gives
an isomorphism
$\DDD\bigl(
 \nbigl^{<D_1\leq D_2}[d_X]
 \bigr)
\simeq
 \nbigl^{\lor\,<D_2\leq D_1}[d_X]$.
It is uniquely determined by
its restriction to $X-D$.
Then, we can deduce that
$\DR_{\Xtilde(D)}^{<D_1\leq D_2}\DDD V
\lrarr 
 \DDD\DR^{<D_2\leq D_1}_{\Xtilde(D)}V$
is an isomorphism.
Thus, the proof of Theorem \ref{thm;09.10.26.11}
is finished.
\hfill\qed

\begin{cor}
For any good meromorphic flat bundle $V$
on $(X,D)$,
we have the following commutative diagram
of the isomorphisms:
\[
 \begin{CD}
 R\pi_{\ast}\DR_{\Xtilde(D)}^{<D_1\leq D_2}\DDD V
@>{\simeq}>>
 R\pi_{\ast}\DDD 
 \DR_{\Xtilde(D)}^{<D_2\leq D_1}V \\
 @V{\simeq}VV @V{\simeq}VV \\
 \DR_X V^{\lor}(!D_1)
 @>{\simeq}>>
 \DDD\DR_X V(!D_2)
 \end{CD}
\]
\end{cor}
\pf
It follows from
Theorem \ref{thm;09.10.3.50},
Proposition \ref{prop;09.10.27.1} and
Theorem \ref{thm;09.10.26.11}.
\hfill\qed

%% file: 5.3.tex
Let $X$ be a complex manifold,
and let $D$ be a normal crossing hypersurface
with a decomposition $D=D_1\cup D_2$.
Let $D_3$ be a hypersurface of $X$.
Let $\varphi:X'\lrarr X$ be a proper birational morphism
such that 
(i) $D':=\varphi^{-1}(D\cup D_3)$
is normal crossing,
(ii) $X'\setminus D'\simeq X\setminus (D\cup D_3)$.
Let $\Xtilde(D)\lrarr X$ and $\Xtilde'(D')\lrarr X'$
be the real blow up.
Both the projections are denoted by $\pi$.
Let $\varphitilde:\Xtilde'(D')\lrarr\Xtilde(D)$
be the induced map.
We put
$D_1':=\varphi^{-1}(D_1)$.
We have $D_2'\subset D'$ such that
$D'=D_1'\cup D_2'$ is a decomposition.
Let $V$ be a meromorphic flat bundle on $(X,D)$.
We set
$V':=\varphi^{\ast}(V)
 \otimes\nbigo_{X'}(\ast D')$.

\begin{thm}
\label{thm;09.10.4.100}
We have a morphism 
$\DR^{<D_1\leq D_2}_{\Xtilde(D)}(V)
\lrarr
 R\varphitilde_{\ast}
 \DR^{<D_1'\leq D_2'}_{\Xtilde'(D')}(V')$
in $D^b_c(\cnum_{\Xtilde(D)})$
such that the following diagram 
of perverse sheaves is commutative:
\begin{equation}
 \label{eq;09.10.27.10}
 \begin{CD}
 R\pi_{\ast}
 \DR_{\Xtilde(D)}^{<D_1\leq D_2}(V)
 @>>>
 R\pi_{\ast}
 R\varphitilde_{\ast}
 \DR^{<D_1'\leq D_2'}_{\Xtilde(D')}(V') \\
 @V{\simeq}VV @V{\simeq}VV \\
 \DR_X\bigl(V(!D_1)\bigr)
 @>>>
 R\varphi_{\ast}\DR_{X'}\bigl(V'(!D_1')\bigr)
 \end{CD}
\end{equation}
Here, the vertical isomorphisms are given
by {\rm(\ref{eq;09.10.4.1})}
and {\rm(\ref{eq;13.4.17.2})},
and the lower horizontal arrow is
induced by the morphism of $\nbigd$-modules
$V(!D_1)\lrarr \varphi_{\dagger}V'(!D_1')$.

Similarly,
we have a morphism
$R\varphitilde_{\ast}
 \DR^{<D_2'\leq D_1'}_{\Xtilde'(D')}(V')
\lrarr
 \DR^{<D_2\leq D_1}_{\Xtilde(D)}(V)$
such that the following diagram 
of perverse sheaves is commutative:
\begin{equation}
 \label{eq;09.10.27.11}
 \begin{CD}
 R\pi_{\ast}R\varphitilde_{\ast}
  \DR^{<D_2'\leq D_1'}_{\Xtilde'(D')}(V')
 @>>>
 R\pi_{\ast}\DR^{<D_2\leq D_1}_{\Xtilde(D)}(V) \\
 @V{\simeq}VV @V{\simeq}VV \\
 R\varphi_{\ast}\DR_X(V'(!D_2'))
 @>>>
 \DR_{X}\bigl(V(!D_2)\bigr)
 \end{CD}
\end{equation}
\end{thm}
\pf
We have a naturally induced morphism:
\begin{equation}
 \label{eq;09.10.27.5}
 \varphitilde^{-1}
 \bigl(
 \Omega_{\Xtilde(D)}^{\bullet,\bullet\,<D_1\leq D_2}
\otimes \pi^{-1}V
 \bigr)
\lrarr
 \Omega_{\Xtilde'(D')}^{\bullet,\bullet\,<D_1'\leq D_2'}
\otimes \pi^{-1}V'.
\end{equation}
It induces a morphism
of cohomologically constructible complexes:
\begin{equation}
\label{eq;09.10.3.56}
 \DR_{\Xtilde(D)}^{<D_1\leq D_2}(V)
\lrarr
 \varphitilde_{\ast}
 \DR_{\Xtilde'(D')}^{<D_1'\leq D_2'}(V')
\end{equation}
We can directly check the commutativity
of the following diagram:
\[
 \begin{CD}
 \Omega_X^{\bullet,\bullet\,<D_1\leq D_2}
 \otimes V @>>>
 \varphi_{\ast}\Bigl(
 \Omega_{X'}^{\bullet,\bullet\,<D_1'\leq D_2'}
 \otimes V'
 \Bigr) \\
 @VVV @VVV \\
 \pi_{\ast}\Bigl(
 \Omega_{\Xtilde(D)}^{\bullet,\bullet\,<D_1\leq D_2}
 \otimes \pi^{-1}V
 \Bigr)
 @>>>
 \pi_{\ast}\Bigl(
 \varphitilde_{\ast}
 \Omega_{\Xtilde'(D')}^{\bullet,\bullet\,
 <D_1'\leq D_2'}\otimes\pi^{-1}V'
 \Bigr) 
 \end{CD}
\]
It implies the commutativity of
the following diagram:
\begin{equation} 
 \label{eq;09.10.3.210}
 \begin{CD}
 R\pi_{\ast}\DR_{\Xtilde(D)}^{<D_1\leq D_2}(V)
@>>>
 R\pi_{\ast}R\varphitilde_{\ast}
 \DR_{\Xtilde'(D')}^{<D_1'\leq D_2'}(V')\\
 @V{\simeq}VV @V{\simeq}VV \\
 \DR_X^{<D_1\leq D_2}(V)
 @>>>
 R\varphi_{\ast}
 \DR_{X'}^{<D_1'\leq D_2'}(V')
 \end{CD}
\end{equation}
Then, we obtain 
the commutativity of (\ref{eq;09.10.27.10})
from Theorem \ref{thm;09.10.3.55}.

Considering the dual of (\ref{eq;09.10.3.56})
with $V^{\lor}$
(see Theorem \ref{thm;09.10.26.11}),
we obtain the following morphism:
\begin{equation}
\label{eq;09.10.3.215}
 R\varphitilde_{\ast}
 \DR_{\Xtilde'(D')}^{\leq D_1'<D_3'}(V')
\lrarr
 \DR_{\Xtilde(D)}^{\leq D_1<D_2}(V)
\end{equation}
Let us prove the commutativity of
the diagram (\ref{eq;09.10.27.11}).
From (\ref{eq;09.10.3.210}) for $V^{\lor}$,
we obtain the following commutative diagram:
\[
 \begin{CD}
 \DDD R\pi_{\ast}R\varphitilde_{\ast}
 \DR_{\Xtilde'(D')}^{<D_1'\leq D_2'}
 (V^{\prime\lor})
@>>>
 \DDD R\pi_{\ast}
 \DR_{\Xtilde(D)}^{<D_1\leq D_2}(V^{\lor})
\\
 @V{\simeq}VV @V{\simeq}VV \\
 \DDD R\varphi_{\ast}
 \DR_{X'}^{<D_1'\leq D_2'}(V^{\prime\lor})
@>>>
 \DDD \DR_X^{<D_1\leq D_2}(V^{\lor})
 \end{CD}
\]
By Proposition \ref{prop;09.10.27.1}
and Theorem \ref{thm;09.10.26.11},
we have the following commutative diagram:
\[
\begin{CD}
 \DDD R\pi_{\ast}
 \DR_{\Xtilde(D)}^{<D_1\leq D_2}(V^{\lor})
@>{\simeq}>>
 R\pi_{\ast}
 \DR_{\Xtilde(D)}^{\leq D_1<D_2}(V) \\
 @V{\simeq}VV @V{\simeq}VV \\
 \DDD \DR^{<D_1\leq D_2}(V^{\lor})
@>{\simeq}>>
 \DR_X^{<D_2\leq D_1}(V) 
\end{CD}
\]
We have a similar diagram for $V'$.
Then, we obtain the commutativity of
(\ref{eq;09.10.27.11})
from the constructions of 
(\ref{eq;09.10.3.215}) and
(\ref{eq;09.10.3.220}).
\hfill\qed

%% file: 5.4.tex
The author originally used
Theorem \ref{thm;09.10.17.1} below
for the functoriality of the Betti structure
by projective morphisms.
After the improvement,
it is now not necessary.
But, it seems interesting to the author,
so we keep it.
The reader can skip this subsection.

\subsection{Statement}
\label{subsection;09.10.27.40}

We set $X:=\Delta^n$ and
$D:=\bigcup_{i=1}^{\ell}\{z_i=0\}$.
Let $V$ be a good meromorphic
flat bundle on $(X,D)$.
Let $\nbigl$ be the associated local system
on $\Xtilde(D)$.
Let $g$ be a holomorphic function on $X$
such that $g^{-1}(0)=D$.
We have the naturally defined morphisms:
\[
 \begin{CD}
 \Xtilde(D) @>{\pi_1}>>
 \Xtilde(g) @>{\pi_0}>>X
 \end{CD}
\]
We put
$\pi_2:=\pi_0\circ\pi_1$.
We set
$\nbigk:=R\pi_{1\ast}\nbigl^{\leq D}$.
In this subsection,
we will work on the derived category
of cohomologically constructible sheaves.

\begin{thm}
\label{thm;09.10.17.1}
The restriction
$\Hom(\nbigk,\nbigk)
\lrarr
 \Hom\bigl( 
 \nbigk_{|\pi_0^{-1}(X-D)},\,
 \nbigk_{|\pi_0^{-1}(X-D)}
\bigr)$
is injective.
\end{thm}
We will give a consequence in
\S\ref{subsection;09.12.5.33}.

\subsection{Reduction}

We put $D^{[m]}:=\bigcup_{
 \substack{
 I\subset\ellsitabar\\ |I|=m}}
 D_I$.
It is easy to see that
\[
 \Hom\bigl(
 \nbigk_{|\pi_0^{-1}(X-D^{[2]})},\,
 \nbigk_{|\pi_0^{-1}(X-D^{[2]})}
 \bigr)
\lrarr
 \Hom\bigl(
 \nbigk_{|\pi_0^{-1}(X-D)},\,
 \nbigk_{|\pi_0^{-1}(X-D)}
 \bigr)
\]
is injective.
Hence, it is enough to show the injectivity of
the following morphisms for $m\geq 2$:
\[
 \Hom\bigl(
 \nbigk_{|\pi_0^{-1}(X-D^{[m+1]})},\,
 \nbigk_{|\pi_0^{-1}(X-D^{[m+1]})}
 \bigr)
\lrarr
 \Hom\bigl(
 \nbigk_{|\pi_0^{-1}(X-D^{[m]})},\,
 \nbigk_{|\pi_0^{-1}(X-D^{[m]})}
 \bigr)
\]
Then, it is easy to observe that
it is enough to consider the case
$\ell=n$ and the following morphism:
\[
 \Hom(\nbigk,\nbigk)
\lrarr
 \Hom\bigl(
 \nbigk_{|\pi_0^{-1}(X-O)},\,
 \nbigk_{|\pi_0^{-1}(X-O)}
 \bigr)
\]
By the adjunction
$\Hom\bigl(
 \pi_1^{\ast}\nbigk,\nbigl^{\leq D}
 \bigr)
\simeq
 \Hom(\nbigk,\nbigk)$,
it is enough to show the injectivity of
the following morphism:
\[
 \Hom\bigl(
 \pi_1^{\ast}\nbigk,\,
 \nbigl^{\leq D}
 \bigr)
\lrarr
 \Hom\bigl(
 \pi_1^{\ast}\nbigk_{|\pi_2^{-1}(X-O)},\,
  \nbigl^{\leq D}_{|\pi_2^{-1}(X-O)}
 \bigr)
\]
We have $R^i\pi_{1\ast}\nbigl^{\leq D}=0$
unless $0\leq i\leq n-1$,
because the real dimension of the fiber is
less than $n-1$.
We set
\[
 \nbigk^i:=\pi_1^{\ast}
 R^i\pi_{1\ast}\nbigl^{\leq D}.
\]
Let $\gminij:\pi_2^{-1}(X-O)\lrarr \Xtilde(D)$
and $\gminii:\pi_2^{-1}(O)\lrarr\Xtilde(D)$.

\begin{lem}
\label{lem;09.10.17.3}
To prove Theorem {\rm\ref{thm;09.10.17.1}},
it is enough to prove
\begin{equation}
\label{eq;09.10.17.2}
 \next^j(\gminii_{\ast}
 \gminii^{\ast}\nbigk^i,\nbigl^{\leq D})
=0,
\quad (i,j\leq n-1).
\end{equation}
\end{lem}
\pf
From the distinguished triangle
$\nbigk^i[-i]\lrarr
 \tau_{\geq i}\pi_1^{\ast}\nbigk\lrarr
 \tau_{\geq i+1}\pi_1^{\ast}\nbigk
 \stackrel{+1}{\lrarr}$,
we obtain the long exact sequence:
\begin{multline}
 \Ext^{i-1}(\nbigk^i,\nbigl^{\leq D})
\lrarr
 \Hom\bigl(
 \tau_{\geq i+1}\pi_1^{\ast}\nbigk,\nbigl^{\leq D}
 \bigr)
\lrarr
 \Hom\bigl(\tau_{\geq i}\pi_1^{\ast}\nbigk,
 \nbigl^{\leq D}\bigr)
 \\
\lrarr
 \Ext^i\bigl(\nbigk^i,\nbigl^{\leq D}\bigr)
\end{multline}
We have the corresponding long exact sequences
for the restrictions to $\pi_2^{-1}(X-O)$.
The injectivity of
$\Hom\bigl(
 \tau_{\geq \,i}\pi_1^{\ast}\nbigk,\,
 \nbigl^{\leq D}\bigr)
\lrarr
 \Hom\bigl(
 \tau_{\geq i}\pi_1^{\ast}\nbigk_{|
 \pi_2^{-1}(X-O)},
 \nbigl^{\leq D}_{|\pi_2^{-1}(X-O)}
 \bigr)$ can follow from the injectivity of
\begin{equation}
\label{eq;09.10.17.5}
 \Ext^i\bigl(
 \nbigk^i,\nbigl^{\leq D}
 \bigr)
\lrarr
 \Ext^i\bigl(
 \nbigk^i_{|\pi_2^{-1}(X-O)},\,
 \nbigl^{\leq D}_{|\pi_2^{-1}(X-O)}
 \bigr),
\end{equation}
\begin{equation}
 \Hom\bigl(
 \tau_{\geq \,i+1}\pi_1^{\ast}\nbigk,\,
 \nbigl^{\leq D}\bigr)
\lrarr
 \Hom\bigl(
 \tau_{\geq i+1}\pi_1^{\ast}\nbigk_{|
 \pi_2^{-1}(X-O)},
 \nbigl^{\leq D}_{|\pi_2^{-1}(X-O)}
 \bigr),
\end{equation}
and the surjectivity of
\begin{equation}
\label{eq;09.10.17.6}
 \Ext^{i-1}(\nbigk^i,\nbigl^{\leq D})
\lrarr
 \Ext^{i-1}(\nbigk^i_{|\pi_2^{-1}(X-O)},
 \nbigl^{\leq D}_{|\pi_2^{-1}(X-O)}).
\end{equation}
By an easy inductive argument,
we can reduce Theorem \ref{thm;09.10.17.1}
to the injectivity of (\ref{eq;09.10.17.5})
and the surjectivity of (\ref{eq;09.10.17.6})
for any $i\leq n-1$.

From the exact sequence
$0\lrarr \gminij_!\gminij^{\ast}\nbigk^i
 \lrarr\nbigk^i\lrarr 
 \gminii_{\ast}\gminii^{\ast}\nbigk^i
 \lrarr 0$
and the adjunction 
$\Ext^i\bigl(
 \gminij_!\gminij^{\ast}\nbigk^i,\nbigl^{\leq D}
 \bigr)
\simeq
 \Ext^i\bigl(\gminij^{\ast}\nbigk^i,\,
 \gminij^{\ast}\nbigl^{\leq D}\bigr)$,
we obtain the following exact sequence:
\begin{multline}
 \Ext^{i-1}\bigl(\nbigk^i,\,\nbigl^{\leq D}\bigr)
\lrarr
 \Ext^{i-1}\bigl(\gminij^{\ast}\nbigk^i,\,
 \gminij^{\ast}\nbigl^{\leq D}\bigr)
\lrarr
 \Ext^i\bigl(\gminii_{\ast}\gminii^{\ast}\nbigk^i,
 \,\nbigl^{\leq D}\bigr)
 \\
\lrarr
 \Ext^i\bigl(\nbigk^i,\nbigl^{\leq D}\bigr)
\lrarr
 \Ext^i\bigl(\gminij^{\ast}\nbigk^i,
 \,\gminij^{\ast}\nbigl^{\leq D}\bigr)
\end{multline}
Hence, the proof of Theorem \ref{thm;09.10.17.1}
is reduced to the vanishing
$\Ext^i\bigl(
 \gminii_{\ast}\gminii^{\ast}\nbigk^i,
 \,\nbigl^{\leq D}
 \bigr)=0$
for any $0\leq i\leq n-1$.
For that purpose,
it is enough to prove
(\ref{eq;09.10.17.2}).
Thus, the proof of Lemma \ref{lem;09.10.17.3}
is finished.
\hfill\qed

\vspace{.1in}

In the following, 
we will prove
$\next^i\bigl(
 \pi_1^{-1}(I),\nbigl^{\leq D}
 \bigr)=0$
($i=0,\ldots,n-1$)
for any constructible sheaf $I$ on
$\pi_0^{-1}(O)\simeq S^1$.

\subsection{Local form of $\pi_1^{-1}(I)$}

Let $(z_1,\ldots,z_n)$ be a coordinate system
with $z_i^{-1}(0)=D_i$.
It induces a coordinate system
$(\theta_1,\ldots,\theta_n)$
of $\pi_2^{-1}(O)$,
which is independent of the choice of $(z_1,\ldots,z_n)$
up to parallel transport.
We take a coordinate system $t$ of $\cnum$,
which induces a coordinate system $\theta$
of $\pi_0^{-1}(O)$.
The induced map
$\pi_2^{-1}(O)\lrarr \pi_0^{-1}(O)$
is affine with respect to
the coordinate systems
$(\theta_1,\ldots,\theta_n)$
and $\theta$.

Let us consider the behaviour of
$\pi_1^{-1}(I)$ around $P\in \pi_2^{-1}(O)$,
where $I$ is a constructible sheaf on $\pi_0^{-1}(O)$.
We may assume $P=(0,\ldots,0)$.
The map $\pi_2^{-1}(O)\lrarr \pi_0^{-1}(O)$
is of the form
$(\theta_1,\ldots,\theta_n)\longmapsto
 \sum \alpha_i\,\theta_i+\beta$,
where $\beta=\pi_1(P)$.
The sheaf $I$ is the direct sum of sheaves
of the following forms:
\begin{itemize}
\item
 The constant sheaf around $\beta$.
\item
 $j_!\cnum_J$ or $j_{\ast}\cnum_J$,
 where $J$ is an open interval
 such that one of the end points is $\beta$,
 and $j$ denotes the inclusion
 $J\lrarr \pi^{-1}(O)$.
\end{itemize}
Hence, $\pi_1^{-1}(I)$ around $P$
is described as the direct sum
of sheaves of the following forms:
\begin{itemize}
\item
 The constant sheaf
 $\cnum_{\pi_0^{-1}(O)}$.
\item
 $j_{\ast}\cnum_H$
 or $j_!\cnum_H$,
 where 
 $H$ is an open half space
 such that $\del H\ni P$,
 and $j:H\lrarr \pi^{-1}_0(O)$.
 They are denoted by
 $\cnum_{H\ast}$
 and 
 $\cnum_{H!}$.
\end{itemize}

\subsection{Local form of $\nbigl^{\leq D}$
 and $\nbigl/\nbigl^{\leq D}$}

Let $P\in \pi_0^{-1}(O)$.
We have a decomposition around $P$:
\[
 \nbigl=\bigoplus_{\gminia\in\Irr(\nabla)}
 \nbigl_{\gminia}
\quad\quad
 \nbigl^{\leq D}=\bigoplus_{\gminia\in\Irr(\nabla)}
 \nbigl^{\leq D}_{\gminia}
\]

Let us describe $\nbigl_{\gminia}$ 
and $\nbigl/\nbigl^{\leq D}_{\gminia}$ around $P$.
For an appropriate coordinate system,
$\gminia=z_1^{-m_1}\cdots z_n^{-m_n}$
for some $m_i\geq 0$.
Let $q_{\gminia}:\Delta^n\lrarr \Delta$
be given by
$(z_1,\ldots,z_n)\longmapsto
 \prod z_i^{m_i}$.
Let $\pi_{\Delta}:\Deltatilde(0)\lrarr \Delta$
be the real blow up.
We have the induced map:
\[
 q_{\gminia}:
 \Xtilde(D)\lrarr \Deltatilde(0),
\quad 
 (r_i,\theta_i)\longmapsto
 \Bigl(
 \prod_{i=1}^nr_i^{m_i},\,
 \sum m_i\theta_i
 \Bigr)
\]
Let $\nbigq$ be the local system 
on $\Deltatilde(0)$ with Stokes structure,
corresponding to the meromorphic flat bundle
$\bigl(
 \nbigo_{\Delta}(\ast 0),\,d+d(1/z)
 \bigr)$.
Note that
$\nbigq/\nbigq^{\leq 0}$ is 
the constructible sheaf 
$j_{\ast}\cnum_J$
on $\pi_{\Delta}^{-1}(0)$,
where $j:J=(-\pi,\pi)\lrarr \pi_{\Delta}^{-1}(0)$.
Let $r(\gminia)$ be the rank of $\nbigl_{\gminia}$.
We have isomorphisms:
\[
 \nbigl_{\gminia}\simeq
 q_{\gminia}^{\ast}\nbigq^{\oplus\, r(\gminia)}
\quad\quad
 \nbigl^{\leq D}_{\gminia}\simeq
 q_{\gminia}^{\ast}\bigl(
 \nbigq^{\leq 0}\bigr)^{\oplus\, r(\gminia)}
\quad\quad
 \nbigl_{\gminia}/\nbigl^{\leq D}_{\gminia}
\simeq
 q_{\gminia}^{\ast}\bigl(
  \nbigq/\nbigq^{\leq 0}
 \bigr)^{\oplus\,r(\gminia)}
\]
Around $P$,
we have an isomorphism
$q_{\gminia}^{\ast}\bigl(
 \nbigq/\nbigq^{\leq 0}\bigr)
\simeq
 \iota_{\ast}\cnum$,
where 
$Z:=q_{\gminia}^{-1}(J)$
and $\iota:Z\lrarr (S^1)^n\times\real_{\geq 0}^n$.
Note that $Z$ is of the form
 $Z_0\times\del\real_{\geq 0}^n$,
 where 
 $Z_0$ is the inverse image of $J$
 via the induced map
$(S^1)^n\times\{0\}
 \lrarr 
 S^1\times\{0\}$.
Hence,
$q_{\gminia}^{\ast}\bigl(
 \nbigq/\nbigq^{\leq 0}\bigr)$
is isomorphic to one of the following,
around $P$:
\begin{itemize}
\item
The constant sheaf
$\cnum_{(S^1)^n\times\del\real^n_{\geq\,0}}$.
\item
$j_{K\ast}\cnum_{K\times\del\real_{\geq 0}^n}$,
where $K$ is an open half space
such that $\del K\ni P$,
and 
$j_K:K\times\del\real_{\geq 0}^n
 \lrarr (S^1)^n\times\real^n_{\geq\,0}$.
It is denoted by
$\cnum_{K\times\del\real^n_{\geq 0}\ast}$.
\end{itemize}

\subsection{Proof of Theorem \ref{thm;09.10.17.1}}

We reduce the proof of the theorem
to the computation of
$ \next^{i}\bigl(
 \pi_1^{-1}I,
 q_{\gminia}^{-1}(\nbigq/\nbigq^{\leq 0})
 \bigr)$
for $i\leq n-2$,
where $I$ is a constructible sheaf on $\pi_0^{-1}(O)$.

\begin{lem}
\label{lem;09.10.17.10}
We have
$\next^i(\pi_1^{-1}I,q_{\gminia}^{-1}\nbigq)=0$
for any $i$.
In particular, we have
isomorphisms:
\[
\next^i\bigl(
 \pi_1^{-1}I,\,
 q_{\gminia}^{-1}\nbigq^{\leq 0}
 \bigr)
\simeq
 \next^{i-1}\bigl(
 \pi_1^{-1}I,
 q_{\gminia}^{-1}(\nbigq/\nbigq^{\leq 0})
 \bigr). 
\]
\end{lem}
\pf
Let $\iota:(S^1)^n\times\{0\}\lrarr
(S^1)^n\times\del\real_{\geq 0}^n$
denote the inclusion.
There exists a constructible sheaf
$\nbigf$ on $(S^1)^n$
such that 
$\pi_1^{-1}I\simeq\iota_{\ast}\nbigf$.
We have the adjunction
$\next^i\bigl(
 \iota_{\ast}\nbigf,\,
 q_{\gminia}^{-1}\nbigq
 \bigr)
=\iota_{\ast}
 \next^i(\nbigf,i^!q_{\gminia}^{-1}\nbigq)$.
Note
$\iota^!q_{\gminia}^{-1}\nbigq
=\DDD \iota^{-1}\DDD
 \bigl(
 q_{\gminia}^{-1}\nbigq
 \bigr)=0$,
because
$\DDD q_{\gminia}^{-1}\nbigq$
is $0$-extension of
a constant sheaf on
$(S^1)^n\times\real_{>0}^n$
by 
$(S^1)^n\times\real_{>0}^n
\lrarr
 (S^1)^n\times\real_{\geq 0}^n$.
Hence, we obtain
$\next^i\bigl(
 \iota_{\ast}\nbigf,
 q_{\gminia}^{-1}\nbigq
 \bigr)=0$,
and the proof of Lemma \ref{lem;09.10.17.10}
is finished.
\hfill\qed

\vspace{.1in}

Now,
let us prove the following vanishing
of the stalks at $P$:
\begin{equation}
\label{eq;09.10.17.15}
\next^{j}\bigl(
 \pi_1^{-1}I,
 q_{\gminia}^{-1}(\nbigq/\nbigq^{\leq 0})
 \bigr)_P=0,
\quad (j\leq n-2)
\end{equation}
It can be computed 
on $(S^1)^n\times\del\real_{\geq 0}^n$.
We have the following cases,
divided by the local forms of
$\pi_1^{-1}(I)$ and
$q_{\gminia}^{-1}(\nbigq/\nbigq^{\leq 0})$
around $P$:
\begin{description}
\item[(I)]
$\pi_1^{-1}I\simeq \cnum_{(S^1)^n}$
and 
$q_{\gminia}^{-1}\bigl(
 \nbigq\big/\nbigq^{\leq 0}
 \bigr)
\simeq
 \cnum_{(S^1)^n\times\del\real_{\geq 0}^n}$.
\item[(II)]
$\pi_1^{-1}I\simeq \cnum_{(S^1)^n}$
and 
$q_{\gminia}^{-1}\bigl(
 \nbigq\big/\nbigq^{\leq 0}\bigr)
\simeq
 \cnum_{K\times\del\real^n_{\geq 0}\,\ast}$.
\item[(III)]
$\pi_1^{-1}I=\cnum_{H\star}$
and 
$q_{\gminia}^{-1}\bigl(
 \nbigq\big/\nbigq^{\leq 0}\bigr)
\simeq
 \cnum_{(S^1)^n\times\del\real_{\geq\,0}^n}$,
where $\star=\ast,!$.
\item[(IV)]
$\pi_1^{-1}I\simeq\cnum_{H\star}$
and 
$q_{\gminia}^{-1}\bigl(
 \nbigq\big/\nbigq^{\leq 0}\bigr)
\simeq
 \cnum_{K\times\del\real_{\geq\,0}^n\ast}$,
where $\star=\ast,!$.
Moreover, this is divided into
three cases
(IV-1) $\del H$ and $\del K$ are transversal,
(IV-2) $H=K$,
(IV-3) $H=-K$.
\end{description}
In the following,
for a given $i:Y_1\subset Y_2$
and $\star=\ast,!$,
let $\cnum_{Y_1\star}:=
 i_{\star}\cnum_{Y_1}$ on $Y_2$.
It is also denoted just by $\cnum_{Y_1}$
if there is no risk of confusion.

\subsubsection{The case (I)}

Instead of $(S^1)^n\times\{0\}\lrarr
 (S^1)^n\times\del\real^n_{\geq 0}$,
it is enough to consider the inclusion
$\{0\}\lrarr\del\real^n_{\geq 0}\simeq
 \real^{n-1}$.
We obtain (\ref{eq;09.10.17.15})
from the following standard result:
\[
 \next^j\bigl(\cnum_0,\cnum_{\real^{n-1}}\bigr)_0
\simeq
 \left\{
 \begin{array}{ll}
 0 & (j\leq n-2)\\
 \cnum & (j=n-1)
 \end{array}
 \right.
\]

\subsubsection{The case (II)}

We have the exact sequence
{\small$0\lrarr \cnum_{(S^1)^n\setminus K!}\lrarr
 \cnum_{(S^1)^n}\lrarr
 \cnum_{K\ast}\lrarr 0$.}
Let $\iota$ denote the inclusion
$\bigl((S^1)^n\setminus K\bigr)\times
 \del\real_{\geq\,0}^n
\lrarr
 (S^1)^n\times \del\real_{\geq\,0}^n$.
Note $\iota^{\ast}=\iota^!$,
and hence
$\iota^{!}\cnum_{K\times\del\real^n_{\geq\,0}\ast}
 =0$.
We have
\[
 \next^j\Bigl(\cnum_{\bigl((S^1)^n\setminus K\bigr)
 \times\{0\}\,!},\,
 \cnum_{K\times\del\real_{\geq\,0}^n\,\ast}
 \Bigr)_P
\simeq
 \iota_{\ast}\next^j\Bigl(
 \cnum_{\bigl((S^1)^n\setminus K\bigr)\times\{0\}},\,
 \iota^{!}\cnum_{K\times\del\real^n_{\geq\,0}\ast}
 \Bigr)_P
=0
\]
Hence, we obtain
\[
 \next^j\Bigl(
 \cnum_{(S^1)^n},
 \cnum_{K\times\del\real_{\geq\,0}^n\,\ast}
 \Bigr)_P
\simeq
 \next^j\Bigl(
 \cnum_{K\ast},
 \cnum_{K\times\del\real^n_{\geq\,0}\ast}
 \Bigr)_P
=\left\{
 \begin{array}{ll}
 0 & (j\leq n-2)\\
 \cnum & (j=n-1)
 \end{array}
 \right.
\]

\subsubsection{The case (III)}

Let us consider the case $\star=\ast$.
We have the exact sequence:
\[
 0\lrarr
 \cnum_{(S^1)^n\times\del\real_{\geq\,0}^n
 \setminus H\times\{0\}\,!}
\lrarr
 \cnum_{(S^1)^n\times\del\real_{\geq\,0}^n}
\lrarr
 \cnum_{H\ast}
\lrarr 0
\]
Let $k_1$ denote the inclusion $H\times\{0\}
 \lrarr (S^1)^n\times\del\real_{\geq\,0}^n$,
and let $k_2$ denote the open embedding
of the complement.
Because 
$k_1^{\ast}
 \cnum_{(S^1)^n\times\del\real^n_{\geq\,0}
 \setminus H\times\{0\}\,!}=0$,
we have the following isomorphisms:
\begin{multline}
 \nrhom\bigl(
 \cnum_{(S^1)^n\times\del\real^n_{\geq\,0}
 \setminus H\times\{0\}\,!},\,
 \cnum_{(S^1)^n\times\del\real^n_{\geq\,0}}
 \bigr)_P 
\simeq \\
 \nrhom\bigl(
 \cnum_{(S^1)^n\times\del\real^n_{\geq\,0}
 \setminus H\times\{0\}!},\,
 \cnum_{(S^1)^n\times\del\real^n_{\geq\,0}
 \setminus H\times\{0\}!}
 \bigr)_P 
\simeq \\
 k_{2\ast}\bigl(
 \cnum_{(S^1)^n\times\del\real^n_{\geq\,0}
 \setminus H\times\{0\}}
 \bigr)_P
\simeq
 \bigl(
 \cnum_{(S^1)^n\times\del\real_{\geq\,0}^n}
 \bigr)_P
\end{multline}
W obtain
{\small $\nrhom\bigl(
 \cnum_{H\ast},\,
 \cnum_{(S^1)^n\times\del\real_{\geq\,0}^n}
 \bigr)_P=0$.}
In particular,
$\next^j(\cnum_{H\ast},\cnum_{(S^1)^n\times
 \del\real^n_{\geq\,0}})_{P}=0$
for any $j$.

Let us consider the case $\star=!$.
We have the exact sequence
$0\lrarr \cnum_{H!}\lrarr
 \cnum_{(S^1)^n}\lrarr
 \cnum_{(S^1)^n\setminus H\,\ast}\lrarr 0$.
Hence, we obtain the following
isomorphisms:
\[
 \next^j\bigl(
 \cnum_{H!},
 \cnum_{(S^1)^n\times\del\real_{\geq\,0}^n}
 \bigr)_{P}
=\next^j\bigl(
 \cnum_{(S^1)^n},\,
 \cnum_{(S^1)^n\times\del\real_{\geq\,0}^n}
 \bigr)_{P}=
 \left\{
 \begin{array}{ll}
 0 & (j\leq n-2)\\
 \cnum & (j=n-1)
 \end{array}
 \right.
\]

\subsubsection{The case (IV-1)}

Let us consider the case $\star=\ast$.
Let $\nbign$ be the kernel of 
$\cnum_{H\ast}\lrarr \cnum_{H\cap K\ast}$.

\begin{lem}
We have 
$ \nrhom\bigl(
 \nbign,\,
 \cnum_{K\times\del\real_{\geq 0}^n*}
 \bigr)_{P}=0$.
\end{lem}
\pf
Let $\iota$ be the inclusion
$\bigl((S^1)^n\setminus K\bigr)\times
 \del\real^n_{\geq\,0}\lrarr
 (S^1)^n\times\del\real^n_{\geq \,0}$.
Then,
$\nbign$ is of the form $\iota_!\nbign_1$.
Then, the claim follows from 
$\iota^!\cnum_{K\times\del\real_{\geq\,0}^n\ast}=0$.
\hfill\qed

\vspace{.1in}

We have the exact sequence:
$0\lrarr \cnum_{
 K\times\del\real_{\geq 0}^n\setminus 
 (H\cap K)\times\{0\}!}
\lrarr
 \cnum_{K\times\del\real_{\geq 0}^n}
\lrarr
 \cnum_{(H\cap K)\times\{0\}\ast}
\lrarr 0$.
Let $k$ denote the inclusion
$K\times\del\real_{\geq\,0}^n
\setminus(H\cap K)\times\{0\}
\lrarr K\times\del\real_{\geq 0}^n$.
We have the following isomorphisms:
\begin{multline}
 \nrhom\bigl(
 \cnum_{K\times\del\real_{\geq 0}^n
 \setminus (H\cap K)\times\{0\}!},\,
 \cnum_{K\times\del\real^n_{\geq\,0}}
 \bigr)_P
\simeq \\
 Rk_{\ast}\nrhom\bigl(
 \cnum_{K\times\del\real_{\geq\,0}^n
 \setminus (H\cap K)\times\{0\}},\,
 \cnum_{K\times\del\real_{\geq\,0}^n
 \setminus (H\cap K)\times\{0\}}
 \bigr)_{P} 
\simeq
 \cnum_{K\times\del\real^n_{\geq\,0},P}
\end{multline}
Hence, we obtain
$\nrhom\bigl(
 \cnum_{(H\cap K)\times \{0\}\,\ast},\,
 \cnum_{K\times\del\real_{\geq\,0}^n\,\ast}
 \bigr)_{P}=0$.
In particular, we have
$\next^j\bigl(
 \cnum_{H\ast},\,
 \cnum_{K\times\del\real_{\geq\,0}^n\,\ast}
 \bigr)_P=0$ for any $j$.

Let us consider the case $\star=!$.
We have an exact sequence
$0\lrarr \cnum_{H!}\lrarr
 \cnum_{(S^1)^n}\lrarr
 \cnum_{(S^1)^n\setminus H\ast}\lrarr 0$
on $(S^1)^n$.
By using the previous results,
we obtain
\[
 \next^j\bigl(
 \cnum_{H!},\,
 \cnum_{K\times\del\real_{\geq 0}^n\ast}
 \bigr)_{P}
=\left\{
 \begin{array}{ll}
  0 & (j\leq n-2)\\
 \cnum & (j=n-1)
 \end{array}
 \right.
\]

\subsubsection{The case (IV-2)}

Let us consider the case $\star=\ast$.
By considering
$0\lrarr\del\real_{\geq 0}^n$,
we obtain
\[
 \next^j\bigl(
 \cnum_{H\ast},
 \cnum_{H\times\del\real_{\geq 0}^n\ast}
 \bigr)_P
\simeq
 \left\{\begin{array}{ll}
 0 & (j\leq n-2)\\
 \cnum & (j=n-1)
 \end{array}
 \right.
\]
Let us consider the case $\star=!$.
We have an exact sequence
$0\lrarr \cnum_{H!}\lrarr \cnum_{H\ast}
 \lrarr \cnum_{\del H\ast}\lrarr 0$.
Let us look at
$ \next^j\bigl(
 \cnum_{\del H\ast},\,
 \cnum_{H\times\del\real_{\geq\,0}^n}
 \bigr)_P$.
For 
$0\lrarr \closedopen{0}{1}\times\real^{n-1}$,
we have
$\next^j\bigl(
 \cnum_0,\cnum_{\closedopen{0}{1}\times
 \real^{n-1}}
 \bigr)=0$ for any $j$.
Hence, we obtain
\[
 \next^j\bigl(
 \cnum_{H!},\,
 \cnum_{H\times\del\real_{\geq 0}^n}
 \bigr)_P=\left\{
 \begin{array}{ll}
 0 & (j\leq n-2)\\
 \cnum & (j=n-1)
 \end{array}
 \right.
\]

\subsubsection{The case (IV-3)}

It is easy to show
$\next^j\bigl(
 \cnum_{H!},\cnum_{K\times\del\real_{\geq 0}^n}
 \bigr)
=0$ for any $j$.
By using the argument in (IV-2),
we can show
$\next^j\bigl(
 \cnum_{H\ast},\cnum_{K\times\del\real^n}
 \bigr)=0$ for any $j$.
Thus, the proof of Theorem
\ref{thm;09.10.17.1} is finished.
\hfill\qed

\subsection{A uniqueness result on the $K$-structure}
\label{subsection;09.12.5.33}

We use the notation in 
\S\ref{subsection;09.10.27.40}.
Let $V$ be a good meromorphic flat bundle on $(X,D)$.
Let $g$ be a holomorphic function on $X$
such that $g^{-1}(0)=D$,
and let $i_g$ be the graph $X\lrarr X\times\cnum$.
We regard 
$\DR^{\nil}_{X\times\cnumtilde}
 (i_{g\dagger}V)$
as a cohomologically constructible sheaf
on $\Xtilde(g)$.

Let $K$ be a subfield of $\cnum$.
A $K$-structure of
$\DR^{\nil}_{X\times\cnumtilde}
 \bigl(i_{g\dagger}V\bigr)$
is defined to be a $K$-cohomologically
constructible complex $\nbigf$ on $\Xtilde(g)$
with an isomorphism
$\alpha:
\nbigf\otimes\cnum\simeq
 \DR^{\nil}_{X\times\cnumtilde}(i_{g\dagger}V)$
in the derived category.
Two $K$-structures
$(\nbigf_i,\alpha_i)$ $(i=1,2)$ are called
equivalent
if there exists an isomorphism
$\beta:\nbigf_1\lrarr\nbigf_2$
for which the following diagram is commutative:
\[
 \begin{CD}
 \nbigf_1\otimes\cnum
 @>{\beta\otimes\cnum}>>
 \nbigf_2\otimes\cnum \\
 @V{\alpha_1}VV @V{\alpha_2}VV \\
 \DR^{\nil}_{X\times\cnumtilde}
 (i_{g\dagger}V)
 @>{=}>>
 \DR^{\nil}_{X\times\cnumtilde}
 (i_{g\dagger}V) 
 \end{CD}
\]
\begin{lem}
\label{lem;09.10.28.100}
Let $(\nbigf_i,\alpha_i)$ $(i=1,2)$ be $K$-structures
of $\DR^{\nil}_{X\times\cnumtilde}\bigl(
 i_{g\dagger}V\bigr)$.
If their restriction to $\pi_1^{-1}(X-D)$
are equivalent,
then they are equivalent on $\Xtilde(g)$.
\end{lem}
\pf
We put
$\nbigf^{\cnum}_i:=
 \nbigf_i\otimes\cnum$.
We have the following commutative diagram:
\[
 \begin{CD}
 \Hom(\nbigf_1,\nbigf_2)\otimes\cnum
 @>>>
 \Hom\bigl(
 \nbigf_{1|\pi_1^{-1}(X-D)},
 \nbigf_{2|\pi_1^{-1}(X-D)}
 \bigr)\otimes\cnum
 \\
 @VV{\simeq}V
 @VV{\simeq}V \\
 \Hom\bigl(
 \nbigf_1^{\cnum},\nbigf_2^{\cnum}
 \bigr)
 @>>>
 \Hom\bigl(
 \nbigf^{\cnum}_{1|\pi_1^{-1}(X-D)},
 \nbigf^{\cnum}_{2|\pi_1^{-1}(X-D)}
 \bigr)
 \end{CD}
\]
According to Theorem \ref{thm;09.10.17.1},
the horizontal arrows are injective.
Hence, we obtain the equality
\[
\Hom(\nbigf_1,\nbigf_2)
=
\Hom\bigl(
 \nbigf_{1|\pi_1^{-1}(X-D)},
 \nbigf_{2|\pi_1^{-1}(X-D)}
 \bigr)
\cap
\Hom\bigl(
 \nbigf_1^{\cnum},\nbigf_2^{\cnum}
 \bigr)
\]
in
$\Hom\bigl(
 \nbigf^{\cnum}_{1|\pi_1^{-1}(X-D)},
 \nbigf^{\cnum}_{2|\pi_1^{-1}(X-D)}
 \bigr)$.
Then, the element
of $\Hom(\nbigf_1^{\cnum},\nbigf_2^{\cnum})$
corresponding to
the identity of
$\DR^{\nil}_{X\times\cnumtilde}
(i_{g\dagger}V)$
comes from 
$\Hom(\nbigf_1,\nbigf_2)$.
\hfill\qed

%% file: 6.1.tex
\subsection{Good $K$-structure 
of good meromorphic flat bundles}
\label{subsection;09.10.28.10}

Let $K\subset \cnum$ be a subfield.
Let $X$ be a complex manifold
with a normal crossing hypersurface $D$.

\begin{df}
\mbox{{}}\label{df;13.4.21.1}
Let $V$ be a good meromorphic flat bundle on $(X,D)$.
\begin{itemize}
\item
A $K$-structure of $V$ is a pre-$K$-Betti structure of
the flat bundle $V_{|X-D}$.
\index{$K$-structure}
\item
A $K$-structure of $V$ is good
if the Stokes structures are defined over $K$.
\index{good $K$-structure}
\hfill\qed
\end{itemize}
\end{df}
Later, we shall extend the definition
to the case where $V$ is not necessarily good.
(See \S\ref{subsection;13.4.25.200}.)

Let $D=D_1\cup D_2$ be a decomposition.
Let $\nbigl$ be the local system with the Stokes structure
on $\Xtilde(D)$ associated to $V$.
Recall that the complex
$\DR^{<D_1\leq D_2}_{\Xtilde(D)}
 \bigl(V\bigr)$ 
is quasi-isomorphic to 
$\nbigl^{<D_1\leq D_2}[\dim X]$.
(See \S\ref{subsection;09.10.28.1}.)
If $V$ has a good $K$-structure,
it is naturally equipped with a $K$-structure
$\nbigl_K^{<D_1\leq D_2}[\dim X]$.
By the isomorphisms (\ref{eq;13.4.17.2})
and (\ref{eq;09.10.4.1}),
we obtain a pre-$K$-Betti structure
\[
 \nbigf_V^{<D_1\leq D_2}:=
 R\pi_{\ast}\nbigl_K^{<D_1\leq D_2}[\dim X]
\]
of the holonomic $\nbigd$-module $V\bigl(!D_1\bigr)$.
\index{complex $\nbigf_V^{<D_1\leq D_2}$}
This pre-$K$-Betti structure is called canonical.
\index{canonical pre-$K$-Betti structure}
Let $D_1'\cup D_2'=D$ be another 
decomposition
such that $D_1\subset D_1'$.
The natural morphism
$V\bigl(!D_1'\bigr)\lrarr
 V\bigl(!D_1\bigr)$ is compatible
with the pre-$K$-Betti structures.
We use the symbols
$\nbigf_{V\ast}$ and $\nbigf_{V!}$ to denote
$\nbigf_V^{\leq D}$ and $\nbigf_{V}^{<D}$,
respectively.
We also use the symbol $\nbigf_V$
to denote $\nbigf_{V\ast}$
for simplicity.
\index{complex $\nbigf_{V\ast}$}
\index{complex $\nbigf_{V\bikkuri}$}

\vspace{.1in}
More generally,
let $\iota:Z\subset X$ be a complex submanifold
with a normal crossing hypersurface $D_Z$.
Let $V_Z$ be a good meromorphic flat bundle
on $(Z,D_Z)$.
We say that $\iota_{\dagger}V_Z$ has 
a good $K$-structure
if $V_Z$ has a good $K$-structure
in the above sense.
The canonical pre-$K$-Betti structures
for $\iota_{\dagger}V_Z(!D_{Z,1})$ are also
defined in a similar way
for a decomposition
$D_Z=D_{Z,1}\cup D_{Z,2}$.

\subsection{Some basic property}

\subsubsection{Some functoriality}

Let $X$ be any complex manifold
with a normal crossing hypersurface $D$.
The following lemma is clear.
\begin{lem}
Let $V_i$ $(i=1,2)$ be good meromorphic flat bundles
on $(X,D)$ with a good $K$-structure.
If $V_1\oplus V_2$ is good,
then the induced $K$-structure is good.
Similar claims hold for
$V_1\otimes V_2$
and $\nhom(V_1,V_2)$.
\hfill\qed
\end{lem}

Let $V$ be a good meromorphic flat bundle on $(X,D)$.
Let $\varphi:X'\lrarr X$ be a morphism of complex manifolds
such that $D':=\varphi^{-1}(D)$ is normal crossing.
We obtain a good meromorphic flat bundle $V':=\varphi^{\ast}V$
on $(X',D')$.
Suppose that $V$ is equipped with a $K$-structure,
which induces a $K$-structure of $V'$.
\begin{lem}
\label{lem;14.1.15.10}
If the $K$-structure of $V$ is good,
the $K$-structure of $V'$ is also good.
Conversely, suppose that
$\varphi$ is surjective
and that the $K$-structure of $V'$ is good.
Then, the $K$-structure of $V$ is good.
\end{lem}
\pf
Let $P'$ be any point of $D'$.
Let $P:=\varphi(P')$.
We take a small neighbourhood $X_P$
with a coordinate $(z_1,\ldots,z_n)$
around $P$ in $X$ such that
$D=\bigcup_{i=1}^{\ell}\{z_i=0\}$,
and a ramified covering
$\kappa_P:
 (X^{(1)}_P,D^{(1)}_{P})
 \lrarr
 (X_P,D\cap X_P)$
such that
$V_P^{(1)}:=\kappa_P^{\ast}(V)$ is unramifiedly good.
Let $e_i$ $(i=1,\ldots,\ell)$
denote the ramification index of $\kappa_P$
along $z_i=0$.
We take a small neighbourhood 
$X_{P'}'$ of $P'$.
Because $(z_i\circ\varphi)^{-1}(0)$ $(i=1,\ldots,\ell)$
are contained in $D'\cap X_P'$,
we can take a ramified covering
$\kappa'_{P'}:(X^{\prime(1)}_{P'},D^{\prime(1)}_{P'})
 \lrarr (X'_{P'},D'\cap X'_{P'})$
such that there exist functions
$(z_i\circ\varphi\circ\kappa'_{P'})^{1/e_i}$
$(i=1,\ldots,\ell)$ on $X^{\prime(1)}_{P'}$.
Then, we have a morphism 
$\rho:X^{\prime(1)}_{P'}\lrarr X^{(1)}_P$
such that 
$\kappa_P\circ\rho
=\varphi\circ\kappa'_{P'}$.
Then, 
$V^{\prime(1)}:=(\kappa'_{P'})^{\ast}V'
=\rho^{\ast}\kappa_P^{\ast}(V)$
is unramifiedly good.
Let $\nbigl$ be the local system on
$\Xtilde_P^{(1)}(D^{(1)}_P)$ associated to 
$V^{(1)}$.
Let $\nbigl'$ be the local system on
$\Xtilde^{\prime(1)}_{P'}(D^{\prime(1)}_P)$
associated to 
$V^{\prime(1)}$.
Let $\rhotilde:
 \Xtilde^{\prime(1)}_P(D^{\prime(1)}_P)
 \lrarr \Xtilde^{(1)}_P(D^{(1)}_P)$
be the map induced by $\rho$.
We have 
$\nbigl'=\varphitilde^{-1}(\nbigl)$.

Let $\pi^{(1)}:
 \Xtilde^{(1)}_P(D^{(1)}_P)
 \lrarr X^{(1)}_P$
and 
$\pi^{\prime(1)}:
 \Xtilde^{\prime(1)}_{P'}(D^{\prime(1)}_{P'})
 \lrarr X^{\prime(1)}_{P'}$
denote the projections.
Let $Q'_1$ be any point of
$(\pi^{\prime(1)})^{-1}(D^{\prime(1)}_P)$.
We set $Q_1:=\rhotilde(Q'_1)$.
Let $P'_1:=\pi^{\prime(1)}(Q')$
and $P_1:=\pi^{(1)}(Q)$.
The set of the irregular values of $V^{\prime(1)}$
at $P'_1$
is the pull back of 
the set of the irregular values of $V^{(1)}$ at $P_1$.
The partial order $\leq_{Q_1'}$ on the set
is equal to $\leq_{Q_1}$.
The Stokes filtration $\nbigf^{Q_1'}$
is obtained as the pull back of
$\nbigf^{Q_1}$.
Hence, $\nbigf^{Q_1}$ is defined over $K$
if and only if $\nbigf^{Q_1'}$ is defined over $K$.
\hfill\qed

\subsubsection{Curve test}

Let us consider the case $X=\Delta^n$,
$D_i:=\{z_i=0\}$
and $D=\bigcup_{i=1}^{\ell} D_i$.
We set
$D_i^{\circ}:=D_i\setminus \bigcup_{j\neq i} D_j$.
Let $p_i:X\lrarr D_i$ denote the projection.

\begin{prop}
\label{prop;14.1.18.30}
Let $V$ be a good meromorphic flat bundle on $(X,D)$
with a $K$-structure with the following property.
\begin{description}
\item[(C1)]
Let $P$ be any point of $D_i^{\circ}$ for $i=1,\ldots,\ell$.
Then, the induced $K$-structure of 
$V_{|p_i^{-1}(P)}$ is good.
\end{description}
Then, the $K$-structure of $V$ is good.
\end{prop}
\pf
We may assume that $V$ is unramifiedly good.
Let $\pi:\Xtilde(D)\lrarr X$ denote the projection.
Let $\nbigl$ be the local system on $\Xtilde(D)$
with the induced $K$-structure.
Let $Q$ be any point of $\pi^{-1}(D)$.
It is enough to prove that the Stokes filtration
$\nbigf^Q(\nbigl_Q)$ is defined over $K$.
It is enough to consider the case 
$\pi(Q)=(0,\ldots,0)$.
We set $S:=\{(\gminia,\gminib)\in\Irr(V)^2\,|\,
\gminia\neq\gminib \}$.
We have $i$ such that
$\ord_{z_i}(\gminia-\gminib)<0$
for any $(\gminia,\gminib)\in S$.
For any $(\gminia,\gminib)\in S$,
let $H(\gminia,\gminib)$ be 
denote the intersection of
$\pi^{-1}(D_i)$
and the closure of
$\bigl\{R\in X\setminus D\,\big|\,
 \Re(\gminia-\gminib)(R)=0\,
 \bigr\}$in $\Xtilde(D)$.
Let $\nbigu$ be a small neighbourhood of
$Q$ in $\pi^{-1}(D_i)$.
Then, 
for $(\gminia,\gminib)\in S$,
we have
$\gminia<_Q\gminib$
if and only if
we have
$\gminia<_{Q'}\gminib$
for any 
$Q'\in \nbigu':=
 \pi^{-1}(D_i^{\circ})\cap
 \nbigu\setminus 
 \bigcup_{(\gminia,\gminib)\in S}
 H(\gminia,\gminib)$.
We have natural identifications
of $\nbigl_Q$ and $\nbigl_{Q'}$
for $Q'\in\nbigu$.
We have
$\nbigf^Q_{\gminia}
=\bigcap_{Q'\in \nbigu'}
 \nbigf^{Q'}_{\gminia}$.
Under the assumption {\bf (C1)},
$\nbigf^{Q'}_{\gminia}$
are defined over $K$
for any $Q'\in \nbigu'$.
Hence, we obtain that
$\nbigf^Q_{\gminia}$ are defined over $K$.
\hfill\qed

\subsubsection{Sub-quotients}

Let $X$ be any complex manifold
with a normal crossing hypersurface $D$.
Let $0\lrarr V_1\lrarr V\lrarr V_2\lrarr 0$
be an exact sequence of good meromorphic flat bundles on $(X,D)$.
Suppose that $V$ and $V_i$ are equipped with $K$-structures
which are compatible with the morphisms.

\begin{lem}
\label{lem;13.4.23.1}
If the $K$-structure of $V$ is good,
then the $K$-structures of $V_i$ $(i=1,2)$
are good.
\end{lem}
\pf
We may assume that $V$ is unramifiedly good.
We may assume that $X=\Delta$
and $D=\{0\}$.
Let $\nbigl_i$ and $\nbigl$
be the local systems on $\Xtilde(D)$
corresponding to
$V_i$ and $V$, respectively.
For any point $P\in \Xtilde(D)$,
the stalks
$\nbigl_{1P}$ and $\nbigl_P$
are equipped with the Stokes filtrations
$\nbigf^P$.
Note that the Stokes filtrations are 
characterized by the growth order.
Hence, 
$\nbigl_{1P}\lrarr\nbigl_P$
is strict with respect to the filtrations,
i.e., 
$\nbigf^P(\nbigl_{1P})$
is equal to the filtration
obtained as the restriction of $\nbigf^P(\nbigl_P)$.
Then, if $\nbigl_{1P}$ and 
$\nbigf^P(\nbigl_P)$
are defined over $K$,
the filtration $\nbigf^P(\nbigl_{1P})$
is also defined over $K$.
\hfill\qed

\begin{lem}
\label{lem;13.5.18.1}
Let $V_i$ $(i=1,2)$ be good meromorphic flat bundles
on $(X,D)$.
Let $f:V_1\lrarr V_2$ be a morphism 
of meromorphic flat bundles.
\begin{itemize}
\item
$\Ker(f)$, $\Image(f)$ and $\Cok(f)$
are also good.
\item
Suppose that $V_i$ are equipped with good $K$-structures,
and that $f$ is compatible with the $K$-structures.
Then, the induced $K$-structures of
$\Ker(f)$, $\Cok(f)$
and $\Image(f)$ are good.
\end{itemize}
\end{lem}
\pf
It is enough to check the claims locally 
around any point of $D$.
We may assume that $V_i$ are unramifiedly good.
Let $P$ be any point of $D$.
Let $f_{|\Phat}$ denote the induced morphism
$V_{1|\Phat}\lrarr V_{2|\Phat}$.
Because the formal completion is exact,
we have 
$\Ker(f)_{|\Phat}
\simeq
 \Ker(f_{|\Phat})$
and similar isomorphisms
for $\Image$ and $\Cok$.
We have the decompositions
$V_{i|\Phat}
=\bigoplus_{\gminia\in\Irr(V_i,P)}
 V_{i,\Phat,\gminia}$.
It is easy to check that $f_{|\Phat}$ is compatible
with the decompositions.
Then, the first claim follows.
The second claim follows from the first claim
and Lemma \ref{lem;13.4.23.1}.
\hfill\qed

\vspace{.1in}
If $V_i$ are unramifiedly good in Lemma \ref{lem;13.5.18.1},
we have
$\Irr(\Ker f,P)\subset\Irr(V_1,P)$,
$\Irr(\Cok f,P)\subset\Irr(V_2,P)$
and 
$\Irr(\Image f,P)\subset\Irr(V_1,P)\cap\Irr(V_2,P)$.

\subsection{Functoriality for projective birational morphisms}
\label{subsection;13.4.21.2}

Let $D_3$ be a hypersurface of $X$.
Let $\varphi:X'\lrarr X$ be a projective birational morphism
such that $D':=\varphi^{-1}(D\cup D_3)$ is normal crossing,
and that 
 $X'\setminus D'\simeq X\setminus(D_3\cup D)$.
Let $V$ be a good meromorphic flat bundle 
on $(X,D)$.
Suppose that 
$V$ is equipped with a good $K$-structure.
We put 
$V':=\varphi^{\ast}V\otimes\nbigo_{X'}(\ast D')$.
The induced $K$-structure of $V'$ is good.
Let $D_1\cup D_2$ be a decomposition of $D$.
We set $D'_1:=\varphi^{-1}(D_1)$.
We take $D'_2\subset D'$ such that
$D_1'\cup D_2'$ is a decomposition of $D'$.

\begin{prop}
\label{prop;09.10.28.12}
The natural morphisms
\[
V(!D_1)\lrarr\varphi_{\dagger}V'(!D_1'),
\quad\quad
\varphi_{\dagger}V'(!D_2')
\lrarr V(!D_2)
\]
are compatible with
the canonical pre-$K$-Betti structures.
\end{prop}
\pf
Let us prove the second claim.
We use the notation 
introduced in \S\ref{subsection;09.10.28.13}.
Let $\varphitilde:\Xtilde'(D')\lrarr \Xtilde(D)$
be the induced map.
By construction,
it is easy to see that
the morphisms
$\DR^{<D_1\leq D_2}_{\Xtilde(D)}(V)
\lrarr
 R\varphitilde_{\ast}
 \DR^{<D_1'\leq D_2'}_{\Xtilde'(D')}(V')$
and 
$R\varphitilde_{\ast}
 \DR^{<D_2'\leq D_1'}_{\Xtilde'(D')}(V')
\lrarr
\DR^{<D_2\leq D_1}_{\Xtilde(D)}(V)$
are compatible with 
the induced $K$-structures.
Then, the second claim follows from
Theorem \ref{thm;09.10.4.100}.
\hfill\qed

\subsection{A characterization of
 compatibility with Stokes filtrations}
\label{subsection;13.4.20.350}

Let $X=\Delta^n$ and $D=\bigcup_{i=1}^{\ell}\{z_i=0\}$.
Let $V$ be an unramifiedly good meromorphic
flat bundle on $(X,D)$.
Its good set of irregular values is denoted by
$\Irr(V)$.
For each $\gminia\in\Irr(V)$,
put $L(-\gminia)=\nbigo_X(\ast D)\,e$
with the meromorphic flat connection
$\nabla e=e\,d(-\gminia)$.
We fix a $K$-structure of $L(-\gminia)$
by the trivialization $\exp(\gminia)\,e$.
We have a constructible sheaf 
$\DR_{\Xtilde(D)}^{\rapid}
 \bigl(V\otimes L(-\gminia)\bigr)$
on $\Xtilde(D)$.
The following lemma will be useful
to check that a $K$-structure is good.

\begin{lem}
\label{lem;13.4.20.250}
Suppose that $V$ has a $K$-structure
with the following property:
\begin{itemize}
\item
For each $\gminia\in\Irr(V)$,
the induced $K$-structure of
$\bigl(V\otimes L(-\gminia)\bigr)_{|X-D}$
is extended to a $K$-structure of
$\DR^{\rapid}_{\Xtilde(D)}
 \bigl(V\otimes L(-\gminia)\bigr)$.
\end{itemize}
Then, the $K$-structure of $V$ is good.
\end{lem}
\pf
Let $\nbigl$ be the local system with the Stokes structure
on $\Xtilde(D)$ associated to $V_{|X\setminus D}$.
It is equipped with the Stokes structure
i.e., for each $P\in \pi^{-1}(D)$,
the stalk $\nbigl_P$ has 
the Stokes filtration $\nbigf^P$.
By the assumption,
the local system $\nbigl$ has a $K$-structure.
Let $O=(0,\ldots,0)\in X$.
Let $\pi$ denote the projection $\Xtilde(D)\lrarr X$.
It is enough to prove that
the Stokes filtrations $\nbigf^P$ of $\nbigl_P$
are defined over $K$
for $P\in\pi^{-1}(O)$.

Let $S$ denote the set of pairs
$(\gminia,\gminib)$ in $\Irr(V)$
with $\gminia\neq\gminib$.
For  any $(\gminia,\gminib)\in S$,
let $H(\gminia,\gminib)$
denote the closure of the set
$\{\Re(\gminia-\gminib)\}$
in $\Xtilde(D)$.
Take a small neighbourhood $U_1$ of $P$
in $\pi^{-1}(O)$
such that
for any $(\gminia,\gminib)\in S$,
we have
$H(\gminia,\gminib)\cap U_1\neq\emptyset$
if and only if
$P\in H(\gminia,\gminib)$.
Let $U_1':=
U_1\setminus \bigcup_{(\gminia,\gminib)\in S}
 H(\gminia,\gminib)$.
We have
$\gminia<_P\gminib$
if and only if
$\gminia<_{P'}\gminib$
for any $P'\in U_1'$.
We have natural identifications
$\nbigl_P\simeq
 \nbigl_{P'}$
for any $P'\in U_1$.
Under the identifications,
we have
$\nbigf^P_{\gminia}
=\bigcap_{P'\in U_1'}
 \nbigf^{P'}_{\gminia}$.
So, if $\nbigf^{P'}_{\gminia}$
are defined over $K$
for any $P'\in U_1'$,
$\nbigf^P_{\gminia}$
is also defined over $K$.
For the points
$P'\in U_1'$,
the order $\leq_{P'}$ is totally ordered.
So, it is enough to prove that
$\nbigf^{P'}_{<\gminia}$
are defined over $K$ for any $\gminia\in\Irr(V)$
and for any $P'\in U_1'$.
But, it follows from the assumption of the lemma.
\hfill\qed

\subsection{The behaviour of the pre-$K$-Betti structure
by the nearby cycle functor and the maximal functor}

We set $X:=\Delta^n$
and $D:=\bigcup_{i=1}^{\ell}\{z_i=0\}$.
Let $V$ be a good meromorphic flat bundle on $(X,D)$
with a good $K$-structure.
For each $I\subset\ellsitabar$,
we set $I_{!i}:=I\cup\{i\}$
and $I_{\ast i}:=I\setminus \{i\}$.
The $\nbigd$-module
\[
 \Pi^{a,b}_{i\star}\bigl(V(!D(I))\bigr)
=\Bigl(V\otimes\gbigi^{a,b}_{z_i}
 \Bigr)(!D(I_{\star i})) 
\]
has the canonical pre-$K$-Betti structure,
where $\star=\ast,!$.
Hence, $\psi^{(a)}_i\bigl(V(!D(I))\bigr)$
and $\Xi^{(a)}_i\bigl(V(!D(I))\bigr)$ have
the induced pre-$K$-Betti structures.

\begin{lem}
\label{lem;09.10.21.1}
The induced $K$ structure of 
$\psi^{(a)}_{i}(V)$ is good,
i.e.,
it is compatible with the Stokes filtrations.
The induced pre-$K$-Betti structure of
$\psi^{(a)}_i\bigl(V(!D(I))\bigr)$ is canonical
for each $I\subset\ellsitabar$.
\end{lem}
\pf
It is enough to consider the case $a=0$ and $i=1$.
We give a preparation.
We set
$\Pi^{-\infty,a}_{f\star}V:=
 \varinjlim_{N}\Pi^{-N,a}_{f\star}(V)$.
By Lemma \ref{lem;13.4.17.1},
we have the following commutative diagram:
\begin{equation}
 \label{eq;09.12.4.50}
 \begin{CD}
 \DR_X\Bigl(
 \Pi^{-\infty,0}_{1!}\bigl(
V(!D(I))\bigr) 
 \Bigr)
 @>>>
 \DR_X\Bigl(
 \Pi^{-\infty,0}_{1\ast}\bigl(
 V(!D(I))
 \bigr)
 \Bigr)\\
 @A{\simeq}AA @A{\simeq}AA \\
 \DR_X^{<D(I_{\ast 1})}\bigl(
 \Pi^{-\infty,0}_{1!}V
 \bigr)
 @>>>
 \DR_X^{<D(I_{\ast 1})}\bigl(
 \Pi^{-\infty,0}_{1\ast}V
 \bigr)
 \\
 @A{\simeq}AA @A{\simeq}AA\\
 \DR_X^{<D(I_{!1})}\bigl(
 V\otimes\gbigi^{-\infty,0}_{z_1}
 \bigr)
 @>>>
 \DR_X^{<D(I_{\ast 1})}\bigl(
 V\otimes\gbigi^{-\infty,0}_{z_1}
 \bigr)
 \end{CD}
\end{equation}
By the upper square,
the induced $K$-structure
of $\DR_X\psi^{(0)}_1\bigl(V(!D(I))\bigr)$
can be identified with
the $K$-structure of the following:
\begin{equation}
\label{eq;09.12.4.41}
 \DR^{<D(I_{\ast 1})}_X\psi^{(0)}_{1}(V)
\simeq 
 \Cone
 \Bigl(
 \DR_X^{<D(I_{\ast 1})}\bigl(
 \Pi^{-\infty,0}_{1!}V
 \bigr)
\lrarr
 \DR^{<D(I_{\ast 1})}_{X}\bigl(
 \Pi^{-\infty,0}_{1\ast}V
 \bigr)
 \Bigr)
\end{equation}
We set $D':=\bigcup_{i=2}^{\ell}D_i$.
Let $\pi_1:\Xtilde(D')\lrarr X$
be the real blow up.
We obtain (\ref{eq;09.12.4.41})
as the push-forward of 
the following on $\Xtilde(D')$:
\begin{multline}
 \label{eq;09.12.4.20}
  \DR_{\Xtilde(D')}^{<D(I_{\ast 1})
  \leq D(\ellsitabar-I_{!1})}\psi^{(0)}_{1}(V)
\simeq \\
 \Cone
 \Bigl(
 \DR^{<D(I_{\ast 1})
 \leq D(\ellsitabar-I_{!1})}_{\Xtilde(D')}
 \bigl(
 \Pi^{-\infty,0}_{1!}V
 \bigr)
\lrarr
 \DR_{\Xtilde(D')}
 ^{<D(I_{\ast 1})\leq D(\ellsitabar-I_{!1})}
 \bigl(\Pi^{-\infty,0}_{1\ast}V\bigr)
\Bigr) 
\end{multline}
We prepare some commutative diagram
in a general setting.
For any holonomic $\nbigd_X$-module $\nbigm$,
we put 
\[
\DR_{\Xtilde(D')}
 ^{<D(I_{!1})\leq D(\ellsitabar-I_{!1})}\nbigm
 :=
 \Tot
 \Omega_{\Xtilde(D')}^{\bullet,\bullet,
 <D(I_{!1})\leq D(\ellsitabar-I_{!1})}
 \otimes_{\pi_1^{-1}\nbigo_X}
 \pi_1^{-1}\nbigm[\dim X],
\]
\[
\DR_{\Xtilde(D')}
 ^{<D(I_{\ast 1})\leq D(\ellsitabar-I_{\ast 1})}\nbigm
 :=
\Tot
 \Omega_{\Xtilde(D')}^{\bullet,\bullet,<D(I_{\ast 1})\leq D(\ellsitabar-I_{!1})}
 (\ast D_1)
 \otimes_{\pi_1^{-1}\nbigo_X}
 \pi_1^{-1}\nbigm[\dim X].
\]
We have the following commutative diagram:
\[
 \begin{CD}
 \DR^{<D(I_{\ast 1})
 \leq D(\ellsitabar-I_{!1})}_{\Xtilde(D')}\nbigm(!D_1)
@>>>
 \DR_{\Xtilde(D')}
 ^{<D(I_{\ast 1})\leq D(\ellsitabar-I_{!1})}
 \nbigm(\ast D_1)\\
 @AAA @A{=}AA \\
 \DR_{\Xtilde(D')}
 ^{<D(I_{!1})\leq D(\ellsitabar-I_{!1})}\nbigm
@>>>
 \DR_{\Xtilde(D')}
 ^{<D(I_{\ast 1})\leq D(\ellsitabar-I_{\ast 1})}\nbigm
 \end{CD}
\]
If $\nbigm$ is a good meromorphic flat bundle,
the left vertical arrow is also a quasi-isomorphism,
which follows from Lemma \ref{lem;13.4.27.1}.

Let $\rho:\Xtilde(D)\lrarr\Xtilde(D')$ be the induced map.
We have the following natural commutative diagram,
where the vertical arrows are quasi-isomorphisms
by Theorem \ref{thm;09.12.4.30}:
\[
 \begin{CD}
 \DR^{<D(I_{! 1})\leq D(\ellsitabar-I_{!1})}
 _{\Xtilde(D')}
 \nbigm
 @>>>
 \DR^{<D(I_{\ast 1})\leq D(\ellsitabar-I_{\ast 1})}
 _{\Xtilde(D')}
 \nbigm \\
 @V{\simeq}VV @V{\simeq}VV \\
 \rho_{\ast}
  \DR^{<D(I_{! 1})\leq D(\ellsitabar-I_{!1})}
 _{\Xtilde(D)}
 \nbigm
@>>>
  \rho_{\ast}
  \DR^{<D(I_{\ast 1})\leq D(\ellsitabar-I_{\ast 1})}
 _{\Xtilde(D)}
 \nbigm
 \end{CD}
\]
Thus, we obtain the following commutative diagram,
in which the vertical arrows are quasi-isomorphisms:
\begin{equation}
 \label{eq;09.12.4.21}
\begin{CD}
 \DR^{<D(I_{\ast 1})
 \leq D(\ellsitabar-I_{!1})}_{\Xtilde(D')}
 \bigl(
 \Pi^{-\infty,0}_{1!}V
 \bigr)
@>>>
  \DR_{\Xtilde(D')}^{
 <D(I_{\ast 1})\leq D(\ellsitabar-I_{!1})}
 \bigl(\Pi^{-\infty,0}_{1\ast}V\bigr)
 \\
 @A{\simeq}AA @A{\simeq}AA \\
 \rho_{\ast}\DR_{\Xtilde(D)}^{
 <D(I_{!1})\leq D(\ellsitabar-I_{! 1})}
 \bigl(V\otimes\gbigi^{-\infty,0}_{z_1}
 \bigr)
@>>>
 \rho_{\ast}\DR_{\Xtilde(D)}^{
 <D(I_{\ast 1})\leq D(\ellsitabar-I_{\ast 1})}
 \bigl(V\otimes\gbigi^{-\infty,0}_{z_1}
 \bigr)
\end{CD}
\end{equation}
Because 
$\DR_{\Xtilde(D)}^{
 <D(I_{!1})\leq D(\ellsitabar-I_{! 1})}
 \bigl(V\otimes\gbigi^{-\infty,0}_{z_1}
 \bigr)$
and
$\DR_{\Xtilde(D)}^{
 <D(I_{\ast 1})\leq D(\ellsitabar-I_{\ast 1})}
 \bigl(V\otimes\gbigi^{-\infty,0}_{z_1}
 \bigr)$ are equipped with $K$-structures
compatible with the morphism,
we obtain a $K$-structure of
$\DR^{<D(I_{\ast 1})\leq 
 D(\ellsitabar-I_{!1})}_{\Xtilde(D')}
 \psi^{(0)}_{1}(V)$
from (\ref{eq;09.12.4.20})
and (\ref{eq;09.12.4.21}).
The lower square in (\ref{eq;09.12.4.50})
is obtained as the push-forward of 
(\ref{eq;09.12.4.21}).
Hence, the $K$-structure of 
$\DR_X\psi^{(0)}_1\bigl(V(!D(I))\bigr)$
is obtained as the push-forward of
the $K$-structure of
$\DR^{<D(I_{\ast 1})\leq 
 D(\ellsitabar-I_{!1})}_{\Xtilde(D')}
 \psi^{(0)}_{1}(V)$.

Let us consider the case $I=\{1,\ldots,\ell\}$.
By the above consideration,
we obtain that 
$\nbigf^{P}_{< 0}$
is compatible with the $K$-structure,
where $\nbigf^P$ denotes
the Stokes filtration of $\psi^{(0)}_{1}(V)$
at each point $P\in \pi_1^{-1}(\del D_1)$.
By considering the tensor product
with meromorphic flat bundles
with rank one,
we can deduce that 
$\nbigf^P$ is defined over $K$,
as in Lemma \ref{lem;13.4.20.250}.
Since the pre-$K$-Betti structure of
$\psi^{(0)}_{1}\bigl(V(!D(I))\bigr)$ comes from
the $K$-structure of
$\DR_{\Xtilde(D')}
 ^{<D(I_{\ast 1})\leq D(\ellsitabar-I_{!1})}
 \psi^{(0)}_{1}(V)$,
it is canonical.
\hfill\qed

%% file: 6.2.tex
\subsection{Definition}

Let $X=\Delta^n$ and $D=\bigcup_{i=1}^{\ell}\{z_i=0\}$.
Set $\ellsitabar:=\{1,\ldots,\ell\}$.
Let $\nbigm$ be a good holonomic
$\nbigd$-module on $(X,D)$.
\begin{df}
We say that $\nbigm$ has a good $K$-structure
if 
(i) for each $I\subset\ellsitabar$,
$\phi_I(\nbigm)\bigl(\ast D(I^c)\bigr)$
is equipped with a good $K$-structure
(put $\phi_{\emptyset}(\nbigm):=\nbigm$),
(ii) for $i\not\in I$,
the induced morphisms
\begin{equation}
 \label{eq;13.4.26.10}
 \psi^{(1)}_{i}\Bigl(
 \phi_I(\nbigm)\bigl(\ast D(I^c)\bigr)
 \Bigr)
\lrarr
 \Bigl(
 \phi_{i}\phi_I(\nbigm)
 \Bigr)
 \bigl(
 \ast D(I_{! i}^c)
 \bigr)
\lrarr
 \psi^{(0)}_{i}\Bigl(
 \phi_I(\nbigm)\bigl(\ast D(I^c)\bigr)
 \Bigr)
\end{equation}
are compatible with the $K$-structures,
where $I_{! i}:=I\sqcup\{i\}$.
\hfill\qed
\end{df}
\index{good $K$-structure}
Morphisms of  good holonomic $\nbigd$-modules
with a good $K$-structure
$f:\nbigm_1\lrarr\nbigm_2$
are morphisms of
$\nbigd$-modules such that
$\phi_I(f)$ are compatible with
$K$-structures for any $I\subset\ellsitabar$.

Let $\Hol^{\good}(X,D,K)$ denote the category of 
good holonomic $\nbigd_X$-modules
with a good $K$-structure on $(X,D)$.
\index{category $\Hol^{\good}(X,D,K)$}

\begin{lem}
Let $f:\nbigm_1\lrarr \nbigm_2$ be a morphism
in $\Hol^{\good}(X,D,K)$.
Then, the $\nbigd$-modules 
$\Ker(f)$, $\Image(f)$ and $\Cok(f)$
are naturally objects in $\Hol^{\good}(X,D,K)$.
\end{lem}
\pf
It follows from Lemma \ref{lem;13.5.18.1}.
(See also the reconstruction  of 
a good holonomic $\nbigd$-module $\nbigm$
from $\phi_I^{(\veczero)}(\nbigm)$
in \S\ref{subsection;13.4.27.21}.)
\hfill\qed

\subsection{Cells}
\label{subsection;09.10.28.50}

Let $V$ be any good meromorphic flat bundle
on $X$ with a good $K$-structure.
Let us observe that 
we have natural objects 
in $\Hol^{\good}(X,D,K)$
associated to $V$.
\begin{lem}
\label{lem;13.4.26.200}
Let $D^{(1)}$ be a hypersurface of $X$
contained in $D$.
\begin{itemize}
\item
We can naturally regard $V(!D^{(1)})$ 
as an object in $\Hol^{\good}(X,D,K)$.
\item
Suppose that we are given an object $\nbigm$
in $\Hol^{\good}(X,D,K)$
such that
(i) the underlying $\nbigd_X$-module is
 isomorphic to $V(!D^{(1)})$,
(ii) the $K$-structure on 
 $X\setminus D$ is equal 
 to that of $V(!D^{(1)})$ under the isomorphism.
Then, $\nbigm$ is isomorphic to
$V(!D^{(1)})$
in $\Hol^{\good}(X,D,K)$.
\end{itemize}
\end{lem}
\pf
We have $I\subset\ellsitabar$
such that $D^{(1)}=D(I)$.
We have a natural isomorphism
\[
 \phi^{(\veczero_J)}_J\bigl(V\bigl(!D(I)\bigr)\bigr)\bigl(\ast D(J^c)\bigr)
\simeq 
 \psi^{(\vecdelta_{J\cap I})}_{J\cap I}
 \psi^{(\veczero_{J\setminus I})}_{J\setminus I}(V)
\]
for any $J\subset\ellsitabar$,
where $\vecdelta_{J\cap I}=(1,\ldots,1)\in\seisuu^{J\cap I}$
and $\veczero_{J\setminus I}=(0,\ldots,0)\in \seisuu^{J\setminus I}$.
They are equipped with good $K$-structures,
satisfying the compatibility condition
(\ref{eq;13.4.26.10}).
Via these $K$-structures,
we regard $V(!D(I))\in\Hol^{\good}(X,D,K)$.
Thus, we obtain the first claim.

Let us prove the second claim.
We are given the isomorphism of $\nbigd_X$-modules
$V(!D^{(1)})\simeq\nbigm$ under which
the $K$-structures on $X\setminus D$ are equal.
Suppose that we have already known that
$\phi^{(\veczero)}_I(V(!D^{(1)}))
\simeq
 \phi^{(\veczero)}_I(\nbigm)$
preserves the $K$-structures.
Set $V_1:=V(!D^{(1)})$ and $V_2:=\nbigm$.
Because one of
$\psi^{(1)}_i\phi^{(\veczero)}_I(V_i)
\lrarr
 \phi^{(0)}_i\phi^{(\veczero)}_I(V_i)$
or 
$\psi^{(1)}_i\phi^{(\veczero)}_I(V_i)
\lrarr
 \phi^{(0)}_i\phi^{(\veczero)}_I(V_i)$
is an isomorphism
compatible with $K$-structures.
Hence, we obtain that
$\phi_i^{(0)}\phi^{\veczero}_I(V_1)
\lrarr
\phi_i^{(0)}\phi^{\veczero}_I(V_2)$
is also compatible with the $K$-structures.
\hfill\qed

\vspace{.1in}

More generally,
take $J\sqcup I\subset\ellsitabar$.
Let $V_J$ be a good meromorphic flat bundle
on $D_J$ with a good $K$-structure.
Then, we can naturally regard
$\iota_{\dagger}V_J\bigl(!D(I)\bigr)$
as an object in
$\Hol^{\good}(X,D,K)$.

\vspace{.1in}

Let $g$ be a meromorphic function on $(X,D)$
such that $g^{-1}(0)\subset D$.
Let $D=D_1\cup D_2$ be a decomposition
such that $D_1\supset g^{-1}(\infty)$
and $D_2\subset g^{-1}(0)$.
(Note that $D_i$ are not necessarily irreducible.)
Because 
$\Xi^{(0)}_g(V,\ast D_1)$ and $\psi^{(0)}_g(V,\ast D_1)$
are the kernel of
$\Bigl(
 V\otimes\gbigi_g^{-\infty,a}(!D_2)
 \Bigr)(\ast D_1)
\lrarr
V\otimes\gbigi_g^{-\infty,0}(\ast D)$
for $a=1,0$,
they are naturally objects
in $\Hol^{\good}(X,D,K)$.

\subsection{Some operations}

Let us observe that some operations on $\Hol(X)$
are naturally lifted on 
$\Hol^{\good}(X,D,K)$.
Let $\Forget$ denote the forgetful functor
from $\Hol^{\good}(X,D,K)$
to $\Hol(X)$.

\begin{lem}
We have a naturally defined dual functor
$\DDD$ on $\Hol^{\good}(X,D,K)$
such that
$\DDD\circ\Forget
=\Forget\circ\DDD$.
\end{lem}
\pf
Let $\nbigm\in\Hol^{\good}(X,D,K)$.
For each $I\subset\{1,\ldots,\ell\}$,
$\phi_I^{(\veca)}(\DDD\nbigm)(\ast D(I^c))
\simeq
 \DDD \phi_I^{(-\veca-\vecdelta)}(\nbigm)(\ast D(I^c))$
has an induced $K$-structure.
For $I_0:=I\sqcup\{i\}$,
the morphisms
\[
\psi_i^{(1)}
 \phi^{(\veczero)}_I(\DDD\nbigm)(\ast D(I_0^c))
\lrarr
 \phi_i^{(0)}
 \phi^{(\veczero)}_I(\DDD\nbigm)(\ast D(I_0^c))
\lrarr
 \psi_i^{(0)}
 \phi^{(\veczero)}_I(\DDD\nbigm)(\ast D(I_0^c))
\]
are obtained as the dual of 
$\psi_i^{(0)}
 \phi^{(-\vecdelta)}_I(\nbigm)(\ast D(I_0^c))
\lrarr
 \phi_i^{(-1)}
 \phi^{(-\vecdelta)}_I(\nbigm)(\ast D(I_0^c))
\lrarr
 \psi_i^{(-1)}
 \phi^{(-\vecdelta)}_I(\nbigm)(\ast D(I_0^c))$,
they are compatible with the $K$-structure.
Hence,
they give a good $K$-structure on $\DDD\nbigm$.
The construction gives a contravariant functor
$\DDD$ on $\Hol^{\good}(X,D,K)$.
\hfill\qed

\begin{lem}
\label{lem;13.4.26.40}
Let $D^{(1)}\subset D$ be a hypersurface of $X$.
We have a functor
$\Phi_{\ast D^{(1)}}:
 \Hol^{\good}(X,D,K)
\lrarr
 \Hol^{good}(X,D,K)$
such that
\[
 \Forget\circ
 \Phi_{\ast D^{(1)}}
 (\nbigm)
=\Forget(\nbigm)(\ast D^{(1)})
\]
for any $\nbigm$ in $\Hol^{\good}(X,D,K)$.
We also have a natural transformation
$\nbigm\lrarr \Phi_{\ast D^{(1)}}(\nbigm)$.
Such a functor is unique.
\end{lem}
\pf
First, let us observe the uniqueness.
Let $\nbigm\in\Hol^{\good}(X,D,K)$.
We have $I\subset\ellsitabar$
such that $D^{(1)}=D(I)$.
For any $J\subset\ellsitabar$,
the following isomorphism is compatible
with the $K$-structure.
\[
\begin{CD}
 \phi^{(\veczero)}_{J\setminus I}(\nbigm)
 \Bigl(\ast D((J\setminus I)^c)\Bigr)
 @>{\alpha}>>
 \phi^{(\veczero)}_{J\setminus I}(\Phi_{\ast D^{(1)}}\nbigm)
 \Bigl(\ast D((J\setminus I)^c)\Bigr) 
\end{CD}
\]
The following induced isomorphism is
compatible with the $K$-structure:
\[
\begin{CD}
 \psi^{(\veczero)}_{J\cap I}
 \phi^{(\veczero)}_{J\setminus I}(\nbigm)
 \bigl(D(J^c)\bigr)
@>{\psi^{(\veczero)}_{J\cap I}(\alpha)}>>
 \psi^{(\veczero)}_{J\cap I}
 \phi^{(\veczero)}_{J\setminus I}(\Phi_{\ast D^{(1)}}\nbigm)
 \bigl(D(J^c)\bigr)
\end{CD}
\]
Note that the following natural morphism
is an isomorphism:
\[
\begin{CD}
 \phi^{(\veczero)}_J(\Phi_{\ast D^{(1)}}\nbigm)
 \bigl(\ast D(J^c)\bigr)
@>>>
 \psi^{\veczero}_{J\cap I}
 \phi^{(\veczero)}_{J\setminus I}(\Phi_{\ast D^{(1)}}\nbigm)
  \bigl(\ast D(J^c)\bigr)
\end{CD}
\]
It is compatible with the $K$-structure
by the condition for $\Phi_{\ast D^{(1)}}\nbigm$.
Hence, 
the good $K$-structure of
\[
 \phi^{(\veczero)}_{J\setminus I}(\nbigm)
 \Bigl(\ast D((J\setminus I)^c)\Bigr)
\]
uniquely determines  
the $K$-structure of
$\phi^{(\veczero)}_J(\Phi_{\ast D^{(1)}}\nbigm)
 \bigl(\ast D(J^c)\bigr)$.
It means the uniqueness of
$\Phi_{\ast D^{(1)}}$.

As for the existence of $\Phi_{\ast D^{(1)}}$,
it is enough to consider the case $I=\{1\}$.
If $i\in J$,
we have
$\phi^{(\veczero)}_J\bigl(\nbigm(\ast D^{(1)})\bigr)
\simeq
 \psi_1^{(0)}\phi^{(\veczero)}_{J\setminus \{1\}}(\nbigm)$.
If $i\not\in J$,
we have
$\phi^{(\veczero)}_J\bigl(\nbigm(\ast D^{(1)})\bigr)
\simeq
 \phi^{(\veczero)}_{J}(\nbigm)(\ast D^{(1)})$.
The induced $K$-structures on 
$\phi_J^{(\veczero)}\bigl(\nbigm(\ast D^{(1)})\!\bigr)
 (\ast D(J^c))$
give a good $K$-structure of 
$\nbigm(\ast D^{(1)})$,
for which
the natural morphism
$\nbigm\lrarr\nbigm(\ast D^{(1)})$
is a morphism in $\Hol^{\good}(X,D,K)$.
\hfill\qed

\begin{lem}
\label{lem;13.4.26.41}
\mbox{{}}
\begin{itemize}
\item
For any hypersurface $D^{(1)}$ of $X$ contained in $D$,
we have a unique functor
$\Phi_{! D^{(1)}}:
 \Hol^{\good}(X,D,K)
\lrarr
 \Hol^{good}(X,D,K)$
such that
\[
 \Forget\circ
 \Phi_{! D^{(1)}}
 (\nbigm)
=\Forget(\nbigm)(! D^{(1)})
\]
for any $\nbigm$ in $\Hol^{\good}(X,D,K)$,
with a natural transformation
$\Phi_{! D^{(1)}}\lrarr \id$.
\item
We have
$\Phi_{\star D^{(1)}}\circ
 \Phi_{\star D^{(2)}}
=\Phi_{\star (D^{(1)}\cup D^{(2)})}$.
\item
If $\dim(D^{(1)}\cap D^{(2)})<n-1$,
then
$\Phi_{!D^{(1)}}\circ\Phi_{\ast D^{(2)}}
=\Phi_{\ast D^{(2)}}\circ\Phi_{!D^{(1)}}$.
\end{itemize}
\end{lem}
\pf
The first claim follows
from Lemma \ref{lem;13.4.26.40}
as the dual.
The second claim follows from
the uniqueness.
For $\nbigm\in\Hol^{\good}(X,D,K)$,
the underlying $\nbigd_X$-modules of
$\Phi_{! D^{(1)}}\circ\Phi_{\ast D^{(2)}}(\nbigm)$
and 
$\Phi_{\ast D^{(2)}}\circ\Phi_{! D^{(1)}}(\nbigm)$
are 
$\nbigm(!D^{(1)}\ast D^{(2)})
=\nbigm(\ast D^{(2)}!D^{(1)})$.
We have the natural morphisms
$\Phi_{! D^{(1)}}(\nbigm)
\lrarr
 \Phi_{! D^{(1)}}\circ\Phi_{\ast D^{(2)}}(\nbigm)$
and
$\Phi_{! D^{(1)}}(\nbigm)
\lrarr
\Phi_{\ast D^{(2)}}\circ\Phi_{! D^{(1)}}(\nbigm)$
in $\Hol^{\good}(X,D,K)$.
Then, by the argument for the uniqueness 
in the proof of Lemma \ref{lem;13.4.26.40},
we obtain that the $K$-structures are the same.
\hfill\qed

\vspace{.1in}
We denote $\Phi_{\star D^{(1)}}(\nbigm)$
by $\nbigm(\star D^{(1)})$
for $\star=\ast,!$.

%% file: 6.3.tex
\subsection{Statements}

Let $X=\Delta^n$ and $D=\bigcup_{i=1}^{\ell}\{z_i=0\}$.
Let $\Hol^{\pre}(X,K)$
denote the category of pre-$K$-holonomic
$\nbigd_X$-modules.

\begin{prop}
\label{prop;13.4.26.100}
We have a naturally defined exact fully faithful functor 
$\Upsilon:\Hol(X,D,K)\lrarr
 \Hol^{\pre}(X,K)$
over $\Hol(X)$.
We have
$\Upsilon\circ\DDD=\DDD\circ\Upsilon$.
The essential image of $\Upsilon$
is independent of the choice of 
a holomorphic coordinate system.
\end{prop}
\index{functor $\Upsilon$}

\begin{df}
\label{df;09.10.22.10}
Any object in the essential image of  $\Upsilon$
is called a good pre-$K$-holonomic $\nbigd$-module
on $(X,D)$.
The pre-$K$-Betti structure is called
a good pre-$K$-Betti structure.
(The definition will be globalized
 in Definition {\rm\ref{df;13.4.27.22}} below.)
\hfill\qed
\end{df}
\index{good pre-$K$-Betti structure}
\index{good pre-$K$-holonomic $\nbigd$-module}

Let $V$ be a good meromorphic flat bundle
on $(X,D)$ with a good $K$-structure.
Let $D^{(1)}\subset D$ be a hypersurface of $X$.

\begin{prop}
\label{prop;13.4.26.120}
The canonical pre-$K$-Betti structure
of $V(!D^{(1)})$
is associated to the good $K$-structure
of $V(!D^{(1)})$
by $\Upsilon$.
\end{prop}

We shall construct the functor
in \S\ref{subsection;09.10.22.2}--\S\ref{subsection;13.4.26.20}.
We shall prove the full faithfulness
in \S\ref{subsection;13.4.26.21}.
The independence from the coordinate system
will be proved in \S\ref{subsection;13.4.26.101}.
Proposition \ref{prop;13.4.26.120} will be proved
in \S\ref{subsection;13.4.26.130}.

\subsection{Some consequences}

Before going to the proof of Proposition \ref{prop;13.4.26.100},
we give some consequences.
The full faithfulness and the independence
on the coordinate system in Proposition \ref{prop;13.4.26.100} 
ensure that 
we can globalize 
the notion of good pre-$K$-holonomic $\nbigd$-modules
in Definition \ref{df;09.10.22.10}.

\begin{df}
\label{df;13.4.27.22}
Let $Y$ be any complex manifold
with a normal crossing hypersurface $D_Y$.
Let $\nbigm$ be a good holonomic $\nbigd$-module
on $(Y,D_Y)$ with a pre-$K$-Betti structure $\nbigf$.
It is called a good pre-$K$-holonomic $\nbigd$-module
if its restriction to any holomorphic coordinate neighbourhood
is a good pre-$K$-holonomic $\nbigd$-module.
In that case,
$\nbigf$ is called a good pre-$K$-Betti structure.
\hfill\qed
\end{df}
\index{good pre-$K$-Betti structure}
\index{good pre-$K$-holonomic $\nbigd$-module}
The category of good pre-$K$-holonomic $\nbigd$-modules
on $(Y,D_Y)$ is not abelian (see \S\ref{subsection;13.5.5.1}).
If we would like to work on abelian categories,
for example,
the full subcategory of 
$\vecnbigi$-good pre-$K$-holonomic $\nbigd$-modules
is abelian,
where $\vecnbigi$ is any good system of ramified irregular values
on $(Y,D_Y)$.

\vspace{.1in}

Let $Y$ be any complex manifold
with a normal crossing hypersurface $D$.
Let $V$ be a good meromorphic flat bundle on $(Y,D)$
with a good $K$-structure.
Let $g$ be any meromorphic function on $(Y,D)$
such that it is invertible on $Y\setminus D$.
We take a hypersurface $D^{(1)}\subset D$
such that $g^{-1}(0)\subset D^{(1)}$.
We obtain a good meromorphic flat bundle
$V\otimes\gbigi_g^{a,b}$
with a good $K$-structure on $(Y,D)$.
It induces pre-$K$-holonomic $\nbigd$-modules
$\Pi^{a,b}_{g\star}(V)(\ast D^{(1)})$,
$\Xi^{(a)}_g(V,\ast D^{(1)})$
and 
$\psi^{(a)}_g(V,\ast D^{(1)})$
with the canonical pre-$K$-Betti structures.
We obtain the following proposition
from Proposition \ref{prop;13.4.26.120}.

\begin{prop}
\label{prop;13.4.27.30}
The holonomic $\nbigd_Y$-modules
$\Pi^{a,b}_{g\star}(V)(\ast D^{(1)})$,
$\Xi^{(a)}_g(V,\ast D^{(1)})$,
$\psi^{(a)}_g(V,\ast D^{(1)})$
and $\phi^{(a)}_g(V,\ast D^{(1)})$
are naturally
good pre-$K$-holonomic $\nbigd$-modules
on $(Y,D)$.
\hfill\qed
\end{prop}

The claims for 
$\psi^{(a)}_g(V,\ast D^{(1)})$
and 
$\phi^{(a)}_g(V,\ast D^{(1)})$
will be particularly useful.

\subsection{Induced pre-$K$-Betti structures
of $\Xi^{(\veca)}_I\psi^{\vecb}_J(\iota_{\dagger}V_I)$}
\label{subsection;09.10.22.2}

In the following,
we shall prove Proposition \ref{prop;13.4.26.100}
and Proposition \ref{prop;13.4.26.120}.

Let $K\sqcup J\sqcup I= L\subset \ellsitabar$.
Let $V_I$ be an $\nbigi$-good meromorphic flat bundle
on $(D_I,\del D_I)$.
Let $\iota:D_I\lrarr X$.
For a map $f:K\sqcup J\lrarr\{0,1\}$,
we set
$K_0(f):=f^{-1}(0)\cap K$.
We put
\[
 \nbigc_f(J,K,\iota_{\dagger}V_I):=
\Bigl(
 \iota_{\dagger}V_I\otimes
 \bigotimes_{k\in K_0(f)}
 \gbigi_{z_k}^{-\infty,1}
 \otimes
 \bigotimes_{k\not\in K_0(f)}
 \gbigi_{z_k}^{-\infty,0}
\Bigr)\bigl(
 !D(f^{-1}(0))
 \bigr).
\]
Let $\veczero$
denote the constant map valued in $\{0\}$.
Let $\vecdelta_i$  
denote the map such that
$\vecdelta_i(j)=0$ $(j\neq i)$
and $\vecdelta_i(i)=1$.
We can identify 
$\Xi^{(\veczero)}_K\psi^{(\veczero)}_J
 \bigl(\iota_{\dagger}V_I\bigr)$
as the kernel of the following morphism:
\begin{equation}
 \label{eq;09.10.22.1}
 \nbigc_{\veczero}
 \bigl(J,K,\iota_{\dagger}V_I\bigr)
\lrarr
 \bigoplus_{i\in K\sqcup J}
 \nbigc_{\vecdelta_i}
 \bigl(J,K,\iota_{\dagger}V_I\bigr)
\end{equation}
If $V_I$ has a good $K$-structure,
we obtain a pre-$K$-Betti structure of
$\Xi^{(\veczero)}_K\psi^{(\veczero)}_J(\iota_{\dagger}V_I)$
by (\ref{eq;09.10.22.1}).
By taking the tensor product
with $\gbigi^{a,a+1}$ appropriately,
we also obtain an induced pre-$K$-Betti structure of 
$\Xi^{(\veca)}_K\psi^{(\vecb)}_J(\iota_{\dagger}V_I)$.

\begin{lem}
\label{lem;09.12.4.100}
The following morphisms
are compatible with the pre-$K$-Betti structures:
\[
 \Xi^{(\veca)}_{K}
 \psi^{(\vecb)}_{J}
 \psi^{(1)}_i
 \bigl(\iota_{\dagger}V_I\bigr)
\lrarr
 \Xi^{(\veca)}_{K}
\psi^{(\vecb)}_{J}
 \Xi^{(0)}_i
 \bigl(\iota_{\dagger}V_I\bigr)
\lrarr
 \Xi^{(\veca)}_{K}
 \psi^{(\vecb)}_{J}
 \psi^{(0)}\bigl(\iota_{\dagger}V_I\bigr)
\]
\end{lem}
\pf
It is clear by construction.
\hfill\qed

\vspace{.1in}
Recall that we have the naturally induced
good $K$-structure on 
$\psi_i^{(0)}\bigl(\iota_{\dagger}V_I\bigr)$
for $i\not\in I$
(Lemma \ref{lem;09.10.21.1}).
\begin{lem}
\label{lem;09.12.4.101}
For any $i\not\in L$,
the natural isomorphism
\[
 \Xi^{(\veczero)}_{K}
 \psi^{(\veczero)}_{Ji}\bigl(\iota_{\dagger}V_I\bigr)
\simeq
 \Xi^{(\veczero)}_K\psi^{(\veczero)}_J\Bigl(
 \psi^{(0)}_i\bigl(\iota_{\dagger}V_I\bigr)
 \Bigr) 
\]
is compatible with the induced $K$-structures.
\end{lem}
\pf
Both the $K$-structures are obtained as the kernel of
the morphism (\ref{eq;09.10.22.1})
for $(Ji,K)$.
\hfill\qed

\subsection{$\ell$-squares of complexes}
\label{subsection;09.10.4.300}

For small categories $A_i$ $(i=1,\ldots,\ell)$,
let $\prod_{i=1}^{\ell} A_i$  
denote their product,
i.e.,
the category whose objects 
and morphisms are given by
$\ob\Bigl(
 \prod_{i=1}^{\ell}A_i
 \Bigr)
=\prod_{i=1}^{\ell}\ob(A_i)$
and
$\Mor(\veca,\vecb)=\prod\Mor(a_i,b_i)$.
Let $\Gamma$ be a small category
given by the following commutative diagram:
\[
 \begin{CD}
 (0,0) @>{a}>> (0,1)\\
 @V{b}VV @V{c}VV \\
 (1,0) @>{d}>> (1,1)
 \end{CD}
\quad\quad
 c\circ a=d\circ b
\]

Let $A$ be an abelian category.
Let $C(A)$ be the category of complexes
in $A$.
A square in $C(A)$ is a functor
$F:\Gamma\lrarr C(A)$.
For a given square $F$,
let $H(F)$ be the total complex
of the following double complex:
\[
\begin{CD}
 F(0,0)[1] @>{F(a)+F(b)}>>
 F(0,1)\oplus F(1,0) @>{F(c)-F(d)}>>
 F(1,1)[-1]
\end{CD}
\]
An $\ell$-square in $C(A)$ is 
a functor $F:\Gamma^{\ell}\lrarr C(A)$.
Let $\pi_i:\Gamma^{\ell}\lrarr\Gamma^{\ell-1}$
be the projection forgetting the $i$-th component.
For a given $\ell$-square $F$,
we obtain an $(\ell-1)$-square $\pi_{i\ast}F$ by 
$\pi_{i\ast}F(\veca)
=H\bigl( F_{|\pi_i^{-1}(\veca)} \bigr)$.

\begin{lem}
For $i<j$,
we have an isomorphism
$\pi_{i\ast}\pi_{j\ast}F
\simeq\pi_{j-1\ast}\pi_{i\ast}F$.
\end{lem}
\pf
It is enough to consider the case
$\ell=2$, $(i,j)=(1,2)$.
The $i$-th terms of the both complexes are given by
\[
 \bigoplus_{a_1+a_2+b_1+b_2=i-2}
 F(a_1,a_2,b_1,b_2).
\]
The multiplication of $-1$ on
$F(0,0,0,0)\oplus F(1,1,0,0)\oplus
 F(0,0,1,1)\oplus F(1,1,1,1)$
interpolates the differentials
for $\pi_{i\ast}\pi_{j\ast}F$ and
$\pi_{j-1\ast}\pi_{i\ast}F$.
\hfill\qed

\vspace{.1in}

More generally,
for any subset $I\subset\ellsitabar$,
$I$-square in $C(A)$ is a functor
$\Gamma^I\lrarr C(A)$.
For the naturally defined projection
$\pi_I:\Gamma^{\ell}\lrarr\Gamma^{I}$,
we take 
$I=I_0\subset I_1\subset\cdots\subset I_m
 =\ellsitabar$,
which induces the factorization
$\pi_I=\pi^{(1)}\circ\pi^{(2)}\circ
 \cdots\circ\pi^{(m)}$,
where 
$\pi^{(i)}:\Gamma^{I_i}
 \lrarr\Gamma^{I_{i-1}}$.
Then, we obtain an $I$-square
$\pi_{I\ast}F:=\pi^{(1)}_{\ast}\circ
 \cdots \circ\pi^{(m)}_{\ast}F$.
It is well defined up to isomorphisms as above.

\subsection{A construction of the functor $\Upsilon$}
\label{subsection;13.4.26.20}

Let $m$ be any positive integer.
Let $\nbigi\subset M(X^{(m)},D^{(m)})/H(X^{(m)})$
be any good set of ramified irregular values
as in \S\ref{subsection;09.10.19.5}.
Let $\nbigm$ be 
any $\nbigi$-good holonomic $\nbigd$-module
on $(X,D)$.

Let $H\subset\ellsitabar$.
Let us construct an $H$-square
in the category of 
$\nbigi$-good holonomic $\nbigd$-modules
on $(X,D)$.
For $(\veci,\vecj)=
 \bigl(
 (i_k,j_k)
 \,\big|\,k\in H
 \bigr)\in\ob\Gamma^H$,
we have the following subsets of $H$:
\[
 I(\veci,\vecj)=\bigl\{
 k\,\big|\,
 (i_k,j_k)=(0,1)
 \bigr\},
\quad
\quad
 K(\veci,\vecj)=\bigl\{
 k\,\big|\,
 (i_k,j_k)=(1,0)
 \bigr\},
\]
\[
 J_0(\veci,\vecj)=\bigl\{
 k\,\big|\,
 (i_k,j_k)=(0,0)
 \bigr\},
\quad\quad
  J_1(\veci,\vecj)=\bigl\{
 k\,\big|\,
 (i_k,j_k)=(1,1)
 \bigr\}.
\]
Then, we put
$\nbigq^H(\nbigm,\veci,\vecj):=
 \Xi^{(\veczero)}_{I(\veci,\vecj)}
 \psi^{(\vecdelta_{J_0(\veci,\vecj)})}_{J_0(\veci,\vecj)}
 \psi^{(\veczero)}_{J_1(\veci,\vecj)}
 \phi^{(\veczero)}_{K(\veci,\vecj)}\nbigm$.
For $k_0\not\in H$,
we have the following naturally induced diagram:
\begin{equation}
 \label{eq;09.9.16.1}
\begin{CD}
 \psi^{(1)}_{k_0}\Xi^{(\veczero)}_I
 \psi^{(\vecdelta_{J_0})}_{J_0}
 \psi^{(\veczero)}_{J_1}
 \phi^{(\veczero)}_K\nbigm
 @>>>
 \Xi^{(0)}_{k_0}
 \Xi^{(\veczero)}_I
 \psi^{(\vecdelta_{J_0})}_{J_0}
 \psi^{(\veczero)}_{J_1}
 \phi^{(\veczero)}_K\nbigm\\
 @VVV @VVV \\
 \phi^{(0)}_{k_0}
 \Xi_I^{(\veczero)}
 \psi^{(\vecdelta_{J_0})}_{J_0}
 \psi^{(\veczero)}_{J_1}
 \phi_K^{(\veczero)}\nbigm
 @>>>
 \psi^{(0)}_{k_0}
 \Xi^{(\veczero)}_I
 \psi^{(\vecdelta_{J_0})}_{J_0}
 \psi^{(\veczero)}_{J_1}
 \phi^{(\veczero)}_K\nbigm
\end{CD}
\end{equation}
For each decomposition
$H=\{h\}\cup (H-\{h\})$,
we have a similar diagram.
Thus, we obtain an $H$-square $\nbigq^H(\nbigm)$
of good holonomic $\nbigd$-modules.
The cohomology of the complex
associated to (\ref{eq;09.9.16.1})
is naturally isomorphic to
$\Xi_I^{(\veczero)}
 \psi^{(\vecdelta_{J_0})}_{J_0}
 \psi^{(\veczero)}_{J_1}
 \phi^{(\veczero)}_K\nbigm$.
Hence, we have a natural quasi-isomorphism
$\pi_{H\ast}\nbigq^{\ellsitabar}(\nbigm)
 \simeq\nbigq^H(\nbigm)$.
In particular, we have a natural quasi-isomorphism
$\pi_{\ellsitabar\ast}\nbigq^{\ellsitabar}(\nbigm)
 \simeq \nbigm$.

\vspace{.1in}

If $\nbigm$ has a good $K$-structure,
each $\nbigq^{\ellsitabar}(\nbigm,\veci,\vecj)$
is equipped with the pre-$K$-Betti structure
$\nbigf_{\nbigm}^{\ellsitabar}(\veci,\vecj)$
given as in \S\ref{subsection;09.10.22.2}.
\begin{lem}
The morphisms in {\rm(\ref{eq;09.9.16.1})}
are compatible with the induced pre-$K$-Betti
structures.
\end{lem}
\pf
The morphisms
{\small
\[
 \psi^{(1)}_{k_0}\Xi^{(\veczero)}_I
 \psi^{(\vecdelta_{J_0})}_{J_0}
 \psi^{(\veczero)}_{J_1}
 \phi^{(\veczero)}_K\nbigm
\lrarr
 \Xi^{(0)}_{k_0}
 \Xi^{(\veczero)}_I
 \psi^{(\vecdelta_{J_0})}_{J_0}
 \psi^{(\veczero)}_{J_1}
 \phi^{(\veczero)}_K
 \nbigm
\lrarr
 \psi^{(0)}_{k_0}
 \Xi^{(\veczero)}_I
 \psi^{(\vecdelta_{J_0})}_{J_0}
 \psi^{(\veczero)}_{J_1}
 \phi^{(\veczero)}_K\nbigm 
\]}
are compatible with the pre-$K$-Betti structures
by construction,
as remarked in Lemma \ref{lem;09.12.4.100}.
Let $K':=\ellsitabar-(K\sqcup{k_0})$.
By definition, the morphisms
\[
 \psi^{(1)}_{k_0}\phi^{(\veczero)}_K\nbigm(\ast D(K'))
 \lrarr
 \phi^{(0)}_{k_0}\phi^{(\veczero)}_K\nbigm(\ast D(K'))
 \lrarr
 \psi^{(0)}_{k_0}\phi^{(\veczero)}_K\nbigm(\ast D(K')) 
\]
are compatible with the $K$-structures.
We remark Lemma \ref{lem;09.12.4.101},
and then it follows that the morphisms
\[
  \psi^{(1)}_{k_0}\Xi^{(\veczero)}_I
 \psi^{(\vecdelta_{J_0})}_{J_0}
 \psi^{(\veczero)}_{J_1}
 \phi^{(\veczero)}_K\nbigm
\lrarr
 \phi^{(0)}_{k_0}
 \Xi^{(\veczero)}_I
 \psi^{(\vecdelta_{J_0})}_{J_0}
 \psi^{(\veczero)}_{J_1}
 \phi^{(\veczero)}_K
 \nbigm
\lrarr
 \psi^{(0)}_{k_0}
 \Xi^{(\veczero)}_I
 \psi^{(\vecdelta_{J_0})}_{J_0}
 \psi^{(\veczero)}_{J_1}
 \phi^{(\veczero)}_K\nbigm 
\]
are compatible with the pre-$K$-Betti structures.
\hfill\qed

\vspace{.1in}

Thus, we obtain a pre-$K$-Betti structure of 
$\pi_{\ellsitabar\ast}\nbigq^{\ellsitabar}(\nbigm)
\simeq \nbigm$,
which is independent of the choice of
a factorization of $\pi_{\ellsitabar}$.
It is called the pre-$K$-Betti structure of $\nbigm$
associated to the good $K$-structure,
and denoted by $\nbigf_{\nbigm}$.
We obtain a pre-$K$-holonomic $\nbigd_X$-module
$\Upsilon(\nbigm):=(\nbigm,\nbigf_{\nbigm})$.
Thus, we obtain the desired exact functor
$\Upsilon:\Hol^{\good}(X,D,K)
\lrarr \Hol^{\pre}(X,K)$.
It is clearly exact.

\subsection{Proof of Proposition \ref{prop;13.4.26.120}}
\label{subsection;13.4.26.130}

If $\nbigm\bigl(\ast D(H^c)\bigr)=\nbigm$, 
any $\nbigq^H(\nbigm,\veci,\vecj)$
are equipped with the pre-$K$-Betti structures,
which induce a pre-$K$-Betti structure of $\nbigm$.
\begin{lem}
The associated pre-$K$-Betti
structures of $\nbigm$ are the same.
\end{lem}
\pf
The naturally defined morphisms
\[
 \Xi^{(\veczero)}_{H^c}
 \Xi^{(\veczero)}_K
 \psi^{(\vecdelta_{J_0})}_{J_0}
 \psi^{(\veczero)}_{J_1}
 \phi^{(\veczero)}_I(\nbigm)
\lrarr
\Xi^{(\veczero)}_K
 \psi^{(\vecdelta_{J_0})}_{J_0}
 \psi^{(\veczero)}_{J_1}
 \phi^{(\veczero)}_I(\nbigm)
\]
induce 
the quasi-isomorphism
$\pi_{\ellsitabar\ast}\nbigq^{\ellsitabar}(\nbigm)
\lrarr
 \pi_{H\ast}\nbigq^{H}(\nbigm)$,
which is compatible with 
the pre-$K$-Betti structures.
\hfill\qed

\vspace{.1in}

Let us prove Proposition \ref{prop;13.4.26.120}.
By the above consideration,
the following isomorphisms
are compatible with the pre-$K$-Betti structures:
\[
 V\bigl(!D(H)\bigr)
\stackrel{\simeq}{\lrarr}
 \nbigq^H(V(!D(H)))
\stackrel{\simeq}{\llarr}
 \nbigq^{\ellsitabar}(V(!D(H)))
\]
Thus,
we obtain Proposition \ref{prop;13.4.26.120}.
\hfill\qed

\subsection{Full faithfulness}
\label{subsection;13.4.26.21}

Let us prove that 
the functor $\Upsilon$ is fully faithful.
We denote $\Upsilon(\nbigm_i)$
by $\nbigm_i$ to simplify the notation.
Let $\nbigm_i\in\Hol^{\good}(X,D,K)$ $(i=1,2)$.
Suppose we are given
a morphism $\varphi:\nbigm_1\lrarr\nbigm_2$
in $\Hol^{\pre}(X,K)$.
We would like to prove that
$\varphi$ gives a morphism
in $\Hol^{\good}(X,D,K)$.

\vspace{.1in}

We use an induction 
on $\rho(\nbigm_1\oplus\nbigm_2)$.
(See \S\ref{subsection;09.10.28.3}
for $\rho$.)
We take a subset $J\subset \ellsitabar$ such that
$|J|=n-\dim\Supp(\nbigm_1\oplus\nbigm_2)$
and $(\nbigm_1\oplus\nbigm_2)
 \bigl(\ast D(J^c)\bigr)\neq 0$.
Let $g$ be a holomorphic function such that
$g^{-1}(0)=D(J^c)$.
Then,
$\nbigm_i(\ast g)$
and 
$\nbigm_i\otimes\gbigi^{a,b}_g$
come from good meromorphic flat bundles
with good $K$-structures on $(D_J,D_J(J^c))$.
We have the following morphisms
in $\Hol^{\good}(X,D,K)$:
\[
 \nbigm_i(!g)
\lrarr
 \Xi^{(0)}_g\bigl(
 \nbigm_i(\ast g)
 \bigr)
\lrarr
 \nbigm_i(\ast g)
\]
They are compatible with
the associated pre-$K$-Betti structures.
By the localization in Lemma \ref{lem;13.4.26.40}
and Lemma \ref{lem;13.4.26.41},
we obtain the following
in $\Hol^{\good}(X,D,K)$:
\[
 \nbigm_i(!g)
\lrarr
\nbigm_i
\lrarr
 \nbigm_i(\ast g)
\]
Note the uniqueness of good $K$-structure
on $\nbigm_i(\star g)$
in Lemma \ref{lem;13.4.26.200}.
We obtain the following diagram 
of the pre-$K$-holonomic $\nbigd$-modules:
\[
 \begin{CD}
 \nbigm_1(!g) @>>>
 \Xi^{(0)}_g(\nbigm_1(\ast g))
\oplus
 \nbigm_1
 @>>>
 \nbigm_1(\ast g)\\
 @VV{\varphi(!g)}V
 @VV{\Xi^{(0)}_g(\varphi)\oplus \varphi}V
 @VV{\varphi(\ast g)}V \\
 \nbigm_2(!g) @>>>
 \Xi^{(0)}_g(\nbigm_2(\ast g))
\oplus
 \nbigm_2
 @>>>
 \nbigm_2(\ast g)
 \end{CD}
\]

We obtain a morphism
$\phi^{(0)}_g(\varphi):
 \phi^{(0)}_g(\nbigm_1)\lrarr
 \phi^{(0)}_g(\nbigm_2)$
in $\Hol^{\pre}(X,K)$.
By using the inductive assumption,
$\phi^{(0)}_g(\varphi)$
is a morphism in $\Hol^{\good}(X,D,K)$.
Then,
$\varphi$ is obtained as the cohomology of
the following:
\begin{equation}
 \label{eq;13.4.26.202}
 \begin{CD}
 \psi^{(1)}_g(\nbigm_1(\ast g)) @>>>
 \Xi^{(0)}_g(\nbigm_1(\ast g))
\oplus
 \phi^{(0)}_g(\nbigm_1)
 @>>>
 \psi^{(0)}_g(\nbigm_1(\ast g))
 \\
 @VV{\psi_g^{(1)}\varphi}V
 @VV{\Xi^{(0)}_g(\varphi)\oplus \phi_g^{(0)}\varphi}V
 @VV{\psi_g^{(0)}\varphi}V \\
 \psi_g^{(1)}(\nbigm_2(\ast g))
 @>>>
 \Xi^{(0)}_g(\nbigm_2(\ast g))
\oplus
 \phi_g^{(0)}(\nbigm_2)
 @>>>
 \psi_g^{(0)}(\nbigm_2(\ast g))
 \end{CD}
\end{equation}
The morphisms in (\ref{eq;13.4.26.202})
are morphisms in $\Hol^{\good}(X,D,K)$.
Therefore, we obtain that $\varphi$
is also a morphism
in $\Hol^{\good}(X,D,K)$.
\hfill\qed

\subsection{Independence from the coordinate system}
\label{subsection;13.4.26.101}

Let us prove that
the essential image of $\Upsilon$
is independent of the choice of a coordinate system.
Let $(w_1,\ldots,w_n)$ be another
holomorphic coordinate system
such that $w_i^{-1}(0)=z_i^{-1}(0)$.
It is enough to prove the following lemma.

\begin{lem}
\label{lem;13.4.26.201}
If $\nbigm$ has a good $K$-structure
with respect to the coordinate system $(z_1,\ldots,z_n)$,
it has an induced good $K$-structure
with respect to $(w_1,\ldots,w_n)$ such that 
the associated pre-$K$-Betti structures
are the same.
\end{lem}
\pf
We use symbols
$\phi^{(\veczero)}_{\vecz,I}$ and $\phi^{(\veczero)}_{\vecw,I}$
to distinguish the dependence on 
the coordinate systems.
As remarked in \S\ref{subsection;09.10.22.5},
we have the natural isomorphisms
(\ref{eq;09.10.22.6}).
They induce isomorphisms
$\phi^{(\veczero)}_{\vecz,I}(\nbigm)
\simeq
 \phi^{(\veczero)}_{\vecw,I}(\nbigm)$
and
$\psi^{(a)}_i\phi^{(\veczero)}_{\vecz,I}(\nbigm)
\simeq
 \psi^{(a)}_i\phi^{(\veczero)}_{\vecz,I}(\nbigm)$.
Hence, we obtain good $K$-structure
of $\nbigm$ with respect to $(w_1,\ldots,w_n)$.
Let $\nbigq^{\ellsitabar}_{\vecz}(\nbigm)$
and $\nbigq^{\ellsitabar}_{\vecw}(\nbigm)$
denote the $\ellsitabar$-square associated to
$\nbigm$ with respect to
the coordinate systems $(z_1,\ldots,z_n)$
and $(w_1,\ldots,w_n)$,
respectively.
It is easy to observe that 
isomorphisms (\ref{eq;09.10.22.6})
induce
$\pi_{\ellsitabar\ast}
 \nbigq^{\ellsitabar}_{\vecz}(\nbigm)\simeq
\pi_{\ellsitabar\ast}
 \nbigq^{\ellsitabar}_{\vecw}(\nbigm)$
compatible with pre-$K$-Betti structures,
and they induce the identity on $\nbigm$.
Hence, the associated pre-$K$-Betti structures
on $\nbigm$ are the same.
Thus, the proof of Lemma \ref{lem;13.4.26.201}
and Proposition \ref{prop;13.4.26.100}
are finished.
\hfill\qed

%% file: 6.4.tex
\subsection{Good $K$-structure of meromorphic flat connections}
\label{subsection;13.4.21.30}

Let $X$ be a complex manifold
with a hypersurface $D$.
Let $V$ be a meromorphic flat connection on $(X,D)$,
i.e., $V$ is a reflexive $\nbigo_X(\ast D)$-coherent sheaf
with a flat connection. 
\index{meromorphic flat connection}
We do not assume that $V$ is good.
\begin{df}
As in the case of good meromorphic flat bundles,
a $K$-structure of $V$
means a pre-$K$-Betti structure of
the flat bundle $V_{|X\setminus D}$.
\index{$K$-structure}
\end{df}

Recall that, according to K. Kedlaya
(\cite{kedlaya}, Theorem 8.2.2 of \cite{kedlaya2}),
for any point $P\in X$,
there exist a neighbourhood $X_P\subset X$
and a projective birational morphism
$\lambda_P:\check{X}_P\lrarr X_P$
such that 
(i) $\lambda_P:\check{X}_P\setminus \lambda_P^{-1}(D)
\simeq X_P\setminus D$,
(ii) $\check{D}_P:=\lambda_P^{-1}(D)$ is normal crossing,
(iii) $\lambda_P^{\ast}V$ is a good meromorphic 
flat bundle.
(See also \cite{mochi6} and 
Theorem 16.2.1 of \cite{mochi7}
for the algebraic case.)
Such $(X_P,\lambda_P)$ is called
a local resolution of $V$ in this paper.
In the situation, we set $D_P:=D\cap X_P$.
\index{local resolution}

\begin{df}
A $K$-structure of $V$ is called good at $P$
if the following holds:
\begin{itemize}
\item
For any local resolution $(X_P,\lambda_P)$ around $P$,
the induced pre-$K$-Betti structure of
$\lambda_P^{\ast}(V_{|X_P\setminus D})$
is a good $K$-structure of $\lambda_P^{\ast}V$.
\end{itemize}
A $K$-structure of $V$ is called good
if it is good at any point of $X$.
\index{good $K$-structure}
\hfill\qed
\end{df}

If a $K$-structure of $V$ is good,
the induced $K$-structure on the dual $V^{\lor}$
is also good.
The following lemma is easy to see.

\begin{lem}
\mbox{{}}\label{lem;13.4.22.20}
Let $V_i$ $(i=1,2)$ be meromorphic flat bundles
on $(X,D)$ with a good $K$-structure.
\begin{itemize}
\item
The naturally induced $K$-structures on 
$V_1\oplus V_2$,
$V_1\otimes V_2$
and $\nhom(V_1,V_2)$
are good.
\item
Let $f:V_1\lrarr V_2$ be a flat morphism
which is compatible with the $K$-structures.
Then, the naturally induced $K$-structures of
$\Ker f$, $\Cok f$ and $\Image(f)$
are good.
\hfill\qed
\end{itemize}
\end{lem}

Let $\varphi:X'\lrarr X$ be a morphism of complex manifolds
such that $D':=\varphi^{-1}(D)$ is normal crossing.
We have the induced good meromorphic flat bundle
$V'=\varphi^{\ast}V$.
A $K$-structure of $V$ induces a $K$-structure of $V'$.
\begin{lem}
\label{lem;13.1.15.20}
If the $K$-structure of $V$ is good,
the $K$-structure of $V'$ is also good.
Conversely,
suppose that the $K$-structure of $V'$ is good
and that $\varphi$ is surjective.
Then, the $K$-structure of $V$ is also good.
\end{lem}
\pf
Let $(X_P,\lambda_P)$ be a local resolution for $V$ 
around $P\in X$.
We take a projective birational morphism
$\lambda:\check{X}_P'\lrarr \check{X}_P\times_{X}X'$
such that
(i) $\check{X}_P'$ is smooth,
(ii) the induced morphism
 $\varphi_P:\check{X}_P'\lrarr \check{X}_P$
 gives
 $\check{X}_P'\setminus \check{D}'_P
 \simeq
 \check{X}_P\setminus \check{D}_P$,
where $\check{D}_P':=\lambda^{-1}(\check{X}_P\times_XD')$.
The induced map
$\lambda_P':\check{X}'_P\lrarr X'$
gives a local resolution for $V'$.
Then, the claim follows from
Lemma \ref{lem;14.1.15.10}.
\hfill\qed

\vspace{.1in}
We obtain the following lemma
from Proposition \ref{prop;14.1.18.30}.
\begin{lem}
Let $V$ be a meromorphic flat connection on $(X,D)$
with a $K$-structure.
Suppose that,
for any morphism $\Delta\lrarr X$ 
with $\varphi(\Delta)\cap D=\{\varphi(0)\}$,
the induced $K$-structure of $\varphi^{\ast}(V)$
is good.
Then, the $K$-structure of $V$ is also good.
\hfill\qed
\end{lem}

We obtain the following lemma
from Lemma \ref{lem;13.4.23.1}.
\begin{lem}
\label{lem;13.4.23.2}
Let $V$ be a meromorphic flat connection
with a good $K$-structure.
Let $V_1\subset V$ be a sub-connection
such that 
$V_{1|X\setminus D}$
is compatible with the $K$-structure.
Then, the induced $K$-structure of $V_1$
is good.
A similar claim holds for quotients of $V$.
\hfill\qed
\end{lem}

\subsection{Canonical pre-$K$-Betti structures}
\label{subsection;13.4.22.1}

\index{canonical pre-$K$-Betti structure}

Let $V$ be a meromorphic flat connection on $(X,D)$ 
with a good $K$-structure.
Let $D=D_1\cup D_2$ be a decomposition,
i.e., $D_i$ are unions of irreducible components of $D$
such that $\codim_X(D_1\cap D_2)>1$.
Let $(X_P,\lambda_P)$ be any local resolution of $V$
around $P\in X$.
Put $D_{P1}=D_1\cap X_P$
and $\check{D}_{P1}:=\lambda_{P}^{-1}(D_{1})$.
We have the decomposition 
$\check{D}_P=\check{D}_{P1}\cup \check{D}_{P2}$.
We set $V_P:=V_{|X_P}$ and $\check{V}_P:=\lambda_P^{\ast}V$.
The canonical pre-$K$-Betti structure 
$\nbigv_{\check{V}_P}^{<\check{D}_{P1}\leq \check{D}_{P2}}$
of $\check{V}_P(!\check{D}_{P1})$
induces a pre-$K$-Betti structure 
$\nbigg$ of $V_P(!D_{P1})$.
Let $(X_P^{(1)},\lambda^{(1)}_P)$ be another local resolution of $V$
around $P\in X$.
It induces a pre-$K$-Betti structure
$\nbigg^{(1)}$ of $V_{|X_P^{(1)}}$.
We have $\nbigg^{(1)}=\nbigg$ on $X_P\cap X_P^{(1)}$.
Indeed, we can find a local resolution 
$(X_P^{(2)},\lambda^{(2)}_P)$ with morphisms
$a:\check{X}_P^{(2)}\lrarr \check{X}_P^{(1)}$
and $b:\check{X}_P^{(2)}\lrarr \check{X}_P$
such that
$\lambda_P^{(2)}=\lambda_P^{(1)}\circ a=\lambda_P\circ b$.
By using $(X_P^{(2)},\lambda^{(2)}_P)$
with Proposition \ref{prop;09.10.28.12},
we can prove that
the pre-$K$-Betti structures are equal.
Therefore,
by gluing the pre-$K$-Betti structures
around any $P\in X$,
we obtain 
a pre-$K$-Betti structure of $V(!D_1)$.
(See Proposition 10.2.9 of \cite{kashiwara-schapira}.)
We denote it by
$\nbigf_{V}^{<D_1}$.
It is called the canonical pre-$K$-Betti structure
of $V(!D_1)$.
By taking the dual of
$\bigl(V^{\lor}\bigr)(!D_1)$,
we obtain a pre-$K$-Betti structure of
$\bigl(V(!D)\bigr)(\ast D_1)$,
denoted by
$\nbigf_{V}^{<D\leq D_1}$.

\vspace{.1in}

\vspace{.1in}

Let $D_3$ be a hypersurface of $X$.
Let $\varphi:X'\lrarr X$ be a projective birational morphism
such that 
(i) $X'\setminus D'\simeq X\setminus (D\cup D_3)$
where $D':=\varphi^{-1}(D\cup D_3)$,
(ii) $D'$ is normal crossing.
We set $D_1':=\varphi^{-1}(D_1)$.
We have $D_2'$ such that
$D'=D_2'\cup D_1'$ is a decomposition.
We set
$V'=\varphi^{\ast}V(\ast D')$.
\begin{prop}
\label{prop;13.4.22.10}
The natural morphisms
\[
V(!D_1)\lrarr\varphi_{\dagger}V'(!D_1'),
\quad
\varphi_{\dagger}\bigl(
 V^{\prime}(!D')(\ast D_1')
 \bigr)
\lrarr
 V(!D)(\ast D_1)
\]
are compatible with the canonical pre-$K$-Betti structures.
\end{prop}
\pf
Let $(X_P,\lambda_P)$ be a local resolution for $V$ 
around $P\in X$.
We take a projective birational morphism
$\lambda:\check{X}_P'\lrarr \check{X}_P\times_{X}X'$
such that
(i) $\check{X}_P'$ is smooth,
(ii) the induced morphism
 $\varphi_P:\check{X}_P'\lrarr \check{X}_P$
 gives
 $\check{X}_P'\setminus \check{D}'_P
 \simeq
 \check{X}_P\setminus \check{D}_P$,
where $\check{D}_P':=\lambda^{-1}(\check{X}_P\times_XD')$.
The induced map
$\lambda_P':\check{X}'_P\lrarr X'$
gives a local resolution for $V'$.
By Proposition \ref{prop;09.10.28.12},
$\lambda_P^{\ast}(V)(!\check{D}_{P1})
\lrarr
 \varphi_{P\dagger}\bigl(
 \lambda_P^{\prime\ast}
 V'(!\check{D}'_{P1})
 \bigr)$
is compatible with the pre-$K$-Betti structures.
Then, we obtain that
$V(!D_1)\lrarr\varphi_{\dagger}V'(!D_1')$
is compatible with the pre-$K$-Betti structures.
We obtain the claim for the other
as the dual.
\hfill\qed

\subsection{Pre-$K$-Betti structure on the real blow up}

Let $X$, $D$ and $V$ be as in
the beginning of \S\ref{subsection;13.4.22.1}.
Let $G:X\lrarr\cnum^{\ell}$ be a holomorphic function
such that
$G^{-1}(D_0)\subset D_1$,
where $D_0=\bigcup_{i=1}^{\ell}\{z_i=0\}$.
We obtain an object $(X,G)$ in $\Cat_{\ell}$.
Let $\pi:\Xtilde(G)\lrarr X$ denote the real blow up.

\begin{lem}
\label{lem;13.4.21.10}
The natural morphism
$R\pi_{\ast}\DR^{\rapid}_{X,G}(V(!D_1))
\lrarr
 \DR_X(V(!D_1))$
is an isomorphism
in $D^b(\cnum_X)$.
\end{lem}
\pf
It is enough to check the claim locally around each 
$P\in X$.
Let $(X_P,\lambda_P)$ be a local resolution of $V$
around $P$.
We set 
$G_P:=G_{|X_P}$ and
$\check{G}_P:=G\circ\lambda_P$.
We obtain a morphism
$\lambda_P:
 (\check{X}_P,\check{G}_P)
\lrarr
 (X_P,G_P)$
in $\Cat_{\ell}$.
We set
$\check{\nbigm}_P:=\check{V}_P(!\check{D}_1)$.
By Corollary \ref{cor;13.4.20.230},
we have the following isomorphism
in $D^b(\cnum_{\Xtilde_P(G_P)})$:
\[
 R\lambdatilde_{P\ast}
 \DR^{\rapid}_{\check{X}_P,\check{G}_P}
 (\check{\nbigm}_P)
\simeq
 \DR^{\rapid}_{X_P,G_P}
 (\lambda_{P\dagger}\check{\nbigm}_P)
=\DR^{\rapid}_{X,G}(V(!D_1))_{|\Xtilde_P(G_P)}
\]
By using
$R\pi_{\check{G}_P\,\ast}
 \DR^{\rapid}_{\check{X}_P,\check{G}_P}
 (\check{\nbigm}_P)
\simeq
 \DR_{\check{X}_P}(\check{\nbigm}_P)$,
we obtain the claim of the lemma.
\hfill\qed

\vspace{.1in}

In the situation of the proof of
Lemma \ref{lem;13.4.21.10},
let $\widetilde{\check{X}}_P(\check{D}_P)$
be the real blow up along $\check{D}_P$.
We have the natural map
$\rho:
 \widetilde{\check{X}}_P(\check{D}_P)
\lrarr
\widetilde{\check{X}}_P(\check{G}_P)$.
As in Lemma \ref{lem;13.4.27.3},
we have the following natural isomorphism:
\[
 R\rho_{\ast}
 \DR^{\leq \check{D}_{P2}<\check{D}_{P1}}
 _{ \widetilde{\check{X}}_P(\check{D}_P)}
 (\check{V}_P)
\simeq
 \DR^{\rapid}_{\check{X}_P,\check{G}_P}
 (\check{V}_P(!\check{D}_{P1}))
\]
In particular, a good $K$-structure of 
$\check{V}_P$
induces a $K$-structure of
$\DR^{\rapid}_{\check{X}_P,\check{G}_P}
 (\check{V}_P(!\check{D}_{P1}))$.
We would like to glue them.

\begin{lem}
\label{lem;13.4.21.31}
Suppose that there exists a finite family
$\{(\nbigu_i,\lambda_i)\,|\,i\in\Lambda\}$ 
$(|\Lambda|<\infty)$
of local resolutions of $V$
such that $X=\bigcup\nbigu_i$.
Then, there exists an object
$\nbigk$ in $D^b(K_{\Xtilde(G)})$
with isomorphisms 
\[
 c_1:\nbigk\otimes\cnum
\simeq
 \DR^{\rapid}_{X,G}(V(!D_1))
\quad\quad
\mbox{\rm in $D^b(\cnum_{\Xtilde(G)})$},
\]
\[
 c_2:
 R\pi_{\ast}\nbigk
\simeq
 \nbigf_V^{<D_1}
\quad
\mbox{\rm in $D^b(K_{X})$},
\]
such that 
$c_2\otimes\cnum$
is equal to
$R\pi_{\ast}c_1$.
\end{lem}
\pf
We shall construct a $K$-complex $\nbigk$
on $\Xtilde(G)$ as follows.
For $I\subset \Lambda$,
we set $\nbigu_I:=\bigcap_{i\in I}\nbigu_i$.
Let $\iota_I:\nbigu_I\lrarr X$ denote the inclusion.
We set $G_I:=G_{|\nbigu_I}$.
Take local resolutions
$\lambda_I:\check{\nbigu}_I\lrarr\nbigu_I$
of $V$.
We may assume to have
$\lambda_{IJ}:\check{\nbigu}_J\lrarr \check{\nbigu}_I$
such that
$\iota_I\circ\lambda_I\circ\lambda_{IJ}=\iota_J\circ\lambda_J$
for any $I\subset J$.
We have $\lambda_{I_1I_2}\circ\lambda_{I_2I_3}=\lambda_{I_1I_3}$.
We put $V_I:=\lambda_I^{\ast}V$.

We set $\check{D}_I:=\lambda_I^{-1}(D)$,
and $\check{D}_{I1}:=\lambda_I^{-1}(D_1)$.
Let $\check{D}_{I2}$ denote the complement of 
$\check{D}_{I1}$ in $\check{D}_I$.
Let $\check{\pi}_I:
 \widetilde{\check{\nbigu}}_I(\check{D}_I)
\lrarr
 \check{\nbigu}_I$
denote the real blow up.
We have the induced morphisms
$\lambdatilde_{IJ}:
 \widetilde{\check{\nbigu}}_J(\check{D}_J)
\lrarr
 \widetilde{\check{\nbigu}}_I(\check{D}_I)$.
We also have the induced morphisms
$\lambdatilde_{I}:
  \widetilde{\check{\nbigu}}_I(\check{D}_I)
\lrarr
 \widetilde{\nbigu}_I(G_I)$.
Let $\iotatilde_I:\widetilde{\nbigu}_I(G_I)
\lrarr \Xtilde(G)$ denote the inclusion.

Let $\nbigl_{K,I}$ denote the $K$-local system 
on $\widetilde{\check{\nbigu}}_I(\check{D}_I)$
with the Stokes structure associated to $V_I$
with good $K$-structure.
We have the constructible sheaves
$\nbigl_{K,I}^{<\check{D}_{I1}\leq \check{D}_{I2}}$
on $\widetilde{\check{\nbigu}}_I(\check{D}_I)$,
and natural morphisms
$\lambdatilde_{IJ}^{-1}\nbigl_{K,I}^{<\check{D}_{I1}\leq \check{D}_{I2}}
\lrarr
 \nbigl_{K,J}^{<\check{D}_{J1}\leq \check{D}_{J2}}$.
For any sheaf $\nbigf$,
let $\Gd(\nbigf)$ denote its Godement resolution.
By the construction,
we have natural morphisms
\begin{equation}
 \label{eq;14.1.11.1}
 \lambdatilde_{IJ}^{-1}
\Gd\bigl(
\nbigl_{K,I}^{<\check{D}_{I1}\leq \check{D}_{I2}}
\bigr)
\lrarr
 \Gd\bigl(
 \lambdatilde_{IJ}^{-1}\nbigl_{K,I}^{<\check{D}_{I1}\leq \check{D}_{I2}}
 \bigr)
\lrarr
 \Gd\bigl(
 \nbigl_{K,J}^{<\check{D}_{J1}\leq \check{D}_{J2}}
 \bigr).
\end{equation}
We set 
$\nbigg^{\bullet}_{K,I}:=
 \iotatilde_{I\ast}
 \lambdatilde_{I\ast}
 \Gd\bigl(
\nbigl_{K,I}^{<\check{D}_{I1}\leq \check{D}_{I2}}
\bigr)[d_X]$
on $\Xtilde(G)$.
The morphisms (\ref{eq;14.1.11.1})
induce 
$\lambda_{JI}:
 \nbigg^{\bullet}_{K,I}
\lrarr
 \nbigg^{\bullet}_{K,J}$.
They satisfy
$\lambda_{I_1I_2}\circ\lambda_{I_2I_3}
=\lambda_{I_1I_3}$.

We take a $K$-vector space $U_K$
with a basis $\{e_i\,|\,i\in\Lambda\}$.
Let $U_{K,I}$ denote the subspace in
$\bigwedge^{\bullet}U_K$
generated by
$e_{i_1}\wedge\cdots\wedge e_{i_m}$
where $I=(i_1,\ldots,i_m)$.
For $m\in\seisuu_{\geq 0}$,
we set
\[
 \nbigk^{m,\bullet}_K:=
 \bigoplus_{|I|=m+1}
 \nbigg^{\bullet}_{I,K}\otimes U_{K,I}
\]
We have the morphism
$\nbigk^{m,\bullet}_K\lrarr\nbigk^{m+1,\bullet}_K$
induced by the morphisms
$\lambda_{I,I\cup\{j\}}\otimes (e_j\wedge\bullet)$.
They give a double complex $\nbigk^{\bullet,\bullet}_K$
of  $K_{\Xtilde(G)}$-modules.
The total complex is denoted by
$\nbigk^{\bullet}$.

\vspace{.1in}
We have the $\cnum$-local systems
$\nbigl_I$ with the Stokes structure
on $\widetilde{\check{\nbigu}}_I(\check{D}_I)$
associated to $V$.
Using $\nbigl_I$ with the same construction,
we obtain complexes $\nbigg_{\cnum,I}^{\bullet}$,
a double complex $\nbigk^{\bullet,\bullet}_{\cnum}$
and a complex $\nbigk^{\bullet}_{\cnum}$.

We have naturally defined isomorphisms
$\nbigl_{K,I}^{<\check{D}_{I1}\leq \check{D}_{I2}}
 \otimes\cnum
\lrarr
\nbigl_{I}^{<\check{D}_{I1}\leq \check{D}_{I2}}$.
The natural morphisms
$\Gd\bigl(
 \nbigl_{K,I}^{<\check{D}_{I1}\leq \check{D}_{I2}}
 \bigr)
 \otimes\cnum
\lrarr
 \Gd\bigl(
\nbigl_{I}^{<\check{D}_{I1}\leq \check{D}_{I2}}
 \bigr)$
are quasi-isomorphisms.
By the projection formula,
we have the following natural isomorphisms
\[
 \iotatilde_{I\ast}
 \lambdatilde_{I\ast}
 \bigl(
 \Gd\bigl(
 \nbigl_{K,I}^{<\check{D}_{I1}\leq \check{D}_{I2}}
 \bigr)[d_X]
 \otimes\cnum
 \bigr)
\simeq
\nbigg_{K,I}\otimes\cnum.
\]
It also implies that the complex
$(\iotatilde\circ\lambdatilde_I)_{\ast}
 \bigl(
 \Gd\bigl(
 \nbigl_{K,I}^{<\check{D}_{I1}\leq \check{D}_{I2}}
 \bigr)
 \otimes\cnum
 \bigr)$
represents
 $R(\iotatilde\circ\lambdatilde_I)_{\ast}
 \bigl(
 \Gd\bigl(
 \nbigl_{K,I}^{<\check{D}_{I1}\leq \check{D}_{I2}}
 \bigr)
 \otimes\cnum
 \bigr)$.
Hence, the natural morphism
$\nbigg_{K,I}\otimes\cnum
\lrarr
 \nbigg_{\cnum,I}$ is a quasi-isomorphism.
Then, it is easy to deduce that the natural morphism
$\nbigk^{\bullet}\otimes_K\cnum
\lrarr
 \nbigk_{\cnum}^{\bullet}$
is a quasi-isomorphism.

\vspace{.1in}

We have the natural quasi-isomorphism
$\nbigl_{I}^{<\check{D}_{I1}\leq \check{D}_{I2}}[d_X]
\lrarr
 \DR^{<\check{D}_{I1}\leq \check{D}_{I2}}
 _{\widetilde{\check{\nbigu}}_I(\check{D}_I)}(\check{V}_I)$.
We have morphisms
\[
 \lambdatilde_{JI}^{-1}
 \DR^{<\check{D}_{I1}\leq \check{D}_{I2}}
 _{\widetilde{\check{\nbigu}}_I(\check{D}_I)}(\check{V}_I)
\lrarr
  \DR^{<\check{D}_{J1}\leq \check{D}_{J2}}
 _{\widetilde{\check{\nbigu}}_J(\check{D}_J)}(\check{V}_J).
\]
By applying the above construction
to $\DR^{<\check{D}_{I1}\leq \check{D}_{I2}}
 _{\widetilde{\check{\nbigu}}_I(\check{D}_I)}(\check{V}_I)$
instead of
$\nbigl_{I}^{<\check{D}_{I1}\leq \check{D}_{I2}}[d_X]$,
we obtain double complexes
$\nbigg_{I,\DR}^{\bullet,\bullet}$
on $\Xtilde(G)$,
and a complex
$\nbigk_{\DR}^{\bullet}$
on $\Xtilde(G)$.
The natural morphism
$\nbigk^{\bullet}_{\cnum}
\lrarr
 \nbigk^{\bullet}_{\DR}$
is a quasi-isomorphism.

\vspace{.1in}
Set $\check{H}_I:=\check{G}_I^{-1}(0)$.
We have the complexes
$\DR^{<\check{D}_{I1}}_{\widetilde{\check{\nbigu}}_I(\check{H}_I)}
 (\check{V}_I)$
and 
$\DR^{<\check{H}_{I}}_{\widetilde{\check{\nbigu}}_I(\check{H}_I)}
 (\check{V}_{I}(!\check{D}_{I1}))$
on $\widetilde{\check{\nbigu}}_I(\check{H}_I)$.
By applying the above construction to them,
we obtain double complexes
$\nbigg^{\bullet,\bullet}_{I,a}$ $(a=1,2)$,
and complexes
$\nbigk_{a}$ $(a=1,2)$ on $\Xtilde(G)$.
We have the following natural quasi-isomorphisms of complexes,
as in Lemma \ref{lem;13.4.27.1}:
\[
 \nbigg^{\bullet}_{I,\DR}
\llarr
 \nbigg^{\bullet}_{I,1}
\lrarr
 \nbigg^{\bullet}_{I,2}
\]
Hence, we have the natural quasi-isomorphisms
of complexes
$\nbigk^{\bullet}_{\DR}
\llarr
\nbigk^{\bullet}_{1}
\lrarr
 \nbigk^{\bullet}_2$.
We set
$\nbigg^{\bullet}_{I,3}:=
  \iotatilde_{I\ast}
 \Gd\bigl(
\iotatilde_I^{-1}
 \DR^{\rapid}_{X,G}\bigl(V(!D_1)\bigr)
 \bigr)$.
As before, by the \v{C}ech construction
we obtain 
a complex $\nbigk^{\bullet}_3$.
We have natural quasi-isomorphism
$\nbigg_{I,3}
\lrarr
 \nbigg_{I,2}$,
which induce
$\nbigk^{\bullet}_{3}
\lrarr
 \nbigk^{\bullet}_2$.
By construction,
we have natural quasi-isomorphisms
$\Gd\bigl(
 \DR^{\rapid}_{X,G}\bigl(V(!D_1)\bigr)
\lrarr
 \nbigk^{\bullet}_{3}$.
(See Proposition 2.8.4 of \cite{kashiwara-schapira}.)
In all, we obtain the following sequence of
quasi-isomorphisms:
\begin{equation}
 \label{eq;14.1.11.10}
 \nbigk^{\bullet}_{K}\otimes\cnum
\lrarr
 \nbigk^{\bullet}_{\DR}
\llarr
 \nbigk^{\bullet}_{1}
\lrarr
 \nbigk^{\bullet}_2
\llarr
 \Gd\DR^{\rapid}_{X,G}\bigl(V(!D_1)\bigr)
\end{equation}
We define $c_1$
as the composite of the morphisms.

The projections
$\varphi_i:
 \nbigk^{\bullet}_{|\nbigutilde_i(G_i)}
 \lrarr \nbigg^{\bullet}_{K,i|\nbigutilde_i(G_i)}$
are quasi-isomorphisms.
It is easy to see that
$\lambda_{\{ij\},i|\nbigutilde_{ij}(G_{ij})}
 \circ
 \varphi_{i|\nbigutilde_{ij}(G_{ij})}$
and 
$\lambda_{\{ij\},j|\nbigutilde_{ij}(G_{ij})}
 \circ
 \varphi_{j|\nbigutilde_{ij}(G_{ij})}$
are chain homotopic.
Hence,
$\pi_{\ast}\nbigk^{\bullet}$
is a $K$-perverse sheaf
obtained as the gluing of
$\pi_{\ast}\nbigg_{K,i|\nbigutilde_i(G_i)}$.
We obtain an isomorphism
of $K$-perverse sheaves
$\nbigf_V^{<D_1}
\simeq
 \pi_{\ast}\nbigk^{\bullet}$,
which is $c_2$.
We can easily compare 
$(c_2\otimes\cnum)_{|\nbigu_i}$
and 
$R\pi_{\ast}(c_1)_{|\nbigu_i}$,
and we obtain
$c_2\otimes\cnum
=R\pi_{\ast}(c_1)$.
\hfill\qed

\subsection{Sequence of hypersurface pairs}
\label{subsection;13.4.27.20}

\index{sequence of hypersurface pairs}

Let $X$ be a complex manifold.
Let $\vecH=(H_{!},H_{\ast})$
be an ordered pair of 
(possibly empty) hypersurfaces of $X$.
Such a pair is called a hypersurface pair
in the following.
For any coherent $\nbigd_X$-module $\nbigm$.
we define
\[
 \gbigp_{\vecH}(\nbigm):=
 \bigl(\nbigm(\ast H_{\ast})
 \bigr)(!H_!),
\quad\quad
 \gbigp'_{\vecH}(\nbigm):=
 \bigl(\nbigm(!H_!)
 \bigr)(\ast H_{\ast}).
\]
We set $\DDD \vecH=(H_{\ast},H_{!})$.
Then, we have natural isomorphisms:
\[
 \DDD\bigl(
 \gbigp_{\vecH}(\nbigm)
 \bigr)
\simeq
 \gbigp'_{\DDD\vecH}(\DDD\nbigm)
\]
If we are given a sequence of hypersurface pairs
$\gbigh=(\vecH_1,\vecH_2,\ldots,\vecH_N)$,
we set
$\gbigp_{\gbigh}:=
 \gbigp_{\vecH_N}\circ\cdots
 \circ\gbigp_{\vecH_2}\circ
 \gbigp_{\vecH_1}$
and 
$\gbigp'_{\gbigh}:=
 \gbigp'_{\vecH_N}\circ\cdots
 \circ\gbigp'_{\vecH_2}\circ
 \gbigp'_{\vecH_1}$.
Clearly,
$\gbigp_{\gbigh}$
can be described as
$\gbigp'_{\gbigh_1}$
for an appropriate $\gbigh_1$.
\index{functor $\gbigp_{\gbigh}$}
\index{functor $\gbigp'_{\gbigh}$}
We shall use a special case of this operation
in \S\ref{subsection;09.10.17.150}.

\subsection{Generalization}
\label{subsection;13.4.23.100}

Let $X$, $D$ and $V$ be as in 
the beginning of \S\ref{subsection;13.4.22.1}.
Suppose that we are given a sequence of 
hypersurface pairs
$\gbigh=(\vecH_1,\ldots,\vecH_N)$
contained in $D$.
Let us observe that
$\gbigp_{\gbigh}(V)$ 
and
$\gbigp'_{\gbigh}(V)$
are naturally equipped with 
pre-$K$-Betti structures.

Let $P$ be any point of $X$.
We take a local resolution $(X_P,\lambda_P)$
of $V$ around $P$.
By taking the pull back,
we obtain a sequence of  hypersurface pairs
$\check{\gbigh}_P:=\lambda_P^{\ast}(\gbigh)$
contained in $\check{D}_P$.
For the irreducible decomposition
$\check{D}_{P}=\bigcup_{j\in\Lambda_P}\check{D}_{Pj}$,
there uniquely exists a subset
$I_P\subset\Lambda_P$
such that
$\gbigp_{\check{\gbigh}_P}(\check{V}_P)
\simeq
 \check{V}_P\bigl(!\check{D}_P(I_P)\bigr)$,
where
$\check{D}_P(I_P)=\bigcup_{j\in I_P}\check{D}_{Pj}$.
Hence, we have the canonical pre-$K$-Betti structure
$\check{V}_P\bigl(!\check{D}_P(I_P)\bigr)$
induced by the good $K$-structure of $\check{V}_P$.
\index{canonical pre-$K$-Betti structure}
By the natural isomorphism
$\lambda_{P\dagger}\gbigp_{\check{\gbigh}_P}(\check{V}_P)
\simeq
 \gbigp_{\gbigh}(V)_{|X_P}$,
we obtain a pre-$K$-Betti structure of
$\gbigp_{\gbigh}(V)_{|X_P}$.

Suppose that we are given other local resolutions
$(X_P^{(i)},\lambda_P^{(i)})$ $(i=1,2)$
as in \S\ref{subsection;13.4.22.1}.
We put
$\check{V}_P^{(2)}:=\lambda_P^{(2)\ast}V$.
We have the expression
$\gbigp_{\check{\gbigh}^{(2)}_P}
 (\check{V}^{(2)}_P)
\simeq
 \check{V}^{(2)}_P
 \bigl(!\check{D}^{(2)}_P(I_P^{(2)})\bigr)$.
For the morphism
$a:\check{X}^{(2)}_P\lrarr\check{X}_P$,
we have 
$\check{D}_P(I_P)
=a\bigl(
 \check{D}^{(2)}_P(I_P^{(2)})
\bigr)$.
We have the natural isomorphisms of
holonomic $\nbigd$-modules
$a_{\dagger}\gbigp_{\check{\gbigh}^{(2)}_P}(\check{V}^{(2)}_P)
\simeq
 a_{\dagger}\bigl(
 \check{V}^{(2)}_P\bigl(
 !a^{-1}(\check{D}_P(I_P))
 \bigr)
\simeq
 \gbigp_{\check{\gbigh}_P}(\check{V}_P)$
which are compatible with 
the pre-$K$-Betti structures.
Therefore,
we obtain the pre-$K$-Betti structures of
$\gbigp_{\gbigh}(V)$
by gluing the locally given pre-$K$-Betti structures.
We obtain 
a pre-$K$-Betti structure of 
$\gbigp'_{\gbigh}(V)$
in the same way.
They are
called the canonical pre-$K$-Betti structure
of $\gbigp_{\gbigh}(V)$
and $\gbigp'_{\gbigh}(V)$,
denoted by $\nbigf_{\gbigh,V}$
and $\nbigf'_{\gbigh,V}$.
\begin{lem}
Let $\gbigh^{\circ}=(\vecH^{\circ}_1,\ldots,\vecH^{\circ}_H)$ 
be a sequence of hypersurface pairs
such that
$H_{i\ast}^{\circ}\subset H_{i\ast}$
and 
$H_{i!}^{\circ}\supset H_{i!}$
for any $i$.
The natural morphisms
$\gbigp_{\gbigh^{\circ}}(V)
\lrarr
 \gbigp_{\gbigh}(V)$
and
$\gbigp'_{\gbigh^{\circ}}(V)
\lrarr
 \gbigp'_{\gbigh}(V)$
are compatible with the $K$-Betti structures.
\end{lem}
\pf
It is reduced to the case where $V$ is good.
Then, it is easy to check.
\hfill\qed

\vspace{.1in}
Let $G:X\lrarr \cnum^{\ell}$ be a holomorphic function.
The following lemma can be shown 
by the same arguments
as those in the proof of 
Lemma \ref{lem;13.4.21.10}
and Lemma \ref{lem;13.4.21.31}.
\begin{prop}
\label{prop;13.4.23.200}
Suppose that,
for $\vecH_N=(H_{N!},H_{N\ast})$,
we have
$G^{-1}(D_0)\subset H_{N!}$.
Then, the natural morphism
\[
 R\pi_{\ast}\DR^{\rapid}_{X,G}(\nbigm_N)
\lrarr
 \DR_X(\nbigm_N)
\]
is an isomorphism.
If we are given a finite family of local resolutions of $V$
as in Lemma {\rm\ref{lem;13.4.21.31}},
then there exists an object
$\nbigk$ in $D^b(K_{\Xtilde(G)})$
with isomorphisms 
$c_1:\nbigk\otimes\cnum
\simeq
 \DR^{\rapid}_{X,G}(\nbigm_N)$
in $D^b(\cnum_{\Xtilde(G)})$,
and
$c_2:
 R\pi_{\ast}\nbigk
\simeq
 \nbigf^{\can}_{\nbigm_N}$
in $D^b(K_{X})$,
such that 
$c_2\otimes\cnum$
is equal to
$R\pi_{\ast}c_1$.
\hfill\qed
\end{prop}

Let $D_3$, $\varphi:X'\lrarr X$
and $V'$ be as in Proposition \ref{prop;13.4.22.10}.
By the pull back,
we obtain a sequence of hypersurface pairs
$\gbigh':=\varphi^{-1}\gbigh$.

\begin{prop}
\label{prop;13.4.27.111}
The natural morphisms
$\varphi_{\dagger}
 \gbigp_{\gbigh'}(V'(!D'))
 \lrarr
 \gbigp_{\gbigh}(V(!D))$
and 
$\gbigp_{\gbigh}(V)
 \lrarr
 \varphi_{\dagger}
 \gbigp_{\gbigh'}(V')$
are compatible 
with the canonical pre-$K$-Betti structures.
The natural morphisms
$\varphi_{\dagger}
 \gbigp'_{\gbigh'}(V'(!D'))
 \lrarr
 \gbigp'_{\gbigh}(V(!D))$
and 
$\gbigp'_{\gbigh}(V)
 \lrarr
 \varphi_{\dagger}
 \gbigp'_{\gbigh'}(V')$
are also compatible 
with the canonical pre-$K$-Betti structures.
\end{prop}
\pf
It is reduced to the case where $V$ is good.
We have 
$\gbigp_{\gbigh}(V)=V(!D^{(1)})$
and 
$\gbigp_{\gbigh'}(V')=V'(!D^{\prime(1)})$
for some 
$D^{(1)}\subset D$
and $D^{\prime(1)}\subset D'$.
We have 
$\varphi(D^{\prime(1)})=D^{(1)}$.
We set 
$L^{(1)}:=\varphi^{-1}(D^{(1)})$.
Then,
the natural morphisms
$ \gbigp_{\gbigh}(V)
\simeq
 \varphi_{\dagger}
 V'(!L^{(1)})
\lrarr
 \varphi_{\dagger}
 \gbigp_{\gbigh'}(V')$
are compatible with
the pre-$K$-Betti structures.
Similarly, we obtain that
$ \gbigp'_{\gbigh}(V)
\lrarr
 \varphi_{\dagger}
 \gbigp'_{\gbigh'}(V')$
is compatible with
the pre-$K$-Betti structure.
We obtain the others by the dual.
\hfill\qed

%% file: 6.5.tex
\label{subsection;13.4.20.251}

Let $Y$ be a complex manifold
with a hypersurface $D_Y$.
Let $G:X\lrarr Y$ be a projective morphism of 
complex manifolds.
We set $D_{X0}:=G^{-1}(D_Y)$.
Let $D_X$ be a hypersurface of $X$ 
with a decomposition
$D_X=D_{X1}\cup D_{X2}$
such that $D_{X0}\subset D_{X2}$.

Let $V$ be a meromorphic flat connection on $(X,D_X)$
with a good $K$-structure.
Put $\nbigm:=V(!D_{X2})$.
Let $\nbigf_{\nbigm}$ be the canonical
pre-$K$-Betti structure.
Assume the following:
\begin{itemize}
\item
 $G^i_{\dagger}\nbigm=0$ for any $i\neq 0$,
 and $V_1:=G^0_{\dagger}(\nbigm)(\ast D_Y)$ is
 a meromorphic flat connection on $(Y,D_Y)$.
\end{itemize}
We put $\nbigg:=RG_{\ast}(\nbigf_{\nbigm})_{|Y-D_Y}$,
which gives a pre-$K$-Betti structure of 
$G_{\dagger}^0(\nbigm)_{|Y-D_Y}$.
The following theorem will be used in
the proof of Theorem \ref{thm;09.10.16.5}.
(See \S\ref{subsection;09.12.5.120}.)

\begin{thm}
\label{thm;13.4.20.410}
The $K$-structure $\nbigg$ of $V_1$
is good,
i.e.,
it is compatible with the Stokes filtrations.
Moreover, $RG_{\ast}\nbigf_{\nbigm}$
is the canonical pre-$K$-Betti structure
of $G_{\dagger}^0(\nbigm)$.
\end{thm}
\pf
It is enough to consider the issues
locally around any point $P$ of $Y$.
Let $(Y_P,\lambda_P)$ be a local resolution of $V_1$.
We take a projective birational morphism
$\lambda:X'\lrarr \check{Y}_P\times_YX$
such that
(i) $X'$ is smooth,
(ii) $D'_X:=\check{X}_P\times_XD_X$ is normal crossing,
(iii) the induced morphism
 $X'\setminus D'_X\lrarr X\setminus D_X$ 
 is an isomorphism.
Let $\mu:X'\lrarr X$ 
and $G':X'\lrarr \check{Y}_P$
be the induced maps.
We obtain a meromorphic flat connection
$V'=\mu^{\ast}V$ with a good $K$-structure.
We set $D_{X2}':=\mu^{-1}(D_{X2})$.
We have
$\mu_{\dagger}(V'(!D'_{X2}))=V(!D_{X2})$,
$G'_{\dagger}(V'(!D'_{X}))(\ast D_{P})
=\lambda_P^{\ast}V_1$
and 
$\lambda_{P\dagger}
 G'_{\dagger}(V'(!D_X'))\simeq\nbigm_{|Y_P}$.
It is enough to prove the claims
on $\check{Y}_P$.
Hence, we may and will assume that
$D_Y$ is normal crossing,
and that
$V_1$ is a good meromorphic flat bundle.

\vspace{.1in}

It is enough to consider the case where
 $Y:=\Delta^n$ and 
$D_Y:=\bigcup_{i=1}^{\ell}\{z_i=0\}$.
We have $G_{\dagger}^0(\nbigm)=V_1(!D_Y)$.
Let $F:Y\lrarr \cnum^{\ell}$ be given by
$(z_1,\ldots,z_{\ell})$.
We set $F_X:=F\circ G$.
We obtain a projective morphism
$G:(X,F_X)\lrarr (Y,F)$ in $\Cat_{\ell}$.
We have $\Ytilde(F)=\Ytilde(D_Y)$.
According to Corollary \ref{cor;13.4.20.400},
we have the following isomorphism in $D^b_c(\Ytilde(D_Y))$:
\[
 R\Gtilde_{\ast}
 \DR^{\rapid}_{\Xtilde(F_X)}(\nbigm)
\simeq
 \DR^{\rapid}_{\Ytilde(D)}(G^0_{\dagger}\nbigm)
\]
The good $K$-structure of $V$ induces 
a $K$-structure of 
$\DR^{\rapid}_{\Xtilde(D_{X2})}(\nbigm)$
on $\Xtilde(D_{X_2})$
(Lemma \ref{lem;13.4.21.31}).
It induces a $K$-structure of
$R\Gtilde_{\ast}\DR^{\rapid}_{\Xtilde(F_X)}(\nbigm)$,
which is compatible with 
the natural $K$-structure of
$G_{\dagger}^0(\nbigm)_{|Y\setminus D_Y}$.

\vspace{.1in}

Let us prove that the $K$-structure of $V_1$
is good.
First, we consider the case where $V_1$ is unramifiedly good.
Take $\gminia\in\Irr(V_1)$.
Let $L(-\gminia)$ be a meromorphic flat bundle 
with a $K$-structure
as in \S\ref{subsection;13.4.20.350}.
Then, $V\otimes G^{\ast}L(-\gminia)$
has a good $K$-structure.
By applying the previous argument,
we obtain that
$\DR^{\rapid}_{\Ytilde(D_Y)}(V_1\otimes L(-\gminia))$
has a $K$-structure,
whose restriction to $Y\setminus D_Y$
is the same as one
induced by the $K$-structure of
$V_1$ and $L(-\gminia)$.
Hence, by Lemma \ref{lem;13.4.20.250},
we obtain that 
the $K$-structure of $V_1$ is good
if $V_1$ is unramifiedly good.

\vspace{.1in}

Let us consider the case where
$V_1$ is not necessarily unramified.
Let $\kappa:Y'\lrarr Y$ be a ramified covering
such that $\kappa^{\ast}V_1$ is unramifiedly good.
We put $D'_Y:=\kappa^{-1}(D_Y)$.
We take a projective birational map
$\mu:X'\lrarr X\times_{Y}Y'$ 
such that 
(i) $X'$ is smooth,
(ii) $X'-\mu^{-1}\bigl(X\times_YD'\bigr)
\simeq X-(X\times_Y D')$.
We set
 $D_{X}':=\mu^{-1}(D_X\times_YY')$.
Let $\mu_1:X'\lrarr X$ and 
 $G':X'\lrarr Y'$ be the induced morphisms.
We have the decomposition
$D_{X}'=D'_{X1}\cup D'_{X2}$
such that $D_{X2}':=\mu_1^{-1}(D_{X2})$.
Let $\nbigm':=\mu_1^{\ast}(V)(!D_{X2}')$.
Applying the previous argument to
$G^{\prime\,0}_{\dagger}(\nbigm')$,
we obtain that the $K$-structure of $V_1$ is good
even in the ramified case.

\vspace{.1in}

Because the pre-$K$-Betti structure $\nbigg$ of
$G^0_{\dagger}\nbigm$
is induced by 
the $K$-structure of
$\DR^{\rapid}_{\Xtilde(D)}(G^0_{\dagger}\nbigm)$,
it is canonical.
Thus, the proof of Theorem \ref{thm;13.4.20.410}
is finished.
\hfill\qed

\begin{cor}
\label{cor;13.4.20.440}
Under the assumption, 
the induced $K$-structure of 
a meromorphic flat connection
$G_{\dagger}^0(\DDD\nbigm)$
is good,
and
$RG_{\ast}\DDD\nbigf_{\nbigm}$ gives 
the canonical pre-$K$-Betti structure of
$G_{\dagger}^0\DDD\nbigm$
\hfill\qed
\end{cor}

We have a variant of Theorem \ref{thm;13.4.20.410}
and Corollary \ref{cor;13.4.20.440}.
Let $\gbigh=(\vecH_1,\ldots,\vecH_N)$ 
be a sequence of hypersurface pairs
of $X$ contained in $D_X$.

\begin{thm}
\label{thm;13.4.23.110}
Suppose either
(i) $D_{X0}\subset H_{N!}$;
or 
(ii) $H_{N!}=\emptyset$
and $D_{X0}\subset H_{N\ast}$.
We also assume that
$G^i_{\dagger}\gbigp_{\gbigh}(V)=0$
unless $i=0$.
Then, the induced $K$-structure of
$G^0_{\dagger}\gbigp_{\gbigh}(V)(\ast D_Y)$ is good,
and the induced pre-$K$-Betti structure
$RG_{\ast}(\nbigf_{\gbigh,V})$
is the canonical pre-$K$-Betti structure
of $G^0_{\dagger}\gbigp_{\gbigh}(V)$.
\end{thm}
\pf
The case (i) can be proved by 
Proposition \ref{prop;13.4.23.200}
and the argument
in the proof of Theorem \ref{thm;13.4.20.410}.
The case (ii) can be obtained as the dual.
\hfill\qed

%% file: 7.1.tex
\subsection{Cells and cell functions}
\label{subsection;09.12.5.130}

Let $X$ be a complex manifold
or a smooth complex algebraic variety.
In the complex analytic case,
we use ordinary topology.
In the algebraic case, we consider Zariski topology.
In the algebraic setting,
$\nbigd$-modules are assumed to be algebraic.
An open subset $\nbigu$ is called principal
if it is the complement of a hypersurface.
\index{principal open subset}
Let $P$ be a point of $X$.
For any closed subvariety $W$ of $X$,
let $\dim_PW$ denote the dimension of 
the germ of $W$ at $P$.
Let $\nbigm$ be a holonomic $\nbigd$-module
on $X$ with $\dim_P\Supp\nbigm\leq n$.
An $n$-dimensional cell of $\nbigm$ at $P$ is a tuple
$\nbigc=(Z,U,\varphi,V)$ as follows:
\index{cell}
\begin{description}
\item[(Cell 1)]
 $\varphi:Z\lrarr X$ is a morphism of
 complex manifolds or smooth complex algebraic varieties,
 such that $P\in\varphi(Z)$ and $\dim Z=n$.
 We assume that
 there exists a neighbourhood of $X_P$
 of $P$ in $X$ such that
 $\varphi:\varphi^{-1}(X_P)\lrarr X_P$ is projective.
 We permit that 
 $Z$ may be non-connected or empty.
\item[(Cell 2)]
 $U\subset Z$ is a principal open subset
with the complementary hypersurface denoted by $D_Z$.
 We assume that 
 the restriction $\varphi_{|U}$ is an immersion,
 and that
 there exists a hypersurface $H$ of $X_P$
 such that $\varphi^{-1}(H)=D_Z\cap\varphi^{-1}(X_P)$.
\item[(Cell 3)]
 $V$ is a meromorphic flat connection
 on $(Z,D_Z)$ with morphisms
\[
 \varphi_{\dagger}(V_!)_{P}
 \lrarr\nbigm_{P}\lrarr 
 \varphi_{\dagger}(V)_{P}
\]
 such that 
 $\nbigm_{P}(\ast H)\simeq
 \varphi_{\dagger}(V)_{P}$
 for the hypersurface $H$ in (Cell 2),
where the subscript ``$P$'' means the restriction to $X_P$.
Note that we have
 $\nbigm_P(!H)\simeq
 \varphi_{\dagger}(V_{!})_{P}$,
where $V_!:=V(!D_Z)$.
 The restriction of $V$ to some connected components of $Z$
 may be $0$.
\end{description}
The cell $\nbigc$ is called good
if (i) $D_Z$ is normal crossing,
(ii) $V$ is good on $(Z,D_Z)$.
For a given holonomic $\nbigd_X$-module
$\nbigm$ and $P\in\Supp\nbigm$,
there always exists a cell for $\nbigm$ at $P$.
If $\dim_P\nbigm=1$, any cell is good.
If $\dim_P\nbigm=2$,
there always exists a good cell for $\nbigm$ at $P$,
due to Kedlaya \cite{kedlaya}.
(See also \cite{mochi6} for the algebraic case.)
In the algebraic case,
there always exists a good cell for $\nbigm$ at $P$
(see \cite{kedlaya2},  \cite{mochi6} and \cite{mochi7}).
\index{good cell}

\begin{rem}
Let $(Z,U,\varphi)$ be a tuple 
satisfying {\rm(Cell 1)} and {\rm(Cell 2)}.
If we are given a meromorphic flat connection
$V$ on $(Z,D_Z)$,
the tuple $(Z,U,\varphi,V)$ is called a cell at $P$.
\hfill\qed
\end{rem}

Let $g$ be a holomorphic or algebraic function
on $X_P$.
It is called a cell function for $\nbigc$
if $U=\varphi\bigl(\Supp\nbigm_P\setminus g^{-1}(0)\bigr)$.
\index{cell function}
For such $g$,
we obtain a description of $\nbigm_P$
as the cohomology of the complex
in the category of analytic or algebraic
holonomic $\nbigd_{X_P}$-modules:
\[
\psi^{(1)}_g\bigl(\varphi_{\dagger}(V)_P\bigr)
\lrarr
 \Xi^{(0)}_g\bigl(\varphi_{\dagger}(V)_P\bigr)
\oplus 
 \phi^{(0)}_g(\nbigm_P)
\lrarr
 \psi^{(0)}_g\bigl(\varphi_{\dagger}(V)_P\bigr) 
\]
For a given cell,
a cell function always exists
after we shrink $X_P$ and $Z$ appropriately.

\begin{rem}
Let $\nbigc$ be a cell of $\nbigm$ at $P$.
If we have a neighbourhood $X_P$ of $P$
for which {\rm (Cell 1--3)} are satisfied,
they are also satisfied
for any neighbourhood $X_P'\subset X_P$.
Hence, we do not have to be careful
with a choice of $X_P$.
\hfill\qed
\end{rem}

\subsection{Refinement and enhancement}

\index{refinement}
\index{enhancement}

Let $\nbigc'=(Z',\varphi',U',V')$
and $\nbigc=(Z,\varphi,U,V)$
be $n$-cells of $\nbigm$ at $P$.
We say that $\nbigc'$ is a refinement
of $\nbigc$, and denote $\nbigc'\mnuleq\nbigc$
if the following holds:
\begin{itemize}
\item
$\varphi'$ factors through $\varphi$
in the sense that 
there exists $\varphi_1:Z'\lrarr Z$
such that
(i) $\varphi'=\varphi\circ\varphi_1$,
(ii) $\varphi_1(U')\subset U$.
\item
$V'=\varphi_1^{\ast}V
 \otimes\nbigo_{Z'}(\ast D_{Z'})$,
where $D_{Z'}:=Z'-U'$.
\end{itemize}
In that situation,
there exist naturally induced morphisms:
\begin{equation}
\label{eq;09.12.5.1}
 \varphi_{\dagger}'(V'_!)_P\lrarr
 \varphi_{\dagger}(V_!)_P\lrarr
 \nbigm_P
 \lrarr \varphi_{\dagger}(V)_P
 \lrarr\varphi'_{\dagger}(V')_P 
\end{equation}
We say that $\nbigc'$ is a dominant refinement
of $\nbigc$
if $U'$ is dense in $U$.

\vspace{.1in}

Let $\nbigc=(Z,U,\varphi,V)$ 
be an $n$-cell
of $\nbigm$ at $P$.
We take an $n$-dimensional closed subvariety
$Z'\subset X$ such that
$\dim\bigl(
 \varphi(Z)\cap Z'\bigr)<n$.
We take a refinement of $\nbigc$
such that $\varphi(U)\cap Z'=\emptyset$.
Let $Z_1$ be a complex manifold
with a projective birational morphism
$\varphi_1:Z_1\lrarr Z'$
and a smooth open subset $U_1\subset Z_1$
such that (i) $\varphi_{1|U_1}$ is an immersion,
(ii) $Z_1-U_1$ is normal crossing
 and the pull back of a hypersurface 
 in $X$ around $P$.
We set $\Ztilde:=Z\sqcup Z_1$
and $\Utilde:=U\sqcup U_1$.
We have the induced map
$\varphitilde:\Ztilde\lrarr X$.
Let $\Vtilde$ be a meromorphic flat connection
on $\Ztilde$
such that $\Vtilde_{|Z}=V$
and $\Vtilde_{|Z_1}=0$.
Then, 
it is easy to observe that
$\nbigctilde:=
 (\Ztilde,\Utilde,\varphitilde,\Vtilde)$
is an $n$-cell of $\nbigm$,
which is called an enhancement of $\nbigc$.

\vspace{.1in}
In the following,
for a cell $\nbigc=(Z,U,\varphi,V)$,
we implicitly assume $\varphi^{-1}(X_P)=Z$
by taking a refinement of $\nbigc$.
So we omit the subscript ``$P$''
in $\varphi_{\dagger}(V_{!})_P$
and $\varphi_{\dagger}(V)_P$.

\subsection{$K$-cells
and the induced pre-$K$-Betti structure
on the nearby cycle sheaves} 

Let $\nbigf$ be a pre-$K$-Betti structure
of $\nbigm$.
Let $\nbigc=(Z,U,\varphi,V)$ 
be an $n$-cell of $\nbigm$ at $P$.
\begin{df}
We say that $\nbigf$ and $\nbigc$
are compatible
if the following holds:
\index{compatible}
\begin{itemize}
\item
The induced $K$-structure of
$V_{|U}$ is good.
(We do not assume that
$V$ is a good meromorphic flat bundle.
See {\rm\S\ref{subsection;13.4.25.200}}.)
\item
The induced morphisms
$\varphi_{\dagger}(V_{!})\lrarr
 \nbigm_P\lrarr
 \varphi_{\dagger}(V)$
are compatible with the pre-$K$-Betti structures.
(See {\rm\S\ref{subsection;13.4.22.1}}
for the canonical pre-$K$-Betti structures of
$V_!$ and $V$.)
\end{itemize}
Such a cell $\nbigc$ is called a $K$-cell
of $(\nbigm,\nbigf)$.
\index{$K$-cell}
\hfill\qed
\end{df}

It is not difficult to construct an example
of a pre-$K$-holonomic $\nbigd$-module,
for which there does not exist a $K$-cell
at some point.

\begin{lem}
Let $\nbigc=(Z,U,\varphi,V)$ be a $K$-cell of 
$(\nbigm,\nbigf)$ at $P$.
Any refinement $\nbigc'=(Z',U',\varphi',V')$
of  $\nbigc$ is also a $K$-cell.
Moreover,
the induced morphisms in {\rm(\ref{eq;09.12.5.1})}
are compatible with pre-$K$-Betti structures.
\end{lem}
\pf
It follows from Proposition \ref{prop;13.4.22.10}.
\hfill\qed

\vspace{.1in}

Let $g$ be any cell function for a $K$-cell $\nbigc$.
We observe that 
$\Xi^{(a)}_g\bigl(\varphi_{\dagger}(V)\bigr)$, 
$\psi^{(a)}_g\bigl(\varphi_{\dagger}(V)\bigr)$
and $\phi^{(a)}_g(\nbigm_P)$ are equipped with
induced pre-$K$-Betti structures.
We set $V^{a,b}_{g\star}:=
 \Pi_{\varphi^{-1}(g)\star}^{a,b}V$
for $\star=\ast,!$.
Note that $\varphi_{\dagger}\bigl(
 V_{g\star}^{a,b}\bigr)$
have the canonical pre-$K$-Betti structures.
Since 
$\Xi^{(a)}_g\bigl(\varphi_{\dagger}V\bigr)$ and 
$\psi^{(a)}_g\bigl(\varphi_{\dagger}V\bigr)$ are 
of the form
$\Ker\Bigl(
 \varphi_{\dagger}\bigl(
 V_{g!}^{a,b}\bigr)
\lrarr
 \varphi_{\dagger}\bigl(
 V_{g\ast}^{a',b'}\bigr)
 \Bigr)$,
they are equipped with induced pre-$K$-Betti
structures,
denoted by $\lefttop{D}\Xi^{(a)}_g(\varphi_{\ast}\nbigf_{V})$
and $\lefttop{D}\psi^{(a)}_g(\varphi_{\ast}\nbigf_{V})$.
We will use the following obvious
lemma implicitly.
\begin{lem}
The natural isomorphisms
\[
 \Xi^{(a)}_g\bigl(\varphi_{\dagger}(V)\bigr)
 \simeq
 \varphi_{\dagger}\bigl(
 \Xi^{(a)}_{g\circ\varphi}(V)\bigr),
\quad
 \psi^{(a)}_g\bigl(\varphi_{\dagger}V\bigr)
 \simeq
 \varphi_{\dagger}\psi^{(a)}_{g\circ\varphi}(V)
\]
are compatible with
the induced pre-$K$-Betti structures.
\hfill\qed
\end{lem}

Since $\phi^{(0)}_g(\nbigm_P)$ is 
the cohomology of the complex
$\varphi_{\dagger}V_!\lrarr
 \Xi^{(0)}_g(\varphi_{\dagger}V)
\oplus
 \nbigm
 \lrarr \varphi_{\dagger}V$,
we obtain a pre-$K$-Betti structure
of $\phi^{(0)}_g(\nbigm_P)$,
denoted by $\lefttop{D}\phi^{(0)}_g(\nbigf)$.
The tuples
$\bigl(
 \Xi^{(a)}_g(\varphi_{\dagger}V),
 \lefttop{D}\Xi^{(a)}_g(\varphi_{\ast}\nbigf_V)
 \bigr)$,
$\bigl(
 \psi^{(a)}_g(\varphi_{\dagger}V),
 \lefttop{D}\psi^{(a)}_g(\varphi_{\ast}\nbigf_V)
 \bigr)$
and 
$\bigl(
 \phi^{(a)}_g(\nbigm),
 \lefttop{D}\phi^{(a)}_g(\nbigf)
 \bigr)$
are also denoted by
$\Xi^{(a)}_g\varphi_{\dagger}(V,\nbigf_{V})$,
$\psi^{(a)}_g\varphi_{\dagger}(V,\nbigf_{V})$
and 
$\phi^{(a)}_g(\nbigm,\nbigf)$.
We will often omit to denote 
the pre-$K$-Betti structures
if there is no risk of confusion.

%% file: 7.2.tex
\subsection{Definition of $K$-Betti structure}

Let $X$ be any complex manifold,
and $P$ be any point of $X$.
Let $(\nbigm,\nbigf)$ be 
a pre-$K$-holonomic $\nbigd$-module on $X$.
Let us define the notion of $K$-Betti structure
of $\nbigm$ at $P$,
inductively on the dimension of $\Supp\nbigm$ at $P$.
\begin{df}
\label{df;09.10.4.10}
In the case $\dim_P\Supp\nbigm=0$,
a $K$-Betti structure is defined to be
a pre-$K$-Betti structure.

Let us consider the case $\dim_P\Supp\nbigm\leq n$.
We say that $\nbigf$ is a $K$-Betti structure of  
$\nbigm$ at $P$
if there exists an $n$-dimensional $K$-cell
$\nbigc_0=(Z_0,\varphi_0,U_0,V_0)$ 
of $(\nbigm,\nbigf)$ at $P$
with the following property:
\begin{itemize}
\item
 $\dim_P\Bigl(
 \bigl(\Supp\nbigm\cap X_P\bigr)
 \setminus
 \varphi_0(Z_0)\Bigr)
 <n$
 for some neighbourhood $X_P$ of $P$ in $X$.
\item
For any dominant refinement
$\nbigc\mnuleq\nbigc_0$
and any cell function $g$ for $\nbigc$,
the induced pre-$K$-Betti structure
$\lefttop{D}\phi^{(0)}_g(\nbigf)$
is a $K$-Betti structure 
of $\phi^{(0)}_g(\nbigm_P)$ at $P$.
Note that $\dim_P\phi^{(0)}_g(\nbigm)<n$.
\end{itemize}
Such an $n$-cell $\nbigc_0$
is called a bounding $n$-cell of $\nbigm$
at $P$.
\hfill\qed
\end{df}

If $\nbigc_0$ is a bounding $n$-cell of $\nbigm$,
any dominant refinement and enhancement
is also bounding $n$-cells of $\nbigm$.

\begin{df}
If $\nbigf$ is a $K$-Betti structure
of $\nbigm$ at any point of $X$,
it is called a $K$-Betti structure of $\nbigm$.
A holonomic $\nbigd$-module with a $K$-Betti structure
is called a $K$-holonomic $\nbigd$-module.
\hfill\qed
\end{df}
\index{$K$-Betti structure}
\index{$K$-holonomic $\nbigd$-module}

Morphisms of $K$-holonomic $\nbigd$-modules
$(\nbigm_1,\nbigf_1)\lrarr(\nbigm_2,\nbigf_2)$
are defined to be morphisms
of pre-$K$-holonomic $\nbigd$-modules.
The category of $K$-holonomic $\nbigd_X$-modules
is denoted by $\Hol(X,K)$.
It is a full subcategory of
the category of pre-$K$-holonomic
$\nbigd_X$-modules $\Hol^{\pre}(X,K)$ by definition.

\begin{rem}
As we will see later in 
{\rm\S\ref{section;09.11.11.3}},
for any $K$-cell
$\nbigc=(Z,U,\varphi,V)$
with a cell function $g$
at $P$,
the pre-$K$-holonomic $\nbigd$-modules
$\varphi_{\dagger}(V)$,
$\varphi_{\dagger}(V_!)$,
$\varphi_{\dagger}\Xi^{(a)}_{g\circ\varphi}(V)$,
and $\varphi_{\dagger}\psi^{(a)}(V)$
on a neighbourhood of $P$ are $K$-holonomic.
We will see that
$\Hol(X,K)$ is an abelian category
in Proposition {\rm\ref{prop;09.10.7.2}} below.
So, we may replace the condition 
in the higher dimensional case
in Definition {\rm\ref{df;09.10.4.10}}
with the following,
which is easier to check:
\begin{itemize}
\item
We say that $\nbigf$ is a $K$-Betti structure
of $\nbigm$ at $P$
if there exists an $n$-dimensional $K$-cell
$\nbigc=(Z,\varphi,V,U)$ with a cell function $g$
at $P$ such that 
the induced pre-$K$-Betti structure
$\lefttop{D}\phi^{(0)}_g(\nbigf)$
is a $K$-Betti structure
of $\phi_g^{(0)}(\nbigm)$ at $P$.
\end{itemize}
It seems convenient for the author
to begin with a stronger condition 
as in Definition {\rm\ref{df;09.10.4.10}}
for the development of the theory.
\hfill\qed
\end{rem}

\subsection{Abelian category}

It is basic to obtain the following.
\begin{prop}
\label{prop;09.10.7.2}
$\Hol(X,K)$ is abelian.
\end{prop}
\pf
Let $P$ be any point of $X$.
We use an induction on 
the dimension of $\Supp_P\nbigm$.
Let $(f_{\nbigd},f_{\nbigp}):
 (\nbigm_1,\nbigf_1)\lrarr
 (\nbigm_2,\nbigf_2)$ be a morphism of
$K$-holonomic $\nbigd$-modules.
Let us prove that
$\Ker(f_{\nbigp})$ is a $K$-Betti structure 
of $\Ker f_{\nbigd}$.

Let $n\geq \max\bigl\{
 \dim\Supp_P\nbigm_i
 \bigr\}$.
Let $\nbigc_{i,0}
=(Z_{i,0},U_{i,0},\varphi_{i,0},V_{i,0})$
$(i=1,2)$ be bounding $n$-cells 
for $\nbigm_i$ at $P$.
By considering refinement and enhancement,
we may assume that
$(Z_{1,0},U_{1,0},\varphi_{1,0})
=(Z_{2,0},U_{2,0},\varphi_{2,0})$,
which is denoted by
$(Z_0,U_0,\varphi_0)$.
We have an induced morphism
$f_{Z_0}:V_{1,0}\lrarr V_{2,0}$.
We obtain a cell
$\nbigc_0(\Ker)=\bigl(Z_0,U_0,\varphi_0,
 \Ker f_{Z_0}\bigr)$
of $\Ker f_{\nbigd}$.
The $K$-structure of $\Ker f_{\nbigd}$
is good by Lemma \ref{lem;13.4.22.20}.

Let $\nbigc(\Ker)=(Z,U,\varphi,K_Z)$ be 
a dominant refinement of
$\nbigc_0(\Ker)$.
We have refinements 
$\nbigc_i=(Z,U,\varphi,V_i)$ of $\nbigc_{i,0}$
with the induced morphism
$f_{Z}:V_1\lrarr V_2$.
We have $\Ker f_Z\simeq K_Z$.
We obtain the following commutative diagram
of pre-$K$-holonomic $\nbigd$-modules:
\[
 \begin{CD}
 \varphi_{\dagger}V_{1!} @>>>
 \nbigm_{1\,P} @>>>
 \varphi_{\dagger}V_{1}\\
 @VVV @VVV @VVV\\
 \varphi_{\dagger}V_{2!} @>>>
 \nbigm_{2\,P} @>>>
 \varphi_{\dagger}V_2
 \end{CD}
\]
Hence, the induced morphisms
$\varphi_{\dagger}K_{Z!}\lrarr
 \Ker (f_{\nbigd})_P\lrarr \varphi_{\dagger}K_Z$
are compatible with
the pre-$K$-Betti structures.
We have the following commutative diagram 
of pre-$K$-holonomic $\nbigd$-modules:
\[
\begin{CD}
 \varphi_{\dagger}\bigl(V^{a,b}_{1,g!}
 \bigr)
 @>>>
 \varphi_{\dagger}(V^{a,b}_{1,g\ast}) \\
 @VVV @VVV \\
 \varphi_{\dagger}(V^{a,b}_{2,g!})
 @>>>
 \varphi_{\dagger}(V^{a,b}_{2,g\ast})
\end{CD} 
\]
Hence, the morphisms
$\Xi^{(0)}_g\bigl(\varphi_{\dagger}V_1\bigr)
 \lrarr 
 \Xi^{(0)}_g\bigl(\varphi_{\dagger}V_2\bigr)$
and 
$\psi^{(0)}_g\bigl(\varphi_{\dagger}V_1\bigr)
 \lrarr
 \psi^{(0)}_g\bigl(
 \varphi_{\dagger}V_2\bigr)$
preserve the pre-$K$-Betti structures.
Therefore, $\phi^{(0)}_g(f_{\nbigd})$
preserves the pre-$K$-Betti structures,
i.e.,
$\lefttop{D}\phi^{(0)}_g(f_{\nbigp}):
 \lefttop{D}\phi^{(0)}_g(\nbigf_1)
\lrarr
 \lefttop{D}\phi^{(0)}_g(\nbigf_2)$
is induced.
By the assumption of the induction,
$\Ker\lefttop{D}\phi^{(0)}_g(f_{\nbigp})$
is a $K$-Betti structure.
It is easy to obtain that
$\lefttop{D}\phi^{(0)}_g\Ker f_{\nbigp}
=\Ker\lefttop{D}\phi^{(0)}_g(f_{\nbigp})$.
Then, we can conclude that 
$\bigl(
 \Ker f_{\nbigd},\Ker f_{\nbigp} 
 \bigr)$ is a $K$-holonomic $\nbigd$-module.
The claims for the cokernel and the image 
can be proved similarly.
\hfill\qed

\subsection{Dual}

\index{dual functor $\DDD$}
\begin{prop}
\label{prop;09.10.19.1}
Let $(\nbigm,\nbigf)$ be a $K$-holonomic $\nbigd_X$-module.
Then, the dual $\DDD(\nbigm,\nbigf):=
 (\DDD\nbigm,\DDD\nbigf)$ is also $K$-holonomic.
\end{prop}
\pf
Let $P$ be any point of $\Supp\nbigm$,
and let $\nbigc_0$ be a bounding $n$-cell at $P$.
Let $\nbigc=(Z,U,\varphi,V)$ be 
any refinement of $\nbigc_0$.
Let $\nbigf_V$ and $\nbigf_{V!}$
be the canonical pre-$K$-Betti structures
of $V$ and $V_!$.
Let $\nbigc^{\lor}:=(Z,U,\varphi,V^{\lor})$.
We have the induced $K$-structure of
$V^{\lor}$.
By using Proposition \ref{prop;09.10.27.1} and
Theorem \ref{thm;09.10.26.11},
we obtain that
$\DDD\nbigf_{V!}$ and $\DDD\nbigf_{V}$
are the canonical pre-$K$-Betti structures of 
$V^{\lor}$ and $V^{\lor}_!$.
Hence, we obtain that
$\nbigc^{\lor}$ and $\DDD\nbigf$
are compatible.
We also obtain that
$\DDD \lefttop{D}\Xi^{(a)}_g\varphi_{\ast}\nbigf_V$
is equal to the canonical pre-$K$-Betti
structure of $\Xi^{(-a-1)}_g\varphi_{\ast}V^{\lor}$.
Moreover,
the induced $K$-structure of
$\phi^{(a)}_g(\DDD\nbigm_P)$ is equal to 
$\DDD\lefttop{D}\phi^{(-a-1)}_g\nbigf$
under the isomorphism
$\phi^{(a)}_g\DDD\nbigm_P\simeq
 \DDD\phi^{(-a-1)}_g\nbigm_P$.
By the inductive assumption,
it is $K$-Betti structure.
Thus, we obtain that
$\DDD(\nbigm,\nbigf)$ is $K$-holonomic.
\hfill\qed

\subsection{Sub-objects and quotient objects} 

Let $(\nbigm,\nbigf)$ be a $K$-holonomic $\nbigd$-module.

\begin{prop}
\label{prop;09.10.4.110}
If $(\nbigm_1,\nbigf_1)$ is a subobject
of $(\nbigm,\nbigf)$ in $\Hol^{\pre}(X,K)$,
it is also $K$-holonomic.
A similar claim holds for quotients.
\end{prop}
\pf
Let $P$ be any point of $X$.
We use an induction on the dimension
of the support of $\nbigm$.
Let $n\geq \dim_P\Supp\nbigm$.
Let $\nbigc=(Z,U,\varphi,V)$ be 
a bounding $n$-cell of $\nbigm$ at $P$.
Let $V_1\subset V$ denote the sub-connection
induced by $\nbigm_1$.
Then, $\nbigc_1=(Z,U,\varphi,V_1)$
is an $n$-cell of $\nbigm_1$ at $P$.
Let us prove that $\nbigc_1$ and $\nbigf_1$
are compatible.
By Lemma \ref{lem;13.4.23.2},
the $K$-structure of $V_1$ is good.
Let $\nbigf_{\ast}$ and $\nbigf_!$
denote the canonical $K$-structures of
$\varphi_{\dagger}V$ and $\varphi_{\dagger}V_{!}$.
Let $\nbigf_{1\ast}$ and $\nbigf_{1!}$ denote
the canonical $K$-structures of
$\varphi_{\dagger}V_1$ and $\varphi_{\dagger}V_{1!}$.
We have the following morphisms:
\[
 \begin{CD}
 \varphi_{\dagger}(V_!)
 @>>>
 \nbigm @>>>
 \varphi_{\dagger}(V)\\
 @AAA @AAA @AAA \\
 \varphi_{\dagger}(V_{1!})
 @>>> \nbigm_1 @>>>
 \varphi_{\dagger}(V_{1})
 \end{CD}
\quad\quad\quad
 \begin{CD}
 \nbigf_{!} @>>>
 \nbigf @>>> \nbigf_{\ast}\\
 @AAA @AAA @AAA \\
 \nbigf_{1!} @.
 \nbigf_1 @.
 \nbigf_{1\ast}
 \end{CD}
\]
Because the morphism
$\varphi_{\dagger}(V_{1!})\lrarr
 \nbigm/\nbigm_1$ is $0$,
the morphism
$\nbigf_{1!}\lrarr 
 \nbigf/\nbigf_1$
is also $0$,
i.e.,
$\nbigf_{1!}\lrarr\nbigf$
factors through $\nbigf_1$.
Similarly, 
we obtain that
$\nbigf_1\lrarr\nbigf_{\ast}$
factors through $\nbigf_{1\ast}$.
Hence, $\nbigc_1$ is compatible
with $\nbigf_1$.

Let $f$ be a cell function for $\nbigc$.
We have 
$\lefttop{D}\Xi^{(a)}_f(\nbigf)
\supset
 \lefttop{D}\Xi^{(a)}_f(\nbigf_1)$
and
$\lefttop{D}\psi^{(a)}_f(\nbigf)
\supset
 \lefttop{D}\psi^{(a)}_f\nbigf_1$.
Hence, we obtain
$\lefttop{D}\phi^{(a)}_f(\nbigf)\supset
 \lefttop{D}\phi^{(a)}_f(\nbigf_1)$,
which are pre-$K$-Betti structures of
$\phi^{(a)}_f\nbigm$ and $\phi^{(a)}_f\nbigm_1$.
By the assumption of the induction,
we obtain that
$\lefttop{D}\phi^{(a)}_f(\nbigf_1)$
is a $K$-Betti structure
of $\phi_f\nbigm_1$.
\hfill\qed

\subsection{Twist}

Let $(\nbigm,\nbigf)$ be a $K$-holonomic
$\nbigd$-module on $X$.
Let $\nbigv$ be a flat bundle on $X$
with a $K$-structure,
i.e.,
we have a $K$-local system $\nbigf_{\nbigv}$
such that $\nbigf_{\nbigv}\otimes\cnum[\dim X]
\simeq \DR_X(V)$.
Then, we obtain a pre-$K$-Betti structure
$\nbigf\otimes\nbigf_{\nbigv}$ of 
$\nbigm\otimes \nbigv$.

\begin{lem}
\label{lem;09.11.2.1}
$\nbigf\otimes\nbigf_{\nbigv}$ is 
a $K$-Betti structure of 
$\nbigm\otimes \nbigv$.
\end{lem}
\pf
Let $P$ be any point of $X$.
We use an induction on $\dim_P\Supp\nbigm$.
Let $\nbigc=(Z,U,\varphi,V)$ be a $K$-cell of $\nbigm$
at $P$.
Then, $\nbigc'=
 \bigl(Z,U,\varphi,V\otimes\varphi^{\ast}\nbigv\bigr)$
is a $K$-cell of $\nbigm\otimes\nbigv$ at $P$.
Let $g$ be a cell function of $\nbigc$.
Then, we have natural isomorphism
of pre-$K$-holonomic $\nbigd_X$-modules
$\psi^{(a)}_{g}\bigl(
 \varphi_{\dagger}(V\otimes\varphi^{\ast}\nbigv)\bigr)
\simeq
\psi^{(a)}_g\bigl(\varphi_{\dagger}(V)\bigr)
 \otimes\nbigv$
and 
$\Xi^{(a)}_{g}\bigl(
 \varphi_{\dagger}(V\otimes\varphi^{\ast}\nbigv)\bigr)
\simeq
 \Xi^{(a)}_g\bigl(\varphi_{\dagger}(V)\bigr)
 \otimes\nbigv$.
Hence, we obtain an isomorphism of 
pre-$K$-holonomic $\nbigd$-modules
$\phi^{(a)}_g(\nbigm\otimes \nbigv)
\simeq
 \phi^{(a)}_g(\nbigm)\otimes\nbigv$.
By using the inductive assumption,
we obtain that 
$\phi^{(a)}_g(\nbigm\otimes\nbigv)$
is $K$-holonomic.
Hence, we obtain that
$\nbigm\otimes\nbigv$ is $K$-holonomic
at $P$.
\hfill\qed

\subsection{$K$-cells}

\begin{prop}
\label{prop;13.4.23.10}
Let $(\nbigm,\nbigf)$ be a $K$-holonomic
$\nbigd$-module.
Then, 
any cell $\nbigc=(Z,U,\varphi,V)$ of $\nbigm$
is a $K$-cell.
\end{prop}
\pf
Let $P$ be any point of $\Supp(\nbigm)$.
Let $\nbigc'_P=(Z'_P,U'_P,\varphi'_P,V'_P)$ be a 
bounding $K$-cell of $\nbigm$ at $P$,
which is a refinement of $\nbigc$.
By Lemma \ref{lem;13.1.15.20},
we obtain that
the induced $K$-structure of $V$ is good
around $\varphi^{-1}(P)$.
By varying $P$,
we obtain that the $K$-structure of $V$
is good.
Moreover, for $P$ and $\nbigc'_P$ as above,
the induced morphisms
$\nbigm_P\lrarr \varphi'_{P\dagger}V'_P$
and $\varphi_{\dagger}(V)_P\lrarr\varphi'_{P\dagger}V'_P$
are compatible with pre-$K$-Betti structures,
where $\varphi_{\dagger}(V)_P$ denotes the restriction to
a small neighbourhood of $P$.
Because $\varphi_{\dagger}(V)_P\lrarr\varphi'_{P\dagger}V'_P$
is a monomorphism,
we obtain that $\nbigm_P\lrarr\varphi_{\dagger}(V)_P$
is also compatible with pre-$K$-Betti structures.
By varying $P$ in $X$,
we obtain that
$\nbigm_P\lrarr\varphi_{\dagger}(V)$
is also compatible with pre-$K$-Betti structures.
We can prove that 
$\varphi_!V\lrarr\nbigm$ is also 
compatible with pre-$K$-Betti structures
with a similar argument.
\hfill\qed

%% file: 7.3.tex
We introduce a variant notion of
$K(\ast D)$-Betti structure
of holonomic $\nbigd_{X(\ast D)}$-modules,
where $D$ is a hypersurface.
It is rather auxiliary.
Indeed, as proved in \S\ref{section;09.11.11.3},
it is equivalent to $K$-Betti structure
for holonomic $\nbigd_{X(\ast D)}$-modules,
although it will be convenient in some arguments.

\subsection{Cells and cell functions
 for holonomic $\nbigd_{X(\ast D)}$-modules}

Let $X$ be any complex manifold
or smooth complex algebraic variety,
and let $D$ be any hypersurface of $X$.
Let $\nbigm$ be any holonomic 
$\nbigd_{X(\ast D)}$-module,
i.e.,
$\nbigm$ is a holonomic $\nbigd_X$-module
such that $\nbigm(\ast D)=\nbigm$.
A cell of a holonomic 
$\nbigd_{X(\ast D)}$-module $\nbigm$ 
is defined to be a cell of 
a holonomic $\nbigd_X$-module $\nbigm$.
\index{cell}
The notions of refinement and enhancement
of a cell of a holonomic $\nbigd_{X(\ast D)}$-module
are defined in the same way.
\index{refinement}
\index{enhancement}
However, 
we will be interested in the morphisms
$\varphi_{\dagger}(V_!)(\ast D)\lrarr
 \nbigm_P\lrarr
\varphi_{\dagger}V$.

The notion of cell functions is modified.
\index{cell function}
Let $\nbigc=(Z,U,\varphi,V)$ be a cell of 
a holonomic $\nbigd_{X(\ast D)}$-module $\nbigm$.
A cell function $g$ of $\nbigc$ is a meromorphic
function on $X$
whose poles are contained in $D$,
such that $U=\Supp\nbigm\setminus 
 \bigl(g^{-1}(0)\cup D\bigr)$.

\subsection{$K(\ast D)$-cell}

Let $\nbigm$ be a holonomic $\nbigd_{X(\ast D)}$-module.
Let $\nbigf$ be a pre-$K$-Betti
structure of $\nbigm$.
Let $\nbigc=(Z,U,\varphi,V)$ be an $n$-cell
of $\nbigm$ at $P$.
We say that $\nbigf$ and $\nbigc$ are compatible
if (i) the induced $K$-structure of $V$ is good,
(ii) the induced morphisms
$\varphi_{\dagger}(V_!)(\ast D)\lrarr
 \nbigm_P\lrarr
 \varphi_{\dagger}(V)$
are compatible with
the pre-$K$-Betti structures.
Such a cell $\nbigc$ is called a $K(\ast D)$-cell
of $(\nbigm,\nbigf)$.
\index{$K(\ast D)$-cell}
Note that 
the condition (i) implies that
$\varphi_{\dagger}(V_!)(\ast D)$
and
$\varphi_{\dagger}(V)$
are equipped with the canonical pre-$K$-Betti structure.

Let $g$ be a cell function for a $K(\ast D)$-cell
$\nbigc$.
We set $V^{a,b}_{g\star}(\ast D)\!\!:=\!\!
 \bigl(
 V\otimes\gbigi_{\varphi^{-1}(g)}^{a,b}
 \bigr)(\ast\varphi^{-1}D)$
for $\star=\ast,!$.
Note that $\varphi_{\dagger}
 \bigl(V_{g\star}^{a,b}(\ast D)
 \bigr)$
have the canonical pre-$K$-Betti structures.
Since 
$\Xi^{(c)}_g\bigl(\varphi_{\dagger}V,\ast D\bigr)$ and 
$\psi^{(c)}_g\bigl(\varphi_{\dagger}V,\ast D\bigr)$ are 
of the form
\[
 \Ker\Bigl(
 \varphi_{\dagger}\bigl(
 V_{g!}^{a,b}(\ast D)\bigr)
\lrarr
 \varphi_{\dagger}
 \bigl(V_{g\ast}^{a',b'}(\ast D)\bigr)
 \Bigr),
\]
they are equipped with the induced
pre-$K$-Betti structures
$\lefttop{D}
 \Xi^{(c)}_g(\varphi_{\ast}\nbigf_{V},\!\ast D)$
and $\lefttop{D}
 \psi^{(c)}_g(\varphi_{\ast}\nbigf_{V},\!\ast D)$.
The tuples
\[
 \bigl(
 \Xi^{(c)}_g(\varphi_{\dagger}V,\ast D),
 \lefttop{D}\Xi_g(\varphi_{\ast}\nbigf_V,\ast D)
 \bigr),
\quad
 \bigl(
 \psi^{(c)}_g(\varphi_{\dagger}V,\ast D),
 \lefttop{D}\psi_g(\varphi_{\ast}\nbigf_V,\ast D)
 \bigr)
\]
are also denoted by
$\Xi^{(c)}_g\varphi_{\dagger}(V,\nbigf_{V},\ast D)$
and
$\psi^{(c)}_g\varphi_{\dagger}(V,\nbigf_{V},\ast D)$.
We will often omit to denote the pre-$K$-Betti structures.
We will use the following obvious lemma implicitly.
\begin{lem}
The natural isomorphisms
\[
 \Xi^{(a)}_g\bigl(\varphi_{\dagger}V,\ast D\bigr)
 \simeq
 \varphi_{\dagger}\Xi^{(a)}_g(V,\ast \varphi^{-1}D),
\quad\quad
\psi^{(a)}_g\bigl(\varphi_{\dagger}V,\ast D\bigr)
 \simeq
 \varphi_{\dagger}\psi^{(a)}_g(V,\ast \varphi^{-1}D)
\]
are compatible with
the induced pre-$K$-Betti structures.
\hfill\qed
\end{lem}

Since $\phi^{(0)}_g(\nbigm_P,\ast D)$
is the cohomology of the complex
in the category of pre-$K$-holonomic $\nbigd_X$-modules
\[
\varphi_{\dagger}(V_!)(\ast D)\lrarr
 \Xi^{(0)}_g(\varphi_{\dagger}V,\ast D)
\oplus
 \nbigm_P
 \lrarr \varphi_{\dagger}(V)(\ast D), 
\]
we obtain a pre-$K$-Betti structure
of $\phi^{(a)}_g(\nbigm_P,\ast D)$
denoted by 
$\lefttop{D}\phi^{(a)}_g(\nbigf,\ast D)$.
Let 
$\phi^{(a)}_g(\nbigm_P,\nbigf,\ast D)$
denote the tuple
$\bigl(
 \phi^{(a)}_g(\nbigm_P,\ast D),
 \lefttop{D}\phi^{(a)}_g(\nbigf,\ast D)
 \bigr)$.
We will often omit to denote the pre-$K$-Betti
structure.

\subsection{Definition of $K(\ast D)$-Betti structure}

Let us define the notion of $K(\ast D)$-Betti structure
at any point of $D$,
inductively on the dimension of the support of
$\nbigd_{X(\ast D)}$-modules.
Let $(\nbigm,\nbigf)$ be a pre-$K$-holonomic
$\nbigd_{X(\ast D)}$-module.
Note that
we have $\nbigm=0$ around $P\in D$
in the case $\dim_P\Supp\nbigm=0$.

\begin{df}
Let $P$ be any point of $D$.
Suppose $\dim_P\Supp\nbigm\leq n$.
We say that $\nbigf$ is a $K(\ast D)$-Betti structure of  
$\nbigm$ at $P$
if there exists an $n$-dimensional
$K(\ast D)$-cell
$\nbigc_0=(Z_0,\varphi_0,U_0,V_0)$ at $P$
with the following property:
\begin{itemize}
\item
 $\dim_P\Bigl(
 \bigl(\Supp\nbigm\cap X_P\bigr)
 \setminus
 \varphi_0(Z_0)\Bigr)
 <n$
 for some neighbourhood $X_P$ of $P$ in $X$.
\item
For any dominant refinement
$\nbigc\mnuleq\nbigc_0$
and any cell function $g$ for $\nbigc$
as a $\nbigd_{X(\ast D)}$-module,
the induced pre-$K$-Betti structure
$\lefttop{D}\phi^{(0)}_g(\nbigf,\ast D)$
is a $K(\ast D)$-Betti structure at $P$.
\end{itemize}
Such an $n$-cell $\nbigc_0$
is called a bounding $n$-cell of $\nbigm$
at $P$.
\hfill\qed
\end{df}

If $\nbigc_0$ is a bounding $n$-cell of $\nbigm$,
its dominant refinements and enhancements
are also bounding $n$-cells of $\nbigm$.

\begin{df}
A pre-$K$-Betti structure $\nbigf$ of $\nbigm$
is called  a $K(\ast D)$-Betti structure
if it is $K$-Betti structure
of $\nbigm$ at any points of $X\setminus D$,
and if it is $K(\ast D)$-Betti structure
of $\nbigm$ at any points of $D$.
A holonomic $\nbigd_{X(\ast D)}$-module with 
a $K(\ast D)$-Betti structure
is called a $K(\ast D)$-holonomic 
$\nbigd_{X(\ast D)}$-module.
\hfill\qed
\end{df}
\index{$K(\ast D)$-Betti structure}
\index{$K(\ast D)$-holonomic $\nbigd_{X(\ast D)}$-module}

Let $\Hol\bigl(X,\ast \!D\!,\!K\bigr)\subset
 \Hol^{\pre}(X,K)$
denote the full subcategory of
$K(\ast D)$-holonomic 
$\nbigd_{X(\ast D)}$-modules.
\index{category $\Hol(X,\ast D,K)$}
The following lemma is similar to 
Proposition \ref{prop;09.10.7.2}.
\begin{lem}
The category $\Hol\bigl(X,\ast D,K\bigr)$
is abelian.
\hfill\qed
\end{lem}

The following lemma is similar to
Proposition \ref{prop;09.10.4.110}.

\begin{lem}
\label{lem;09.10.18.511}
Let $(\nbigm,\nbigf)$ be any $K(\ast D)$-holonomic
$\nbigd_X$-module.
Any subobject of $(\nbigm,\nbigf)$
in $\Hol^{\pre}(X,K)$ is also $K(\ast D)$-holonomic.
A similar claim holds for quotients.
\hfill\qed
\end{lem}

The following lemma is analogue of 
Proposition \ref{prop;13.4.23.10}.
\begin{lem}
\label{lem;13.4.23.11}
Let $(\nbigm,\nbigf)$ be a $K(\ast D)$-holonomic
$\nbigd_{X(\ast D)}$-module.
Then, 
any cell $\nbigc=(Z,U,\varphi,V)$ of $\nbigm$
is a $K(\ast D)$-cell.
\hfill\qed
\end{lem}

\subsection{Uniqueness}
We have the following uniqueness.

\begin{prop}
\label{prop;09.10.18.250}
Let $\nbigm$ be a holonomic 
$\nbigd_{X(\ast D)}$-module
with 
$K(\ast D)$-Betti structures 
 $\nbigf_i$ $(i=1,2)$.
If $\nbigf_{1|X-D}=\nbigf_{2|X-D}$,
then we have $\nbigf_1=\nbigf_2$.
\end{prop}
\pf
It is enough to consider the issue
locally around any point $P\in D$.
We use an induction on 
$\dim_P\Supp\nbigm$.
In the case $\dim_P\Supp\nbigm=0$,
the claim is clear.
Suppose $\dim_P\Supp\nbigm\leq n$.
Let $\nbigc$ be any bounding cell at $P$,
and let $g$ be any cell function of $\nbigc$.
Let $\lefttop{D}\phi^{(0)}_g(\nbigf_i,\ast D)$ be
the induced pre-$K(\ast D)$-Betti structures
of $\phi^{(0)}_g(\nbigm,\ast D)$.
By the assumption of the induction,
we have 
$\lefttop{D}\phi^{(0)}_g(\nbigf_1,\ast D)
=\lefttop{D}\phi^{(0)}_g(\nbigf_2,\ast D)$.
Because $\nbigf_i$ can be reconstructed from
$\lefttop{D}\phi^{(0)}_g(\nbigf_i,\ast D)$
and the canonical pre-$K(\ast D)$-Betti structures
of $\psi^{(a)}_g(\varphi_{\ast}V,\ast D)$ and 
$\Xi^{(a)}_g(\varphi_{\ast}V,\ast D)$,
we obtain $\nbigf_1=\nbigf_2$.
\hfill\qed

\subsection{Independence from a compactification}

Let $F:X'\lrarr X$ be a projective birational morphism
of complex manifolds
such that $X'-D'\simeq X-D$,
where $D':=F^{-1}(D)$.
Recall that $F_{\dagger}$ denotes
the push-forward of pre-$K$-holonomic $\nbigd$-modules.

\begin{prop}
\label{prop;09.10.19.2}
\mbox{{}}
The functor $F_{\dagger}$
induces an equivalence of
the categories 
$\Hol\bigl(X,\ast D,K\bigr)$ 
and $\Hol\bigl(X',\ast D',K\bigr)$.
\end{prop}
\pf
It is enough to check the claims
locally around any $P\in D$.
We begin with a remark.
Let $\nbigm'$ be a holonomic $\nbigd_{X'(\ast D')}$-module.
We set $\nbigm:=F_{\dagger}\nbigm$.
Let $\nbigc=(Z,U,\varphi,V)$ be a cell
of $\nbigm$ at $P$.
By taking a refinement,
we may assume that $\varphi$
factors through $F$,
i.e.,
$\varphi=F\circ\varphi'$,
and that
$\nbigc'=(Z,U,\varphi',V)$ is 
a cell of $\nbigm'$.
Let $g$ be a cell function for $\nbigc$
as a $\nbigd_{X(\ast D)}$-module.
Note that $g'=g\circ F$ is a cell function for $\nbigc'$.
We have a description of $\nbigm'$
as the cohomology of the following complex:
\begin{equation}
 \label{eq;09.10.18.100}
 \psi^{(1)}_{g'}(\varphi'_{\dagger}V,\ast D')
\lrarr
 \Xi^{(0)}_{g'}(\varphi'_{\dagger}V,\ast D')
\oplus 
 \phi^{(0)}_{g'}(\nbigm',\ast D')
\lrarr
 \psi^{(0)}_{g'}(\varphi'_{\dagger}V,\ast D')
\end{equation}
By the push-forward $F_{\dagger}$,
it induces a description of $\nbigm$
as the cohomology of the following complex:
\begin{equation}
\label{eq;13.4.23.12}
 \psi^{(1)}_{g}(\varphi_{\dagger}V,\ast D)
\lrarr
 \Xi^{(0)}_{g}(\varphi_{\dagger}V,\ast D)
\oplus 
 \phi^{(0)}_{g}(\nbigm,\ast D)
\lrarr
 \psi^{(0)}_{g}(\varphi_{\dagger}V,\ast D)
\end{equation}

Suppose that $\nbigf'$ is a $K(\ast D)$-Betti structure
of $\nbigm'$.
Let us prove that $F_{\dagger}\nbigf'$
is a $K(\ast D)$-Betti structure of $\nbigm$.
By Lemma \ref{lem;13.4.23.11},
$\nbigc'$ is a $K(\ast D)$-cell of $\nbigm'$.
We obtain that $\nbigc$ is a $K(\ast D)$-cell
of $\nbigm$.
Because the pre-$K$-holonomic $\nbigd$-module
$\phi^{(0)}_g(\nbigm,\ast D)$ is 
obtained as 
$F_{\dagger}\phi^{(0)}_{g'}(\nbigm',\ast D)$,
we obtain that
$\phi^{(0)}_g(\nbigm,\ast D)$ is 
$K(\ast D)$-holonomic
by the inductive assumption.
Hence, $\nbigf$ is also a $K(\ast D)$-Betti structure.
Thus, $F_{\dagger}$ induces a functor
$\Hol(X',\ast D',K)\lrarr \Hol(X,\ast D,K)$.
It is clearly faithful.

Let us prove that it is full.
We use an induction on the dimensions 
of the supports of the holonomic $\nbigd$-modules.
Let $(\nbigm'_i,\nbigf'_i)$ $(i=1,2)$
be objects in $\Hol(X',\ast D',K)$.
Let $f:F_{\dagger}(\nbigm_1',\nbigf_1')\lrarr 
 F_{\dagger}(\nbigm_2',\nbigf_2')$
be a morphism in 
$\Hol(X,\ast D,K)$.
We have a morphism
$f':\nbigm'_1\lrarr\nbigm_2'$
of holonomic $\nbigd_{X'(\ast D')}$-modules.
It is enough to show that it is compatible
with the $K(\ast D)$-Betti structures.
For the cohomological descriptions 
(\ref{eq;09.10.18.100}) for $\nbigm_i'$,
$\psi^{(a)}_{g'}(f')$and $\Xi^{(a)}_{g'}(f')$
are compatible with the pre-$K$-Betti structures.
Because
$\phi^{(a)}_{g}(f)$ is compatible with 
the $K(\ast D)$-Betti structures,
we obtain that $\phi^{(a)}_{g'}(f')$ is compatible with
the $K(\ast D')$-Betti structures.
Thus, we obtain that $f'$ is compatible
with the $K(\ast D')$-Betti structures.

Let us prove the essential surjectivity.
We use an induction on the dimension of the support.
Let $\nbigm$ and $\nbigm'$ be as above.
Let $\nbigf$ be a $K(\ast D)$-Betti structure of $\nbigm$.
By the inductive assumption,
the $K(\ast D)$-Betti structure of
$\psi^{(a)}_g\bigl(\varphi_{\dagger}(V),\ast D\bigr)$
and $\phi^{(a)}_g(\nbigm,\ast D)$
induce $K(\ast D)$-Betti structures of
$\psi^{(a)}_{g'}\bigl(\varphi'_{\dagger}(V),\ast D'\bigr)$
and $\phi^{(a)}_{g'}(\nbigm',\ast D')$,
which are compatible with the natural morphisms.
We also have the canonical $K$-Betti structures
of $\psi^{(a)}_{g'}\bigl(\varphi'_{\dagger}(V),\ast D'\bigr)$
and $\Xi^{(a)}_{g'}\bigl(\varphi'_{\dagger}V,\ast D'\bigr)$.
By Proposition \ref{prop;09.10.18.250},
the induced $K(\ast D)$-Betti structures on
$\psi^{(a)}_{g'}\bigl(\varphi'_{\dagger}(V),\ast D'\bigr)$
are the same.
Hence, (\ref{eq;09.10.18.100})
is a complex of $K(\ast D)$-holonomic 
$\nbigd(\ast D)$-modules.
Hence, we have an induced $K(\ast D)$-Betti structure
of $\nbigm'$.
The functoriality is clear from the above construction.
\hfill\qed

%% file: 8.1.tex
We give several statements.

\begin{thm}
\label{thm;09.10.16.5}
Let $F:X\lrarr Y$ be any projective morphism
of complex manifolds.
For any $K$-holonomic $\nbigd_X$-module
$(\nbigm,\nbigf)$,
the push-forward $F_{\dagger}^i(\nbigm,\nbigf)
:=\bigl(F^i_{\dagger}\nbigm,F^i_{\dagger}\nbigf\bigr)$
are also $K$-holonomic for any $i$.
\end{thm}
Here, $F^i_{\dagger}\nbigf$
denotes the $i$-th cohomology of $RF_{\ast}\nbigf$
with respect to the middle perversity.

\begin{thm}
\label{thm;09.10.17.151}
Let $X$ be any complex manifold
with a normal crossing hypersurface $D$.
Any good pre-$K$-holonomic $\nbigd$-module
on $(X,D)$
is $K$-holonomic.
\end{thm}
See Definition {\rm\ref{df;13.4.27.22}}
for good pre-$K$-holonomic $\nbigd$-modules.

\begin{thm}
\label{thm;13.4.23.20}
Let $X$ be a complex manifold
with a hypersurface $D$.
Let $\gbigh$ be a sequence of hypersurface pairs
contained in $D$.
Let $V$ be any meromorphic flat connection on $(X,D)$
with a good $K$-structure.
Then, the pre-$K$-holonomic $\nbigd$-module
$\gbigp_{\gbigh}(V)$ is $K$-holonomic.
\end{thm}
See \S\ref{subsection;13.4.25.200}
for hypersurface pairs and $\gbigp_{\gbigh}(V)$.

\begin{thm}
\label{thm;09.10.18.400}
Let $X$ be any complex manifold with a hypersurface $D$.
We have a unique functor
$\Hol(X,K)\lrarr\Hol(X,\ast D,K)$
with the following properties:
\begin{itemize}
\item
It is an enhancement of
the functor $\Hol(X)\lrarr\Hol(X,\ast D)$
given by $\nbigm\longmapsto\nbigm(\ast D)$.
\item
For any $(\nbigm,\nbigf)\in\Hol(X,K)$,
the natural morphism
$\nbigm\lrarr\nbigm(\ast D)$
is compatible with the induced pre-$K$-Betti structures.
\end{itemize}
\end{thm}

\subsection{Auxiliary statements}

We will use an induction on the dimension
of the supports of $\nbigd$-modules
for the proof.
Let $SI(\leq\!n)$ 
denote the statement of Theorem \ref{thm;09.10.16.5}
in the case $\dim\Supp\nbigm\leq n$.
Let $GOOD(\leq n)$ means the following:
\begin{itemize}
\item
 The claim of Theorem \ref{thm;09.10.17.151}
 holds if $\dim\Supp\nbigm\leq n$.
\item
 The claim of Theorem \ref{thm;13.4.23.20}
 holds if $\dim X\leq n$.
\end{itemize}
 For any complex manifold $X$
 with a hypersurface $D$,
 let $\Hol_{\leq n}(X,K)\subset\Hol(X,K)$
 denote the full subcategory
 of $K$-holonomic $\nbigd_X$-modules
 $(\nbigm,\nbigf)$ with $\dim\Supp\nbigm\leq n$.
 We use the symbols
 $\Hol_{\leq n}(X)$, $\Hol_{\leq n}(X,\ast D)$
 and $\Hol_{\leq n}(X,\ast D,K)$
with a similar meaning.
Let $LOC(\leq \!n)$ means the following:
\begin{itemize}
\item
 The claim of Theorem \ref{thm;09.10.18.400}
 holds
 if we replace $\Hol(X,K)$, $\Hol(X,\ast D,K)$, etc.,
 by $\Hol_{\leq n}(X,K)$, $\Hol_{\leq n}(X,\ast D,K)$, etc..
\end{itemize}

Our induction will proceed as follows:
\begin{itemize}
\item
 $SI(<n)+GOOD(<n)\Longrightarrow
 GOOD(\leq n)$
 (\S\ref{subsection;09.12.17.1} and
 \S\ref{subsection;13.4.27.100}).
\item
 $SI(<n)+GOOD(\leq n)+LOC(<n)\Longrightarrow
 LOC(\leq n)$
 (\S\ref{subsection;09.12.5.10}).
\item
 $SI(<n)+GOOD(\leq n)+LOC(\leq n)\Longrightarrow
 SI(\leq n)$
 (\S\ref{subsection;09.10.17.150}).
\end{itemize}

\begin{rem}
In the proof,
we will observe the equivalence of
$K(\ast D)$-Betti structure
and $K$-Betti structure.
(See Lemma {\rm\ref{lem;09.10.18.510}}.)
\hfill\qed
\end{rem}

%% file: 8.2.tex
\subsection{$K$-cell}
\label{subsection;09.12.5.20}

Let $\varphi:Z\lrarr X$ be a projective morphism
of complex manifolds such that $\dim Z=n$.
Let $D_Z$ be a hypersurface of  $Z$.
Assume that  $\varphi_{|Z-D_Z}$ is an immersion.
Let $V$ be a meromorphic flat connection on $(Z,D_Z)$
with a good $K$-structure.
We have the canonical pre-$K$-Betti structures 
$\nbigf_V$ and $\nbigf_{V!}$
of $V$ and $V(!D_Z)$, respectively.
More generally,
for any sequence of hypersurface pairs
$\gbigh$ contained in $D_Z$,
we obtain the canonical pre-$K$-holonomic 
$\nbigd_Z$-modules
$\gbigp_{\gbigh}(V)$.
Note that
the natural morphisms
$V(!D_Z)\lrarr
\gbigp_{\gbigh}(V)
\lrarr
 V$
are compatible with
the pre-$K$-Betti structures.
Hence, we can regard
$(Z,U,\id,V)$
as a $K$-cell of
$\gbigp_{\gbigh}(V)$.

\begin{lem}
\label{lem;13.4.23.20}
Suppose $SI(<n)$ and $GOOD(<n)$.
Let $g$ be any cell function for 
$\nbigc_0=(Z,U,\varphi,V)$.
We set $g_Z:=g\circ\varphi$.
The pre-$K$-holonomic 
$\phi^{(a)}_{g_Z}\bigl(\gbigp_{\gbigh}(V)\bigr)$
and 
$\phi^{(a)}_{g}\bigl(
 \varphi_{\dagger}\gbigp_{\gbigh}(V)
 \bigr)$
are $K$-holonomic.
In particular,
$\psi^{(a)}_{g_Z}(V)$ and
$\psi^{(a)}_{g}\varphi_{\dagger}(V)$
are $K$-holonomic.
\end{lem}
\pf
By $SI(<n)$,
it is enough to prove that
$\phi^{(0)}_{g_Z}\bigl(\gbigp_{\gbigh}(V)\bigr)$ 
is $K$-holonomic. 
It is enough to consider the issue  
locally around any point $P\in D_Z$. 
We take a local resolution 
$(Z_P,\lambda_P)$ of $V$.
We put 
$\check{g}_{P}:=g_Z\circ\lambda_P$.
We set $\check{\gbigh}_P:=\lambda_P^{-1}(\gbigh)$
and $\check{V}_P:=\lambda_P^{\ast}V$.
We have the good pre-$K$-holonomic
 $\nbigd_{\check{Z}_P}$-module
$\phi^{(0)}_{\check{g}_P}
 \bigl(
  \gbigp_{\check{\gbigh}_P}(\check{V}_P)
\bigr)$
(Proposition \ref{prop;13.4.27.30}).
By $GOOD(<n)$,
it is $K$-holonomic.
By $SI(<n)$,
$\lambda_{P\dagger}
 \phi^{(0)}_{\check{g}_P}
 \bigl(
  \gbigp_{\check{\gbigh}_P}(\check{V}_P)
\bigr)$
is $K$-holonomic,
which means that
$\phi^{(0)}_{g_Z}\bigl(\gbigp_{\gbigh}(V)\bigr)$ 
is $K$-holonomic at $P$.
\hfill\qed

\begin{prop}
\label{prop;09.10.18.200}
Suppose that $SI(<n)$ and $GOOD(<n)$ hold.
Then, the pre-$K$-holonomic $\nbigd$-modules
$\varphi_{\dagger}(V,\nbigf_V)$
and 
$\varphi_{\dagger}(V_!,\nbigf_{V!})$
are $K$-holonomic.
\end{prop}
\pf
Let us prove the claim for 
$\varphi_{\dagger}(V,\nbigf_V)$.
The other can be proved as the dual.
Let us prove that 
$\nbigc_0=(Z,U,\varphi,V)$ 
is a bounding $n$-cell for $\varphi_{\dagger}(V,\nbigf_V)$.
Let $P$ be any point of $X$.
Let $\nbigc'=(Z',U',\varphi',V')$
be a dominant refinement at $P$
with a cell function $g$.
We have a factorization
$\varphi'=\varphi\circ\varphi_1$,
where $\varphi_1:Z'\lrarr Z$.
We put $g':=g\circ\varphi$.
We have 
$V'=\varphi_1^{-1}V\otimes\nbigo_{Z'}(\ast g')$.
We have the canonical pre-$K$-Betti structures
$\nbigf_{V'}$ and $\nbigf_{V'!}$
of $V'$ and $V'_!$, respectively. 
According to Proposition \ref{prop;13.4.22.10},
the morphisms
$\varphi_{1\dagger}V'_!\lrarr 
 \varphi_{\dagger}V\lrarr
 \varphi_{1\dagger}V'$ are compatible
with pre-$K$-Betti structures.
Hence, $\nbigc'$ is a $K$-cell.
We obtain a monomorphism
$\phi^{(0)}_g(\varphi_{\dagger}V)
\lrarr
 \phi^{(0)}_{g}(\varphi'_{\dagger}V')$
of pre-$K$-holonomic $\nbigd_X$-modules.
By Lemma \ref{lem;13.4.23.20},
$\phi^{(0)}_{g}(\varphi'_{\dagger}V')$
is $K$-holonomic.
Then,
we obtain that $\phi^{(0)}_g(\varphi_{\dagger}V)$
is $K$-holonomic by Proposition \ref{prop;09.10.4.110}.
\hfill\qed

\begin{cor}
\label{cor;09.10.4.301}
Assume that $SI(<\!n)$ and $GOOD(<\!n)$.
Let $f$ be a cell function of 
$\nbigc=(Z,U,\varphi,V)$.
Then,
$\Xi^{(a)}_f(\varphi_{\dagger}V)$
with the canonical pre-$K$-Betti structures
are $K$-holonomic.
\end{cor}
\pf
Applying the previous results to
$\varphi_{\dagger}\bigl(
 \Pi^{a,b}_{f\star}V\bigr)$ 
$(\star=!,\ast)$,
we obtain that they are $K$-holonomic.
Then, we obtain the corollary.
\hfill\qed

\subsection{Gluing}

By Lemma \ref{lem;13.4.23.20}
and Corollary \ref{cor;09.10.4.301},
we have
a gluing construction of 
$K$-holonomic $\nbigd$-modules.
Let $X$ be a complex manifold,
$\nbigc=(Z,U,\varphi,V)$ be a $K$-cell
as in \S\ref{subsection;09.12.5.20}.
Let $f$ be a cell function for $\nbigc$ on $X$.
Let $\nbigq$ be a $K$-holonomic
$\nbigd$-module whose support is contained
in $f^{-1}(0)$.
Assume that we are given morphisms of
$K$-holonomic $\nbigd$-modules
\[
 \psi^{(1)}_f(\varphi_{\dagger}V)\lrarr  \nbigq
 \lrarr \psi^{(0)}_f(\varphi_{\dagger}V),
\]
such that the composite is equal to
the canonical map
$\psi^{(1)}_f(\varphi_{\dagger}V)\lrarr
 \psi^{(0)}_f(\varphi_{\dagger}V)$.
Then, we obtain a $K$-holonomic $\nbigd$-module
as the cohomology of the following complex:
\[
 \psi^{(1)}_f(\varphi_{\dagger}V)
\lrarr
 \Xi^{(0)}_f(\varphi_{\dagger}V)\oplus \nbigq
\lrarr
 \psi^{(0)}_f(\varphi_{\dagger}V)
\]

\subsection{Good holonomic $\nbigd$-module
with good $K$-structure}
\label{subsection;09.12.17.1}

Suppose $SI(<\!n)$ and $GOOD(<\!n)$.
Let $X$ be a complex manifold
with a simply normal crossing hypersurface $D$.
Let $\nbigm$ be a good pre-$K$-holonomic $\nbigd$-module
on $(X,D)$
such that $\dim\Supp\nbigm=n$.
Let us prove that $\nbigm$ is $K$-holonomic.
We may assume that
$X=\Delta^N$ and $D=\bigcup_{i=1}^{\ell}\{z_i=0\}$.
Let $\rho(\nbigm)\in\seisuu_{\geq\,0}\times
 \seisuu_{>0}$ denote the pair of
$\dim\Supp\nbigm$ and
the number of the irreducible components of
$\Supp\nbigm$ with the maximal dimension.
We use the lexicographic order on 
$\seisuu_{\geq \,0}\times\seisuu_{>0}$.
For a good holonomic $\nbigd$-module $\nbigm$
on $(X,D)$,
there exists $J\subset\ellsitabar$
with $n=N-|J|$
such that
$\nbigm\bigl(\ast g\bigr)\neq 0$
comes from a meromorphic flat bundle $V$ on $D_J$,
where $g:=\prod_{\substack{j\not\in J\\ j\leq \ell}}z_j$.
Let $\iota:D_J\lrarr X$ denote the inclusion.
We have a description of $\nbigm$
as the cohomology of the complex
of pre-$K$-holonomic $\nbigd$-modules
$\psi^{(1)}_g(\iota_{\dagger}V)
\lrarr
 \Xi^{(0)}_g(\iota_{\dagger}V)\oplus\phi^{(0)}_g(\nbigm)
\lrarr
 \psi^{(0)}_g(\iota_{\dagger}V)$.
They are good pre-$K$-holonomic $\nbigd$-modules.
By Lemma \ref{lem;13.4.23.20}
and Corollary \ref{cor;09.10.4.301},
$\psi^{(a)}_g(V)$ and $\Xi^{(a)}_g(V)$
are $K$-holonomic.
Because $\rho(\phi^{(0)}_g(\nbigm))<\rho(\nbigm)$,
we obtain that $\phi^{(0)}_g(\nbigm)$
is $K$-holonomic.
Hence, we obtain that
$\nbigm$ is also $K$-holonomic.

\subsection{Generalization}
\label{subsection;13.4.27.100}

We use the notation introduced in \S\ref{subsection;09.12.5.20}.
\begin{prop}
\label{prop;13.4.23.220}
Suppose that $SI(<n)$ and $GOOD(<n)$.
Then, the pre-$K$-holonomic $\nbigd_X$-module
$\varphi_{\dagger}\gbigp_{\gbigh}(V)$
is $K$-holonomic.
\end{prop}
\pf
It is enough to consider the issue
locally around any point $P\in X$.
We will shrink $X$ around $P$
without mention.
Let $\nbigc'=(Z',U',\varphi',V')$ 
be a dominant refinement 
of $\nbigc$
with a cell function $g$ for $\nbigc'$.
We set
$\gbigh':=(\varphi^{\prime})^{-1}(\gbigh)$.

\begin{lem}
\label{lem;13.4.27.110}
Under the assumptions
$SI(<n)$ and $GOOD(<n)$,
$\varphi'_{\dagger}\bigl(
 \gbigp_{\gbigh'}(V')
 \bigr)$
is $K$-holonomic.
\end{lem}
\pf
We have the expression of 
$\varphi'_{\dagger}\bigl(
 \gbigp_{\gbigh'}(V')
 \bigr)$
as the cohomology of the following
complex of pre-$K$-holonomic $\nbigd$-modules:
\[
 \psi^{(1)}_{g}\bigl(
 \varphi'_{\dagger}(V')
 \bigr)
\lrarr
 \phi^{(0)}_{g}
 \varphi'_{\dagger}\bigl(
 \gbigp_{\gbigh'}(V')
 \bigr)
 \oplus
 \Xi^{(0)}_g\varphi'_{\dagger}(V')
\lrarr
  \psi^{(0)}_{g}\bigl(
 \varphi'_{\dagger}(V')
 \bigr)
\]
By Lemma \ref{lem;13.4.23.20}
and Corollary \ref{cor;09.10.4.301},
we obtain that
$\psi^{(a)}_{g}\bigl(
 \varphi'_{\dagger}(V')
 \bigr)$
and 
$\Xi^{(a)}_g\varphi'_{\dagger}(V')$
are $K$-holonomic.
By Lemma \ref{lem;13.4.23.20},
$\phi^{(0)}_{g}
 \varphi'_{\dagger}\bigl(
 \gbigp_{\gbigh'}(V')
 \bigr)$ 
is $K$-holonomic.
Hence, we obtain that
$\varphi'_{\dagger}\bigl(
 \gbigp_{\gbigh'}(V')
 \bigr)$ 
is $K$-holonomic.
Thus, we obtain Lemma \ref{lem;13.4.27.110}.
\hfill\qed

\vspace{.1in}
We have a natural monomorphism
of pre-$K$-holonomic $\nbigd$-modules
$\varphi_{\dagger}\bigl(
 \gbigp_{\gbigh}(V)
\bigr)
\lrarr
 \varphi'_{\dagger}\bigl(
 \gbigp_{\gbigh'}(V')
 \bigr)$,
as remarked in Proposition \ref{prop;13.4.27.111}.
Then, by Proposition \ref{prop;09.10.4.110},
we obtain that 
$\varphi_{\dagger}\bigl(\gbigp_{\gbigh}(V)\bigr)$
is $K$-holonomic.
\hfill\qed

\subsection{$K(\ast D)$-cell}

We use the notation introduced
in \S\ref{subsection;09.12.5.20}.
Let $D$ be a hypersurface of $X$
such that
$D_{Z1}:=\varphi^{-1}(D)\subset D_Z$.
We have the pre-$K$-holonomic
$\nbigd_{Z}$-module $V(!D_{Z1})$.
We obtain the following proposition
as a special case of 
Proposition \ref{prop;13.4.23.220}.
\begin{prop}
\label{prop;13.4.23.30}
$\varphi_{\dagger}\bigl(V(!D_{Z1})\bigr)$ 
is $K$-holonomic.
\hfill\qed
\end{prop}

%% file: 8.3.tex
\subsection{Equivalence of
$K(\ast D)$-Betti structure 
and $K$-Betti structure}

Let $X$ be any complex manifold
with a hypersurface $D$.
Let $(\nbigm,\nbigf)$ be any pre-$K$-holonomic
$\nbigd_{X(\ast D)}$-module
with $\dim\Supp\nbigm\leq n$.

\begin{lem}
\mbox{{}}\label{lem;09.10.18.510}
\begin{itemize}
\item
Assume $SI(<\!n)$ and $GOOD(<\!n)$.
If $\nbigf$ is a $K(\ast D)$-Betti structure,
then it is a $K$-Betti structure.
\item
Assume $LOC(\leq n)$.
If $\nbigf$ is a $K$-Betti structure,
then it is a $K(\ast D)$-Betti structure.
\end{itemize}
\end{lem}
\pf
Let us prove the first claim.
We use an induction on
the dimension of the support.
Let $P$ be any point of $D\cap\Supp\nbigm$.
We take a bounding cell $\nbigc=(Z,U,\varphi,V)$
of $(\nbigm,\nbigf)$ at $P$,
and a cell function $g$ of $\nbigc$
as $\nbigd_{X(\ast D)}$-module.
We have a description of $\nbigm$ 
as the cohomology of the following complex
of $K(\ast D)$-holonomic 
$\nbigd_{X(\ast D)}$-modules:
\[
 \psi^{(1)}_g\bigl(
 \varphi_{\dagger}(V),\,\ast D
 \bigr)
\lrarr
 \Xi^{(0)}_g(\varphi_{\dagger}V,\ast D)
\oplus
 \phi^{(0)}_g(\nbigm,\ast D)
\lrarr
 \psi^{(0)}_g\bigl(\varphi_{\dagger}(V),\ast D\bigr)
\]
By the inductive assumption,
$\phi^{(0)}_g(\nbigm,\ast D)$ is $K$-holonomic.
By Proposition \ref{prop;13.4.23.30},
$\psi^{(a)}_g\bigl(
 \varphi_{\dagger}(V),\ast D
 \bigr)$
and $\Xi^{(a)}_g(\nbigm,\ast D)$
are $K$-holonomic.
Hence, we obtain that
$\nbigm$ is also $K$-holonomic.

Let us prove the second claim.
By the assumption $LOC(\leq n)$,
we obtain a $K(\ast D)$-holonomic
$\nbigd_{X(\ast D)}$-module
$\bigl(\nbigm(\ast D),\nbigf(\ast D)\bigr)$
with a morphism
\[
 (\nbigm,\nbigf)\lrarr
 \bigl(\nbigm(\ast D),\nbigf(\ast D)\bigr)
\]
of pre-$K$-holonomic $\nbigd$-modules.
Because $\nbigm=\nbigm(\ast D)$,
we obtain $\nbigf=\nbigf(\ast D)$,
and hence $\nbigf$ is a $K(\ast D)$-Betti structure.
\hfill\qed

\vspace{.1in}
We reformulate the uniqueness 
(Proposition \ref{prop;09.10.18.250})
as follows.

\begin{cor}
\label{cor;13.4.27.300}
Let $\star$ be $\ast$ or $!$.
Assume $SI(<n)$, $GOOD(<n)$ and $LOC(\leq n)$.
Let $\nbigm$ be a holonomic $\nbigd$-module on $X$
such that $\nbigm(\star D)=\nbigm$.
Let $\nbigf_i$ $(i=1,2)$ be $K$-Betti structures
on $\nbigm$.
If $\nbigf_{1|X-D}=\nbigf_{2|X-D}$,
then $\nbigf_1=\nbigf_2$.
\end{cor}
\pf
The claim for $\star=\ast$ follows from
Lemma \ref{lem;09.10.18.510} and
Proposition \ref{prop;09.10.18.250}.
We obtain the claim for $\star=!$ 
by using the dual with Proposition \ref{prop;09.10.19.1}.
\hfill\qed

\begin{cor}
Suppose $SI(<n)$ and $GOOD(<n)$ and $LOC(\leq n)$.
Let $\nbigm$ be a holonomic $\nbigd_X$-module.
Assume that one of the following holds;
(i) $\nbigm(!D)\lrarr \nbigm$ is surjective,
(ii) $\nbigm\lrarr \nbigm(\ast D)$ is injective.
Let $\nbigf_i$ $(i=1,2)$ be $K$-Betti structures
on $\nbigm$.
If $\nbigf_{1|X-D}=\nbigf_{2|X-D}$,
then $\nbigf_1=\nbigf_2$.
\hfill\qed
\end{cor}

We reformulate the independence
from a compactification (Proposition \ref{prop;09.10.19.2}).
Let $F:X'\lrarr X$ be a projective birational morphism
of complex manifolds.
Let $D$ be a hypersurface,
and we put $D':=F^{-1}(D)$.
Assume $X'-D'\simeq X-D$.

\begin{prop}
\label{prop;09.10.5.50}
\mbox{{}}
Assume $SI(<n)$, $GOOD(<n)$ and $LOC(\leq n)$.
Let $\nbigm'$ be a holonomic 
$\nbigd_{X'(\ast D')}$-module.
We set $\nbigm:=F_{\dagger}\nbigm'$.
\begin{itemize}
\item
If $\nbigf'$ is a $K$-Betti structure of
$\nbigm'$,
then $F_{\ast}\nbigf'$ is a $K$-Betti structure of
$\nbigm$.
\item
If $\nbigf$ is a $K$-Betti structure of $\nbigm$,
then $\nbigm'$ is equipped with a $K$-Betti structure
$\nbigf'$ such that
$\nbigf'_{|X'-D'}=\nbigf_{|X-D}$
under the isomorphism
$\nbigm'_{|X'-D'}\simeq \nbigm_{|X-D}$.
It is functorial.
\hfill\qed
\end{itemize}
\end{prop}

\subsection{$K(\ast D)$-cell}

Let $\varphi:Z\lrarr X$ be a projective morphism
of complex manifolds
such that $\dim Z=n$.
Let $D_Z$ be a normal crossing hypersurface of $Z$.
Assume that $\varphi_{|Z-D_Z}$ is an immersion.
We suppose $D_1:=\varphi^{-1}(D)\subset D_Z$.
Let $V$ be a meromorphic
flat connection on $(Z,D_Z)$
with a good $K$-structure.
We obtain the pre-$K$-holonomic
$\nbigd_Z$-modules
$V$ and $V_!(\ast D_1)$.

\begin{prop}
\label{prop;09.10.18.500}
Assume that $SI(<\!\!n)$, $GOOD(<\!\!n)$
and $LOC(<\!\!n)$ hold.
Then, 
$\varphi_{\dagger}
 V_!(\ast D_1)$
and 
$\varphi_{\dagger}V$
are $K(\ast D)$-holonomic.
\end{prop}
\pf
Let us prove that
$\nbigc_0=(Z,U,\varphi,V)$ is 
a bounding $n$-cell at any $P\in D\cap\varphi(Z)$.
It is enough to consider the issue locally.
We shall shrink $X$ without mention.
Let $\nbigc'=(Z',U',\varphi',V')$
be a dominant refinement at $P$
with a cell function $g$
as $\nbigd_X(\ast D)$-modules.
We have a factorization
$\varphi'=\varphi\circ\varphi_1$,
where $\varphi_1:Z'\lrarr Z$.
We put $g':=g\circ\varphi'$
and $D_1':=(\varphi')^{-1}D$.
We have 
$V'=\varphi_1^{-1}V\otimes\nbigo_{Z'}(\ast g')$.
According to Proposition \ref{prop;13.4.22.10},
the morphisms
$\varphi'_{\dagger}(V'_!)(\ast D)
 \lrarr
 \varphi_{\dagger}(V_!)(\ast D)
 \lrarr 
 \varphi_{\dagger}V\lrarr
 \varphi'_{\dagger}V'$
are compatible
with the canonical pre-$K$-Betti structures.
We obtain the induced pre-$K$-Betti
structures of
$\phi^{(a)}_g\bigl(
 \varphi_{\dagger}(V),\ast D
 \bigr)$
and 
$\phi^{(a)}_g\bigl(
 \varphi_{\dagger}(V_!),\ast D
 \bigr)$.

We obtain pre-$K$-holonomic $\nbigd$-modules
$\phi^{(a)}_{g'}\bigl(
 V'_!,\ast D_1'
 \bigr)$
and $\phi^{(a)}_{g'}\bigl(V',\ast D_1'\bigr)$ 
on $Z'$.
They are $K$-holonomic,
which can be proved by the argument
in the proof of Lemma \ref{lem;13.4.23.20}.
We obtain that 
$\phi^{(a)}_{g}\bigl(
 \varphi'_{\dagger}V',\ast D
 \bigr)$
and
$\phi^{(a)}_{g}\bigl(
 \varphi'_{\dagger}\bigl(V'_!),\ast D
 \bigr)$ 
are $K$-holonomic
by the assumption $SI(<\!n)$.
By Lemma \ref{lem;09.10.18.510}
and the assumption $LOC(<n)$,
$\phi^{(a)}_g(\varphi'_{\dagger}V',\ast D)$ and
$\phi^{(a)}_{g}(\varphi'_{\dagger}V'_!,\ast D)$ 
are $K(\ast D)$-holonomic.
Because
$\phi^{(a)}_g(\varphi_{\dagger}V,\ast D)
\subset 
 \phi^{(a)}_{g}\bigl(\varphi'_{\dagger}V',\ast D\bigr)$ 
is compatible with the pre-$K$-Betti structures,
$\phi^{(a)}_g(\varphi_{\dagger}V,\ast D)$
is also a $K(\ast D)$-holonomic
by Lemma \ref{lem;09.10.18.511}.
Since the surjection 
$\phi^{(a)}_{g}\bigl(
\varphi'_{\dagger}V'_!,\ast D
 \bigr)
\lrarr
 \phi^{(a)}_g\bigl(
 \varphi_{\dagger}V_!,\ast D
 \bigr)$ 
is compatible with the pre-$K$-Betti structures,
$\phi^{(a)}_g\bigl(
 \varphi_{\dagger}V_!,\ast D
 \bigr)$ 
is also $K(\ast D)$-holonomic
by Lemma \ref{lem;09.10.18.511}.
\hfill\qed

\begin{cor}
\label{cor;09.10.18.501}
Assume that
$SI(<\!n)$, $GOOD(<\!n)$ and $LOC(<\!n)$ hold.
Let $f$ be a cell function of 
an $n$-dimensional cell
$\nbigc=(Z,U,\varphi,V)$ as 
$\nbigd_{X(\ast D)}$-module.
Then,
$\psi^{(a)}_f(\varphi_{\dagger}V,\ast D)$ and 
$\Xi^{(a)}_f(\varphi_{\dagger}V,\ast D)$
with the canonical pre-$K$-Betti structures
are $K(\ast D)$-holonomic.
\end{cor}
\pf
Applying the previous results to
$\Pi^{a,b}_{f\star}\bigl(
 \varphi_{\dagger}V,\ast D
 \bigr)$ for $\star=\ast,!$,
we obtain that they are $K(\ast D)$-holonomic.
Then, we obtain the corollary.
\hfill\qed

\subsection{Localization}
\label{subsection;09.12.5.10}

Let us prove $LOC(\leq n)$ by assuming
$SI(<n)$, $GOOD(<n)$ and $LOC(<n)$.
By Proposition \ref{prop;09.10.18.250},
the problem is local.
Let $\nbigm$ be a $K$-holonomic
$\nbigd_X$-module
with $\dim \Supp\nbigm\leq n$.

Let $P$ be any point of $D$.
Let $(Z,U,\varphi,V)$ be a bounding cell of $\nbigm$
at $P$
with a cell function $g$ as $K$-holonomic $\nbigd$-modules.
By taking a refinement,
we may assume $U\cap D=\emptyset$.
We put $g_1:=\varphi^{-1}(g)$
and $D_1:=\varphi^{-1}(D)$.
We have the expression of $\nbigm$
as the cohomology of
the following complex of 
the $K$-holonomic $\nbigd$-modules:
\begin{equation}
 \label{eq;09.10.5.6}
\psi^{(1)}_g\varphi_{\dagger}(V_!)
\lrarr
 \Xi^{(0)}_g\varphi_{\dagger}(V)
\oplus
 \phi^{(0)}_g(\nbigm)
\lrarr
 \psi^{(0)}_g\varphi_{\dagger}(V)
\end{equation}
By the assumption of the induction,
$\psi^{(a)}_g\bigl(
 \varphi_{\dagger}V_!,\ast D\bigr)$
and $\phi^{(a)}_g(\nbigm,\ast D)$
are equipped with
the induced $K(\ast D)$-Betti structures.
We also have the following 
commutative diagram of
pre-$K$-holonomic $\nbigd$-modules:
\[
\begin{CD}
 \psi^{(1)}_g(\varphi_{\dagger}V)
 @>>>
 \phi^{(0)}_g(\nbigm) 
 @>>>
 \psi^{(0)}_g(\varphi_{\dagger}V)\\
 @VVV @VVV @VVV \\
 \psi^{(1)}_g\bigl(
 \varphi_{\dagger}V_!,\ast D\bigr)
 @>>>
\phi^{(0)}_g(\nbigm,\ast D)
 @>>>
 \psi^{(0)}_g\bigl(\varphi_{\dagger}V_!,
 \ast D\bigr)
\end{CD}
\]

We have the canonical 
pre-$K$-Betti structures of
$\psi^{(a)}_{g_1}\bigl(
 V,\ast D_1
 \bigr)$
and $\Xi^{(a)}_{g_1}\bigl(V,\ast D_1\bigr)$.
According to Corollary \ref{cor;09.10.18.501},
their push-forward 
$\varphi_{\dagger}\psi^{(a)}_{g_1}\bigl(
 V,\ast D_1\bigr)$
and $\varphi_{\dagger}\Xi^{(a)}_{g_1}\bigl(
 V,\ast D_1\bigr)$ 
are $K(\ast D)$-holonomic.
We also have the following commutative diagram
of pre-$K$-holonomic $\nbigd$-modules:
\[
 \begin{CD}
 \varphi_{\dagger}\psi^{(1)}_{g_1}(V)
@>>>
 \varphi_{\dagger}\Xi^{(0)}_{g_1}(V)
@>>> 
 \varphi_{\dagger}\psi^{(0)}_{g_1}(V)
 \\ @VVV @VVV @VVV \\
 \varphi_{\dagger}\psi^{(1)}_{g_1}\bigl(
 V,\ast D_1 \bigr)
@>>>
 \varphi_{\dagger}\Xi^{(0)}_{g_1}\bigl(
 V,\ast D_1 \bigr)
@>>> 
 \varphi_{\dagger}\psi^{(0)}_{g_1}\bigl(
 V,\ast D_1
 \bigr)
 \end{CD}
\]
By Proposition \ref{prop;09.10.18.250},
the identification
$\varphi_{\dagger}\psi^{(a)}_{g_1}\bigl(
 V,\ast D_1\bigr)
 \simeq
 \psi^{(a)}_g\bigl(\varphi_{\dagger}V,\ast D\bigr)$ 
is compatible with the pre-$K$-Betti structures.
Hence, we obtain a $K(\ast D)$-Betti structure of
$\nbigm(\ast D)$
with a morphism
of pre-$K$-holonomic $\nbigd$-modules
$\nbigm\lrarr\nbigm(\ast D)$
whose restriction to $X-D$ is an isomorphism.
The functoriality is clear from the above construction.
\hfill\qed

\subsection{Twist}

Let $(\nbigm,\nbigf)$ be any $K(\ast D)$-holonomic
$\nbigd(\ast D)$-module
such that $\dim\Supp\nbigm\leq n$.
Let $\nbigv$ be a meromorphic flat connection on 
$(X,D)$ with a good $K$-structure
$\nbigf_{\nbigv}$.
According to Lemma \ref{lem;09.11.2.1},
$\nbigf_{\nbigm|X-D}\otimes
 \nbigf_{\nbigv|X-D}$
is a $K$-Betti structure of
$(\nbigm\otimes\nbigv)_{|X-D}$.

\begin{prop}
\label{prop;09.12.6.1}
Assume $SI(<n)$, $GOOD(<n)$
and $LOC(<n)$.
There exists a $K(\ast D)$-Betti structure 
$\nbigf_{\nbigm\otimes\nbigv}$ of
$\nbigm\otimes\nbigv$
such that
$\nbigf_{\nbigm\otimes\nbigv|X-D}
\simeq
 \nbigf_{\nbigm|X-D}
\otimes
 \nbigf_{\nbigv|X-D}$.
It is functorial with respect to
$\nbigm$ and $\nbigv$.
\end{prop}
\pf
Let $P\in D$.
It is enough to consider the issue locally 
around $P$.
We use an induction on 
$\dim_P\Supp\nbigm$.
Let $\nbigc=(Z,U,\varphi,V)$ be 
a dominating cell of $\nbigm$ at $P$.
Let $g$ be a cell function for $\nbigc$
as $\nbigd_{X(\ast D)}$-module.
By the inductive assumption,
we have the $K(\ast D)$-Betti structures of
$\psi^{(a)}_g(\varphi_{\dagger}V,\ast D)\otimes\nbigv$
and $\phi^{(a)}_g(\varphi_{\dagger}V,\ast D)\otimes\nbigv$.
According to Corollary \ref{cor;09.10.18.501},
we have the $K(\ast D)$-Betti structures
of $\psi^{(a)}_g(\varphi_{\dagger}V,\ast D)\otimes\nbigv$
and $\Xi^{(a)}_g(\varphi_{\dagger}V,\ast D)\otimes\nbigv$
induced by the isomorphisms
$\psi^{(a)}_g(\nbigm,\ast D)\otimes\nbigv
\simeq
\psi^{(a)}_{g}(\nbigm\otimes\nbigv,\ast D)$
and 
$\Xi^{(a)}_g(\nbigm,\ast D)\otimes\nbigv\simeq
\Xi^{(a)}_{g}(\nbigm\otimes\nbigv,\ast D)$.
By the uniqueness,
the induced $K(\ast D)$-Betti structures
on $\psi^{(a)}_g(\nbigm,\ast D)\otimes\nbigv$
are equal.
Because $\nbigm\otimes\nbigv$
is expressed as the cohomology of
the complex
\[
 \psi^{(1)}_g(\nbigm,\ast D)\otimes\nbigv
\lrarr
 \Xi^{(0)}_g(\nbigm,\ast D)\otimes\nbigv
\oplus
 \phi^{(0)}_g(\nbigm,\ast D)\otimes\nbigv
\lrarr
 \psi^{(0)}_g(\nbigm,\ast D)\otimes\nbigv,
\]
we obtain a $K(\ast D)$-Betti structure
on $\nbigm\otimes\nbigv$
with the desired property.
\hfill\qed

\subsection{Nearby, vanishing and maximal functors}
\label{subsection;13.4.23.31}

Suppose that
$SI(<n)$, $GOOD(\leq n)$ and $LOC(< n)$ hold.
Let $(\nbigm,\nbigf)$ be a $K$-holonomic 
$\nbigd_X$-module
with $\dim\Supp\nbigm\leq n$.
Let $f$ be any holomorphic function on $X$.
As proved in \S\ref{subsection;09.12.5.10},
we obtain a morphism of 
$K$-holonomic $\nbigd_X$-modules
$\nbigm\lrarr \nbigm(\ast f)$.
By considering the dual,
we also obtain 
a morphism of $K$-holonomic
$\nbigd_X$-modules 
$\nbigm(!f)\lrarr\nbigm$.
By Proposition \ref{prop;09.12.6.1},
for any $a\leq b$,
we have $K$-holonomic $\nbigd_X$-modules
$\Pi^{a,b}_{f\star}(\nbigm)$ $(\star=\ast,!)$.
Hence, we obtain $K$-holonomic $\nbigd_X$-modules
$\Pi^{a,b}_{f\star !}(\nbigm)$.
In particular,
we obtain $K$-holonomic $\nbigd_X$-modules
$\Xi^{(a)}_f(\nbigm)$
and 
$\psi^{(a)}_f(\nbigm)$
with morphisms
$\nbigm(!f)\lrarr
 \Xi^{(0)}_f(\nbigm)\lrarr \nbigm(\ast f)$
and
$\psi^{(1)}_f(\nbigm)\lrarr
 \Xi^{(0)}_f(\nbigm)
\lrarr \psi^{(0)}_f(\nbigm)$
in $\Hol(X,K)$.
We obtain a $K$-holonomic
$\nbigd_X$-module
$\phi^{(0)}_f(\nbigm)$
as the cohomology of
the complex 
$\nbigm(!f)
\lrarr
 \Xi^{(0)}_f(\nbigm)
\oplus\nbigm
\lrarr \nbigm(\ast f)$
in $\Hol(X,K)$.
We can recover $\nbigm$
as the cohomology of the complex
$\psi^{(1)}_f(\nbigm)
\lrarr
 \Xi^{(0)}_f(\nbigm)\oplus\phi^{(0)}_f(\nbigm)
\lrarr
 \psi^{(0)}_f(\nbigm)$
in $\Hol(X,K)$.

%% file: 8.4.tex
This subsection is a preliminary for the proof of
Theorem \ref{thm;09.10.16.5}.

\subsection{Non-characteristic condition}

Let $\nbigm$ be a holonomic $\nbigd$-module on 
a complex  manifold $X$.
There exists a stratification 
$\Supp(\nbigm)=\coprod_{i\in\Lambda} Z_i$
such that
(i) each $Z_i$ is a smooth locally closed analytic subset
of $X$,
(ii) $\Ch(\nbigm)=\coprod_{i\in\Lambda}
 T_{Z_i}^{\ast}X$.

\begin{lem}
A complex submanifold $W\subset X$ 
is non-characteristic with respect to $\nbigm$
if and only if $W$ and $Z_i$ are transversal
for any $i\in\Lambda$.
In that case,
for the inclusion $\iota:W\lrarr X$,
we have
$\Ch(\iota_{\ast}\iota^{\ast}\nbigm)
=\coprod_{i\in\Lambda}
 T_{Z_i\cap W}^{\ast}X$.
\end{lem}
\pf
For $P\in W\cap Z_i$,
we have subspaces
$(T_{Z_i}^{\ast}X)_{|P}$
and $(T_{W}^{\ast}X)_{|P}$
of $(T^{\ast}X)_{|P}$.
Then,
$W$ and $Z_i$ are transversal at $P$
if and only if
$(T_{W}^{\ast}X)_{|P}\cap
(T_{Z_i}^{\ast}X)_{|P}=\{0\}$.
Then, the first claim of the lemma is clear.
The second claim follows from
general formulas of the characteristic varieties
for the pull back by a non-characteristic closed immersion
and the push-forward by a closed immersion.
\hfill\qed

\begin{lem}
\label{lem;13.4.15.21}
Let $D$ be a smooth hypersurface of $X$.
If $D$ is non-characteristic with respect to $\nbigm$,
the natural morphism
$\nbigm(!D)\lrarr
 \nbigm\otimes\nbigo(!D)$
is an isomorphism
\end{lem}
\pf
Let $i:D\lrarr X$ be the closed immersion.
Because $D$ is non-characteristic with respect to $\nbigm$,
we have
the exact sequence
$0\lrarr 
i_{\ast}i^{\ast}\nbigm
\lrarr\nbigm(!D)
\lrarr\nbigm\lrarr 0$.
We have
$0\lrarr i_{\ast}i^{\ast}\nbigo_X
\lrarr \nbigo_X(!D)
\lrarr \nbigo_X\lrarr 0$.
By the non-characteristic condition
and the projection formula,
we obtain
$0\lrarr i_{\ast}i^{\ast}\nbigm
\lrarr \nbigm\otimes\nbigo_X(!D)
\lrarr \nbigm\lrarr 0$.
Then, we obtain the claim of the lemma.
\hfill\qed

\begin{lem}
\label{lem;13.4.15.22}
Let $D_i$ $(i=1,2)$ be smooth hypersurfaces of $X$
such that 
(i) $D_1$ and $D_2$ are transversal,
(ii) $D_i$ $(i=1,2)$ and $D_1\cap D_2$
are non-characteristic with respect to $\nbigm$.
Then, 
$D_2$ is non-characteristic with respect to
$\nbigm(\ast D_1)$,
and we have
\begin{equation}
 \label{eq;13.4.15.20}
 \bigl(
 \nbigm(\ast D_1)
 \bigr)(!D_2)
\simeq
  \bigl(
 \nbigm(! D_2)
 \bigr)(\ast D_1)
\simeq
  \nbigm\otimes\nbigo(!D_2)\otimes\nbigo(\ast D_2).
\end{equation}
\end{lem}
\pf
By the assumption,
$D_i$ $(i=1,2)$
and $D_1\cap D_2$ are transversal
to $Z_j$ for $j\in\Lambda$.
It is elementary to check that
$D_2$ is transversal to $D_1\cap Z_j$
$(j\in\Lambda)$.
We obtain that $D_2$ is non-characteristic with respect to
$\nbigm(\ast D_1)$.
We obtain the isomorphisms
(\ref{eq;13.4.15.20})
from Lemma \ref{lem;13.4.15.21}.
\hfill\qed

\subsection{Non-characteristic tuple of
hyperplane subbundles}
\label{subsection;09.10.16.1}

Let $\nbige$ be a locally free sheaf on 
any complex manifold $Y$.
Let $X$ be its projectivization 
with the projection $G:X\lrarr Y$.
If a section $s$ of $\nbigo_{\proj(\nbige)/Y}(1)$
gives a nowhere vanishing section
of $G_{\ast}\bigl(
  \nbigo_{\proj(\nbige)/Y}(1)\bigr)$,
the zero set of $s$ is called a hyperplane subbundle of $X$.
For any hyperplane subbundle $H$ of $X$
and $P\in Y$,
let $H_{|P}$ denote the fiber over $P$.

Let $\nbigm$ be any holonomic $\nbigd_X$-module.
Let $\vecH:=(H_1,\ldots,H_N)$ be a tuple of
hyperplane subbundles of $X$
such that,
for each $P\in Y$,
the tuple of hyperplanes
$(H_{1|P},H_{2|P},\ldots,H_{N|P})$
is of general position.
We say that $\vecH$ is non-characteristic 
with respect to $\nbigm$
if $H_I:=\bigcap_{i\in I}H_i$ are
non-characteristic with respect to $\nbigm$
for any $I\subset \{1,\ldots,N\}$.
We can prove the following lemma
by a standard argument of genericity.

\begin{lem}
Suppose that $(H_1,\ldots,H_N)$ is
non-characteristic with respect to $\nbigm$.
Let $P$ be any point of $Y$.
Then, if we shrink $Y$ around $P$,
we can take a hyperplane subbundle $H_{N+1}$
such that 
$(H_1,\ldots,H_{N},H_{N+1})$ is also
non-characteristic with respect to $\nbigm$.
\hfill\qed
\end{lem}

Recall the following general lemma.
\begin{lem}
\label{lem;14.1.15.1}
Let $(H_1,H_2)$ be a tuple of hyperplane bundles
of $X$, which is non-characteristic with respect to $\nbigm$.
Then, 
$G_{\dagger}^i\bigl(
 \nbigm(\ast H_1!H_2)
 \bigr)=0$
for any $i\neq 0$.
\end{lem}
\pf
Let $\nbigm_i$ $(i=1,2)$
be holonomic $\nbigd_X$-modules,
and let $H_i$ be hypersurfaces 
which is non-characteristic with respect to $\nbigm_i$.
Because $\nbigm_i$ has a global good filtration
according to \cite{Malgrange-filtration},
we have an exhaustive filtration 
$\nbigg_a$ $(a=1,2,\ldots)$
by coherent $\nbigo_X$-submodules 
of $\nbigm_1$.
We have
$R^bG_{\ast}(\nbigg_a(\ast H_1)\otimes\Omega^j_{X/Y})=0$
for any $b>0$.
Hence, we have 
$R^bG_{\ast}\nbigm_1(\ast H_1)\otimes\Omega^j_{X/Y}=0$.
Then, we obtain
$G_{\dagger}^i\nbigm_1(\ast H_1)=0$
for any $i>0$.
By using the duality,
we obtain that
$G_{\dagger}^i\bigl(
 \nbigm_2(! H_2)
 \bigr)=0$ for any $i<0$.
Then, the claim follows from
Lemma \ref{lem;13.4.15.22}.
\hfill\qed

\subsection{Resolutions}
\label{subsection;09.10.16.2}

Let $X$, $Y$ and $\nbigm$ be as 
in \S\ref{subsection;09.10.16.1}.
Let $\vecH=(H_1,\ldots,H_N)$ be a tuple of 
hyperplane subbundles of $X$,
non-characteristic with respect to $\nbigm$.
Let $\ibar:=\{1,\ldots,i\}$,
and let  $\iota_{H_{\ibar}}$ denote
the inclusion $H_{\ibar}\subset X$.
We put 
$\nbign_0:=\nbigm(\ast H_1)$.
We also put $\nbigc_i:=
 \iota_{H_{\ibar}\dagger}
 \iota_{H_{\ibar}}^{\ast}\nbigm$,
and 
$\nbign_i:=\nbigc_i(\ast H_{i+1})$.
We have the natural exact sequences:
\begin{equation}
 0\lrarr \nbigm\lrarr \nbign_0\lrarr \nbigc_{1}
 \lrarr 0,
\quad\quad
0\lrarr \nbigc_i\lrarr \nbign_i
\lrarr \nbigc_{i+1}\lrarr 0
\end{equation}
Hence, we obtain the following 
exact sequence:
\begin{equation}
 \label{eq;13.4.27.200}
0\lrarr
 \nbigm\lrarr \nbign_{0}\lrarr
 \nbign_{1}\lrarr\cdots\lrarr \nbign_n\lrarr
\cdots
\end{equation}
Let $\vecH'=(H_j'|\,j=1,\ldots,N')$ be a tuple of
hyperplane subbundles of $X$
such that $\vecH\sqcup\vecH'$ is non-characteristic with respect to
$\nbigm$.
We set 
$\nbigq_{i,0}:=\nbign_i(!H_1')$.
We also put
$\nbigk_{i,-j}:=\iota_{H'_{\jbar}\dagger}
 \iota_{H'_{\jbar}}^{\ast}\nbign_i$
and 
$\nbigq_{i,-j}:=\nbigk_{i,-j}(!H_{j+1})$.
We have the natural exact sequences:
\[
 0\lrarr \nbigk_{i,-1}\lrarr
 \nbigq_{i,0}\lrarr \nbign_i\lrarr 0,
\quad\quad
 0\lrarr \nbigk_{i,-j-1}\lrarr 
 \nbigq_{i,-j}\lrarr\nbigk_{i,-j}\lrarr 0
\]
Hence, we obtain the following exact sequences:
\[
 0\llarr \nbign_{i}\llarr
 \nbigq_{i,0}\llarr \nbigq_{i,-1}  \llarr
 \nbigq_{i,-2}\llarr \cdots
\]
By construction, we have the naturally defined
morphisms $\nbigq_{i,-j}\lrarr \nbigq_{i+1,-j}$
and the commutative diagrams:
\[
 \begin{CD}
 \nbigq_{i,-j} @>>> \nbigq_{i+1,-j}\\
 @VVV @VVV \\
 \nbigq_{i,-j+1} @>>> \nbigq_{i+1,-j+1}
 \end{CD}
\]
Let $\Tot\bigl(\nbigq_{\bullet,\bullet}\bigr)$
denote the total complex of the double complex
$\nbigq_{\bullet,\bullet}$.
We have natural quasi-isomorphisms
$\Tot\bigl(\nbigq_{\bullet,\bullet}\bigr)
 \stackrel{\simeq}{\lrarr}
 \nbign_{\bullet}
 \stackrel{\simeq}{\llarr}
 \nbigm$.

By the construction,
for each $\nbigq_{i,-j}$,
there exists a holonomic $D$-module
$\nbigp_{i,-j}$ such that
(i) $(H_{i+1},H'_{j+1})$ is non-characteristic with respect to
$\nbigp_{i,-j}$,
(ii) $\nbigq_{i,-j}=\nbigp_{i,-j}(\ast H_{i+1}!H'_{j+1})$.

%% file: 8.5.tex
Let us prove that
$SI(<\!n)$, $GOOD(\leq \!n)$
and $LOC(\leq\!n)$ imply
$SI(\leq \!n)$.
The following argument is 
inspired by \cite{beilinson1}.

\subsection{Special case I}
\label{subsection;09.12.5.120}

Let $G:X\lrarr Y$ be any projective morphism
of complex manifolds
with $\dim X\leq n$.
Let $D$ be a hypersurface of $X$.
Let $V$ be a meromorphic flat connection on $(X,D)$
with a good $K$-structure.
Suppose that we are given a sequence of
hypersurface pairs $\gbigh$
contained in $D$.
We obtain a $K$-holonomic $\nbigd_X$-module
$\nbigm:=\gbigp_{\gbigh}(V)$
with the canonical $K$-Betti structure $\nbigf$.

\begin{prop}
\label{prop;09.10.16.21}
If $G_{\dagger}^i\nbigm=0$ for $i\neq 0$,
then
$RG_{\ast}\nbigf$ is a $K$-Betti structure of
$G_{\dagger}^0\nbigm$.
\end{prop}
\pf
It is enough to argue the issue locally around
any points of $Y$.
Let us consider the case
 $\Supp G_{\dagger}^{0}\nbigm\subsetneq G(X)$.
We take a holomorphic function $f$ such that
$\Supp G_{\dagger}^0\nbigm\subset f^{-1}(0)$
and $G(X)\not\subset f^{-1}(0)$.
We set $f_X:=f\circ G$.
As remarked in \S\ref{subsection;13.4.23.31},
we have a description of
the $K$-holonomic $\nbigd$-module
$\phi^{(0)}_{f_X}\nbigm$ as the cohomology of the following:
\[
 \nbigm(!f_X)
 \lrarr
 \Xi^{(0)}_{f_X}\nbigm(\ast f_X)\oplus 
 \nbigm
\lrarr \nbigm(\ast f_X)
\]
By the assumption,
$G_{\dagger}\nbigm(!f_X)
=G_{\dagger}\nbigm(\ast f_X)=
 G_{\dagger}\Xi^{(0)}_{f_X}\nbigm(\ast f_X)=0$.
Hence, we obtain
$G_{\dagger}(\nbigm,\nbigf)
\simeq
 G_{\dagger}\phi^{(0)}_{f_X}(\nbigm,\nbigf)$
as pre-$K$-holonomic $\nbigd$-modules.
By the assumption $SI(<n)$,
we obtain that $RG_{\ast}\nbigf$
is a $K$-Betti structure of $G_{\dagger}^0\nbigm$.

\vspace{.1in}

Let us consider the case 
$G(X)=\Supp G_{\dagger}^0\nbigm$.
Let $P\in\Supp G_{\dagger}^0\nbigm$.
Let $\nbigc=(Z,U,\varphi,E)$ be a cell
of $G_{\dagger}^0\nbigm$ at $P$
with a cell function $g$.
We set $g_Z:=\varphi^{-1}g$
and $g_X:=G^{-1}g$.
We have the $K$-Betti structures $\nbigf(\ast g_X)$
of $\nbigm(\ast g_X)$
by $LOC(\leq n)$.
By considering the dual,
we obtain the $K$-Betti structure 
$\nbigf(!g_X)$ of $\nbigm(!g_X)$.

\begin{lem}
\label{lem;09.10.16.13}
The $K$-structure of $E$ is good,
and the natural isomorphisms
\[
 \varphi_{\dagger}E(\star g_Z)\simeq 
 G_{\dagger}(\nbigm)(\star g)
\]
are compatible with the pre-$K$-Betti structures
for $\star=\ast,!$.
\end{lem}
\pf
We argue the case $\star=!$.
The case $\star=\ast$
can be argued similarly.
We take a projective birational morphism
$\kappa:X'\lrarr X$ such that
(i) $X'$ is smooth,
(ii) $X'-(g_X\circ\kappa)^{-1}(0)\simeq X-g_X^{-1}(0)$,
(iii) the induced morphism $X'\lrarr Y$ factors into
 $X'\stackrel{G_Z}{\lrarr}Z
 \stackrel{\varphi}{\lrarr} Y$.

We set $g_{X'}:=g_X\circ\kappa$
and
$\gbigh':=\varphi^{-1}(\gbigh)$.
We set $V':=\kappa^{\ast}V\otimes\nbigo(\ast g_{X'})$
and $\nbigm':=\gbigp_{\gbigh'}(V')(!g_{X'})$.
Note that
$\kappa_{\dagger}\nbigm'
  \simeq \nbigm(! g_X)$
and $G_{Z\dagger}\nbigm'=E(!g_Z)$.

We have 
the canonical pre-$K$-Betti structure
$\nbigf'$ of $\nbigm'$.
We have
$R\kappa_{\ast}\nbigf'=\nbigf(! g_X)$.
By Theorem \ref{thm;13.4.20.410},
we obtain that 
the $K$-structure of $E$ is compatible with
the Stokes structures,
and that 
$RG_{Z\ast}\nbigf'$ is the canonical
$K$-Betti structure of 
$G_{Z\dagger}\nbigm'$.
Hence, we obtain that
$RG_{\ast}\nbigf(! g_X)$
is the canonical $K$-Betti structure
of $G_{\dagger}(\nbigm)(! g)=
 \varphi_{\dagger}E(!g_Z)$.
Thus, we obtain Lemma \ref{lem;09.10.16.13}.
\hfill\qed

\begin{lem}
\label{lem;09.10.16.20}
The natural isomorphisms
$G_{\dagger}\Xi^{(a)}_{g_X}\bigl(\nbigm(\ast g_X)\bigr)
\simeq
 \Xi^{(a)}_g\bigl(\varphi_{\dagger}E\bigr)$
and 
$G_{\dagger}\psi^{(a)}_{g_X}\bigl(\nbigm(\ast g_X)\bigr)
\simeq
 \psi^{(a)}_g\bigl(\varphi_{\dagger}E\bigr)$
are compatible with the induced pre-$K$-Betti structures.
\end{lem}
\pf
By Lemma \ref{lem;09.10.16.13},
we obtain that the natural isomorphisms
$G_{\dagger}\bigl(
 \nbigm(\ast g_X)\otimes\gbigi_{g_X}^{a,b}
 \bigr)(\star g_X)
\simeq
 \varphi_{\dagger}E\otimes
 \gbigi_{g_Z}^{a,b}(\star g_Z)$
are compatible with the induced pre-$K$-Betti structures.
Hence, we obtain Lemma \ref{lem;09.10.16.20}.
\hfill\qed

\vspace{.1in}

By Lemma \ref{lem;09.10.16.13},
the morphisms
$\varphi_{\dagger}E_!\lrarr
 G_{\dagger}\nbigm\lrarr
 \varphi_{\dagger}E$
are compatible with the induced pre-$K$-Betti
structures,
i.e.,
$\nbigc$ is a $K$-cell.
Hence, we have an induced pre-$K$-Betti structure
$\lefttop{D}\phi^{(0)}_g(RG_{\ast}\nbigf)$
of $\phi^{(0)}_g(G^0_{\dagger}\nbigm)$.
We also have the induced $K$-Betti structure
$\lefttop{D}\phi^{(0)}_{g_X}(\nbigf)$
of $\phi^{(0)}_{g_X}\nbigm$.
By using Lemma \ref{lem;09.10.16.20},
we obtain
$\lefttop{D}\phi^{(0)}_g(RG_{\ast}\nbigf)
=RG_{\ast}\lefttop{D}\phi^{(0)}_{g_X}(\nbigf)$
under the isomorphism
$\phi^{(0)}_g(G^0_{\dagger}\nbigm)
\simeq
 G^0_{\dagger}\phi^{(0)}_{g_X}\nbigm$.
By the assumption $SI(<\dim X)$,
we obtain that 
$\lefttop{D}\phi^{(0)}_g\bigl(
 RG_{\ast}\nbigf\bigr)$
is a $K$-Betti structure of
$\phi^{(0)}_g\bigl(G_{\dagger}\nbigm\bigr)$.
Thus, we obtain 
Proposition \ref{prop;09.10.16.21}.
\hfill\qed

\subsection{Special case II}

Let $G:X\lrarr Y$ be a projective morphism
of complex manifolds.
Let $\varphi:Z\lrarr X$ be a projective morphism.
Let $D_Z$ be a hypersurface of $Z$.
Assume that
$\varphi_{|Z-D_Z}$ is an immersion.
Let $V$ be a meromorphic flat connection
on $(Z,D_Z)$ with a good $K$-Betti structure.

Suppose that we are given a sequence of
hypersurface pairs $\gbigh_Z$  of $Z$
contained in $D_Z$.
We obtain the $K$-holonomic $\nbigd_Z$-modules
$\nbigm:=\varphi_{\dagger}\gbigp_{\gbigh_Z}(V)$.

\begin{lem}
\label{lem;09.10.16.21}
Suppose $G^i_{\dagger}\nbigm=0$
unless $i=0$.
Then, 
the pre-$K$-holonomic $\nbigd_Y$-module
$G_{\dagger}^0\nbigm$
is $K$-holonomic.
\end{lem}
\pf
It follows from Proposition \ref{prop;09.10.16.21}.
\hfill\qed

\subsection{Special case III}

Let $\nbige$ be a locally free sheaf 
on a complex manifold $Y$.
Let $X$ be its projectivization.
Let $H_i$ $(i=0,1,2)$ be hyperplane subbundles.
Let $\nbign$ be a $K$-holonomic $\nbigd$-module on $X$
such that $\nbign(\ast H_0)=\nbign$.
By shrinking $Y$,
we may assume that $X=Y\times\proj^n$ for some $n$.
\begin{lem}
Let $A\subsetneq X$ be
any closed complex analytic subset.
If we shrink $Y$ appropriately,
there exists a meromorphic function $g$
on $X$ such that
(i) the poles of $g$ are contained in $H_0$,
(ii) $A$ is contained in 
$H_0\cup g^{-1}(0)$.
\end{lem}
\pf
Let $\nbigi_A$ denote the ideal sheaf of $A$ on $X$.
If $m$ is sufficiently large,
we have a non-zero section of
$\nbigi_A(mH_0)$ for $m$.
\hfill\qed

\begin{lem}
We can take a meromorphic function $g$ on $X$
such that (i) the poles of $g$ are contained in $H_0$,
(ii) $\nbign(\ast g)$ is obtained as
$\varphi_{\dagger}V$
for a cell
$\nbigc=(Z,U,\varphi,V)$.
(Note that we do not assume that
$V$ is a good meromorphic flat bundle
on $Z$.)
\end{lem}
\pf
We have a decomposition
of $\Supp(\nbign)$
into the locally closed complex analytic subsets 
$\coprod A_i$
such that the characteristic variety of $\nbign$
is $\coprod T^{\ast}_{A_i}X$.
Applying the previous lemma
to the lower dimensional strata,
we find a meromorphic function $g$ on $X$
such that
(i) the poles are contained in $H_0$,
(ii) $A_i\subset H_0\cup g^{-1}(0)$
 if $\dim A_i<\dim\Supp(\nbign)$.
By using the resolution of singularity
to the irreducible components of
$\Supp(\nbign)$ with the maximal dimension,
we obtain the cell.
\hfill\qed

\vspace{.1in}
Suppose that
$\vecH=(H_1,H_2)$ is non-characteristic with respect to
$\nbign$,
$\nbign(\ast g)$,
$\nbign(!g)(\ast H_0)$,
$\psi^{(a)}_g(\nbign,\ast H_0)$,
$\Xi^{(a)}_g(\nbign,\ast H_0)$ and
$\phi^{(a)}_g(\nbign,\ast H_0)$.
In this case,
$\vecH$ is non-characteristic with respect to
$\Pi^{a,b}_{g!}(\nbign,\ast H_0)$
and $\Pi^{a,b}_{g\ast}(\nbign)$
for any $a,b$.

\begin{lem}
\label{lem;09.10.16.25}
The induced pre-$K$-Betti structure of
$G_{\dagger}^0\gbigp_{\vecH}\nbign$
is a $K$-Betti structure.
\end{lem}
\pf
By $LOC(\leq n)$,
$\gbigp_{\vecH}
 \bigl(
 \Pi^{a,b}_{g!}(\nbign,\ast H_0)
 \bigr)$
and
$\gbigp_{\vecH}
 \bigl(
\Pi^{a,b}_{g\ast}\nbign
 \bigr)$
are naturally $K$-holonomic $\nbigd$-modules.
By Lemma \ref{lem;14.1.15.1},
we have
\[
 G_{\dagger}^i
 \gbigp_{\vecH}
 \bigl(
 \Pi^{a,b}_{g!}(\nbign,\ast H_0)
 \bigr)=0,
 \quad
 G_{\dagger}^i
 \gbigp_{\vecH}
\bigl(
\Pi^{a,b}_{g\ast}\nbign
\bigr)
=0
\]
unless $i=0$.
According to Lemma \ref{lem;09.10.16.21},
$G_{\dagger}^0
 \gbigp_{\vecH}
 \bigl(
 \Pi^{a,b}_{g!}(\nbign,\ast H_0)
 \bigr)$
and
$G_{\dagger}^0
 \gbigp_{\vecH}
\bigl(
\Pi^{a,b}_{g\ast}\nbign
\bigr)$
are $K$-holonomic.
Hence, we obtain that
$G^0_{\dagger}
 \gbigp_{\vecH}
 \Xi^{(a)}_g(\nbign,\ast H_0)$
and
$G^0_{\dagger}
 \gbigp_{\vecH}
 \psi^{(a)}_g(\nbign,\ast H_0)$
are $K$-holonomic.
We have the description of
$G_{\dagger}^0
 \gbigp_{\vecH}\nbign$ 
as the cohomology of the following complex
of pre-$K$-holonomic $\nbigd_Y$-modules:
\begin{multline*}
 G_{\dagger}^0
 \gbigp_{\vecH}
 \psi^{(1)}_g\bigl(\nbign,\ast H_0\bigr)
\lrarr 
 G_{\dagger}^0
 \gbigp_{\vecH}
 \Xi^{(0)}_g\bigl(\nbign,\ast H_0\bigr)
\oplus
 G_{\dagger}^0
 \gbigp_{\vecH}
 \phi^{(0)}_g(\nbign,\ast H_0)
 \\
\lrarr
 G_{\dagger}^0
 \gbigp_{\vecH}
 \psi^{(0)}_g\bigl(\nbign,\ast H_0\bigr).
\end{multline*}
By $SI(<n)$,
we obtain that 
$G_{\dagger}^0
 \gbigp_{\vecH}
 \phi^{(0)}_g(\nbign,\ast H_0)$
is $K$-holonomic.
Then, we obtain Lemma \ref{lem;09.10.16.25}.
\hfill\qed

\subsection{Proof of Theorem 
\ref{thm;09.10.16.5}}
\label{subsection;09.10.16.6}

It is enough to consider the case 
$X=\proj(\nbige)$ for some locally free sheaf
$\nbige$ on $Y$.
Let $(\nbigm,\nbigf)$ be a $K$-holonomic $\nbigd_X$-module
with $\dim\Supp\nbigm\leq n$.
Let us prove that
$F^i_{\dagger}(\nbigm,\nbigf)$
are $K$-holonomic.

We take a resolution
$\nbign_{\bullet}$ of $\nbigm$
as in (\ref{eq;13.4.27.200}) of \S\ref{subsection;09.10.16.2}.
Then, by applying the construction $\nbigq_{\bullet,\bullet}$
in \S\ref{subsection;09.10.16.1} to each $\nbign_i$,
we take a resolution 
$\Tot(\nbigq(\nbign_{\bullet})_{\bullet,\bullet})$
of $\nbigm$.
It is naturally equipped with
the $K$-Betti structure
$\Tot\bigl(\nbigf^{\nbigq}_{\bullet,\bullet,\bullet}\bigr)$.
Then, $F_{\dagger}^i(\nbigm,\nbigf)$
is described as the $i$-th cohomology of 
$\Tot\Bigl(
 F_{\dagger}^0\bigl(
 \nbigq(\nbign_{\bullet})_{\bullet,\bullet},
 \nbigf^{\nbigq}_{\bullet,\bullet,\bullet}
 \bigr)
 \Bigr)$. 
Hence, it is enough to show that
$F_{\dagger}^0\bigl(
 \nbigq(\nbign_{\bullet})_{\bullet,\bullet},
 \nbigf^{\nbigq}_{\bullet,\bullet,\bullet}
 \bigr)$
are $K$-holonomic.
By the construction,
we have
$\dim\Supp\nbigq(\nbign_k)_{i,j}
<\dim\Supp\nbigm$
for $(k,i,j)\neq (0,0,0)$,
to which we can apply the inductive assumption.
Hence, it is enough to show that
$F_{\dagger}^0\bigl(
 \nbigq(\nbign_0)_{0,0},
 \nbigf^{\nbigq}_{0,0,0}
 \bigr)$
is $K$-holonomic,
which follows from Lemma \ref{lem;09.10.16.25}.
\hfill\qed

%% file: 9.tex
We study the standard functors
on the derived category of 
algebraic $K$-holonomic $\nbigd$-modules.
It is enough to follow very closely 
the arguments 
in \cite{beilinson1}, \cite{beilinson2},
\cite{bbd}
and \cite{saito2}, \cite{saito4}.
This section is included
for a rather expository purpose.

%% file: 9.1.tex
Let $X$ be a smooth complex quasi-projective variety.
We take a smooth projective completion
$X\subset \Xbar$ such that 
$D=\Xbar-X$ is a hypersurface.
We set
$\Hol(X,K):=\Hol\bigl(\Xbar, \ast D,K\bigr)$,
which is independent of 
the choice of a completion $\Xbar$
(Proposition \ref{prop;09.10.5.50}).
\index{category $\Hol(X,K)$}
Let $D^b\bigl(\Hol(X,K)\bigr)$ denote
the derived category of $\Hol(X,K)$.
We will implicitly use the following 
obvious lemma.
(Later, we will prove a stronger version
in Theorem \ref{thm;09.11.13.21}.)
\begin{lem}
The forgetful functor
$\Hol(X,K)\lrarr\Hol(X)$ is faithful.
\hfill\qed
\end{lem}

\subsection{Dual}
\index{dual functor $\DDD$}
For any
$\nbigm \in \Hol\bigl(\Xbar,\ast D,K\bigr)$,
we have the $K$-holonomic
$\nbigd_{\Xbar(\ast D)}$-module
$\DDD_X\nbigm:=
\DDD_{\Xbar}(\nbigm)(\ast D)$.

\begin{lem}
$\DDD_X(\nbigm)$
is well defined in 
$\Hol(X,K)$.
\end{lem}
\pf
Let $\Xbar'$ be another smooth projective
compactification of $X$.
Put $D':=\Xbar'-X$.
We may assume to have a projective morphism
$\varphi:\Xbar'\lrarr \Xbar$ such that
$\varphi_{|X}=\id_X$.
We have a $K$-holonomic $\nbigd_{\Xbar'(\ast D')}$-module
$\nbigm'$
such that 
$\varphi_{\dagger}\nbigm'=\nbigm$,
which is unique up to canonical isomorphisms.
Then,
the natural isomorphism
$\varphi_{\dagger}(\DDD\nbigm')(\ast D')
\simeq
 \DDD(\nbigm)(\ast D)$
preserves the $K$-Betti structure
by the uniqueness (Corollary \ref{cor;13.4.27.300}).
It implies the claim of the lemma.
\hfill\qed

\begin{cor}
There exists a functor
$\DDD_X$ on $\Hol(X,K)$
which is compatible with
the standard duality functors
on $\Hol(X)$ and the category
of $K$-perverse sheaves.
We also have a functor
$\DDD_X$ on $D^b\bigl(\Hol(X,K)\bigr)$,
compatible with
the standard duality functors
on $D^b_{\hol}(X)$
and $D^b_c(K_X)$.
They are unique up to natural equivalences.
\hfill\qed
\end{cor}

We use the symbol $\lefttop{K}\DDD_X$
if we would like to emphasize
that it is a functor for $K$-holonomic $\nbigd$-modules.

\begin{lem}
For $\nbigm,\nbign\in\Hol(X,K)$,
we have a natural isomorphism:
\[
 \Ext^{i}_{\Hol(X,K)}(\nbigm,\nbign)
\simeq
 \Ext^i_{\Hol(X,K)}\bigl(
 \lefttop{K}\DDD_X\nbign,
 \lefttop{K}\DDD_X\nbigm
 \bigr)
\]
\end{lem}
\pf
It follows from the comparison of
Yoneda extensions.
\hfill\qed

\subsection{Localization}

\index{localization}

Let $H$ be a hypersurface of $X$.
As is shown in Theorem \ref{thm;09.10.18.400}
and Proposition \ref{prop;09.10.5.50},
we have the localization:
\[
\ast H:\Hol(X,K)\lrarr\Hol(X,K),
\quad
\nbigm\longmapsto \nbigm(\ast H)
\]
It is an exact functor.
By considering the dual,
we obtain an exact functor:
\[
 !H:
 \Hol(X,K)\longmapsto
 \Hol(X,K),
\quad
\nbigm\longmapsto \nbigm(!H)
\]
They induce exact functors
$\ast H$ and $!H$ on
$D^b\bigl(\Hol(X,K)\bigr)$.

\begin{lem}
\label{lem;10.1.13.1}
For $\nbigm,\nbign\in \Hol(X,K)$,
we have the following natural isomorphisms:
\[
 \Ext^i_{\Hol(X,K)}
 \bigl(\nbigm,\nbign(\ast D)\bigr)
\simeq
 \Ext^i_{\Hol(X,K)}
 \bigl(\nbigm(\ast D),\nbign(\ast D)\bigr)
\]
\[
  \Ext^i_{\Hol(X,K)}
 \bigl(\nbigm(!D),\nbign\bigr)
\simeq
 \Ext^i_{\Hol(X,K)}
 \bigl(\nbigm(!D),\nbign(!D)\bigr)
\]
\end{lem}
\pf
It follows from comparisons
of Yoneda extensions.
\hfill\qed

\subsection{Nearby cycle,
vanishing cycle and maximal functors}

\index{nearby cycle functor}
\index{vanishing cycle functor}
\index{maximal functor}

Let $g$ be an algebraic function on $X$.
By Proposition \ref{prop;09.12.6.1},
we have the exact functors
$\Pi^{a,b}_{g\star}$
($\star=\ast,!$)
on $\Hol(X,K)$ 
given by
$\Pi^{a,b}_{g\star}(\nbigm):=
 \bigl(\nbigm\otimes\gbigi_g^{a,b}\bigr)
 (\star g)$
and $a,b\in\seisuu$.
Hence, we obtain the exact functors
$\Xi^{(a)}_g$, $\psi^{(a)}_g$ 
and $\phi^{(a)}_g$ on $\Hol(X,K)$.
They induce the corresponding exact functors
on $D^b\bigl(\Hol(X,K)\bigr)$.
We use the symbols
$\lefttop{K}\Xi^{(a)}_g$,
$\lefttop{K}\psi^{(a)}_g$
and $\lefttop{K}\phi^{(a)}_g$,
when we would like to emphasize
that they are functors for
$K$-holonomic $\nbigd$-modules.
We remark that the functors are not 
compatible with the forgetful functor
$D^b\bigl(\Hol(X,K)\bigr)\lrarr D^b_c(K_X)$.
The $K$-Betti structure of
$\lefttop{K}\psi^{(a)}_g(\nbigm,\nbigf)$
is denoted by $\lefttop{D}\psi^{(a)}_g(\nbigf)$
for the distinction,
when we would like to emphasize it.
Similar notations such as
$\lefttop{D}\Xi_g^{(a)}$
and 
$\lefttop{D}\phi_g^{(a)}$
are used.
\index{functor $\lefttop{D}\psi^{(a)}_g$}
\index{functor $\lefttop{D}\phi^{(a)}_g$}
\index{functor $\lefttop{D}\Xi^{(a)}_g$}

%% file: 9.2.tex
\subsection{Statements}

Let $f:X\lrarr Y$ be an algebraic morphism
of quasi-projective varieties.
We take a commutative diagram
\[
 \begin{CD}
 X @>{f}>> Y \\
 @V{a_1}VV @V{a_2}VV \\
 \Xbar @>{\fbar}>> \Ybar
 \end{CD}
\]
where (i) $a_i$ are open immersions,
(ii) $\Xbar$ and $\Ybar$ are smooth projective,
(iii) $H_X=\Xbar-X$ and $H_Y:=\Ybar-Y$
 are hypersurfaces.
We have a natural equivalence
between $\Hol\bigl(\Xbar,\ast H_X,K\bigr)$
and $\Hol(X,K)$.
Let $\nbigmbar\in
 \Hol\bigl(\Xbar,\ast H_X,K\bigr)$
correspond to $\nbigm\in\Hol(X,K)$.
According to Theorem \ref{thm;09.10.16.5},
we obtain the following objects in $\Hol(Y,K)$:
\[
\lefttop{K}f^i_{\ast}(\nbigm)
:=f^i_{\dagger}\nbigmbar,
\quad
\lefttop{K}f^i_{!}(\nbigm)
:=f^i_{\dagger}\bigl(
 \nbigmbar(!H_X)\bigr)(\ast H_Y).
\]
They are independent of the choice of $\Xbar$
up to natural isomorphisms.
Thus, we obtain cohomological functors
$\lefttop{K}f^i_{\ast},
 \lefttop{K}f^i_{!}
 :\Hol(X,K)\lrarr\Hol(Y,K)$
for $i\in\seisuu$.

\begin{prop}
\label{prop;09.11.10.10}
For $\star=!,\ast$,
there exists a functor of
triangulated categories
\[
 \lefttop{K}f_{\star}:
 D^b\bigl(\Hol(X,K)\bigr)\lrarr D^b(\Hol(Y,K))
\]
such that
(i) it is compatible with the standard functor
$f_{\star}:D^b_{\hol}(X)\lrarr D^b_{\hol}(Y)$,
(ii) the induced functor 
$H^i(\lefttop{K}f_{\star}):
 \Hol(X,K)\lrarr \Hol(Y,K)$
is isomorphic to $\lefttop{K}f_{\star}^i$.
It is characterized by the property (i) and (ii)
up to natural equivalences.
\end{prop}
\index{functor $\lefttop{K}f_{\ast}$}
\index{functor $\lefttop{K}f_{\bikkuri}$}

As in \S4 of \cite{saito2},
the pull back is defined to be
the adjoint of the push-forward.

\begin{prop}
\label{prop;09.11.10.12}
$\lefttop{K}f_{!}$
has the right adjoint
$\lefttop{K}f^{!}$,
and
$\lefttop{K}f_{\ast}$
has the left adjoint
$\lefttop{K}f^{\ast}$.
Thus, we obtain the following functors:
\[
 \lefttop{K}f^{\star}:
 D^b(\Hol(Y,K))\lrarr
 D^b\bigl(\Hol(X,K)\bigr)
\quad\quad
 (\star=!,\ast)
\]
They are compatible with
the corresponding functors of
holonomic $\nbigd$-modules
with respect to the forgetful functor.
\end{prop}
\index{functor $\lefttop{K}f^{\bikkuri}$}
\index{functor $\lefttop{K}f^{\ast}$}

Let us consider the case where
$f$ is a closed immersion,
via which $X$ is regarded as 
a submanifold of $Y$.
Let $D^{b}_{X}(\Hol(Y,K))$ be
the full subcategory of $D^b(\Hol(Y,K))$
which consists of the objects $\nbigm^{\bullet}$
such that the supports of the cohomology
$\bigoplus_i\nbigh^i\nbigm^{\bullet}$
are contained in $X$.
\index{category $D^{b}_{X}(\Hol(Y,K))$}
\begin{prop}
\label{prop;10.1.12.10}
The natural functor
$\lefttop{K}f_{!}:
 D^b\Hol(X,K)\!\lrarr\!
 D^b_{X}\Hol(Y,K)$
is an equivalence.
\end{prop}

\begin{rem}
\footnote{This remark is due to the referee.}
It is a deep theorem of Z. Mebkhout
that the irregularity sheaf of any holonomic $\nbigd$-module
$\nbigm$ is a perverse sheaf.
See {\rm\cite{Mebkhout-positivity}}.
By using the above functors,
in the algebraic case,
we obtain that
the irregularity sheaf of a $K$-holonomic $\nbigd$-module
is equipped with an induced $K$-structure
which is clear by the definition of the irregularity sheaf.
We may apply the argument
even in the analytic case.
\hfill\qed
\end{rem}

\subsection{Preliminary}

Let $X$ be a smooth complex projective variety
with a hypersurface $D$.
Let $D^b\bigl(\Hol(X,\ast D,K)\bigr)$
denote the derived category of
$\Hol\bigl(X,\ast D,K\bigr)$.
Similarly,
let $D^b\bigl(\Hol(X,\ast D)\bigr)$
denote the derived category of
$\Hol\bigl(X,\ast D\bigr)$.

Let $f:X\lrarr Y$ be a morphism
of smooth projective varieties.
Let $D_X$ and $D_Y$ be hypersurfaces
of $X$ and $Y$ respectively,
such that $D_X\supset f^{-1}(D_Y)$.
We have the functor
$\lefttop{K}f^i_{\ast}:
 \Hol\bigl(X,\ast D_X,K\bigr)\lrarr
 \Hol\bigl(Y,\ast D_Y,K\bigr)$,
naturally given by $f^i_{\dagger}$.
We have a decomposition
$D_{X}=D_{X1}\cup D_{X2}$
such that $D_{X2}=f^{-1}(D_Y)$.
We set $\vecD_X:=(D_{X1},D_{X2})$.
We have the functor
$\lefttop{K}f^i_{!}:
 \Hol\bigl(X,\ast D_X,K\bigr)\lrarr
 \Hol\bigl(Y,\ast D_Y,K\bigr)$
given by
$\lefttop{K}f^i_!(\nbigm,\nbigf)
=f^i_{\dagger}
 \gbigp'_{\vecD_X}
 \nbigm$.

\begin{lem}
\label{lem;09.11.10.20}
For $\star=\ast,!$,
there exist functors
$\lefttop{K}f_{\star}:
 D^b\bigl(\Hol(X,\ast D_X,K)\bigr)
 \lrarr D^b\bigl(\Hol(Y,\ast D_Y,K)\bigr)$
such that
(i) they are compatible with the standard functors
$f_{\star}:
 D^b\bigl(\Hol(X,\ast D_X)\bigr)
\lrarr
 D^b\bigl(\Hol(Y,\ast D_Y)\bigr)$
by the forgetful functors,
(ii) the induced functor
$H^i(\lefttop{K}f_{\star}):
 \Hol\bigl(X,\ast D_X,K\bigr)\lrarr
 \Hol\bigl(Y,\ast D_Y,K\bigr)$
are isomorphic to $\lefttop{K}f_{\star}^i$.
It is characterized by (i) and (ii)
up to natural equivalences.
\end{lem}
\pf
Let us consider the case $\star=\ast$.
Let $\nbigm$ be a $K$-holonomic 
$\nbigd_{X(\ast D_X)}$-module.
Let $\vecH=(H_1,\ldots,H_M)$ 
be a tuple of hypersurfaces of $X$.
We put $H_I:=\bigcup_{i\in I}H_i$.
We take a $K$-vector space $U$
with a base $(e_1,\ldots,e_M)$.
For $I=(i_1,\ldots,i_m)\subset \{1,\ldots,M\}$,
let $U_I$ denote the subspace of
$\bigwedge^{\bullet}U$
generated by $e_{i_1}\wedge\cdots\wedge e_{i_m}$.
For $m\geq 0$,
we set
\[
 \nbigc^m_{\ast\vecH}(\nbigm):=
 \bigoplus_{|I|=m+1}
 \nbigm(\ast H_I)
 \otimes U_I.
\]
For $Ii:=I\sqcup\{i\}\subset \{1,\ldots,M\}$,
the natural morphism
$\nbigm(\ast H_I)
\lrarr
 \nbigm(\ast H_{Ii})$
and the multiplication of $e_i$
induce
$\nbigm(\ast H_I)
 \otimes U_I
\lrarr
 \nbigm(\ast H_{Ii})
 \otimes U_{Ii}$.
They give a complex
$\bigl(
 \nbigc^{\bullet}_{\ast\vecH}(\nbigm),\del_{\ast\vecH}
\bigr)$.
We have a natural morphism of complexes
$\nbigm\lrarr \nbigc^{\bullet}_{\ast\vecH}(\nbigm)$.
If $\bigcap H_i=\emptyset$,
it is a quasi-isomorphism.
Suppose we are given a tuple of hypersurfaces
$\vecL=(L_1,\ldots,L_{N})$.
We put $\vecH\vecL=(H_1,\ldots,H_M,L_1,\ldots,L_N)$.
The natural projection
$\nbigc^{\bullet}_{\ast\vecH\vecL}(\nbigm)
\lrarr
 \nbigc^{\bullet}_{\ast\vecH}(\nbigm)$
gives a complex of morphisms.

Let $\vecH'=(H_1',\ldots,H_{N}')$
be a tuple of hypersurfaces on $X$.
We take a $K$-vector space $U'$
with a base $(e'_1,\ldots,e'_N)$.
For $J=(j_1,\ldots,j_n)\subset\{1,\ldots,N\}$,
let $U_J'$ be the subspace of
$\bigwedge U'$
generated by
$e'_{j_1}\wedge\cdots\wedge e_{j_n}'$.
For $n\leq 0$,
we set
\[
 \nbigc^{n}_{!\vecH'}(\nbigm):=
 \bigoplus_{|J|=-n+1}\nbigm(!H_J')\otimes U_J'.
\]
Let $e^{\prime\lor}_j$ denote the dual base.
For $Jj=J\sqcup\{j\}\subset\{1,\ldots,N\}$,
the natural morphism
$\nbigm(H_{Jj}')
\lrarr
 \nbigm(H_J')$
and the inner product of $e_j^{\prime\lor}$
induce
$\nbigm(H_{Jj}')\otimes U'_{Jj}
\lrarr
 \nbigm(H_J')\otimes U_J'$.
They give a complex
$\bigl(
 \nbigc^{\bullet}_{!\vecH'}(\nbigm),
 \del_{!\vecH'}
 \bigr)$.
We have a natural morphism of complexes
$\nbigc_{!\vecH'}(\nbigm)
\lrarr
 \nbigm$.
If $\bigcap H_i'=\emptyset$,
it is a quasi-isomorphism.
Suppose that we are given a tuple of hypersurfaces
$\vecL'=(L_1',\ldots,L'_M)$.
We put
$\vecH'\vecL'=(H_1',\ldots,H'_N,L_1',\ldots,L'_M)$.
The natural inclusion
$\nbigc^{\bullet}_{!\vecH'}(\nbigm)
\lrarr
 \nbigc^{\bullet}_{!\vecH'\vecL'}(\nbigm)$
gives a quasi-isomorphism.

Let $\nbigm^{\bullet}$ be a complex of
$K$-holonomic $\nbigd_{X(\ast D_X)}$-modules.
Let $\vecH$ and $\vecH'$
be tuples of hypersurfaces.
The total complex of
$\nbigc^{\bullet}_{\ast\vecH}
 \nbigc^{\bullet}_{!\vecH'}
 (\nbigm^{\bullet})$
is denoted by
$\nbigc^{\bullet}_{\ast\vecH!\vecH'}(\nbigm^{\bullet})$.
The total complexes of
$\nbigc^{\bullet}_{\ast\vecH}(\nbigm^{\bullet})$
and 
$\nbigc^{\bullet}_{!\vecH'}(\nbigm^{\bullet})$
are also denoted by the same notation.
We assume $\bigcap H_i=\bigcap H_j'=\emptyset$.
We have the following
natural quasi-isomorphisms of complexes:
\[
\begin{CD}
 \nbigc^{\bullet}_{\ast\vecH!\vecH'}(\nbigm^{\bullet})
@>>>
 \nbigc^{\bullet}_{\ast\vecH}(\nbigm^{\bullet})
@<<<
 \nbigm^{\bullet}
\end{CD}
\]

Let $(\vecH_i,\vecH_i')$ $(i=1,2)$
be tuples of hypersurfaces as above.
We say that we have a morphism
$(\vecH_1,\vecH_1')\lrarr
 (\vecH_2,\vecH_2')$
if $\vecH_1\supset \vecH_2$
and $\vecH_1'\subset\vecH_2'$
are satisfied.
Then, we have a naturally defined quasi-isomorphism
of complexes:
\[
 \nbigc^{\bullet}_{\ast\vecH_1!\vecH'_1}(\nbigm^{\bullet})
\lrarr
 \nbigc^{\bullet}_{\ast\vecH_2!\vecH'_2}(\nbigm^{\bullet})
\]

For a tuple of ample hypersurfaces
$(\vecH,\vecH')$ which is
non-characteristic with respect to $\nbigm^{\bullet}$
(\S\ref{subsection;09.10.16.1}),
we have
$f^i_{\dagger}\nbigm^p(\ast H_I!H_J)=0$
unless $i=0$.
For each $\nbigm^{\bullet}$,
we choose such 
$(\vecH(\nbigm^{\bullet}),\vecH'(\nbigm^{\bullet}))$.
We obtain a complex of
$K$-holonomic $\nbigd_{Y(\ast D_Y)}$-modules:
\[
 \lefttop{K}f_{\ast}(\nbigm^{\bullet}):=
 f_{\dagger}^0 
 \nbigc^{\bullet}_{\ast\vecH(\nbigm^{\bullet})!\vecH'(\nbigm^{\bullet})}
 (\nbigm^{\bullet})
\]

Let 
$\nbigm_1^{\bullet}
 \stackrel{a}{\llarr}
 \nbigm_1^{\prime\bullet}
 \stackrel{b}{\lrarr}
 \nbigm_2^{\bullet}$
be morphisms,
where $a$ is a quasi-isomorphism.
We take a tuple of ample hypersurfaces
$(\vecH,\vecH')$
such that 
(i) the tuple $(\vecH,\vecH')$ is non-characteristic with respect to 
$\nbigm_i^{\bullet}$ and $\nbigm_1^{\prime\bullet}$,
(ii) the tuple
$\bigl(\vecH,\vecH(\nbigm_i^{\bullet}),
 \vecH',\vecH'(\nbigm_i')\bigr)$
is non-characteristic with respect to $\nbigm^{\bullet}_i$.
We have the following morphism of complexes
\[
\begin{CD}
\nbigc^{\bullet}_{\ast\vecH!\vecH'}(\nbigm_1^{\bullet})
@<{a_0}<<
\nbigc^{\bullet}_{\ast\vecH!\vecH'}(\nbigm_1^{\prime\bullet})
@>>>
\nbigc^{\bullet}_{\ast\vecH!\vecH'}(\nbigm_2^{\bullet})
\end{CD}
\]
Here, $a_0$ is a quasi-isomorphism.
We set $\vecH_i=\vecH(\nbigm_i^{\bullet})$
and $\vecH'_i=\vecH'(\nbigm_i^{\bullet})$.
We have the following quasi-isomorphisms:
\[
 \begin{CD}
  \nbigc^{\bullet}
 _{\ast\vecH!\vecH'}
 (\nbigm_i^{\bullet})
 @>{a_{i1}}>>
 \nbigc^{\bullet}_{\ast\vecH!\vecH'}
 \nbigc^{\bullet}_{\ast\vecH_i}
 (\nbigm_i^{\bullet})
 @<{a_{i2}}<<
 \nbigc^{\bullet}_{\ast\vecH,\vecH'}
 \nbigc^{\bullet}_{\ast\vecH_i!\vecH'_i}
 (\nbigm_i^{\bullet})
 \end{CD}
\]
Note that 
$\nbigc^{\bullet}_{\ast\vecH,\vecH'}
 \nbigc^{\bullet}_{\ast\vecH_i!\vecH'_i}
 (\nbigm_i^{\bullet})$
and 
$ \nbigc^{\bullet}_{\ast\vecH_i!\vecH'_i}
 \nbigc^{\bullet}_{\ast\vecH,\vecH'}
 (\nbigm_i^{\bullet})$
are naturally isomorphic.
We also have the following quasi-isomorphisms:
\[
 \begin{CD}
  \nbigc^{\bullet}_{\ast\vecH_i!\vecH'_i}
 \nbigc^{\bullet}_{\ast\vecH!\vecH'}
 (\nbigm_i^{\bullet})
@>{a_{i3}}>>
  \nbigc^{\bullet}_{\ast\vecH_i!\vecH'_i}
 \nbigc^{\bullet}_{\ast\vecH}
 (\nbigm_i^{\bullet})
@<{a_{i4}}<<
  \nbigc^{\bullet}_{\ast\vecH_i\vecH'_i}
 (\nbigm_i^{\bullet})
 \end{CD}
\]
Note that $f_{\dagger}^0(a_0)$
and $f_{\dagger}^0(a_{ij})$
are quasi-isomorphisms.
They induce a morphism
in $D^b(\Hol(Y,\ast D_Y,K))$:
\begin{equation}
 \label{eq;14.1.12.1}
 \lefttop{K}f_{\ast}^0(\nbigm_1^{\bullet})
\lrarr
 \lefttop{K}f_{\ast}^0(\nbigm_2^{\bullet})
\end{equation}
If we are given morphisms
$\nbigm_1^{\bullet}
 \stackrel{a}{\llarr}
 \nbigm_1^{\prime\bullet}
 \stackrel{b}{\lrarr}
 \nbigm_2^{\bullet}$
such that $a'$ and $b'$
are chain homotopic to $a$ and $b$
respectively,
it is easy to check 
that the induced morphisms
(\ref{eq;14.1.12.1}) 
in $D^b(\Hol(Y,\ast D_Y,K))$
are the same.

Let us check that (\ref{eq;14.1.12.1})
is independent from the choice of
$(\vecH,\vecH')$.
Let $(\vecL,\vecL')$ be other choice.
Take a sequence of
sufficiently generic ample hypersurfaces
$(\vecH^{(j)},\vecH^{\prime(j)})$
$(j=1,\ldots,2L)$
satisfying the above conditions,
such that
(i)
$(\vecH^{(1)},\vecH^{\prime(1)})
=(\vecH,\vecH')$ and
 $(\vecH^{(2L)},\vecH^{\prime(2L)})
=(\vecL,\vecL')$,
(ii) we have morphisms
\[
 (\vecH^{(2m-1)},\vecH^{\prime(2m-1)})
\llarr
(\vecH^{(2m)},\vecH^{\prime(2m)})
\lrarr
 (\vecH^{(2m+1)},\vecH^{\prime(2m+1)}).
\]
Then, it is easy to check that
$(\vecH,\vecH')$  and $(\vecL,\vecL')$
induce 
the same morphism (\ref{eq;14.1.12.1})
in $D^b(\Hol(Y,\ast D_Y,K))$.
Hence, the morphism (\ref{eq;14.1.12.1})
depends only on the morphism
in $D^b(\Hol(X,\ast D_X,K))$
determined by $(a,b)$,
i.e.,
we obtain a morphism
\[
\Hom_{D^b\bigl(\Hol(X,K)\bigr)}\bigl(
 \nbigm^{\bullet}_1,\nbigm^{\bullet}_2
 \bigr)
\lrarr
 \Hom_{D^b(\Hol(Y,K))}\bigl(
 \lefttop{K}f_{\ast}\nbigm^{\bullet}_1,\,
 \lefttop{K}f_{\ast}\nbigm^{\bullet}_2
 \bigr).
\]
Thus, we obtain a functor
$D^b(\Hol(X,\ast D_X,K))
\lrarr
 D^b(\Hol(Y,\ast D_Y,K))$.
We set
$\lefttop{K}f_{!}:=
 \lefttop{K}\DDD_Y\circ
 \lefttop{K}f_{\ast}\circ 
 \lefttop{K}\DDD_X$.
By the construction,
they satisfy the conditions (i) and (ii).
The uniqueness follows from
the existence of a resolution 
by $K$-holonomic $\nbigd$-modules
$\nbign$ such that
$f_{\dagger}^i\nbign=0$
unless $i=0$.
\hfill\qed

\subsection{Proof of Proposition
\ref{prop;09.11.10.10}}

We take projective completions
$X\subset \Xbar$
and $Y\subset \Ybar$
with the following commutative diagram:
\begin{equation}
 \label{eq;09.11.10.11}
\begin{CD}
X @>{\subset}>> \Xbar \\
 @V{f}VV @V{\fbar}VV \\
Y @>{\subset}>> \Ybar 
\end{CD}
\end{equation}
Set $D_X:=\Xbar-X$ and $D_Y:=\Ybar-Y$.
The functor
$\lefttop{K}f_{\star}:
 D^b\bigl(\Hol(\Xbar,\ast D_X,K)\bigr)
\lrarr
 D^b\bigl(\Hol(\Ybar,\ast D_Y,K)\bigr)$
induces
$\lefttop{K}f_{\star}:
 D^b\bigl(\Hol(X,K)\bigr)
\lrarr
 D^b(\Hol(Y,K))$.

Let $X\subset \Xbar'$
and $Y\subset \Ybar'$ be 
other projective completions
with a commutative diagram
as in (\ref{eq;09.11.10.11}).
We set $D_X':=\Xbar'-X$
and $D_Y':=\Ybar-Y$.
Let us prove that the induced morphisms
$\lefttop{K}f_{\star}:
 D^b\bigl(\Hol(X,K)\bigr)
\lrarr
 D^b(\Hol(Y,K))$
are equal up to natural equivalences.
It is enough to consider the case 
where we have the following commutative diagram:
\[
 \begin{CD}
 \Xbar' @>{\fbar'}>> \Ybar' \\
 @V{\varphi_X}VV @V{\varphi_Y}VV \\
 \Xbar @>{\fbar}>> \Ybar
 \end{CD}
\]
Here, $\varphi_X$ and $\varphi_Y$
are projective and birational
such that
$\varphi_X^{-1}(D_X)=D_X'$
and $\varphi_Y^{-1}(D_Y)=D_Y'$.
We have the following diagrams
which are commutative up to equivalences:
\[
 \begin{CD}
 D^b\bigl(\Hol(\Xbar',\ast D_X',K)\bigr)
 @>{\lefttop{K}f_{\star}}>>
 D^b\bigl(\Hol(\Ybar',\ast D_Y',K)\bigr) \\
 @V{\lefttop{K}\varphi_{X\star}}VV 
 @V{\lefttop{K}\varphi_{Y\star}}VV \\
 D^b\bigl(\Hol(\Xbar,\ast D_X,K)\bigr)
 @>{\lefttop{K}f_{\star}}>>
 D^b\bigl(\Hol(\Ybar,\ast D_Y,K)\bigr)
 \end{CD}
\]
It implies that 
$\lefttop{K}f_{\star}:
 D^b\bigl(\Hol(X,K)\bigr)
\lrarr
 D^b(\Hol(Y,K))$
are independent of the choice of
projective completions
up to equivalences.
Thus, the proof of Proposition 
\ref{prop;09.11.10.10}
is finished.
\hfill\qed

\subsection{Proof of Proposition
\ref{prop;10.1.12.10}}

Let $\nbigm,\nbign\in\Hol(X,K)$.
According to Proposition 3.1.16 of \cite{bbd},
it is enough to check the following effaceability:
\begin{itemize}
\item
For any $f\in \Ext^i_{\Hol(Y,K)}(\nbigm,\nbign)$,
there exists a monomorphism $\nbign\lrarr \nbign'$
in $\Hol(X,K)$ such that 
the image of $f$ 
in $\Ext^i_{\Hol(Y,K)}(\nbigm,\nbign')$ is $0$.
\end{itemize}
We can prove it by using the arguments
in \S2.2.1 and \S2.2.2 in \cite{beilinson1}.
\hfill\qed

\subsection{Proof of Proposition
\ref{prop;09.11.10.12}}
\label{subsection;13.4.24.100}

It is enough to consider the cases
(i) $f$ is a closed immersion,
(ii) $f$ is a projection $X\times Y\lrarr Y$.
We closely follow the arguments
in \S2.19 and \S4.4 of \cite{saito2}.

\subsubsection{Closed immersion}

Let $f:X\lrarr Y$ be a closed immersion.
Let $\nbigm^{\bullet}$ 
be a complex of $K$-holonomic $\nbigd_Y$-modules.
Let $H_i$ $(i=1,\ldots,N)$
be sufficiently general ample hypersurfaces of $Y$
such that 
(i) $H_i\supset X$,
(ii) $\nbigm^{\bullet}
 \lrarr \nbigm^{\bullet}(\ast H_i)$ are monomorphisms, 
(iii) $\bigcap_{i=1}^N H_i=X$.
For any subset $I=(i_1,\ldots,i_m)
 \subset \{1,\ldots,N\}$,
let $\cnum_I$ be the subspace of
$\bigwedge^{m}\cnum^N$
generated by
$e_{i_1}\wedge\cdots \wedge e_{i_m}$,
where $e_i\in\cnum^N$ denotes
an element whose $k$-th entry is
$1$ $(k=i)$ or $0$ $(k\neq i)$.
For $I=I_0\sqcup\{i\}$,
we set $H_I=\bigcup_{i\in I}H_i$.
The inclusion
$\nbigm^p(\ast H_{I_0})
\lrarr \nbigm^p(\ast H_{I})$
and the multiplication of
$e_i$ induces
$\nbigm^p(\ast H_{I_0})\otimes\cnum_{I_0}
\lrarr
 \nbigm^p(\ast H_I)\otimes\cnum_I$.
For $m\geq 0$,
we put
$\nbigc^{m}
 (\nbigm^{p},\ast\vecH):=
 \bigoplus_{|I|=m}
 \nbigm^p(\ast H_I)\otimes\cnum_I$,
and we obtain the double complex
$\nbigc^{\bullet}
 (\nbigm^{\bullet},\ast \vecH)$.
The total complex is denoted by
$\Tot\nbigc^{\bullet}(\nbigm^{\bullet},\ast\vecH)$.
It is easy to observe that
the support of the cohomology of
$\Tot\nbigc^{\bullet}(\nbigm^{\bullet},\ast\vecH)$
is contained in $X$.
According to Proposition \ref{prop;10.1.12.10},
we obtain
$\lefttop{K}f^{!}\nbigm^{\bullet}:=
 \Tot\nbigc^{\bullet}
 (\nbigm^{\bullet},\ast \vecH)$
in $D^b(\Hol(X,K))$.
We obtain a functor
$\lefttop{K}f^{!}:
 D^b(\Hol(Y,K))\lrarr
 D^b\bigl(\Hol(X,K)\bigr)$
as in Lemma \ref{lem;09.11.10.20}.
Note that the underlying 
$\nbigd_Y$-complex
is naturally quasi-isomorphic to
$f^{!}\nbigm^{\bullet}$,
where $f^{!}$ is the left adjoint of
$f_{\dagger}:D^b_{\hol}(X)\lrarr D^b_{\hol}(Y)$.

We have the naturally defined morphism
$\alpha:\Tot\nbigc^{\bullet}
 (\nbigm^{\bullet},\ast\vecH)
\lrarr \nbigm^{\bullet}$.
We put $\nbigk^{\bullet}:=\Cone(\alpha)$.
We have another description.
For $m\geq 0$,
we put 
$\nbigcbar^{m}
 (\nbigm^{p},\ast\vecH):=
 \bigoplus_{|I|=m+1}
 \nbigm^p(\ast H_I)\otimes\cnum_I$,
and we obtain the double complex
$\nbigcbar^{\bullet}
 (\nbigm^{\bullet},\ast \vecH)$.
We have a natural quasi-isomorphism
$\nbigk^{\bullet}\simeq 
 \Tot\nbigcbar^{\bullet}(\nbigm^{\bullet},\ast\vecH)$.
By using the second description and
Lemma \ref{lem;10.1.13.1},
we obtain the following vanishing
for any $\nbign^{\bullet}\in D^b\bigl(\Hol(X,K)\bigr)$:
\[
 \Hom_{D^b(\Hol(Y,K))}\bigl(
 \lefttop{K}f_{!}\nbign^{\bullet},
 \nbigk^{\bullet}
 \bigr)=0 
\]
Hence, we have the following isomorphisms
for any $K$-holonomic $\nbigd_X$-complex
$\nbign^{\bullet}$:
\begin{multline}
 \Hom_{D^b(\Hol(Y,K))}
 \bigl(
 \lefttop{K}f_{!}\nbign^{\bullet},\,
 \nbigm^{\bullet}
 \bigr)
\simeq
 \Hom_{D^b(\Hol(Y,K))}
 \bigl(
 \lefttop{K}f_{!}\nbign^{\bullet},\,
 \lefttop{K}f_{!}
 \lefttop{K}f^{!}
 \nbigm^{\bullet}
 \bigr)
 \\
\simeq
  \Hom_{D^b\bigl(\Hol(X,K)\bigr)}
 \bigl(
 \nbign^{\bullet},\,
 \lefttop{K}f^{!}\nbigm^{\bullet}
 \bigr)
\end{multline}
Hence, we obtain that
the above functor
$\lefttop{K}f^{!}$
is the right adjoint of
$\lefttop{K}f_{!}$.
By taking the dual,
we obtain the left adjoint
$\lefttop{K}f^{\ast}$
of $\lefttop{K}f_{\ast}$.

\subsubsection{Projection}

Let $f:Z\times Y\lrarr Y$ be the natural projection.
Let $(\nbigm,\nbigf)$ be a $K$-holonomic
$\nbigd_Y$-module.
We put
$\lefttop{K}f^{\ast}(\nbigm,\nbigf):=
 \bigl(
 \nbigo_{Z}\boxtimes\nbigm[-\dim Z],
 K_Z\boxtimes\nbigf
 \bigr)$.
It is easy to check that
$\lefttop{K}f^{\ast}(\nbigm,\nbigf)$
is $K$-holonomic.
Thus, we obtain the exact functor
$\lefttop{K}f^{\ast}:
 D^b(\Hol(Y,K))\lrarr
 D^b(\Hol(Z\times Y,K))$.
Let us prove that
$\lefttop{K}f^{\ast}$ is the left adjoint of
$\lefttop{K}f_{\ast}$.
It is enough to repeat the argument
in \S4.4 of \cite{saito2},
which we include for the convenience of readers.
It is enough to construct
natural transformations
$\alpha:\id\lrarr  
 \lefttop{K}f_{\ast}\lefttop{K}f^{\ast}$
and
$\beta:\lefttop{K}f^{\ast}
 \lefttop{K}f_{\ast}\lrarr \id$
such that
\[
 \beta\circ \lefttop{K}f^{\ast}\alpha:
  \lefttop{K}f^{\ast}\nbigm^{\bullet}
\lrarr
 \lefttop{K}f^{\ast}\lefttop{K}f_{\ast}
 \lefttop{K}f^{\ast}\nbigm^{\bullet}
\lrarr
 \lefttop{K}f^{\ast}\nbigm^{\bullet},
\]
\[
 \lefttop{K}f_{\ast}\beta\circ\alpha:
 \lefttop{K}f_{\ast}\nbign^{\bullet}
\lrarr
 \lefttop{K}f_{\ast}\lefttop{K}f^{\ast}\lefttop{K}f_{\ast}
 \nbign^{\bullet}
\lrarr
 \lefttop{K}f_{\ast}\nbign^{\bullet}
\]
are the identities.
We define $\alpha$ as
the external tensor product with
the natural map
$(\cnum,K)\lrarr 
 \bigl(H^0_{DR}(Z),H^0(Z,K)\bigr)$.
For the construction of $\beta$,
the following diagram is used:
\[
 \begin{CD}
 Z\times Y @>{i}>>
 Z\times Z\times Y @>{q_1}>> Z\times Y\\
 @. @V{q_2}VV @V{p_1}VV \\
 @. Z\times Y @>{p_2}>> Y
 \end{CD}
\]
Here, $i$ is induced by the diagonal 
$Z\lrarr Z\times Z$,
$q_j$ are induced by the projection
$Z\times Z\lrarr Z$ onto the $j$-th component,
and $p_j$ are the projections.
We have the following morphisms
of $K$-holonomic $\nbigd$-complexes:
\begin{equation}
\label{eq;09.11.10.21}
 \lefttop{K}f^{\ast}
 \lefttop{K}f_{\ast}\nbigm^{\bullet}
=
 \lefttop{K}p_2^{\ast}
 \lefttop{K}p_{1\ast}\nbigm^{\bullet}
\simeq
 \lefttop{K}q_{2\ast}
 \lefttop{K}q_1^{\ast}\nbigm^{\bullet}
\lrarr
 \lefttop{K}q_{2\ast}
 \bigl(\lefttop{K}i_{\ast}\lefttop{K}i^{\ast}
 \lefttop{K}q_1^{\ast}\nbigm^{\bullet}
 \bigr)
\simeq
 \lefttop{K}i^{\ast}
 \lefttop{K}q_1^{\ast}
 \nbigm^{\bullet}
\end{equation}
\begin{lem}
\label{lem;09.11.10.22}
We have a natural isomorphism
$\lefttop{K}i^{\ast}
 \lefttop{K}q_1^{\ast}\nbigm^{\bullet}
\simeq \nbigm^{\bullet}$
in $D^b(\Hol(Z\times Y,K))$.
\end{lem}
\pf
We have the following morphism of
$K$-holonomic $\nbigd$-complexes
\[
 \nbigm^{\bullet}
\stackrel{\alpha}{\lrarr}
 \lefttop{K}q_{1\ast}
 \lefttop{K}q_1^{\ast}\nbigm^{\bullet}
\lrarr
 \lefttop{K}q_{1\ast}
 \lefttop{K}i_{1\ast}
 \lefttop{K}i_1^{\ast}
 \lefttop{K}q_1^{\ast}\nbigm^{\bullet}
\simeq
 \lefttop{K}i_1^{\ast}
 \lefttop{K}q_1^{\ast}\nbigm^{\bullet}
\]
It is enough to check that
the composite of the morphisms is an isomorphism
for the underlying $\nbigd_Y$-modules.
It is enough to consider the issue
locally around any point of $Z\times Y$.
Then, it can be checked by a direct computation.
\hfill\qed

\vspace{.1in}
We define $\beta$ as the composite of
(\ref{eq;09.11.10.21})
with the isomorphism in 
Lemma \ref{lem;09.11.10.22}.
Let us look at
$\lefttop{K}f_{\ast}\beta\circ\alpha$,
which is the composite of the following
morphisms:
\begin{multline}
 \label{eq;13.4.24.1}
 \lefttop{K}f_{\ast}\nbigm^{\bullet}
=\lefttop{K}p_{1\ast}\nbigm^{\bullet}
\lrarr
 \lefttop{K}p_{2\ast}
 \lefttop{K}p_2^{\ast}
 \lefttop{K}p_{1\ast}\nbigm^{\bullet}
\lrarr
 \lefttop{K}p_{2\ast}
 \lefttop{K}q_{2\ast}
 \lefttop{K}q_1^{\ast}\nbigm^{\bullet}
 \\
\lrarr
 \lefttop{K}p_{2\ast}
 \lefttop{K}q_{2\ast}
 \lefttop{K}i_{\ast}
 \lefttop{K}i^{\ast}
 \lefttop{K}q_1^{\ast}\nbigm^{\bullet}
\lrarr
 \lefttop{K}f_{\ast}
 \lefttop{K}i^{\ast}
 \lefttop{K}q_1^{\ast}\nbigm^{\bullet}
\simeq
 \lefttop{K}f_{\ast}\nbigm^{\bullet}
\end{multline}
We have a natural identification
$p_{2\ast}q_2^{\ast}q_1^{\ast}
\simeq
 p_{1\ast}q_{1\ast}q_1^{\ast}$,
and $p_{1\ast}\lrarr p_{2\ast}q_2^{\ast}q_1^{\ast}$
in (\ref{eq;13.4.24.1})
is induced by $\alpha$ for $q_1$
under the identification.
Then, it is easy to see that
the composite is the identity
by the construction.
As for $\beta\circ \lefttop{K}f^{\ast}\alpha$,
it is expressed as follows:
\begin{multline}
 \label{eq;13.4.24.2}
 \lefttop{K}f^{\ast}\nbign^{\bullet}
=\lefttop{K}p_2^{\ast}\nbign^{\bullet}
\lrarr
 \lefttop{K}p_2^{\ast}
 \lefttop{K}p_{1\ast}
 \lefttop{K}p_1^{\ast}
 \nbign^{\bullet}
\lrarr
 \lefttop{K}q_{2\ast}
 \lefttop{K}q_1^{\ast}
 \lefttop{K}p_1^{\ast}\nbign^{\bullet}
 \\
\lrarr
 \lefttop{K}q_{2\ast}
 \lefttop{K}i_{\ast}
 \lefttop{K}i^{\ast}
 \lefttop{K}q_{2}^{\ast}
 \lefttop{K}p_2^{\ast}\nbign^{\bullet}
\simeq
 \lefttop{K}p_2^{\ast}\nbign^{\bullet}
=\lefttop{K}f^{\ast}\nbign^{\bullet}
\end{multline}
We have a natural identification
$p_2^{\ast}p_1^{\ast}p_{1\ast}
\simeq
 q_{2\ast}q_2^{\ast}p_{2}^{\ast}$,
and 
$p_{2}^{\ast}\lrarr
 p_2^{\ast}p_1^{\ast}p_{1\ast}$
in (\ref{eq;13.4.24.2})
is induced by
$\alpha$ for $q_2$.
Then, it is easy to observe that the composite
is the identity.
Thus, the proof of
Proposition \ref{prop;09.11.10.12}
is finished.
\hfill\qed

%% file: 9.3.tex
\subsection{Statement}

Let $(\nbigm_i,\nbigf_i)$ $(i=1,2)$
be $K$-holonomic $\nbigd$-modules on $X_i$.

\begin{prop}
\label{prop;09.11.10.30}
$\nbigf_1\boxtimes\nbigf_2$
is a $K$-Betti structure of
$\nbigm_1\boxtimes\nbigm_2$.
As a result,
we obtain a natural functor
$\boxtimes:
 \Hol(X_1,K)\times\Hol(X_2,K)
\lrarr 
 \Hol(X_1\times X_2,K)$,
compatible with 
the standard external products
$\boxtimes:
 \Hol(X_1)\times\Hol(X_2)
\lrarr
 \Hol(X_1\times X_2)$
and 
$D^b_c(K_{X_1})\times
 D^b_c(K_{X_2})
\lrarr
 D^b_c(K_{X_1\times X_2})$.
\end{prop}
\index{functor $\boxtimes$}

Before going into the proof of Proposition
\ref{prop;09.11.10.30},
We give a standard consequence.
Let $X$ be an algebraic variety.
Let $\delta_X:X\lrarr X\times X$
be the diagonal morphism.
We obtain the functors
$\otimes$ and $\nrhom$
on $D^b\bigl(\Hol(X,K)\bigr)$
in standard ways:
\[
 \nbigm\otimes\nbign:=
 \lefttop{K}\delta_X^{\ast}\bigl(
 \nbigm \boxtimes\nbign
 \bigr),
\quad
 \nrhom(\nbigm,\nbign):=
 \lefttop{K}\delta_X^!\bigl(
 \DDD_X\nbigm\boxtimes\nbign
 \bigr)
\]
They are compatible with the corresponding
functors on $D^b_{\hol}(X)$.
\index{functor $\otimes$}
\index{functor $\nrhom$}

\subsection{Preliminary}

Let $(\nbigm,\nbigf_{\nbigm})$ be a $K$-holonomic
$\nbigd_X$-module.
Let $\nbigv$ be a meromorphic flat connection
on $(Y,D_Y)$ with 
a good $K$-structure.
Let $\nbigf_{\nbigv}$ and 
$\nbigf_{\nbigv!}$ denote the canonical
$K$-Betti structures of 
$\nbigv$ and $\nbigv_!$,
respectively.

\begin{lem}
\label{lem;09.11.2.2}
$\nbigf_{\nbigv}\boxtimes
\nbigf_{\nbigm}$
and $\nbigf_{\nbigv!}\boxtimes\nbigf_{\nbigm}$
are $K$-Betti structures
of $\nbigv\boxtimes\nbigm$
and $\nbigv_!\boxtimes\nbigm$,
respectively.
\end{lem}
\pf
We use an induction on the dimension
of the support of $\nbigm$.
Let $P$ be any point of $X$.
It is enough to consider locally around
$Y\times \{P\}$.
Let $\nbigc=(Z,U,\varphi,V)$
be a $K$-cell of $\nbigm$ at $P$
with a cell function $g$.
The pre-$K$-holonomic $\nbigd$-module
$\nbigv\otimes\nbigm$ is expressed
as the cohomology of the following complex
of pre-$K$-holonomic $\nbigd$-modules:
\[
 \nbigv\boxtimes \psi_g\bigl(
 \varphi_{\dagger}V\bigr)
\lrarr
 \nbigv\boxtimes\Xi_g\bigl(
 \varphi_{\dagger}V
 \bigr)
\oplus
 \nbigv\boxtimes\phi_g(\nbigm)
\lrarr
 \nbigv\boxtimes\psi_g\bigl(
 \varphi_{\dagger}V
 \bigr)
\]
By the inductive assumption,
$\nbigf_{\nbigv}\boxtimes
 \lefttop{D}\psi_g(\varphi_{\ast}\nbigf_V)$
and
$\nbigf_{\nbigv}\boxtimes
 \lefttop{D}\phi_g(\varphi_{\ast}\nbigf_V)$
are $K$-Betti structures of
$\nbigv\boxtimes\psi_g(\varphi_{\dagger}V)$
and 
$\nbigv\boxtimes\phi_g(\varphi_{\dagger}V)$,
respectively.
We put $g_Z:=\varphi^{\ast}g$.
By using Theorem \ref{thm;09.10.17.151},
we obtain that 
$\nbigf_{\nbigv}\boxtimes
 \lefttop{D}\Xi_{g_Z}\bigl(
 \nbigf_V\bigr)$
and 
$\nbigf_{\nbigv}\boxtimes
 \lefttop{D}\psi_{g_Z}\bigl(
 \nbigf_V\bigr)$
are $K$-Betti structures of
$\nbigv\boxtimes \Xi_{g_Z}(V)$ 
and $\nbigv\boxtimes\psi_{g_Z}(V)$,
respectively.
By construction,
the isomorphism
$\nbigv\boxtimes
\varphi_{\dagger}\bigl(
\psi_{g_Z}(V) \bigr)
\simeq
 \nbigv\boxtimes
 \psi_g\bigl(\varphi_{\dagger}V\bigr)$
preserves $K$-Betti structures.
Hence, we obtain that
$\nbigf_{\nbigm}\boxtimes
 \nbigf_{\nbigv}$ is a $K$-Betti structure.
Thus, we obtain the first claim.
By considering the dual,
we obtain the second claim.
\hfill\qed

\vspace{.1in}
Let $g$ be a holomorphic function on $Y$
such that $g^{-1}(0)=D_Y$.
We obtain the following corollary from
Lemma \ref{lem;09.11.2.2}.
\begin{cor}
\label{cor;09.11.2.3}
$\lefttop{D}\psi_g(\nbigf_{\nbigv})
 \boxtimes \nbigf_{\nbigm}$
and 
$\lefttop{D}\Xi_g(\nbigf_{\nbigv})
 \boxtimes\nbigf_{\nbigm}$
are $K$-Betti structures of
$\psi_g(\nbigv)\boxtimes
 \nbigm$ and
 $\Xi_g(\nbigv)\boxtimes
 \nbigm$,
respectively.
\hfill\qed
\end{cor}

\subsection{Proof of Proposition
\ref{prop;09.11.10.30}}

Let $P$ be any point of $X_1$.
It is enough to consider locally around
$\{P\}\times X_2$.
We use an induction on 
$\dim_P\Supp\nbigm_1$.
Let $\nbigc=(Z,U,\varphi,V)$ be a $K$-cell
of $\nbigm_1$.
The pre-$K$-holonomic $\nbigd$-module
$\nbigm_1\boxtimes\nbigm_2$
is expressed as the cohomology of 
the following complex:
\[
 \psi_{g}(\varphi_{\dagger}V)
\boxtimes
 \nbigm_2
\lrarr
 \Xi_g(\varphi_{\dagger}V)
 \boxtimes\nbigm_2
\oplus
 \phi_g(\nbigm_1)
 \boxtimes\nbigm_2
\lrarr
 \psi_g(\varphi_{\dagger}V)
 \boxtimes\nbigm_2
\]
By the inductive assumption,
$\psi_g(\varphi_{\dagger}V)\boxtimes
 \nbigm_2$ and
$\phi_g(\varphi_{\dagger}V)\boxtimes
 \nbigm_2$ are $K$-holonomic.
According to 
Theorem \ref{thm;09.10.16.5}
and Corollary \ref{cor;09.11.2.3},
$\Xi_g(\varphi_{\dagger}V)
 \boxtimes
 \nbigm_2$ is $K$-holonomic.
Hence, we obtain that
$\nbigm_1\boxtimes\nbigm_2$
is also $K$-holonomic.
Thus, we obtain Proposition
\ref{prop;09.11.10.30}.
\hfill\qed

%% file: 9.4.tex
\subsection{Statements}

\begin{thm}
\label{thm;09.11.13.21}
For $M^{\bullet},N^{\bullet}
\in D^b\bigl(\Hol(X,K)\bigr)$,
the induced morphism
\begin{equation}
 \label{eq;13.4.24.110}
\Hom_{D^b\bigl(\Hol(X,K)\bigr)}(M^{\bullet},N^{\bullet})
 \otimes\cnum
\lrarr
 \Hom_{D^b_{\hol}(X)}(M^{\bullet},N^{\bullet}) 
\end{equation}
is an isomorphism.
In other words,
the forgetful functor
$D^b\bigl(\Hol(X,K)\bigr)\otimes\cnum\lrarr
 D^b_{\hol}(X)$ is fully faithful.
\end{thm}
We closely follow Beilinson's argument
in \cite{beilinson1} for the proof.

\begin{thm}
\label{thm;09.11.13.20}
We have the following natural isomorphism
\[
 \Hom_{D_{\hol}^b(X,K)}
 (M^{\bullet},N^{\bullet})
\simeq
 \Hom_{D^b\bigl(\Hol(X,K)\bigr)}
 \bigl(
 \nbigo_X,\,
 \nrhom(M^{\bullet},N^{\bullet})[d_X]
 \bigr).
\]
\end{thm}
We essentially use a commutative diagram
due to Saito in \cite{saito4}.

\subsection{Homomorphisms
and extensions for meromorphic flat connections
with a good $K$-structure}
\label{subsection;09.11.14.1}

Let $X$ be a smooth complex projective variety
with a hypersurface $D$.

\begin{lem}
Let $V$ be a meromorphic flat connection on $(X,D)$
with a good $K$-structure.
Let $\nbigf_V$ be the canonical $K$-Betti structure of $V$.
We have the following natural isomorphisms
for $i=0,1$:
\[
 \Ext^i_{\Hol(X,K)}\bigl(
 \nbigo_X(\ast D),V
 \bigr)
\simeq
 H^i(X,\nbigf_V[-d_X])
\]
\end{lem}
\pf
By taking a global resolution of turning points 
in the algebraic situation
(\cite{kedlaya2}, \cite{mochi7}),
we may assume that $V$ is a good meromorphic flat bundle.
Let $\nbigl(V)$ be the associated local
system with the Stokes structure on $\Xtilde(D)$.
It is naturally equipped with a $K$-structure
$\nbigl_K(V)$.
If we are given an extension
$0\lrarr V\lrarr P\lrarr\nbigo_X(\ast D)\lrarr 0$
as $K$-holonomic $\nbigd_X$-modules,
$P$ is also a good meromorphic flat bundle
with a good $K$-structure,
and it induces an extension
$0\lrarr \nbigl_K(V)^{\leq D}
\lrarr \nbigl_K(P)^{\leq D}
\lrarr K_{\Xtilde(D)}\lrarr 0$
of $K$-constructible sheaves.
Conversely,
assume that we are given 
an extension of $K$-constructible sheaves
$0\lrarr\nbigl_K(V)^{\leq D}\lrarr
 \nbigg_K\lrarr K_{\Xtilde(D)}\lrarr 0$.
We obtain a $K$-local system
$\nbiggtilde_K:=
 \iotatilde_{\ast}\nbigg_{|X\setminus D}$,
where $\iota:X\setminus D\lrarr X$.
The $\cnum$-local system
$\nbiggtilde_{K}\otimes \cnum$
is naturally equipped with
a Stokes structure compatible with
the $K$-structure.
Hence, we obtain an extension of
$K$-holonomic $\nbigd_X$-modules
$0\lrarr V\lrarr P\lrarr\nbigo_X(\ast D)\lrarr 0$.
The above procedures are mutually inverse.
Thus, we obtain a bijection
$\Ext^1_{\Hol(X,K)}\bigl(\nbigo_X(\ast D),\,V\bigr)
\simeq
 \Ext^1_{K_{\Xtilde(D)}}\bigl(
 K_{\Xtilde(D)},\,\nbigl_K(V)^{\leq D}\bigr)
\simeq
 H^1\bigl(X,\nbigf_V[-d_X]\bigr)$. 
Similarly,
we have a natural isomorphism
$\Ext^0_{\Hol(X,K)}\bigl(\nbigo_X(\ast D),V\bigr)
\simeq
 H^0(X,\nbigf_V[-d_X])$.
\hfill\qed

\vspace{.1in}

Let $V,W$ be meromorphic flat connections
on $(X,D)$ with good $K$-structures.
We have a natural bijection
$\Ext^i_{\Hol(X,K)}(W,V)
\simeq
 \Ext^i_{\Hol(X,K)}\bigl(\nbigo_{X}(\ast D),
W^{\lor}\otimes V\bigr)$
for any $i$.
We obtain the natural isomorphisms
$\Ext^i_{\Hol(X,K)}\bigl(W,V\bigr)
\simeq
 H^i\bigl(
 X,\nbigf_{W^{\lor}\otimes V}[-d_X]
 \bigr)$ for $i=0,1$.
Because
\[
H^i\bigl(
 X,\nbigf_{W^{\lor}\otimes V}[-d_X]
 \bigr)\otimes_K\cnum
\simeq
H^i\bigl(X,\DR_X(W^{\lor}\otimes V)[-d_X]\bigr)
=:H^i_{\DR}(X,W^{\lor}\otimes V),
\]
the vector spaces
$H^i_{\DR}(X,W^{\lor}\otimes V)$
have the natural $K$-structure.
We say that an element 
$f\in H^i_{\DR}(X,W^{\lor}\otimes V)$
is compatible with $K$-structure
if it comes from
$H^i\bigl(
 X,\nbigf_{W^{\lor}\otimes V}[-d_X]
 \bigr)$.
An element $f\in H^1_{\DR}(X,W^{\lor}\otimes V)$
induces an extension
$0\lrarr V\lrarr P\lrarr W\lrarr 0$
in $\Hol(X,K)$
as observed above.

\subsection{Some extensions}

Let $X$ be a smooth complex 
quasi-projective variety.
Let $V_i$ $(i=1,2)$ be algebraic flat bundles on $X$ 
with a good $K$-structure,
i.e.,
there exists a projective variety $\Xbar\supset X$
such that 
(i) $D:=\Xbar-X$ is normal crossing,
(ii) $V_i$ are good meromorphic flat bundles on
 $(\Xbar,D)$ with a good $K$-structure.
According to \cite{beilinson1},
we have
$\Ext^i_{\Hol(X)}(V_1,V_2)\simeq
 H^i_{DR}\bigl(X,V_1^{\lor}\otimes V_2\bigr)$.

\begin{lem}
\label{lem;09.11.13.1}
There exist a Zariski open subset $U\subset X$
and an extension 
$V_3\supset V_{2|U}$ on $U$
of algebraic flat bundles with a good $K$-structure,
such that the induced morphisms
$\Ext^i_{\Hol(X)}(V_1,V_2)
\lrarr
 \Ext^i_{\Hol(U)}(V_{1|U},V_3)$
are $0$ for $i>0$.
\end{lem}
\pf
We use an induction on $\dim X$.
In the case $\dim X=0$,
the claim is trivial.
Let us consider the case $\dim X>0$.
We take a Zariski open subset 
$X_1\subset X$ with a smooth affine fibration
$\rho:X_1\lrarr Z_1$ such that 
the relative dimension is $1$.
For any algebraic flat bundle $\nbigv$ on
$X_1$,
we put $\rho^q_{\ast}(\nbigv):=
 R^q\rho_{\ast}\bigl(
 \nbigv\otimes\Omega^{\bullet}_{X_1/Z_1}
 \bigr)$.
For a Zariski open subset
$Z_1'\subset Z_1$,
the induced morphism
$\rho^{-1}(Z_1')\lrarr Z_1'$
is also denoted by $\rho$.

We may assume that
$L_q:=
 \rho^q_{\ast}(V_1^{\lor}\otimes V_2)$
are algebraic flat bundles on $Z_1$,
which is equipped with 
the induced good $K$-structure.
We have $L_q=0$ unless $q=0,1$.
By the argument in \S2.1 of \cite{beilinson1},
we can reduce Lemma \ref{lem;09.11.13.1}
to Lemma \ref{lem;09.11.13.2} below
which is Lemma 2.1.2 of \cite{beilinson1}
with a minor enhancement.

\begin{lem}
\mbox{{}}\label{lem;09.11.13.2}
\begin{description}
\item[(a)]
 There exist a Zariski open subset
 $Z_2\subset Z_1$ and
 an extension 
 $P\supset V_{2|X_2}$ of
 algebraic flat bundles
 with good $K$-structure
 on $X_2:=\rho^{-1}(Z_2)$,
 such that
 the induced morphism
 $\rho^1_{\ast}(V_1^{\lor}\otimes
 V_{2|X_2})
\lrarr
 \rho^1_{\ast}(V_1^{\lor}\otimes P)$ is $0$.
\item[(b)]
There exists a Zariski open subset
$Z_3\subset Z_1$ and
an extension $Q\supset V_{2|X_3}$
of algebraic flat bundles with good $K$-structure
on $X_3:=\rho^{-1}(Z_3)$, such that
the induced maps
\[
 H^p_{\DR}\bigl(Z_3,
 \rho^0_{\ast}(V_1^{\lor}\otimes
 V_{2|X_3})
 \bigr)
\lrarr
 H^p_{\DR}\bigl(Z_3,
 \rho^0_{\ast}(V_1^{\lor}\otimes Q)
 \bigr)
\]
are $0$ for any $p>0$.
\end{description}
\end{lem}
\pf
It is enough to use the argument
in the proof of Lemma 2.1.2 of \cite{beilinson1}.
We give only an indication.
Let 
$\alpha\in H^0_{\DR}\bigl(
 Z_1,L_1^{\lor}\otimes L_1
 \bigr)
=
 H^0_{\DR}\bigl(Z_1,\rho^1_{\ast}(
 (\rho^{\ast}L_1\otimes V_1)^{\lor}\otimes
 V_2)
 \bigr)$
be the element corresponding to
the identity of $L_1$,
which is compatible with $K$-structure.
We have the following exact sequence
compatible with $K$-structures:
\begin{multline*}
 H^1_{\DR}\Bigl(
 X_1,
\bigl(\rho^{\ast}L_1\otimes V_1\bigr)^{\lor}
 \otimes V_2
 \Bigr)
\lrarr
 H^0_{\DR}\Bigl(
 Z_1, \rho^1_{\ast}\bigl(
 (\rho^{\ast}L_1\otimes V_1)^{\lor}\otimes V_2
 \bigr)
 \Bigr) \\
\stackrel{\del}{\lrarr}
 H^2_{\DR}\Bigl(Z_1,
 \rho^0_{\ast}\bigl(
 (\rho^{\ast}L_1\otimes V_1)^{\lor}\otimes
 V_2 \bigr)\Bigr)
=H^2_{\DR}(Z_1,L_1^{\lor}\otimes L_0)
\end{multline*}
Applying the inductive assumption
to $L_0^{\lor}$ and $L_1^{\lor}$,
we have a Zariski open subset
$Z_2\subset Z_1$ and
an extension
$\varphi:L_1^{\lor}\subset R$
of algebraic flat bundles
with a good $K$-structures
on $Z_2$, such that
the induced morphism
$H^2\bigl(
 Z,L_1^{\lor}\otimes L_0
 \bigr)
\lrarr
 H^2(Z_1,R\otimes L_0)$
is $0$.
In particular,
$\varphi(\del\alpha)=0$.
We obtain the element
\[
 \varphi(\alpha)\in
 H^0_{\DR}\bigl(Z_1,R\otimes L_1
\bigr)
=H^0_{\DR}\Bigl(Z_1,
 \rho^1_{\ast}\bigl(
 (\rho^{\ast}R^{\lor}\otimes V_1)^{\lor}
 \otimes V_2 \bigr) \Bigr) 
\]
which is compatible with $K$-structure.
By construction,
we have a lift
\[
 \widetilde{\varphi(\alpha)}\in
 H^1_{\DR}\Bigl(X,
 (\rho^{\ast}R^{\lor}\otimes V_1)^{\lor}
 \otimes V_2\Bigr)
\]
compatible with $K$-structure.
It induces an extension
$0\lrarr
 V_{2|X_2}\lrarr P\lrarr
 \rho^{\ast}R^{\lor}
 \otimes V_{1|X_2}\lrarr 0$
of algebraic flat bundles
with good $K$-structure on $X_2$.
(See \S\ref{subsection;09.11.14.1}.)
It is easy to observe that
$P$ is the desired one.
Thus, we obtain the claim (a).
The claim (b) can also be proved
by the argument in \cite{beilinson1}.
\hfill\qed

\subsection{Vanishing and lifting}

Let $X$ be a smooth quasi-projective variety.
We put $C_1(X):=\Hol(X)$
and $C_2(X):=\Hol(X,K)\otimes\cnum$.
Let $V_i$ $(i=1,2)$ be algebraic flat bundles on $X$
with good $K$-structure.
Let us consider the natural morphism:
\[
g_X: \Ext^i_{C_2(X)}(V_1,V_2)
\lrarr
 \Ext^i_{C_1(X)}(V_1,V_2)
\]
They are isomorphisms in the cases $i=0,1$
(\S\ref{subsection;09.11.14.1}).

\begin{lem}
\mbox{{}}\label{lem;09.11.13.10}
Let $i>0$.
\begin{itemize}
\item Let $a\in \Ext^i_{C_2(X)}(V_1,V_2)$
such that $g_X(a)=0$.
There exists $U\subset X$ such that
$a=0$ in $\Ext^i_{C_2(U)}(V_{1|U},V_{2|U})$.
\item
Let $a\in \Ext^i_{C_1(X)}(V_1,V_2)$.
There exist $U\subset X$
and $b\in \Ext^i_{C_2(U)}(V_{1|U},V_{2|U})$
such that $a_{|U}=g_U(b)$.
\end{itemize}
\end{lem}
\pf
We give only an outline.
We use an induction on $i$.
We have already known the case $i=1$
Let $a\in \Ext^i_{C_2(X)}(V_1,V_2)$
such that $g_X(a)=0$.
We have an extension
$V_2\subset V_3$ of a meromorphic
flat bundle with a good $K$-structure
such that the image of $a$ is mapped to $0$
via $\Ext^i_{C_2(X)}(V_1,V_2)
\lrarr \Ext^{i}_{C_2(X)}(V_1,V_3)$.
Let $\nbigk:=V_3/V_2$.
We have $c\in \Ext^{i-1}_{C_2(X)}(V_1,\nbigk)$
which is mapped to $a$ via
$\Ext^{i-1}_{C_2(X)}(V_1,\nbigk)
 \lrarr
 \Ext^i_{C_2(X)}(V_1,V_2)$.
We have $d\in \Ext^{i-1}_{C_1(X)}(V_1,V_3)$
which is mapped to $g_X(c)$
via $\Ext^{i-1}_{C_1(X)}(V_1,V_3)
\lrarr \Ext^{i-1}_{C_1(X)}(V_1,\nbigk)$.
By using the inductive assumption,
we can find $U\subset X$
and $e\in \Ext^{i-1}_{C_2(U)}(V_1,V_3)$
such that $g_U(e)=d_{|U}$.
By using the inductive assumption,
and by shrinking $U$,
we may assume that $e$ is mapped to $c_{|U}$
via $\Ext^{i-1}_{C_2(X)}(V_1,V_3)
\lrarr \Ext^{i-1}_{C_2(X)}(V_1,\nbigk)$.
Hence, we obtain 
$a_{|U}=0$.

Let $a\in \Ext^i_{C_1(X)}(V_1,V_2)$.
According to Lemma \ref{lem;09.11.13.1},
we can find $U\subset X$ and
an extension $V_{2|U}\subset V_3$
of meromorphic flat bundles with good $K$-structures
such that the induced map
$\Ext^j_{C_1(U)}(V_{1|U},V_{2|U})\lrarr
 \Ext^j_{C_1(U)}(V_{1|U},V_3)$ is $0$
for any $j>0$.
We put $\nbigk:=V_3/V_{2|U}$
We can find $c\in \Ext^{i-1}_{C_1(U)}(V_{1|U},\nbigk)$
which is mapped to $a$
via $\Ext^{i-1}_{C_1(U)}(V_{1|U},\nbigk)\lrarr
 \Ext^i_{C_1(U)}(V_{1|U},V_{2|U})$.
By using the inductive assumption
and by shrinking $U$,
we can find $d\in \Ext^{i-1}_{C_2(U)}(V_{1|U},\nbigk)$
such that $g_{U}(d)=c$.
Let $b$ be the image of $d$ via 
$\Ext^{i-1}_{C_2(U)}(V_{1|U},\nbigk)\lrarr
 \Ext^i_{C_2(U)}(V_{1|U},V_{2|U})$.
Then, it has the desired property.
\hfill\qed

\subsection{Support}

Let $X$ be a smooth quasi-projective variety.
For any subvariety $Z\subset X$,
let $D^b_{j,Z}(X)$ $(j=1,2)$ denote the derived category of 
bounded complexes $M^{\bullet}$ in $C_j(X)$
such that the supports of $\nbigh^{\bullet}(M^{\bullet})$
are contained in $Z$.
For any $M^{\bullet},N^{\bullet}$ in $D^b_{j,Z}(X)$,
we set
\[
 \Hom^k_{j,Z}\bigl(M^{\bullet},N^{\bullet}\bigr)
:=\Hom_{D^b_{j,Z}(X)}\bigl(
 M^{\bullet},N^{\bullet}[k]
 \bigr).
\]
If $Z=X$, we omit to denote $Z$.
If $Z$ is smooth,
then $D^b_{j,Z}(X)$ is equivalent
to the derived category of $C_j(Z)$.
(See Proposition \ref{prop;10.1.12.10}.)

Let $i:Z\lrarr X$ denote the inclusion.
The natural exact functor
$D^b_{j,Z}(X)\lrarr D^b_{j}(X)$
is denoted by $i_{\ast}$.
As in \S\ref{subsection;13.4.24.100},
we have a functor
$i^!:D^b_{j}(X)\lrarr D^b_{j,Z}(X)$.
We set
$i^{\ast}:=\DDD_X\circ i^{!}\circ \DDD_X$.

\subsection{Proof of
Theorem \ref{thm;09.11.13.21}}

Let $X$ be a smooth quasi-projective variety.
Let $M^{\bullet},N^{\bullet}\in D^{b}_2(X)$.
Let us prove that (\ref{eq;13.4.24.110})
is an isomorphism.
We use an induction on $\dim X$.
It is enough to prove that
(\ref{eq;13.4.24.110})
is an isomorphism 
in the case $M,N\in C_2(X)$.
Take any hypersurface 
$D\subset X$.
Let $i:D\lrarr X$ denote the inclusion.
We have the distinguished triangles
$i_{\ast}i^{!}N\lrarr N \lrarr
 N(\ast D)\stackrel{+1}{\lrarr}$
and
$M(!D)\lrarr M\lrarr
 i_{\ast}i^{\ast}M
 \stackrel{+1}{\lrarr}$.
For $j=1,2$,
we obtain the following exact sequence:
\begin{multline}
\label{eq;13.4.24.120}
 \Ext^{k-1}_{C_j(X)}\bigl(M(!D), N(\ast D)\bigr)
\lrarr
 \Hom^k_{j,D}(
 i_{\ast}i^{\ast}M,\,
 i_{\ast}i^{!}N)
\lrarr
 \Ext^k_{C_j(X)}(M,N) \\
\lrarr
  \Ext^{k}_{C_j(X)}\bigl(M(!D), N(\ast D)\bigr)
\lrarr
 \Hom^{k+1}_{j,D}(
 i_{\ast}i^{\ast}M,\,
 i_{\ast}i^{!}N)
\end{multline}
We naturally have
$\Ext^i_{C_j(X)}(M(!D),N(\ast D))
\simeq
 \Ext^i_{C_j(X)}(M(\ast D),N(\ast D))$,
as remarked in Lemma \ref{lem;10.1.13.1}.

By using the exact sequences (\ref{eq;13.4.24.120})
in the case where $D$ is smooth,
and by using the inductive assumption,
we can reduce the issue
to the case where $X$ is affine,
which we will assume in the following.

We use an induction on the dimension
of the support of $M\oplus N$.
We take a projective birational morphism
$\varphi:Z\lrarr \Supp(M\oplus N)$
such that $Z$ is smooth.
There exist an open subset $U\subset Z$,
flat bundles $V_N$ and $V_M$ on $U$
with morphisms
$M\lrarr \varphi_{\dagger}V_M$
and $N\lrarr\varphi_{\dagger}V_N$
which is an isomorphism on generic points of 
$\Supp(M\oplus N)$.
If we shrink $U$ appropriately,
there exists a hypersurface $D\subset X$
such that
$\varphi^{-1}(D)=Z\setminus U$.
In that case,
we have
$M(\ast D)=\varphi_{\dagger}V_M$
and 
$N(\ast D)=\varphi_{\dagger}V_N$.
In the exact sequence (\ref{eq;13.4.24.120}),
the dimension of the supports 
of the cohomology sheaves of
$i_{\ast}i^{\ast}M$
and $i_{\ast}i^!N$
are strictly smaller than
$\dim \Supp(M\oplus N)$.
Then, it is easy to obtain that
(\ref{eq;13.4.24.110})
for $i_{\ast}i^{\ast}M$
and $i_{\ast}i^!N$
is an isomorphism.
By using Proposition \ref{prop;10.1.12.10},
we obtain
\[
 \Ext^k_{C_j(X)}\bigl(M(!D),N(\ast D)\bigr)
\simeq
 \Ext^{k}_{C_j(X)}\bigl(M(\ast D),N(\ast D)\bigr)
\simeq
 \Ext^k_{C_j(U)}\bigl(V_M,V_N\bigr).
\]

For $D_1\subset D_2$,
we have the following commutative diagram:
\[
 \begin{CD}
 M(!D_1) @>>> M\\
 @AAA @A{=}AA \\
 M(!D_2)@>>> M
 \end{CD}
\quad\quad\quad
 \begin{CD}
 N @>>> N(\ast D_1)\\
 @V{=}VV @VVV \\
 N@>>> N(\ast D_2)
 \end{CD}
\]
Let $i_a:D_a\lrarr X$ denote the inclusions.
We set $U_a:=Z\setminus \varphi^{-1}(D_a)$.
Hence, we have the following commutative
diagram:
\[
 \begin{CD}
 \Hom^i_{j,D_1}(i_{1\ast}i_1^{\ast}M,
 i_{1\ast}i_1^!N)
 @>>>
 \Ext^i_{C_j(X)}(M,N)
 @>>>
 \Ext^i_{C_j(U_1)}(V_M,V_N)\\
 @VVV @V{=}VV @VVV \\
 \Hom^i_{j,D_2}(
 i_{2\ast}i_2^{\ast}M,
 i_{2\ast}i_2^!N)
 @>>>
 \Ext^i_{C_j(X)}(M,N)
 @>>>
 \Ext^i_{C_j(U_2)}(V_M,V_N)\\
 \end{CD}
\]
Then, it is easy to prove that
$\Ext^i_{C_2(X)}(M,N)
\lrarr \Ext^i_{C_1(X)}(M,N)$
is an isomorphism
by using Lemma \ref{lem;09.11.13.10}.
\hfill\qed

\subsection{Proof of
Theorem \ref{thm;09.11.13.20}}

Recall a commutative diagram in 
Proposition 4.6 of \cite{saito4}.
For $M^{\bullet},N^{\bullet}\in D^b(\nbigd_X)$,
we have the following commutative diagram:
{\small
\begin{equation}
 \label{eq;09.12.6.10}
  \begin{CD}
 \Hom_{D(\nbigd_X)}(M^{\bullet},N^{\bullet})
 @>{\simeq}>>
 \Hom_{D(\nbigd_{X\times X})}\bigl(
 M^{\bullet}\boxtimes \DDD N^{\bullet},\,
 \delta_{\dagger}\nbigo_X[d_X]
 \bigr)\\
 @VVV @VVV \\
 \Hom_{D(\cnum_X)}\bigl(
 \DR_XM^{\bullet},\,\DR_XN^{\bullet}\bigr)
 @>{\simeq}>>
  \Hom_{D(\cnum_X)}\bigl(
 \DR_XM^{\bullet}\otimes\DDD\DR_XN^{\bullet},\,
 \delta_{\ast}\cnum_X[2d_X]
 \bigr)
 \end{CD}
\end{equation}
}
Let $M$ be a holonomic $\nbigd_X$-module
with a $K$-Betti structure $\nbigf$.
We have 
\[
\Hom_{D(\nbigd_X)}(M,M)
\simeq
 \Hom_{\Hol(X)}(M,M)
\simeq
 \Hom_{\Hol(X,K)}(M,M)\otimes\cnum 
\]
We have similar isomorphisms for 
$\Hom_{D(\nbigd_X)}\bigl(M\boxtimes \DDD M,\,
 \delta_{\dagger}\nbigo_X[d_X]\bigr)$.
Hence, we obtain the following diagram
from (\ref{eq;09.12.6.10}):
{\small
\[
 \begin{CD}
 \Hom_{\Hol(X,K)}(M,M)\otimes\cnum
 @>{c}>{\simeq}>
 \Hom_{\Hol(X\times X,K)}\bigl(
 M\boxtimes \DDD M,\,
 \delta_{\dagger}\nbigo_X[d_X]
 \bigr)\otimes\cnum\\
 @V{a}VV @V{b}VV \\
 \Hom_{D(\cnum_X)}\bigl(
 \DR_XM,\,\DR_XM\bigr)
 @>{\simeq}>>
  \Hom_{D(\cnum_X)}\bigl(
 \DR_XM\otimes\DDD\DR_XM,\,
 \delta_{\ast}\cnum_X[2d_X]
 \bigr)\\
 @A{\simeq}AA @A{\simeq}AA \\
  \Hom_{D(K_X)}\bigl(
 \nbigf,\nbigf \bigr)\otimes\cnum
 @>{\simeq}>>
  \Hom_{D(K_X)}\bigl(
 \nbigf\boxtimes\DDD\nbigf,\,
 \delta_{\ast}K_X[2d_X]
 \bigr)\otimes\cnum
 \end{CD}
\]
}
Note that $a$ is injective.
Hence, $b$ is also injective.
Since $a$ and $b$ are compatible with
$K$-structures,
$c$ is also compatible with $K$-structures.
Let $C:M\boxtimes\DDD M\lrarr\delta_{\ast}\nbigo_X[d_X]$
correspond to $1:M\lrarr M$.
It is compatible with $K$-Betti structures.

For $M^{\bullet}\in D^b\bigl(\Hol(X,K)\bigr)$,
let $C:M^{\bullet}\boxtimes \DDD M^{\bullet}
\lrarr
 \delta_{\dagger}\nbigo_{X}[d_X]$ 
correspond to $1:M^{\bullet}\lrarr M^{\bullet}$.
We obtain that
$C$ is compatible with $K$-Betti structures.
Then, we obtain that
the isomorphism
\[
 \Hom_{D(\nbigd_X)}(M^{\bullet},N^{\bullet})
\lrarr
 \Hom_{D(\nbigd_{X\times X})}
 \bigl(M^{\bullet}\boxtimes
 \DDD N^{\bullet},\,
 \delta_{\dagger}\nbigo_X[d_X]\bigr)
\]
is compatible with $K$-Betti structures
for any $M^{\bullet},N^{\bullet}
 \in D_{\hol}(X,K)$.
By taking the dual,
we obtain Theorem \ref{thm;09.11.13.20}.
\hfill\qed